\documentclass[11pt]{article}
\pdfoutput=1

\usepackage{amsmath, amsthm, amssymb}
\usepackage{booktabs}
\usepackage{tikz}
\usepackage{pgfplots}
\pgfplotsset{compat=1.18}
\usetikzlibrary{fpu}
\usepackage{tikz-cd}
\usepackage[toc,page]{appendix}
\usepackage{algorithm}
\usepackage{algpseudocode}

\usetikzlibrary{cd, arrows, arrows.meta, positioning, fit, calc, backgrounds, shadows, decorations.pathmorphing, patterns, shapes.geometric}
\usepackage{geometry}
\usepackage{enumitem}
\usepackage{array} 
\usepackage{xcolor}
\usepackage{courier}
\usepackage{listings}

\usepackage[utf8]{inputenc}

\usepackage{float}
\floatstyle{boxed}
\restylefloat{figure}

\usepackage[colorlinks=true, urlcolor=blue, linkcolor=red, citecolor=purple]{hyperref}

\usepackage[T1]{fontenc}
\IfFileExists{comicneue.sty}{\usepackage{comicneue}}{}

\tikzset{
    sketch/.style={
        decorate,
        decoration={random steps, segment length=4pt, amplitude=0.6pt}, 
        thick,
        line cap=round,
        font=\sffamily\bfseries 
    }
}

\definecolor{coreBlue}{RGB}{230, 240, 255}
\definecolor{coreBorder}{RGB}{100, 149, 237}
\definecolor{semGreen}{RGB}{235, 250, 235}
\definecolor{semBorder}{RGB}{60, 179, 113}
\definecolor{plugOrange}{RGB}{255, 245, 230}
\definecolor{plugBorder}{RGB}{210, 105, 30}
\definecolor{dataGray}{RGB}{240, 240, 240}

\definecolor{synBlue}{RGB}{230, 240, 255}
\definecolor{synBorder}{RGB}{70, 130, 180}
\definecolor{userGray}{RGB}{245, 245, 245}
\definecolor{auditRed}{RGB}{255, 235, 235}

\newtheorem{theorem}{Theorem}[section]
\newtheorem{lemma}[theorem]{Lemma}
\newtheorem{proposition}[theorem]{Proposition}
\newtheorem{corollary}[theorem]{Corollary}
\theoremstyle{definition}
\newtheorem{definition}[theorem]{Definition}
\newtheorem{example}[theorem]{Example}
\theoremstyle{remark}
\newtheorem{remark}[theorem]{Remark}

\DeclareMathOperator{\dom}{dom}
\DeclareMathOperator{\cl}{cl}
\DeclareMathOperator{\Graph}{Graph}

\DeclareMathOperator{\Image}{Im}

\newcommand{\HSSimp}{\mathbb{HS}} 
\newcommand{\HSPr}{\mathrm{HSP\text{-}r}} 
\newcommand{\A}{\mathcal{A}}
\DeclareMathOperator{\supp}{supp}

\DeclareMathOperator{\DOTS}{DOTS}
\DeclareMathOperator{\Fee}{Fee}\newcommand\R{\mathbb{R}}

\newcommand{\into}{\dashrightarrow} 
\newcommand{\id}{\Delta} 
\newcommand{\vcomp}{\circ} 
\newcommand{\meet}{\wedge} 
\newcommand{\adj}{\dashv} 
\newcommand{\pullback}[1]{#1^*}

\definecolor{codegreen}{rgb}{0,0.6,0}
\definecolor{codegray}{rgb}{0.5,0.5,0.5}
\definecolor{codepurple}{rgb}{0.58,0,0.82}
\definecolor{backcolour}{rgb}{0.95,0.95,0.92}

\title{A Double Categorical Framework for\\
Multi-Stage Portfolio Construction and Alignment}
\author{Wesley Phoa\\[2pt]\texttt{wesley@topos.institute}}

\begin{document}

\maketitle

\begin{abstract}
We construct a thin double category $\mathbb{HS}$ (Hub-and-Spoke) whose objects are closed subsets of standard simplices, horizontal morphisms are continuous maps representing portfolio re-implementation processes, and vertical morphisms are closed relations representing alignment constraints. This framework models industrial portfolio construction pipelines---hierarchical structures in which a single investment strategy is translated through multiple stages into thousands of client portfolios. We establish four structural theorems: compositionality of alignment (functoriality), a pre-trade safety guarantee (adjunction), an order-independence result for compliance checking (lax Beck--Chevalley), and a filter-commutation law (Frobenius reciprocity). The topological requirement that permissible portfolio spaces be closed and compact---ruling out ``phantom portfolios'' that arise from open constraint specifications---is shown to be essential for coherence. Extensions to set-valued re-implementations via the Double Operadic Theory of Systems, stochastic re-implementations via Markov kernels on Polish spaces, and transport-based safety metrics via Wasserstein distances are developed. An abstract axiomatic treatment identifies the equipment axioms sufficient for the main results. The mathematical content is elementary---no novel category theory is required. The contribution is the modelling claim: that these particular objects and morphisms formalise portfolio re-implementation correctly.
\end{abstract}

\medskip

\noindent\textbf{arXiv subject areas:} math.CT (primary); q-fin.PM (cross-list)

\smallskip


\smallskip


\clearpage

\tableofcontents

\clearpage

\listoffigures

\clearpage

\section*{Preface and Guide for the Reader}
\addcontentsline{toc}{section}{Preface and Guide for the Reader}

This work was motivated by the following diagram, which the author drew during an \emph{ad hoc} business discussion. It depicts a hierarchical portfolio construction process where ``hub'' portfolios are used as the basis for constructing ``spoke'' portfolios, which incorporate new assets, different vehicle types or customized client constraints, while having the same essential characteristics. This gives rise to a tree-like structure. In addition, alignment requirements (e.g. exposure, risk) are imposed between nodes in different branches. 

The diagram was originally intended simply to clarify the nature of an investment process, and its interpretation was purely informal. However, it is natural to wonder whether the diagram can be given a mathematical meaning, reflecting formal properties that an investment process should have.

\vspace{0.2cm}
 

\begin{center}
    
\begin{tikzcd}[
    cells={nodes={draw, rectangle, sketch, minimum height=2em}}, 
    arrows={sketch, >=stealth}, 
    math mode=false, 
    font=\sffamily\bfseries,
    row sep=1.2cm,
    column sep=1.5cm
  ]
& FoF & & \\
& ETF 
  \arrow[r]                 
  \arrow[ddr, dashed, -{}]  
  \arrow[ddrr]              
& AP-ETF \arrow[r] 
& AP-ETF+ \arrow[dd, dashed, -{}] 
\\
MF 
  \arrow[uur]      
  \arrow[ur]       
  \arrow[r]        
  \arrow[dr]      
  \arrow[dd, dashed, -{}] 
& VIFoF & & \\
& Client1 
  \arrow[r] 
  \arrow[d, draw=black!40, dashed, text=black!40, -{}]
& Client1-ETF 
& ETF+ 
\\
Client2 
& |[rectangle, draw=black!40, dashed, text=black!40, sketch]| Client1-CIO 
& &
\end{tikzcd}

\end{center}

\begin{center}

\vspace{0.25cm} 

\begin{tikzcd}[
    cells={nodes={draw, rectangle, sketch, minimum height=2em}}, 
    arrows={sketch, >=stealth},
    every label/.append style={font=\sffamily\bfseries},
    math mode=false,
    column sep=4.5cm 
  ]
\text{Portfolio 1} \arrow[r, "directly influences"] & \text{Portfolio 2}
\end{tikzcd}
    
\vspace{0.25cm} 

\begin{tikzcd}[
    cells={nodes={draw, rectangle, sketch, minimum height=2em}}, 
    arrows={sketch, >=stealth},
    every label/.append style={font=\sffamily\bfseries},
    math mode=false,
    column sep=4.5cm 
  ]
\text{Portfolio 1} \arrow[r, dashed, -{}, "should be aligned"] & \text{Portfolio 2}
\end{tikzcd}

\end{center}

\vspace{0.3cm} 

We describe a category-theoretic framework in which this diagram has a formal interpretation. As mass customization gains pace, portfolio construction must become increasingly automated in order to scale, and automated systems need formal rules. Traditional portfolio theory describes the components of a portfolio but does not naturally handle the \emph{composition of constraints} across multiple stages of construction and transformation. The double categorical framework is necessary to handle both the composition of transformations (solid arrows in the diagram) and the composition of alignment constraints (dotted lines, though the simple diagram above contains no composites). The framework is intended to be useful in practice, especially in industrial-scale portfolio construction environments---customization or replication of thousands of portfolios---where manual oversight becomes impracticable.

\subsection*{Guide for Readers: Applied Category Theorists}

The framework, denoted $\mathbb{HS}$ (Hub-and-Spoke), is a thin double category whose objects are closed subsets $K \subseteq \Delta^n$ of standard simplices, horizontal morphisms are continuous maps (representing portfolio re-implementations), vertical morphisms are closed relations (representing alignment constraints), and 2-cells are inclusions. If you are familiar with equipments and framed bicategories \cite{Shulman}, this will be recognizable territory. The financial terminology can be safely ignored on a first pass; the categorical content is self-contained.

The genuine interest, for a categorist, lies in three places. First, the restriction to compact objects makes all horizontal maps automatically proper, which is what allows the pushforward $f_!$ to preserve closedness---the topological hypothesis does categorical work. Second, the main coherence theorems (Beck--Chevalley, Frobenius) turn out to formalize standard audit questions in finance: ``Does filtering for ESG compliance before optimization yield the same result as filtering after?'' is an instance of Frobenius reciprocity. Third, Section~\ref{sec:DOTS} connects the framework to the Double Operadic Theory of Systems \cite{LibkindMyers}, interpreting the set of valid portfolio re-implementations as an algebra over an operad of wiring diagrams. Familiarity with \cite{Markowitz} or \cite{GrinoldKahn} will help with the financial context of the optimization problems that define the morphisms, but is not strictly necessary.

\clearpage

\subsection*{Guide for Readers: Asset Management Professionals}

If you work in portfolio construction, you already know the objects of this paper. A permissible portfolio space ($K$) is the set of portfolios your mandate allows. A re-implementation ($f$) is what your optimizer does when it converts a model portfolio into client accounts. An alignment ($R$) is what your compliance team checks---tracking error within bounds, sector exposures within tolerance.

The contribution of the paper is showing that these three things form a mathematical structure---a double category---whose internal coherence laws answer questions you currently answer by ad~hoc means. Standard optimizer pipelines often fail at the edges, producing portfolios that look correct numerically but have no valid reference portfolio (``phantom portfolios''). The paper proves three things that matter in practice: (a)~if your constraints are specified as closed, compact sets, the optimization process is stable---small changes in market data will not cause wild jumps in the output; (b)~properness ensures that every client portfolio can be traced back to a valid hub decision, with no orphaned accounts; and (c)~verifying compliance early in the pipeline is a conservative strategy that never produces false positives, and interchanging the order of filters is mathematically consistent. The financial prerequisites are \cite{GrinoldKahn} or \cite{Meucci}; the categorical machinery is developed from scratch.

\clearpage

\section*{Non-Technical Summary}
\addcontentsline{toc}{section}{Non-Technical Summary}

\subsection*{The Problem}

A portfolio manager changes her mind about a single stock. That decision must now propagate---through an ETF, through hundreds of tax-managed separately managed accounts, through a target-date retirement suite---to a hundred thousand client portfolios, each with its own constraints. If the propagation is incoherent, some of those portfolios will violate their mandates. The firm may not know until the next audit.

This is the problem of \emph{mass customization}, and it breaks traditional workflows. Consider the lifecycle of a modern investment strategy:
\begin{itemize}
    \item A Chief Investment Officer (CIO) defines a ``Global Growth'' model portfolio.
    \item This model is replicated into an ETF for retail distribution.
    \item It is simultaneously adapted into hundreds of tax-managed SMAs, each with unique exclusion lists (e.g., ``no tobacco,'' ``no Magnificent 7'').
    \item It is further transformed into a ``Glide Path'' for a target-date retirement suite, where the risk profile shifts automatically based on the investor's age.
\end{itemize}

These steps were traditionally handled by disparate teams using different software systems, linked together by manual processes, spreadsheets, and ad-hoc scripts. The fragmentation is not merely inconvenient---it is a source of systemic operational risk. A change in the CIO's model might not propagate correctly to the tax-managed accounts. A compliance rule checked at the beginning of the process might be violated by a re-optimization step at the end. The problem is not that any single step is wrong; it is that no one can prove the chain is right.

The central question of this paper is: \emph{How do we mathematically guarantee that the final client portfolio effectively reflects the original investment intent, regardless of how many transformations it has passed through?}

The analogy to double-entry bookkeeping is not casual. Double-entry bookkeeping is itself a compositional system: every transaction is simultaneously a debit and a credit, and the two sides must balance at every node in the ledger. What bookkeeping does for individual transactions, category theory does for the relationships between processes---it provides a logic for verifying that a chain of transformations preserves what it should preserve. The relevant question is not what a single portfolio looks like (traditional portfolio theory handles that), but whether optimizing for taxes and then filtering for ESG produces the same result as filtering first and then optimizing. This is a question about composition, and composition is what category theory was designed to handle.

  \subsection*{Alignments and Re-implementations: Concrete Examples}

  The abstract dictionary above needs populating. The two primitive notions---\emph{alignment} (a relationship
  between portfolios) and \emph{re-implementation} (a process that transforms one portfolio into another)---each cover a broad range of practical situations.

  \medskip
  \noindent\textbf{Alignments} (vertical arrows) express the many ways in which two portfolios $P$ and $Q$ can be
  ``close'' or ``compatible.'' Each of the following is a distinct alignment relation, and a firm may impose several
  simultaneously:

  \begin{itemize}\itemsep2pt
  \item \emph{Active risk.} The tracking error between $P$ and $Q$ is below a given threshold---the canonical compliance
   check in delegated portfolio management.
  \item \emph{Asset class exposures.} The equity, fixed income, and alternatives weights of $P$ and $Q$ differ by no
  more than prescribed tolerances (e.g., ``equity allocation within 5\% of the model'').
  \item \emph{Factor exposures.} The loadings of $P$ and $Q$ on a specified factor model (value, momentum, quality,
  size, duration, \ldots) differ by no more than given thresholds---a finer-grained version of the asset class check.
  \item \emph{Strategy identity.} $P$ and $Q$ employ the same underlying strategies, possibly with different weights or
  vehicle choices. This is a combinatorial rather than a metric condition.
  \item \emph{Scenario dominance.} Under a given set of stress scenarios, $Q$'s returns exceed $P$'s with a specified
  level of confidence---an alignment that is forward-looking and probabilistic rather than static.
  \end{itemize}

  \noindent The critical observation is that each of these defines a \emph{relation}, not a function: for a given $P$,
  many portfolios $Q$ may satisfy the condition, and vice versa. All are modeled as closed
  subsets of the product of two portfolio spaces.

  \medskip
  \noindent\textbf{Re-implementations} (horizontal arrows) are deterministic processes that convert a portfolio in one
  universe into a portfolio in another:

  \begin{itemize}\itemsep2pt
  \item \emph{Multi-objective optimization.} Starting from a hub portfolio $P$, solve for the spoke portfolio that best
  balances tracking $P$'s exposures against a new objective (lower fees, higher yield, ESG score improvement), over a
  new security universe or with a new strategy overlay.
  \item \emph{Vehicle replacement (minimum tracking error).} Replicate $P$'s exposures using a different set of instruments---e.g. ETFs
   rather than mutual funds---minimizing tracking error against $P$.
  \item \emph{Vehicle replacement (maximum return).} The same vehicle substitution, but maximizing expected return
  subject to a tracking error budget. A change in objective leads to a different map.
  \item \emph{Mechanical substitution.} Replace foreign-listed securities with their ADR equivalents, or swap individual
   bonds for a bond ladder with matching duration and credit quality. These are simpler maps---often piecewise
  linear---but they are maps nonetheless, and they compose with the others.
  \end{itemize}

An extension of the theory allows for non-deterministic re-implementations, which may arise from the use of optimizers in practice.

  \medskip
  \noindent\textbf{From alignment to re-implementation.}
  Alignment and re-implementation are not independent notions. An alignment relation between two portfolio spaces
  defines, implicitly, a \emph{family} of candidate re-implementations: given a hub portfolio $P$, find the spoke
  portfolio $Q$ that is optimal---in some specified sense---among all $Q$ satisfying the alignment. The re-implementation
  is derived from the alignment by solving a constrained optimization problem. Much of the formal
  machinery developed in this paper---the interplay between vertical and horizontal morphisms, the Beck--Chevalley and
  Frobenius laws---exists precisely to ensure that this passage is coherent: that
  checking alignment before optimizing gives the same answer as optimizing and then checking.

\subsection*{The Core Framework: Hubs, Spokes, and Phantom Portfolios}

The paper formalizes the ``Hub-and-Spoke'' model familiar to operations managers. The Hub is the source of investment wisdom; the Spoke is the locus of execution. The mathematical challenge lies in the translation between them.

A key step is the rigorous definition of the ``permissible space.'' In standard quantitative finance, constraints are often specified loosely---e.g., ``target weight $> 0$.'' We show that such loose definitions are a root cause of system instability.

Consider the ``phantom portfolio'' problem. A Hub strategy requires strict diversification: every asset must have weight strictly greater than zero ($x > 0$). A Spoke optimizer, minimizing fees, drives the weight of an expensive asset to $0.0001\%$, then $10^{-9}\%$, then to exactly $0\%$.

The sequence of valid portfolios converges to a limit ($0\%$) that is \emph{outside} the Hub's permissible space. The result is a ``phantom''---a Spoke portfolio that looks compliant, has no valid parent, and was generated by the firm's own system. The phantom satisfies every local check; it fails only the global one. On an industrial SMA platform, these orphaned accounts accumulate quietly.

The fix is a topological requirement: all permissible spaces must be \textbf{closed sets} (containing their boundaries, i.e., $\ge 0$ rather than $> 0$). This ensures that the re-implementation map is \textbf{proper} (preimages of compact sets are compact)---no portfolio can escape the logic of the system, and every client account has a valid parent strategy.

\subsection*{The Main Results}

Using this topological foundation, the paper establishes four structural properties---each a theorem with a proof---that answer questions any portfolio construction team eventually confronts.

\emph{Can we verify a pipeline stage by stage?} Yes. Multi-stage portfolio construction preserves alignment: if each stage individually maintains compliance, so does the composed pipeline. No separate end-to-end verification is needed. (Compositionality, Theorem~\ref{thm:coherence}.)

\emph{Can we prevent errors before they happen?} Yes. The pushforward--pullback adjunction (Theorem~\ref{thm:adjunction}) proves an equivalence between checking a constraint \emph{after} construction and restricting the inputs \emph{before} construction. This defines ``Safe Hubs''---regions of the hub space where \emph{any} strategy will automatically produce compliant spokes. The shift is from reactive fixing of broken accounts to proactive prevention.

\begin{figure}[h]
\centering
\begin{tikzpicture}[
    node distance=2.0cm and 2.0cm,
    >=Latex,
    font=\small,
    block/.style={
        draw=gray!80, 
        fill=gray!5, 
        rectangle, 
        rounded corners=2pt,
        minimum width=3.5cm, 
        minimum height=1.2cm, 
        align=center,
        thick
    },
    safe_block/.style={
        draw=green!60!black, 
        fill=green!5, 
        rectangle, 
        rounded corners=2pt,
        minimum width=3.5cm, 
        minimum height=1.2cm, 
        align=center,
        thick
    },
    risky_block/.style={
        draw=orange!60!black, 
        fill=orange!5, 
        rectangle, 
        rounded corners=2pt,
        minimum width=3.5cm, 
        minimum height=1.2cm, 
        align=center,
        thick
    },
    arrow/.style={->, very thick, gray!70},
    label_text/.style={midway, fill=white, inner sep=2pt, font=\sffamily\footnotesize}
]

    \node[block] (hub) {
        \textbf{Hub} ($K_A$)\\
        Full Universe of Strategies
    };

    \node[safe_block, below left=1.5cm and 1.0cm of hub] (FilterEarly) {
        \textbf{Step 1: Early Audit}\\
        Filter by Guidance $R$\\
        (hub-side Check)
    };
    \node[safe_block, below=1.5cm of FilterEarly] (ResultSafe) {
        \textbf{Safe Spoke Set}\\
        $f'_!(g^*R)$\\
        (Guaranteed Compliant)
    };

    \node[block, below right=1.5cm and 1.0cm of hub] (Reimpl) {
        \textbf{Step 1: Re-implement}\\
        Map to Spoke Space\\
        (No constraints yet)
    };
    \node[risky_block, below=1.5cm of Reimpl] (ResultRisky) {
        \textbf{Apparent Spoke Set}\\
        $h^*(f_!R)$\\
        (Late Audit)
    };

    \draw[arrow] (hub) -| node[label_text, above=0.5cm] {Path A: Filter First} (FilterEarly);
    \draw[arrow] (FilterEarly) -- node[label_text] {Re-implement ($f'$)} (ResultSafe);

    \draw[arrow] (hub) -| node[label_text, above=0.5cm] {Path B: Map First} (Reimpl);
    \draw[arrow] (Reimpl) -- node[label_text] {Late Audit ($h^*$)} (ResultRisky);

    \node[below=1.5cm of hub, yshift=-5.5cm] (VennCenter) {};

    \begin{scope}[shift={(VennCenter)}]
        \filldraw[draw=orange!60!black, fill=orange!10, thick] (0,0) ellipse (4.5cm and 2.2cm);
        \node[orange!60!black, font=\bfseries] at (0, 1.6) {Apparent Set (Late Audit)};
        
        \filldraw[draw=green!60!black, fill=green!10, thick] (-1.2, -0.5) ellipse (2.5cm and 1.0cm);
        \node[green!40!black, font=\bfseries] at (-1.2, -0.5) {Safe Set (Early Audit)};
    \end{scope}

    \draw[->, dashed, green!60!black] (ResultSafe.south) to[out=270, in=120] ($(VennCenter)+(-2, -0.5)$);
    \draw[->, dashed, orange!60!black] (ResultRisky.south) to[out=270, in=60] ($(VennCenter)+(2, 1.0)$);

    \node[below=2.4cm of VennCenter, draw=gray, dashed, rounded corners, fill=white, inner sep=5pt] {
        \textbf{Lax Beck--Chevalley Guarantee:} $\quad f'_!(g^*R) \subseteq h^*(f_!R)$
    };

\end{tikzpicture}
\caption{The ``Audit Safety'' Flowchart. Visualizing the proposition that hub-side verification is conservative. The lax Beck--Chevalley inclusion means that the set of portfolios generated by Path A (Early Audit) is contained within the set generated by Path B (Late Audit).}
\label{fig:audit-safety}
\end{figure}
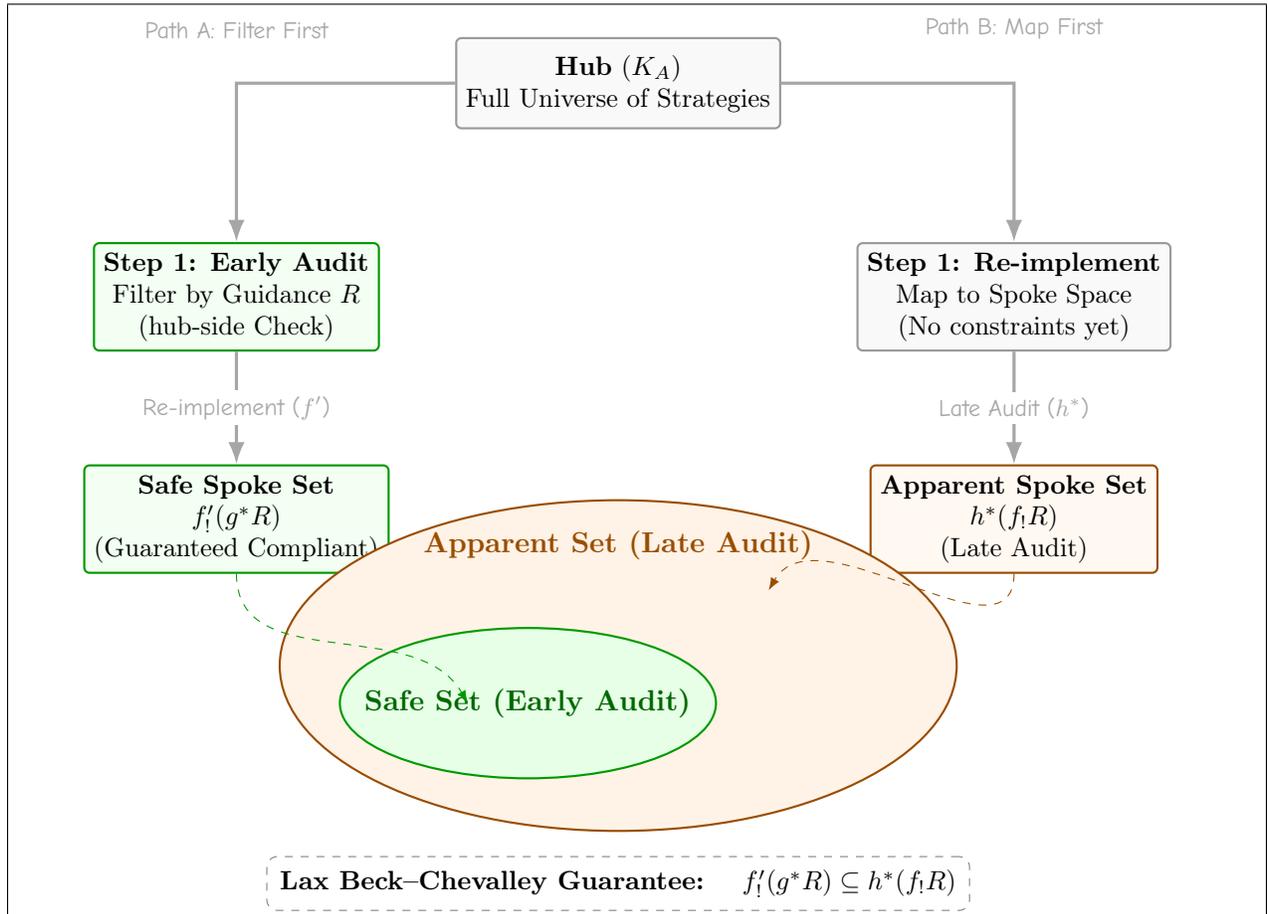

\emph{Can we centralize compliance?} Yes, conservatively. In large firms, a portfolio passes through multiple teams: Asset Allocation $\rightarrow$ Security Selection $\rightarrow$ Tax Overlay $\rightarrow$ Trading. The Beck--Chevalley condition (Theorem~\ref{thm:Beck--Chevalley}, Figure~\ref{fig:audit-safety}) guarantees that checking compliance early---at the hub level---is conservative: if a strategy passes upstream, its re-implementation will pass downstream. Compliance can be centralized rather than repeated at every desk.

\emph{Does the order of filtering matter?} No. Frobenius reciprocity (Theorem~\ref{thm:frobenius}) governs how filtering (e.g., ESG screens) interacts with optimization. Consider two simultaneous requirements: a hub constraint $R$ (``the mutual funds must align with a Growth asset allocation within 5\% tolerance'') and a spoke constraint $S$ (``the resulting ETF portfolio must track a passive basket within 50\,bps''). Frobenius proves that filtering the hub for Growth and then optimizing into ETFs yields exactly the same set of valid client portfolios as optimizing first and filtering after. This is a genuine degree of freedom: filter first for speed, or filter last for reporting transparency---the fiduciary outcome is identical.

\subsection*{Extension 1: Mass Customization and Personas}

The basic theory assumes a static relationship between Hub and Spoke. But strategies evolve---in the retirement market, the strategy \emph{must} evolve, dynamically, based on the client's age, risk tolerance, or tax bracket.

Section~\ref{sec:personas} introduces \textbf{Persona-Indexed Portfolios}: the ``Object'' is not a single portfolio space but a family of spaces indexed by a demographic variable.

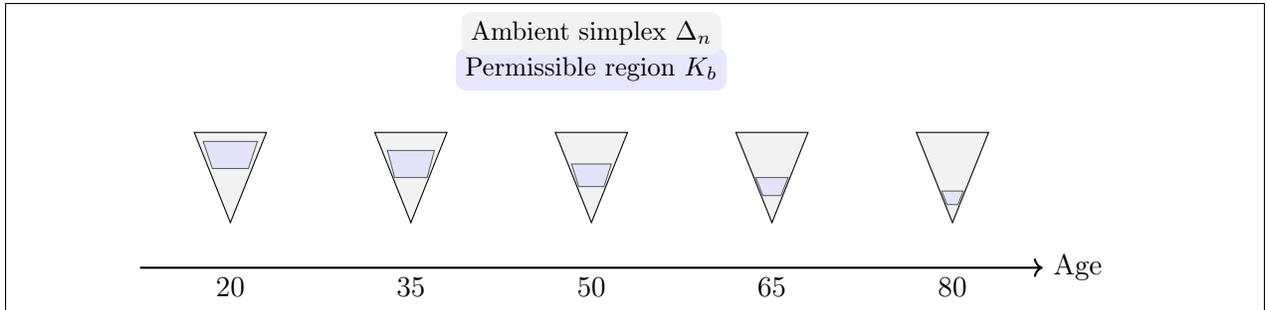
\begin{figure}[h]
    \centering
\begin{tikzpicture}[scale=1.2]
  \draw[->, thick] (0,0) -- (10,0) node[right] {Age};
  
  \foreach \x/\age/\bw/\tw/\by/\ty in {
    1/20/0.20/0.30/1.1/1.4,   
    3/35/0.18/0.26/1.0/1.3,   
    5/50/0.14/0.22/0.9/1.15,  
    7/65/0.10/0.18/0.8/1.0,   
    9/80/0.06/0.12/0.7/0.85   
  } {
    \draw[fill=gray!10] (\x,0.5) -- (\x-0.4,1.5) -- (\x+0.4,1.5) -- cycle;

    \draw[fill=blue!15, opacity=0.6] 
        (\x-\bw, \by) --   
        (\x-\tw, \ty) --   
        (\x+\tw, \ty) --   
        (\x+\bw, \by) --   
        cycle;

    \node[below] at (\x,0) {\age};
  }
  \node[fill=blue!10, rounded corners, font=\small] at (5,2.2) {Permissible region $K_b$};
  \node[fill=gray!10, rounded corners, font=\small] at (5,2.6) {Ambient simplex $\Delta_n$};
\end{tikzpicture}
\caption{The setting for target date glide path construction: permissible spaces indexed by age cohort of plan participant. The permissible regions $K_b$ shrink and shift downward (de-risking) as age increases, while remaining strictly inside the ambient simplex.}
    \label{fig:spoke-date}
\end{figure}

This allows us to model a Glide Path not as a discrete set of funds (2030, 2035, 2040), but as a continuous mathematical object (see Figure~\ref{fig:spoke-date}). We can apply ``Cross-Vintage Constraints,'' such as requiring that the equity allocation decreases smoothly as age increases. The framework ensures that a re-implementation decision made for the ``Age 30'' cohort propagates coherently to the ``Age 31'' cohort, preventing jagged or inconsistent transition paths that can arise when vintages are managed in silos.

\subsection*{Extension 2: Designing Strategy Menus (DOTS)}

Often, a CIO does not want to dictate a single portfolio but to define a ``Menu'' of permissible options---``You may hold any combination of these 5 approved managers, provided the total fee is at most 50 bps.''

Section~\ref{sec:DOTS} introduces the \textbf{DOTS (Double Operadic Theory of Systems)} model. An alignment relation acts as a generator: we take a Hub Space and ``act'' on it with a constraint to produce a ``Menu'' of permissible Spokes (see Figure~\ref{fig:operad-wiring}).
The central office defines the sandbox---the boundaries of risk and fee---while the specific selection is left to local portfolio managers or automated tax-loss harvesting algorithms.

\subsection*{Extension 3: Managing Imperfect Solvers (Probabilistic Risk)}

In practice, optimization is messy. Solvers have numerical noise, jitter, and local optima; market liquidity is finite; a theoretical trade might not be executable.

Section~\ref{sec:Polish} extends the framework to \textbf{HSP-r (Hub-and-Spoke Risk)}. Here, a re-implementation is not a single point, but a probability distribution---a ``cloud'' of likely outcomes.
We introduce two methods for managing this risk:

\textbf{The Safety Radius:} We calculate a ``buffer zone'' around the target portfolio (e.g., defined by the likely error of the solver). The system verifies that this entire buffer zone lies within the compliance limits. If the solver has a 1\% margin of error, we shrink our compliance boundaries by 1\% to ensure safety.

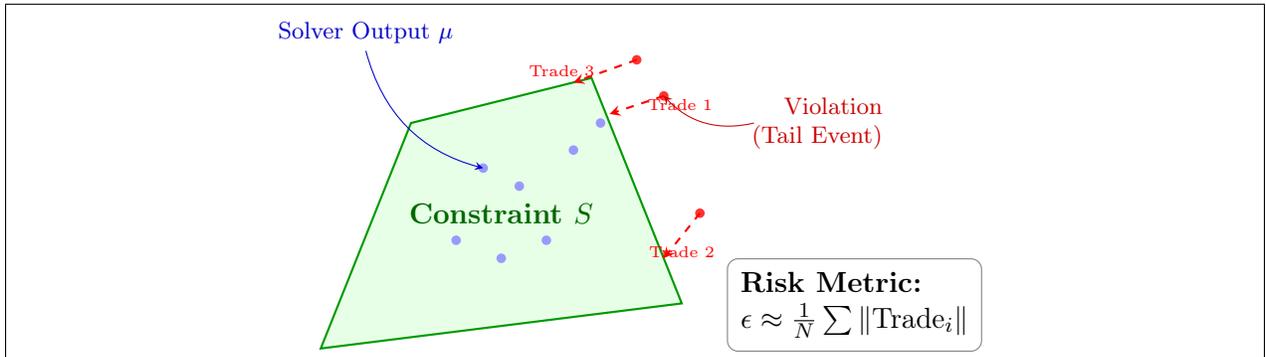
\begin{figure}[htbp]
\centering
\begin{tikzpicture}[scale=1.2, >=stealth]

    \fill[green!10] (0,0) -- (4,0.5) -- (3,3) -- (1,2.5) -- cycle;
    \draw[green!60!black, thick] (0,0) -- (4,0.5) -- (3,3) -- (1,2.5) -- cycle;
    \node[green!40!black, font=\bfseries] at (2, 1.5) {Constraint $S$};


    \foreach \x/\y in {1.5/1.2, 2.2/1.8, 2.5/1.2, 1.8/2.0, 2.8/2.2, 3.1/2.5, 2.0/1.0} {
        \fill[blue!40] (\x, \y) circle (1.5pt);
    }

    \coordinate (Out1) at (3.8, 2.8); 
    \coordinate (Out2) at (4.2, 1.5); 
    \coordinate (Out3) at (3.5, 3.2); 

    \fill[red!80] (Out1) circle (1.5pt);
    \fill[red!80] (Out2) circle (1.5pt);
    \fill[red!80] (Out3) circle (1.5pt);


    \draw[->, red, thick, dashed] (Out1) -- node[midway, right, font=\tiny] {Trade 1} (3.2, 2.6);

    \draw[->, red, thick, dashed] (Out2) -- node[midway, below, font=\tiny] {Trade 2} (3.8, 1.0); 

    \draw[->, red, thick, dashed] (Out3) -- node[midway, left, font=\tiny] {Trade 3} (2.8, 2.95);

    \node[align=left, font=\footnotesize, blue!80!black] at (0.5, 3.5) {Solver Output $\mu$};
    \draw[->, blue!80!black] (0.5, 3.3) to[bend right] (1.8, 2.0);

    \node[align=right, font=\footnotesize, red!80!black] at (5.5, 2.5) {Violation\\(Tail Event)};
    \draw[->, red!80!black] (4.8, 2.5) to[bend left] (Out1);

    \node[draw=gray, fill=white, rounded corners, inner sep=5pt, anchor=north west, align=left] at (4.5, 1.0) {
        \textbf{Risk Metric:}\\
        $\epsilon \approx \frac{1}{N} \sum \| \text{Trade}_i \|$
    };

\end{tikzpicture}
\caption{Visualizing Wasserstein Safety. The constraint $S$ is the green region. The solver produces a cloud of portfolios (dots). Some land outside (red). The Wasserstein risk $\epsilon$ is the expected cost (average length of the red dashed arrows) to ``cure'' the violations.}
\label{fig:wasserstein_cure}
\end{figure}

\textbf{Wasserstein Safety (The ``Cost of Cure"):} Sometimes, a portfolio is technically non-compliant by a fraction of a penny. Rejecting it causes unnecessary churn---the cure is worse than the disease. Using Optimal Transport theory (Section~\ref{sec:wasserstein}), we define alignment not as a binary ``True/False'' but as a cost: ``What is the expected trading cost to fix this portfolio?''
If the ``Cure Cost'' is less than a de minimis threshold (e.g., 1 basis point of turnover), the system accepts the portfolio (see Figure~\ref{fig:wasserstein_cure}). The rigid logic of set membership gives way to the economic logic of trading costs.

\subsection*{System Architecture and Implementation}

Appendix~\ref{app:architecture} translates this theory into a concrete blueprint for IT systems. In current practice, the business logic governing portfolio transformations typically resides in Python scripts or Excel macros, where it is opaque to compliance and difficult to update.

We propose a \textbf{Categorical Orchestrator} (Figure~\ref{fig:system_arch}) that sits at the center of the investment platform.
\begin{itemize}
    \item \emph{Semantic Core:} Encodes the categorical structure---objects, morphisms, composition laws---and enforces ordering constraints (e.g., tax optimization after ESG filtering, when those operations do not commute).
    \item \emph{Solver Interface:} Delegates numerical computation to specialized solvers (Gurobi, Axioma, etc.), making the system solver-agnostic: the optimizer can be swapped without rewriting the compliance logic.
    \item \emph{Evidence Ledger:} Stores not just the final portfolio but the morphisms---the specific audit trail (2-cell) linking Hub to Spoke. This provides full provenance: the system can reconstruct \emph{why} a specific trade was made for a specific client.
\end{itemize}

\subsection*{Conclusion}

The transition to industrial-scale portfolio customization requires a theory of \emph{relationships between portfolios}---not just individual portfolios. Existing portfolio theory was not designed for this; it optimizes one portfolio at a time and is silent on how transformations compose. This paper proposes a mathematical foundation. Whether the framework survives contact with production systems---with their messy solver tolerances, organizational politics, and data that never quite matches the model---remains to be seen. The theorems guarantee coherence under their hypotheses; the hypotheses themselves are idealizations. The first to fail in production will likely be determinism: real optimizers jitter, and the probabilistic extensions in Part III exist for this reason. The second will be the assumption that permissible spaces are known in advance; in practice they are negotiated, revised, and occasionally contradictory. But the idealizations are, we believe, the right ones to aim for---they identify what coherence \emph{means} before asking whether any particular system achieves it.

\clearpage

\section*{Acknowledgments}
\addcontentsline{toc}{section}{Acknowledgments}

The human author thanks Kevin Davenport, Ella Sinfield and Kelly Densmore for comments, Greg Zitelli and David Jaz Myers for the suggestions that led to Sections \ref{sec:Polish} and \ref{sec:DOTS} respectively, and Brendan Fong for support and encouragement. An earlier version of this work was presented informally at Topos Institute, and improvements were made based on suggestions from the audience, particularly the axiomatization (Section~\ref{sec:abstract}). This work was also presented at Capital Group's Theory Seminar, and at the Capital Solutions Group Long-Form Research Call. 

\subsection*{Methods: Use of Large Language Models}

I used large language models---Claude (Sonnet 4.5, Opus 4.5, Opus 4.6), Gemini (Pro 2.5, 3), and GPT-5---as technical assistants: for \LaTeX{} formatting, literature identification, routine verification of standard derivations, and editorial cleanup. None of the conceptual framework, definitions, or proofs originated with a model. The models were treated as fallible drafting tools, and they were. I verified all mathematical proofs, diagrams, and financial logic manually, and bear full responsibility for the final accuracy and originality of this work.

\subsection*{Disclaimers}

The views expressed here are those of the human author and do not necessarily reflect the views or policies of Capital Group. 
This work is for informational purposes only. The examples are generic, and selected for illustrative purposes; they do not necessarily reflect the actual investment process or other business processes at Capital Group. The conceptual system architecture and workflows in Appendix~\ref{app:architecture} are included for expositional purposes only.

The human author and Capital Group make no representations or warranties, express or implied, as to the accuracy or completeness of the information contained herein. The author and Capital Group assume no liability for any errors, omissions, or losses arising from the use of this information.

\clearpage  \section{Introduction}
\label{sec:introduction}

\subsection{Portfolio Spaces and Re-implementations}

Given a finite set of assets $\A = \{a_0,
a_1, \ldots, a_n\}$, a \textbf{long-only portfolio} assigns a non-negative
weight $x_i \geq 0$ to each asset $a_i$ such that the weights sum to one:
$\sum_{i=0}^n x_i = 1$.

The space of all such portfolios forms the standard $n$-simplex:
\[
\Delta^n = \left\{(x_0, \ldots, x_n) \in \mathbb{R}^{n+1} \,\Big|\, x_i \geq 0,
\sum_{i=0}^n x_i = 1\right\}
\]

\begin{remark}[Relaxing the long-only restriction]
If short positions and leverage are permitted, the space of all portfolios is $\mathbb{R}^n$. This setup is too simplistic to address the additional issues related to funding short or leveraged positions; we do not discuss it further.
\end{remark}

In practice, investors are rarely allowed to hold
\textit{any} portfolio in $\Delta^n$. Portfolios are almost always subject to
baseline constraints: position limits, sector bounds, regulatory
requirements, ESG screens.

This motivates our core definition: a \textbf{permissible portfolio space} is a
\textbf{closed subset $K \subseteq \Delta^n$}. The geometric object $K$, not
the full simplex $\Delta^n$, is the true space of allowable portfolios.

\begin{remark}[Asset universes]
Different portfolio spaces may have different
sets of assets available---mutual funds versus ETFs, active versus passive
funds---but represent the same asset classes. This motivates mapping portfolios to a common
``space of quantitative attributes'' to compare their investment characteristics (Section
\ref{sec:optimal-reimpl}).
\end{remark}

\textbf{Portfolio re-implementations} arise when we transform portfolios from
one permissible space to another. A re-implementation $f: K_1 \to K_2$ maps a
``hub'' portfolio (e.g., an active strategy in $K_1 \subseteq \Delta^n$) to a
``spoke'' portfolio (e.g., a passive ETF implementation in $K_2 \subseteq
\Delta^m$). Examples include:
\par
\begin{itemize}
\item \textbf{Portfolio Replication}: Replicating a strategy from $K_1$ in a
new vehicle $K_2$.
\item \textbf{Asset Replacement}: Replacing a portfolio of active mutual funds
($K_1$) with a customized, low-fee basket of ETFs ($K_2$).
\item \textbf{Restriction}: Re-implementing a general portfolio from $K_1$ using
only a restricted set of securities in $K_2$ (e.g., an ESG-screened universe).
\item \textbf{Aggregation}: Converting portfolios over individual stocks
($K_1$) into sector allocations ($K_2$).
\item \textbf{Factor Decomposition}: Transforming security weights ($K_1$) into
factor risk exposures ($K_2$).
\end{itemize}

Mathematically, a portfolio re-implementation from permissible space $K_1
\subseteq \Delta^n$ to $K_2 \subseteq \Delta^m$ is a \textbf{continuous map} $f:
K_1 \to K_2$. We allow $f$ to be nonlinear, as the most practical
re-implementations (e.g., finding a tax-optimal or low-fee replacement) are
constructed via constrained optimization (Section~\ref{sec:optimal-reimpl}) and
are generally not simple affine maps.

\subsection{The Compositional Challenge}

In practice, portfolio management processes often involve chains of re-implementations:
\begin{center}
\begin{tikzcd}[column sep=large]
K_1 (\text{Funds}) \arrow[r, "f_1"] & K_2 (\text{ETFs}) \arrow[r,
"f_2"] & K_3 (\text{Factors})
\end{tikzcd}
\end{center}
At each stage, we impose alignments, which we model as \textbf{closed
relations} $R \subseteq K_i \times K_j$.
\par
\begin{itemize}
\item Tracking error alignment at the ETF level (must track the mutual funds)
\item Factor risk budgets at the factor level
\item Alignment defined against a more abstract portfolio space (e.g., a space of asset classes)
\end{itemize}

\textbf{Critical questions}:
\begin{enumerate}
\item How do alignments propagate through re-implementations?
\item Can we verify alignment at intermediate stages, or must we trace back to
original securities?
\item Does the order of verification matter?
\item When can we decompose a complex multi-stage strategy into simpler pieces?
\end{enumerate}

These questions become urgent at identifiable inflection points: (a)~\emph{volume}---the number of hub-to-spoke relationships outgrows ad~hoc management, with dozens of spokes tracked in spreadsheets; (b)~\emph{complexity}---overlapping regulatory, ESG, and tax screens make manual verification of constraint interactions unreliable; (c)~\emph{preparation}---anticipated structural changes (regulatory reform, benchmark transitions, new asset-class mandates) require confidence that modifications will propagate correctly before they are made.

Classical portfolio theory \cite{Markowitz, Sharpe, GrinoldKahn, Meucci} was not designed to answer these questions---it was designed to answer a different and equally important question, namely how to construct a single portfolio optimally. As Markowitz put it, ``the process of selecting a portfolio may be divided into two stages''---forming beliefs and choosing a portfolio \cite{Markowitz}. We are concerned with a third stage that Markowitz had no reason to consider: the coherent propagation of that choice through a chain of transformations. The gap is not a deficiency of classical theory but a consequence of its scope.

More precisely, we need a \emph{compositional} theory where:
\par
\begin{itemize}
\item Portfolio re-implementations compose correctly
\item Alignment relations compose transitively
\item The two compositions interact coherently
\end{itemize}
This is precisely the structure provided by a \textbf{double category} \cite{Benabou, GrandisPare, Leinster}.

\subsection{The Double Categorical Framework}

A \textbf{double category} $\mathbb{D}$ provides a framework with \textbf{objects}, two distinct types of morphisms---\textbf{horizontal morphisms} (e.g., $f: A \to B$) and \textbf{vertical morphisms} (e.g., $R: A \dashrightarrow C$)---and \textbf{2-cells} $\alpha$ that fill squares formed by these morphisms, like:
\[
    \begin{tikzcd}
    A \arrow[r, "f"] \arrow[d, "R"', dashed] & B \arrow[d, "S", dashed] \\
    C \arrow[r, "g"'] & D
    \arrow[from=1-1, to=2-2, phantom, "\Downarrow \alpha"]
    \end{tikzcd}
\]

Horizontal morphisms compose horizontally ($\circ_h$), vertical morphisms compose vertically ($\circ_v$), and 2-cells compose both horizontally and vertically. These compositions must satisfy associativity, identity laws, and the \textbf{interchange law}, which ensures that the order of composing 2-cells horizontally and vertically does not matter.

We construct a double category $\HSSimp$ (Hub-and-Spoke):

\begin{center}
\begin{tikzcd}[column sep=large, row sep=large]
K_1 \arrow[r, "f (\text{re-impl.})"] \arrow[d, dashed, "R (\text{align})"'] & K_2 \arrow[d, dashed, "S (\text{align})"] \\ 
K_3 \arrow[r, "g (\text{re-impl.})"'] & K_4
\end{tikzcd}
\end{center}
\begin{description}
\item[Objects:] permissible portfolio spaces $K$ (closed subsets of simplices
$\Delta^n$)

\item[Horizontal morphisms:] Continuous portfolio re-implementations $f: K_1
\to K_2$

\item[Vertical morphisms:] Closed alignment relations $R \subseteq K_1 \times
K_3$

\item[2-cells:] Inclusions ensuring re-implementations preserve alignment
\end{description}

A double category in which every square (2-cell) is determined by an inclusion or equality, with no nontrivial 2-cell structure, is called \emph{thin}. In a thin double category, the interchange law holds automatically. For example, the double category $\HSSimp$ is thin. (The thinness assumption is dropped in Section~\ref{sec:2-cells}, which describes a richer framework.)

\subsection{Why Closed Objects}

The requirement that permissible spaces be \emph{closed} subsets of the simplex is not a technical nicety. If objects are open (e.g., requiring strictly positive weights, $x_i > 0$), continuous maps from them need not be proper, and the pushforward $f_!R$ of a closed alignment can fail to be closed---producing ``phantom portfolios'' that appear compliant but have no valid hub pre-image. Section~\ref{sec:proper} proves that closedness (hence compactness) of objects makes all horizontal morphisms automatically proper, and Section~\ref{sec:failure} demonstrates the precise failures that arise without it.

\subsection{Main Results}

The framework yields the following theorems. The definitions and proofs are straightforward---no novel mathematics is required. The novelty lies entirely in the model itself: in the claim that these particular objects and morphisms are the right ones for formalizing portfolio re-implementation.

\begin{enumerate}[label=(\Roman*)]
\item \textbf{Closedness Preservation (Theorem~\ref{thm:proper-closed})}:
Portfolio re-implementations are automatically proper, so they
preserve closedness of alignment relations under forward
propagation (pushforward).

\item \textbf{Compositional Coherence (Theorem~\ref{thm:coherence})}: The
double category $\HSSimp$ (of permissible spaces and proper re-implementations)
satisfies all double category axioms with closed composition operations.

\item \textbf{Adjunction (Theorem~\ref{thm:adjunction})}: The pushforward
$(f_!)$ and pullback $(f^*)$ operations form an adjoint pair:
\[
R \subseteq f^*S \iff f_!R \subseteq S
\]
This states a formal relationship between hub-space and spoke-space
formulations of alignments.

\item \textbf{Beck--Chevalley (Theorem~\ref{thm:Beck--Chevalley})}: In Section~\ref{sec:Beck--Chevalley} we prove a Beck--Chevalley law for
\emph{existing} commuting squares of horizontal morphisms that are
\emph{pointwise cartesian} (Def.~\ref{def:pointwise-cartesian})---a rather restrictive condition. Fortunately, the ``easy'' half of the law holds for general commuting squares, and it is this half that does the work in practice: it ensures that upstream compliance checks are conservative.

\item \textbf{Frobenius Reciprocity (Theorem~\ref{thm:frobenius})}: Intersection of alignment relations commutes with re-implementation:
\[
f_!(R \cap f^*S) = f_!R \cap S
\]
Portfolios can be filtered before or after re-implementation equivalently.

\item \textbf{Construction of Optimal Re-implementations (Theorem
\ref{thm:berge-optimal})}: We \emph{construct}
re-implementations as solutions to optimization problems (e.g., minimize
fees subject to tracking error). Using Berge's Maximum Theorem, we prove this
optimization constructs a well-defined, continuous map $f$ whose domain is a
\textbf{closed set $K$}---precisely the object of our double category.

\item \textbf{Propagation of Portfolio Changes (Section~\ref{sec:propagation})}: Portfolio changes can be propagated coherently. Properness ensures pushforwards of closed alignments remain closed. The adjunction, Beck--Chevalley, and Frobenius laws then ensure that filtering and re-implementation commute.

\item \textbf{Failure Without Closed Objects (Theorem
\ref{thm:closure-fails})}: Attempting to use non-closed permissible spaces (e.g.,
open interiors) breaks properness. The ``fix'' of taking topological closure
of pushforwards breaks adjunction, Beck--Chevalley, and Frobenius
properties, introducing ``phantom'' portfolios. This justifies our restrictive
definition.
\end{enumerate}

Taken together, these results translate into concrete operational properties:
\par
\begin{itemize}
\item \textbf{Path independence}: Multi-stage strategies can be verified in \textit{any}
order
\item \textbf{Decomposition}: Complex optimizations can be broken into simpler
stages
\item \textbf{Boundary stability}: Limit portfolios behave predictably
\item \textbf{Alignment commutation}: Filtering and re-implementation can be
reordered
\item \textbf{Compositional reasoning}: Large strategies can be built from
verified components
\item \textbf{Modeling nonlinear replication or customization}: The framework explicitly
supports optimization-derived re-implementations.
\end{itemize}

\subsection*{Organization}
\enlargethispage{2\baselineskip}

\textbf{Part I: The Basic Theory} develops the central $\HSSimp$ framework. Section~\ref{sec:spaces} defines the objects of our category as permissible regions within ambient simplices. Section~\ref{sec:reimplement} introduces horizontal morphisms (re-implementations), and Section~\ref{sec:constraints} introduces vertical morphisms (alignment relations). Section~\ref{sec:double} assembles these components into the $\HSSimp$ double category. Section~\ref{sec:operations} defines the pullback and pushforward operations, while Section~\ref{sec:proper} establishes the critical result that our re-implementations are automatically proper maps. Sections \ref{sec:Beck--Chevalley} and \ref{sec:frobenius} prove the main coherence theorems regarding path independence and filtering. Section~\ref{sec:failure} analyses the failure modes of non-closed models. Section~\ref{sec:optimal-reimpl} provides a constructive method for deriving proper re-implementations via optimization. Finally, Sections \ref{sec:propagation} and \ref{sec:computation} discuss practical portfolio management applications and computational implementation.

\smallskip

\noindent \textbf{Part II: Axiomatization and Extensions} presents the broader theoretical context and complex use-cases. We begin with an axiomatic framework in Section~\ref{sec:abstract} that serves as the abstract ``syntax'' of the theory. We then extend the framework to include ``evidence of alignment'' via spans in Section~\ref{sec:2-cells}. Subsequent sections develop two alternative semantic models motivated by specific industrial problems: persona-indexed portfolios for glide paths (Section~\ref{sec:personas}) and the DOTS framework for handling menus of permissible spokes (Section~\ref{sec:DOTS}).

\smallskip

\noindent \textbf{Part III: Probabilistic Theory} develops the stochastic extension of the framework.  Section~\ref{sec:Polish} introduces the probabilistic hub-and-spoke category $\mathrm{HSP\text{-}r}$, modeling re-implementations as tight Feller kernels on Polish spaces.  Sections~\ref{sec:safety-radius} and~\ref{sec:HDR-section} develop two approaches to probabilistic compliance---safety radius and highest density regions---while Section~\ref{sec:wasserstein} introduces transport-based safety via Wasserstein distances, together with the associated composition, adjunction, and coherence theorems.  Subsequent sections address the estimation and storage of stochastic kernels in production settings and the interaction of continuous probabilistic models with discrete portfolio constraints (Section~\ref{sec:discrete}).  Finally, Section~\ref{sec:liquidity} applies all three layers of the framework---$\HSSimp$, DOTS, and Wasserstein---to the problem of liquidity-aware portfolio construction with illiquid and semi-liquid assets.

\clearpage

\section*{Logical Structure}
\addcontentsline{toc}{section}{Logical Structure}

\vspace{0.5cm}

\begin{center}
\resizebox{!}{0.85\textheight}{%
\begin{tikzpicture}[
    node distance=0.6cm and 0.4cm, 
    >=stealth,
    core/.style={
        rectangle, 
        draw=blue!60!black, 
        top color=blue!5, 
        bottom color=blue!10, 
        rounded corners=3pt, 
        minimum height=1.0cm, 
        minimum width=3.5cm, 
        text width=3.3cm,
        align=center, 
        drop shadow,
        inner sep=3pt
    },
    input/.style={
        rectangle, 
        draw=orange!60!black, 
        top color=orange!5, 
        bottom color=orange!10, 
        rounded corners=3pt, 
        minimum height=1.0cm, 
        minimum width=3.2cm, 
        text width=3.0cm, 
        align=center, 
        drop shadow
    },
    extension/.style={
        rectangle, 
        draw=green!50!black, 
        top color=green!5, 
        bottom color=green!10, 
        dashed,
        rounded corners=3pt, 
        minimum height=0.9cm, 
        minimum width=2.6cm, 
        text width=2.4cm,
        align=center, 
        font=\sffamily\scriptsize
    },
    flow/.style={->, very thick, draw=gray!60},
    extend/.style={->, dashed, thick, draw=gray!50},
    labeltext/.style={font=\bfseries\sffamily\tiny, text=gray!60, fill=white, inner sep=1pt}
]

    \node[input] (sec11) {
        \textbf{I.Sec 11: Optimal Re-impl.}\\
        Berge's Max Thm\\
        (Minimizing Distance)
    };

    \node[core, below=0.6cm of sec11] (sec3) {
        \textbf{I.Sec 3: Horizontal}\\
        Continuous Maps\\
        $f: K_1 \to K_2$
    };
    \node[core, left=0.3cm of sec3] (sec2) {
        \textbf{I.Sec 2: Objects}\\
        Permissible Spaces\\
        (Closed $K \subseteq \Delta^n$)
    };
    \node[core, right=0.3cm of sec3] (sec4) {
        \textbf{I.Sec 4: Vertical}\\
        Alignments\\
        (Closed $R \subseteq K \times L$)
    };

    \node[core, below=0.8cm of sec3, minimum width=7cm, text width=6.8cm] (sec5) {
        \textbf{I.Sec 5: The $\mathbb{HS}$ Double Category}\\
        Thin Double Category where\\
        2-cells are inclusions.
    };

    \node[core, below=0.8cm of sec5] (sec7) {
        \textbf{I.Sec 7: Properness}\\
        Compact Domain $\to$ Proper.\\
        Guarantees closed images.
    };

    \node[core, below=0.6cm of sec7] (sec6) {
        \textbf{I.Sec 6: Operations}\\
        Pushforward ($f_!R$) \& Pullback ($f^*S$).\\
        \textit{Key:} $f_!$ preserves closedness.
    };

    \node[core, below=0.6cm of sec6] (sec89) {
        \textbf{I.Sec 8 \& 9: Coherence}\\
        \textbf{Adjunction:} $f_! \dashv f^*$\\
        \textbf{Beck--Chevalley:} Path Indep.\\
        \textbf{Frobenius:} Filter Commutes
    };

    \node[core, below=0.8cm of sec89, fill=blue!10] (sec12) {
        \textbf{I.Sec 12: Propagation}\\
        Audit Trails, Stability,\\
        Multi-stage Checks.
    };


    \node[extension, right=0.3cm of sec11] (sec18) {
        \textbf{II.Sec 17: DOTS}\\
        Menu Actions\\
        ($K \odot R$)
    };

    \node[extension, right=0.3cm of sec4] (sec15) {
        \textbf{II.Sec 15: Evidence}\\
        Spans of Evidence
    };

    \node[extension, left=0.3cm of sec5, yshift=0.3cm] (sec14) {
        \textbf{II.Sec 14: Axioms}\\
        Abstract Syntax
    };

    \node[extension, right=0.3cm of sec5, yshift=0.3cm] (sec16) {
        \textbf{II.Sec 16: Personas}\\
        Slice Cat $\mathbb{HS}/B$
    };

    \node[extension, right=0.5cm of sec7] (sec17) {
        \textbf{III.Sec 18 \& 19: Polish}\\
        Feller Kernels\\
        (Alt. Topology)
    };

    \node[extension, below=0.3cm of sec17, fill=green!15] (sec17_OT) {
        \textbf{III.Sec 21: Transport}\\
        Wasserstein $W_1$\\
        (Cure Costs)
    };


    \draw[flow] (sec11) -- (sec3);

    \draw[flow] (sec2.south) -- (sec5.150);
    \draw[flow] (sec3) -- (sec5);
    \draw[flow] (sec4.south) -- (sec5.30);
    
    \draw[flow] (sec5) -- (sec7);
    \draw[flow] (sec7) -- node[right, labeltext] {Enables} (sec6);
    \draw[flow] (sec6) -- (sec89);
    \draw[flow] (sec89) -- (sec12);

    \draw[flow, dashed, gray!40] (sec11.west) -- ++(-4.0,0) |- (sec12.west);

    \draw[extend] (sec18) -- (sec3); 
    \draw[extend] (sec15) -- (sec4); 
    \draw[extend] (sec14) -- (sec5); 
    \draw[extend] (sec16) -- (sec5); 
    \draw[extend] (sec17) -- (sec7); 
    \draw[extend] (sec17_OT) -- (sec17); 

\end{tikzpicture}
}
\end{center}

\section*{Notation and Terminology}
\addcontentsline{toc}{section}{Notation and Terminology}

\begin{center}
\renewcommand{\arraystretch}{1.4}
\begin{tabular}{p{0.15\textwidth} p{0.35\textwidth} p{0.40\textwidth}}
\toprule
\textbf{Symbol} & \textbf{Mathematical Definition} & \textbf{Financial Interpretation} \\
\midrule
\multicolumn{3}{l}{\textit{\textbf{Core Objects and Morphisms}}} \\
$\Delta^n$ & The standard $n$-simplex in $\mathbb{R}^{n+1}$. & The space of all theoretically possible long-only portfolios using $n+1$ assets. \\
$K$ & A closed subset $K \subseteq \Delta^n$. & A \textbf{permissible portfolio space}. The set of portfolios satisfying baseline constraints (e.g., position limits, ESG screens). \\
$f: K_1 \to K_2$ & A continuous map between permissible spaces. & \textbf{Re-implementation}. A deterministic rule mapping every valid hub portfolio to a unique spoke portfolio (e.g., optimizing a model portfolio of funds into an ETF portfolio). \\
$R \subseteq K_1 \times K_2$ & A closed relation. & \textbf{Alignment}. A constraint linking two portfolios. \\
$\mathbb{HS}$ & The Hub-and-Spoke double category. & The compositional framework organizing portfolios, re-implementations, and alignments. \\
\midrule
\multicolumn{3}{l}{\textit{\textbf{Operations}}} \\
$f_! R$ & Pushforward of $R$ along $f$. & \textbf{Forward propagation}. The set of valid spoke portfolios derived from valid hub portfolios. \\
$f^* S$ & Pullback of $S$ along $f$. & \textbf{Pre-alignment}. Identifying which hub portfolios will result in a compliant spoke portfolio. \\
$S \circ R$ & Relational composition. & \textbf{Transitive alignment}. Checking alignment across multiple steps (e.g., Fund $\to$ Benchmark $\to$ Factor Model). \\
$Graph(f)$ & The graph $\{(x, f(x))\}$ of the map $f$. & Viewing a deterministic re-implementation as a specific type of alignment relation. \\
\bottomrule
\end{tabular}
\end{center}

\clearpage

\begin{center}
\renewcommand{\arraystretch}{1.4}
\begin{tabular}{p{0.15\textwidth} p{0.35\textwidth} p{0.40\textwidth}}
\toprule
\textbf{Symbol} & \textbf{Mathematical Definition} & \textbf{Financial Interpretation} \\
\midrule
\multicolumn{3}{l}{\textit{\textbf{Extensions (Parts II and III)}}} \\
$p: K \to B$ & Display map over base space $B$. & \textbf{Persona indexing}. Portfolios constrained by a variable parameter (e.g., investor age in a glide path). \\
$(X, \mathcal{B}, \mu)$ & Polish measure space. & A portfolio space equipped with a notion of volume/probability. \\
$P: X \rightsquigarrow Y$ & Tight Feller kernel. & \textbf{Stochastic Re-implementation}. A rule mapping a hub portfolio to a probability distribution of spoke portfolios (e.g., noisy optimization). \\
$P^* S$ & Probabilistic Pullback. & The set of hubs that map to compliant spokes with probability 1. \\
$K \odot R$ & Action of relation $R$ on set $K$. & \textbf{Menu generation}. The set of all permissible spokes available for a given hub under constraint $R$ (DOTS framework). \\
\bottomrule
\end{tabular}
\end{center}

\clearpage

\part{The Basic Theory}

\clearpage

\definecolor{obj_color}{RGB}{50, 50, 50}      
\definecolor{do_color}{RGB}{0, 100, 180}      
\definecolor{check_color}{RGB}{180, 50, 50}   

\begin{figure}[ht]
    \centering
    \begin{tikzpicture}[
        node distance=4.0cm and 3.5cm,
        annotation/.style={
            rectangle,
            draw=gray!40,
            fill=gray!5,
            rounded corners,
            font=\footnotesize, 
            align=left,
            inner sep=5pt,
            text width=4.5cm
        },
        leader/.style={
            draw=gray!40,
            thick,
            dashed
        },
        obj/.style={
            text=obj_color,
            font=\large
        }
    ]

    \node[obj] (Hub)                     {$K_{\text{Hub}}$};
    \node[obj] (Spoke)  [right=of Hub]   {$K_{\text{Spoke}}$};
    \node[obj] (Bench)  [below=of Hub]   {$K_{\text{Bench}}$};
    \node[obj] (spoke) [right=of Bench] {$K_{\text{Spoke}}$};

    \draw[->, thick, color=do_color] (Hub) -- node[above] {$f$} node[below] {\scriptsize Re-implement} (Spoke);
    \draw[->, thick, color=do_color] (Bench) -- node[above] {$g$} node[below] {\scriptsize Transform} (spoke);

    \draw[->, thick, dashed, color=check_color] (Hub) -- node[left] {$R$} node[right] {\scriptsize Align} (Bench);
    \draw[->, thick, dashed, color=check_color] (Spoke) -- node[left] {$S$} node[right] {\scriptsize Align} (spoke);

    \node (TwoCell) at ($(Hub)!0.5!(spoke)$) {\textcolor{black!70}{$\Downarrow \alpha$}};


    \node[annotation, anchor=east] (NoteObj) at ($(Hub) + (-2.0, 0)$) {
        \textbf{\textcolor{obj_color}{Objects: Permissible Spaces}}\\
        Closed subsets $K \subseteq \Delta^n$.\\
        \textit{Intuition:} The sandbox of valid portfolios (e.g., "Mutual Funds").
    };
    \draw[leader] (NoteObj.east) -- (Hub.west);

    \node[annotation, anchor=south] (NoteHoriz) at ($(Hub)!0.5!(Spoke) + (0, 1.0)$) {
        \textbf{\textcolor{do_color}{Horizontal: Re-implementations}}\\
        Continuous Maps $f: K_1 \to K_2$.\\
        \textit{Intuition:} The \textbf{"DO"} operation. Optimizing or re-shuffling assets.
    };
    \draw[leader] (NoteHoriz.south) -- ($(Hub)!0.5!(Spoke) + (0, 0.2)$);

    \node[annotation, anchor=east] (NoteVert) at ($(Bench) + (-2.0, 0)$) {
        \textbf{\textcolor{check_color}{Vertical: Alignments}}\\
        Closed Relations $R \subseteq K_1 \times K_3$.\\
        \textit{Intuition:} The \textbf{"CHECK"} operation. Constraints like Tracking Error limits.
    };
    \draw[leader] (NoteVert.east) -- ($(Hub)!0.5!(Bench)$);

    \node[annotation, anchor=north] (NoteCell) at ($(Bench)!0.5!(spoke) + (1.5, -1.0)$) {
        \textbf{The 2-Cell: Consistency}\\
        $\Graph(g) \circ R \subseteq S \circ \Graph(f)$.\\
        \textit{Intuition:} \textbf{Audit Safety.} The operations commute/agree.
    };
    \draw[leader] (NoteCell.north) -- (TwoCell.south);

    \end{tikzpicture}
    \caption{Visualizing a 2-Cell in the $\mathbb{HS}$ Double Category.}
    \label{fig:2-cell}
\end{figure}
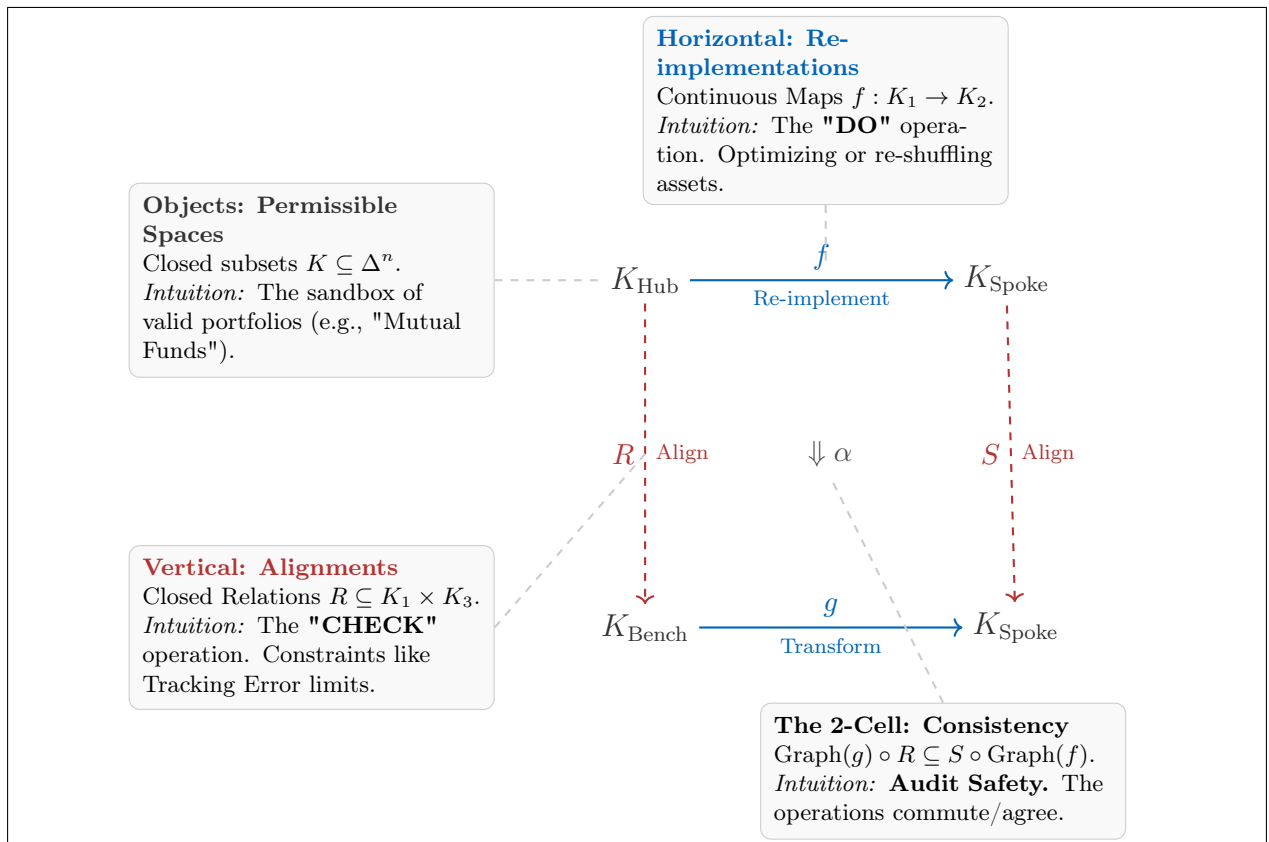

\clearpage

\section{Portfolio Spaces: Ambient and Permissible}
\label{sec:spaces}

\subsection{Ambient Portfolio Spaces (Simplices)}

\begin{definition}[Asset Set]
An \textbf{asset set} is a finite set $\A = \{a_0, a_1, \ldots, a_n\}$ of
distinct investment opportunities.
\end{definition}

For example, the elements of an asset set might be individual securities, mutual funds or ETFs, or more abstract entities such as asset classes or factor exposures.

\begin{definition}[Ambient Portfolio Space]
\label{def:ambient-portfolio-space}
The \textbf{ambient portfolio space} on asset set $\A$ with $|\A| = n+1$ is the
standard $n$-simplex:
\[
\Delta^n(\A) = \left\{(x_0, \ldots, x_n) \in \mathbb{R}^{n+1} \,\Big|\, x_i
\geq 0, \sum_{i=0}^n x_i = 1\right\}
\]
equipped with the subspace topology from $\mathbb{R}^{n+1}$. We write $\Delta^n$
when $\A$ is clear.
\end{definition}

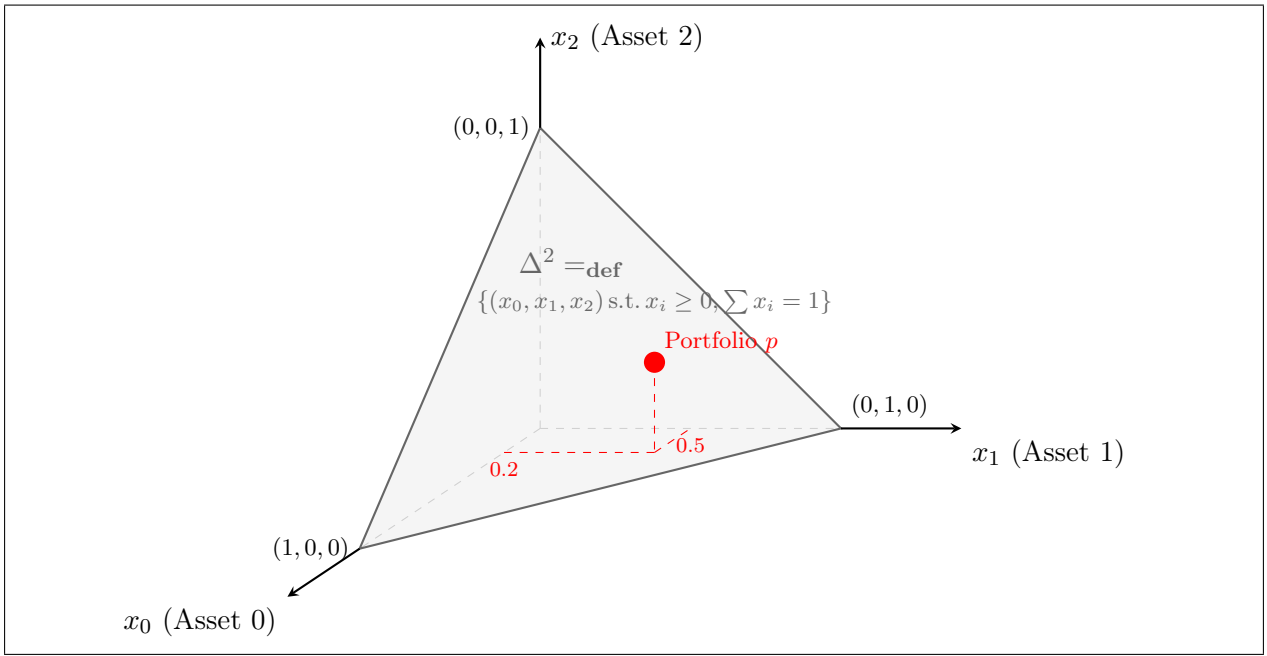
\begin{figure}[h]
\centering
\begin{tikzpicture}[
    scale=4,
    x={(-0.6cm,-0.4cm)}, 
    y={(1cm,0cm)}, 
    z={(0cm,1cm)},
    >=stealth
]

    \coordinate (O) at (0,0,0);
    \coordinate (A) at (1,0,0); 
    \coordinate (B) at (0,1,0); 
    \coordinate (C) at (0,0,1); 
    
    \draw[gray, dashed] (O) -- (A);
    \draw[gray, dashed] (O) -- (B);
    \draw[gray, dashed] (O) -- (C);

    \filldraw[
        fill=gray!10, 
        draw=gray!80!black, 
        thick, 
        fill opacity=0.7
    ] (A) -- (B) -- (C) -- cycle;

    \draw[->, thick] (A) -- (1.4,0,0) node[anchor=north east] {$x_0$ (Asset 0)};
    \draw[->, thick] (B) -- (0,1.4,0) node[anchor=north west] {$x_1$ (Asset 1)};
    \draw[->, thick] (C) -- (0,0,1.3) node[anchor=west] {$x_2$ (Asset 2)};

    
    \node[anchor=east, font=\footnotesize] at (A) {$(1,0,0)$};
    \node[anchor=south west, font=\footnotesize] at (B) {$(0,1,0)$};
    \node[anchor=east, font=\footnotesize] at (C) {$(0,0,1)$};
    
    \node[gray!80!black, font=\bfseries] at (0.25, 0.25, 0.65) {$\Delta^2=_\text{def}$};
    \node[gray!80!black, font=\footnotesize] at (0.2, 0.5, 0.5) {$\{(x_0, x_1, x_2)\text{\,s.t.\,}x_i\geq 0, \sum x_i = 1\}$};

    \coordinate (P) at (0.2, 0.5, 0.3);
    \fill[red] (P) circle (1pt);
    \node[anchor=south west, text=red, font=\footnotesize] at (P) {Portfolio $p$};
    
    \draw[dashed, red, thin] (P) -- (0.2, 0.5, 0); 
    \draw[dashed, red, thin] (0.2, 0.5, 0) -- (0.2, 0, 0); 
    \draw[dashed, red, thin] (0.2, 0.5, 0) -- (0, 0.5, 0); 
    
    \node[text=red, font=\scriptsize, anchor=north] at (0.2,0,0) {$0.2$};
    \node[text=red, font=\scriptsize, anchor=north] at (0,0.5,0) {$0.5$};

\end{tikzpicture}
\caption{Visualization of the (compact) 2-simplex $\Delta^2$ as a subset of the (non-compact) space $\mathbb{R}^3$. Any valid portfolio $p=(x_0, x_1, x_2)$ lies on the triangular plane segment where weights sum to 1 and are non-negative.}
\label{fig:simplex-r3}
\end{figure}

\begin{proposition}[Simplex Properties]
\label{prop:simplex-properties}
The ambient portfolio space $\Delta^n$ is compact, convex, Hausdorff, path-connected and locally compact.
\end{proposition}
The \textbf{compactness} of $\Delta^n$ is the most crucial property for our
framework.

\subsection{Permissible Portfolio Spaces (Objects)}

\begin{definition}[Permissible Portfolio Space (Object)]
\label{def:permissible-portfolio-space}
A \textbf{permissible portfolio space} is a \textbf{closed subset $K \subseteq \Delta^n$} for some $n$. These will serve as the objects of the double category $\mathbb{HS}$ defined in Section~\ref{sec:double}.
We interpret $K$ as the permissible region under baseline alignments (e.g.,
position limits, sector bounds).
\end{definition}

\begin{remark}
Strictly speaking, we need to keep track of what assets are being referred to: a permissible portfolio space should be defined to be a \emph{pair} $(\A, K)$ where $\A$ is a list of assets and $K$ is a geometrical object encoding portfolios which hold those assets. However, for simplicity of notation, we will write $K$ instead of $(\A, K)$, leaving $\A$ implicit.
\end{remark}

\begin{example}[Permissible Region in $\Delta^2$]
\label{ex:permissible}
\par
Let $\A = \{\text{Stocks, Bonds, Cash}\}$. An investor might have constraints:
``no more than 50\% in stocks'' ($x_0 \leq 0.5$) and ``at least 10\% in cash''
($x_2 \geq 0.1$). The permissible space is:
\[
K = \{x \in \Delta^2 \mid x_0 \leq 0.5, x_2 \geq 0.1\}
\]
This is a closed, convex polygon (a trapezoid) \emph{inside} the triangle
$\Delta^2$. This $K$ is a valid object in our category.
\end{example}

\begin{figure}[h]
\centering
\begin{tikzpicture}[
    scale=4,
    x={(-0.6cm,-0.4cm)}, 
    y={(1cm,0cm)}, 
    z={(0cm,1cm)},
    >=stealth
]

    \coordinate (O) at (0,0,0);
    \coordinate (A) at (1,0,0); 
    \coordinate (B) at (0,1,0); 
    \coordinate (C) at (0,0,1); 
    
    \draw[gray, dashed] (O) -- (A);
    \draw[gray, dashed] (O) -- (B);
    \draw[gray, dashed] (O) -- (C);

    \filldraw[
        fill=gray!10, 
        draw=gray!60!black, 
        thin, 
        fill opacity=0.5
    ] (A) -- (B) -- (C) -- cycle;

    \coordinate (K1) at (0.5, 0.4, 0.1);
    \coordinate (K2) at (0.5, 0, 0.5);
    \coordinate (K3) at (0, 0, 1);
    \coordinate (K4) at (0, 0.9, 0.1);

    \filldraw[
        fill=blue!15, 
        draw=blue!60!black, 
        thick, 
        fill opacity=0.8
    ] (K1) -- (K2) -- (K3) -- (K4) -- cycle;

    \draw[->, thick] (A) -- (1.4,0,0) node[anchor=north east] {$x_0$ (Stocks)};
    \draw[->, thick] (B) -- (0,1.4,0) node[anchor=north west] {$x_1$ (Bonds)};
    \draw[->, thick] (C) -- (0,0,1.3) node[anchor=west] {$x_2$ (Cash)};


    \node[gray!80!black, font=\small, anchor=west] at (0, 1.3, 0.5) {Ambient $\Delta^2$};

    \draw[<-, thin, blue!60!black] (0.5, 0.05, 0.45) -- (1.3, -0.3, 0.5) 
        node[anchor=east, font=\scriptsize, align=right] {Stock Cap\\($x_0 \le 0.5$)};

    \draw[<-, thin, blue!60!black] (0.2, 0.7, 0.1) -- (-0.2, 0.7, 1.2) 
        node[anchor=south, font=\scriptsize, align=center] {Cash Floor\\($x_2 \ge 0.1$)};

    \node[anchor=east, font=\scriptsize, gray] at (A) {$(1,0,0)$};
    \node[anchor=south west, font=\scriptsize, gray] at (B) {$(0,1,0)$};
    \node[anchor=east, font=\scriptsize, gray] at (C) {$(0,0,1)$};

    \node[blue!60!black, font=\bfseries] at (0.25, 0.35, 0.4) {Region $K$};

\end{tikzpicture}
\caption{Visualization of the permissible portfolio space $K \subseteq \Delta^2$ from Example~\ref{ex:permissible}. The ambient simplex (gray) represents all possible portfolios. The permissible region $K$ (blue polygon) is defined by the intersection of the constraints $x_0 \le 0.5$ (Stock Cap) and $x_2 \ge 0.1$ (Cash Floor).}
\label{fig:simplex-k-region}
\end{figure}
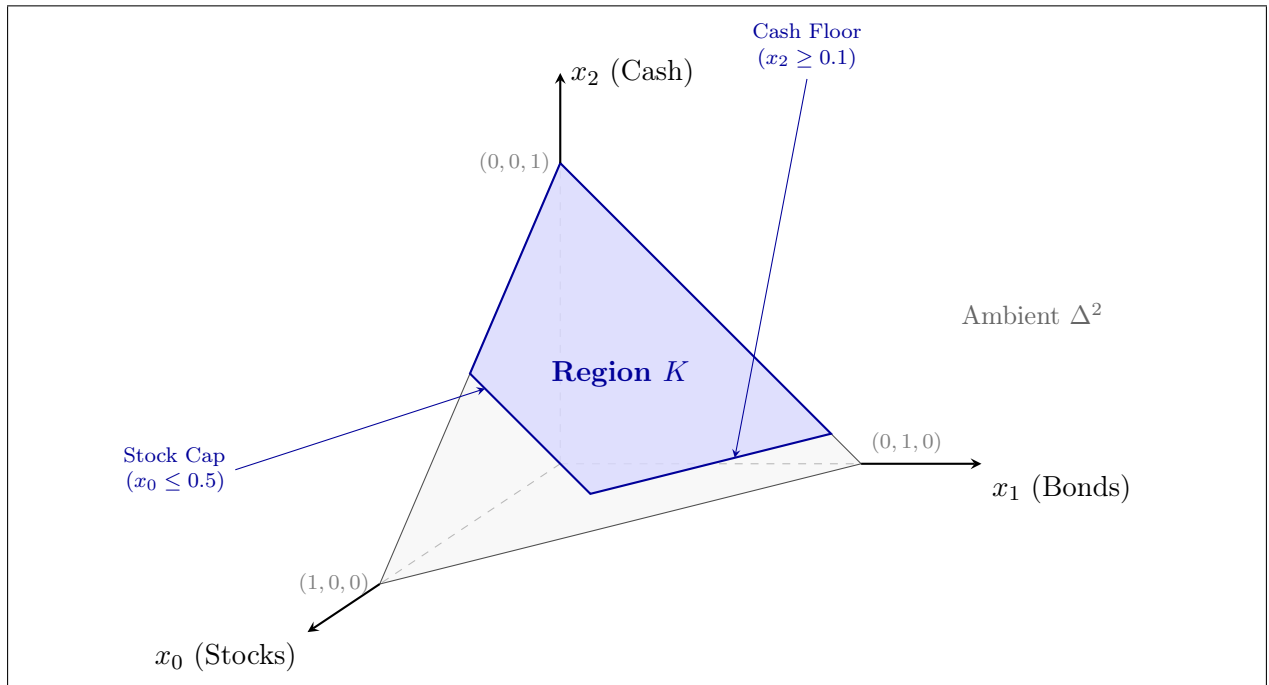

Note that every permissible portfolio space $K$ is a compact topological space.

\clearpage  \section{Portfolio Re-implementations (Horizontal Morphisms)}
\label{sec:reimplement}

\subsection{Definition and Interpretation}

\begin{definition}[Re-implementation (Horizontal Morphism)]
A \textbf{horizontal morphism} (or portfolio re-implementation) $f: K_1 \to
K_2$ from a permissible space $K_1 \subseteq \Delta^n$ to $K_2 \subseteq \Delta^m$
is a \textbf{continuous map} from the topological space $K_1$ to $K_2$.
\end{definition}

\subsection{Examples of Portfolio Re-implementations}

\begin{example}[Sector Aggregation]
\label{ex:sector-agg}
\par
Consider the set of sectors $\mathcal{B} = \{\text{Tech, Finance, Energy}\}$. Let $\mathcal{A}=\{s_{1},...,s_{100}\}$ be a set of 100 stocks partitioned into sectors:

\[
\text{Tech} = \{s_{1}, \dots, s_{30}\}, \quad
\text{Finance} = \{s_{31}, \dots, s_{60}\}, \quad
\text{Energy} = \{s_{61}, \dots, s_{100}\}
\]

Let $K_1 = \Delta^{99}(\A)$ and $K_2 = \Delta^2(\mathcal{B})$
be the full simplices on these asset sets (with coordinates indexed by the assets $s_1,\ldots,s_{100}$ and the sectors respectively). The re-implementation $f: K_1 \to K_2$ is:
\[
f(x_{s_1}, \ldots, x_{s_{100}}) = \left(\sum_{i=1}^{30} x_{s_i}, \sum_{i=31}^{60} x_{s_i},
\sum_{i=61}^{100} x_{s_i}\right)
\]
This is a linear map, which is continuous. It is a valid horizontal morphism.
\end{example}

\begin{example}[Re-implementation with Restricted Universe]
\label{ex:restricted-set}
\par
Let $\A = \{\text{Stock 1, ..., 500}\}$ (S\&P 500) and $\mathcal{B} =
\{\text{Stock 1, ..., 100}\}$ (an ESG-screened subset).
\par
\begin{itemize}
\item The ambient spoke space is $\Delta^{99}(\mathcal{B})$. Let $K_2 =
\Delta^{99}(\mathcal{B})$.
\item The ambient hub space is $\Delta^{499}(\A)$. The permissible region $K_1$
consists only of portfolios that hold the ESG assets:
\[
K_1 = \{x \in \Delta^{499} \mid x_i = 0 \text{ for } i = 101, \ldots, 500\}
\]
\item $K_1$ is a closed face of $\Delta^{499}$ (it is a sub-simplex), so it is a
valid object.
\item The trivial re-implementation $f: K_1 \to K_2$ is the projection:
\[
f(x_1, \ldots, x_{100}, 0, \ldots, 0) = (x_1, \ldots, x_{100})
\]
This map is linear and continuous, so it is a valid horizontal morphism.
\end{itemize}
We may be interested in extending this to a re-implementation $g: \Delta^{499}(\A) \to K_2$ that takes an arbitrary portfolio to its ``nearest ESG-compliant portfolio''. This problem is not well-defined unless we specify what ``nearest'' means.
\end{example}

  \begin{example}[Factor Decomposition]
  \label{ex:factor}
  \par
  Let $\A = \{a_1, \ldots, a_n\}$ (assets) and factors (see \cite{FamaFrench})
  \[\mathcal{B} = \{\text{Market,
  Value, Momentum, Quality}\}\]
  Given a factor loading matrix $F \in \mathbb{R}^{4 \times n}$ where $F_{ji}$ is
  the exposure of asset $a_i$ to factor $j$, we wish to decompose portfolios
  into factor allocations, restricting to portfolios whose factor exposures
  are all non-negative.
  \par
  \begin{itemize}
  \item Let $K_1 = \{x \in \Delta^{n-1} \mid (Fx)_j \geq 0 \text{ for all } j\}$.
  This is an intersection of closed half-spaces with $\Delta^{n-1}$, so $K_1$ is
  a closed, compact set (a valid object).
  \item Let $K_2 = \Delta^3(\mathcal{B})$.
  \item Define the normalization map $f: K_1 \to K_2$ by
  \[
  f(x) = \frac{Fx}{\|Fx\|_1}.
  \]
  This is well-defined and continuous wherever $\|Fx\|_1 > 0$.  However,
  $K_1$ may contain portfolios $x$ with $Fx = 0$ (zero total factor exposure),
  at which $f$ is undefined, and to which it cannot in general be continuously
  extended (the limit of $Fx/\|Fx\|_1$ as $Fx \to 0$ depends on the direction
  of approach).
  \end{itemize}

  To obtain a valid horizontal morphism, we impose one of the following
  conditions ensuring $\|Fx\|_1 > 0$ on $K_1$:
  \begin{enumerate}
  \item[\textbf{(a)}] \textbf{Non-negative loadings with full coverage.}
    If $F$ has non-negative entries, then $(Fx)_j \geq 0$ automatically
    for $x \in \Delta^n$, and $\|Fx\|_1 = (\mathbf{1}^T F)\, x$.  This
    vanishes only for portfolios supported on assets $i$ with
    $\sum_j F_{ji} = 0$ (no factor exposure whatsoever).
    Excluding such factor-inert assets from $\A$---or equivalently
    assuming each asset loads on at least one factor---gives
    $\|Fx\|_1 > 0$ on all of~$\Delta^{n-1}$, and $K_1 = \Delta^{n-1}$.

  \item[\textbf{(b)}] \textbf{Genericity.}  For general (possibly
    negative) loading matrices, assume $\ker(F) \cap K_1 = \emptyset$,
    i.e., no portfolio in $K_1$ has zero exposure to every factor.
    Since $\ker(F)$ has codimension at most $4$ in $\mathbb{R}^n$,
    this condition holds for generic $F$ when $n \gg 4$.
    Then $\|Fx\|_1$ is continuous and strictly positive on the compact
    set $K_1$, so $f$ is continuous on $K_1$.

  \item[\textbf{(c)}] \textbf{Minimum exposure threshold.}
    Replace $K_1$ with
    $K_1^{(c)} = \{x \in \Delta^n \mid (Fx)_j \geq 0 \;\text{for all}\; j,
    \;\|Fx\|_1 \geq c\}$
    for a threshold $c > 0$ representing the minimum meaningful factor tilt.
    This is still a closed, compact subset of $\Delta^n$ (an intersection of
    closed half-spaces and a closed half-space), and $f$ is continuous
    on $K_1^{(c)}$.
  \end{enumerate}

  Under any of these conditions, $f: K_1 \to K_2$ is a continuous map between
  compact sets and hence a valid horizontal morphism.
  \end{example}

\begin{example}[Optimal ETF Replacement (Nonlinear)]
\label{ex:etf-replacement}
\par
This is a key example, developed fully in Section~\ref{sec:optimal-reimpl}.
A map $f$ is constructed that takes a ``hub'' portfolio $x$ (e.g., of mutual
funds) and finds the ``spoke'' portfolio $f(x)$ (e.g., of ETFs) that
\emph{minimizes fees} subject to \emph{tracking $x$'s exposures} within some
error $\epsilon$.

As we will prove in Theorem~\ref{thm:berge-optimal}, this optimization problem
naturally defines:
\begin{enumerate}
\item A closed permissible ``hub'' space $K_1 = \{x \mid \text{tracking is
possible}\}$.
\item A continuous (and generally nonlinear) map $f: K_1 \to K_2$.
\end{enumerate}
This provides a practical source of valid, nonlinear horizontal morphisms.
\end{example}

\clearpage  \section{Alignment Relations (Vertical Morphisms)}
\label{sec:constraints}

\subsection{Closed Relations as Alignments}

\begin{definition}[Alignment Relation (Vertical Morphism)]
A \textbf{vertical morphism} (or alignment relation) $R$ from a permissible space
$K_1 \subseteq \Delta^n$ to $K_3 \subseteq \Delta^p$ is a \textbf{closed subset}:
\[
R \subseteq K_1 \times K_3
\]
We write $(x, z) \in R$ to mean ``portfolio $x \in K_1$ aligns with portfolio
$z \in K_3$.''
\end{definition}

Note that alignment is regarded as a simple binary relation: either $x$ and $z$ are aligned, or they are not. We do not (yet) have a notion of ``how aligned'' $x$ and $z$ are.

\begin{figure}[htbp]
\centering
\begin{tikzpicture}[scale=1.2, >=stealth]
    \draw[->, thick] (0,0) -- (5,0) node[right] {Permissible Space $K_1$ (Hub)};
    \draw[->, thick] (0,0) -- (0,4) node[above] {Permissible Space $K_3$ (spoke)};

    \draw[ultra thick, blue!60!black] (0.5, -0.1) -- (4.5, -0.1) node[midway, below] {Dom($R$)};
    \draw[ultra thick, blue!60!black] (-0.1, 0.5) -- (-0.1, 3.5) node[midway, left] {Cod($R$)};

    \filldraw[fill=purple!20, draw=purple!80!black, thick] 
        (1.0, 1.0) -- (3.5, 1.5) -- (4.0, 3.0) -- (1.5, 2.5) -- cycle;

    \node[purple!80!black, font=\footnotesize] at (2.5, 2.0) {Alignment $R$};
    \node[purple!60!black, font=\footnotesize] at (2.5, 1.7) {(closed subset)};

    \fill[black] (2.8, 2.2) circle (1.5pt);
    \draw[dashed, thin] (2.8, 2.2) -- (2.8, 0) node[below, font=\small] {$x$};
    \draw[dashed, thin] (2.8, 2.2) -- (0, 2.2) node[left, font=\small] {$z$};

    \node[right, align=left, font=\footnotesize, gray] at (4.2, 3) 
        {Points $(x,z) \in R$\\signify that portfolio $x$\\aligns with spoke $z$.};

\end{tikzpicture}
\caption{Visualization of an alignment relation $R \subseteq K_1 \times K_3$ as a closed subset (purple) of the product space. The closedness guarantees that limits of aligned portfolios remain aligned.}
\label{fig:alignment-relation}
\end{figure}
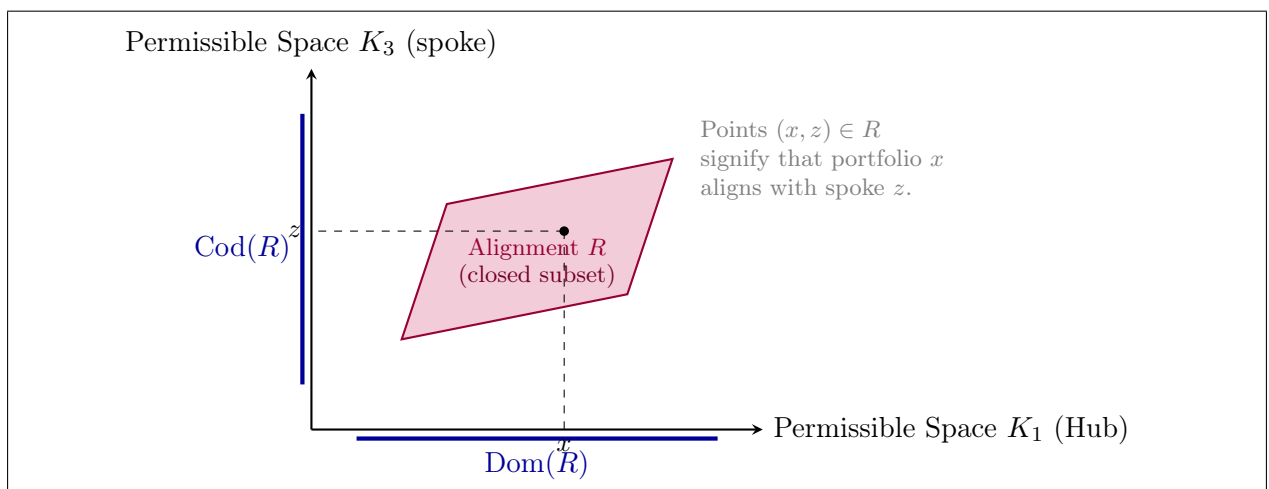

\begin{remark}[Interpretation]
Alignment relations encode various portfolio requirements:
\par
\begin{itemize}
\item \textbf{Benchmark tracking}: $R \subseteq K_1 \times K_{\text{bench}}$
relates a portfolio $x \in K_1$ to a benchmark $z \in K_{\text{bench}}$ if the
tracking error of $x$ with respect to $z$ is $\leq \epsilon$. Note that $K_1$ and $K_{\text{bench}}$ are often different portfolio spaces: e.g. the portfolios in $K_{\text{bench}}$ may be combinations of asset classes or indices, while the portfolios in $K_1$ may be built from a menu of mutual funds.
\item \textbf{Risk budgets}: $R \subseteq K_1 \times K_1$ relates two
portfolios if they have sufficiently close risk profiles.
\item \textbf{Factor exposures}: $R$ relates a portfolio to a vector of
factor loadings.
\end{itemize}
Closedness ensures stability: if a sequence $(x_k, z_k) \in
R$ converges to $(x, z)$, then $(x, z) \in R$. Alignments are preserved under
limits.
\end{remark}

\begin{remark}[Compliance/eligibility/permissibility as alignment]
\label{rem:compliance}
Consider the special case $K_3=\textbf{1}$, where $\textbf{1}$ is the one-point set. In this case, an alignment relation $R \subseteq K_1 \times \textbf{1}$ is simply a closed subset of $K_1$. It can be regarded as carving out a smaller permissible set of portfolios that satisfy some additional compliance requirement.
\end{remark}

\subsection{Vertical Composition}

\begin{definition}[Vertical Composition]
\label{def:vertical-composition}
Given $R: K_1 \to K_2$ and $S: K_2 \to K_3$ (i.e., $R \subseteq K_1 \times K_2$
and $S \subseteq K_2 \times K_3$), the \textbf{vertical composition} is
standard relational composition:
\[
S \circ R = \{(x, z) \in K_1 \times K_3 \mid \exists y \in K_2:
(x, y) \in R \text{ and } (y, z) \in S\}
\]
\end{definition}

\begin{center}
\begin{tikzcd}
K_1 \arrow[d, dashed, "R"'] \\ 
K_2 \arrow[d, dashed, "S"'] \\ 
K_3
\end{tikzcd}
\end{center}

\begin{figure}[htbp]
\centering
\begin{tikzpicture}[
    scale=1.5,
    x={(1cm,0cm)},    
    y={(0.6cm,0.5cm)},  
    z={(0cm,1cm)},    
    >=stealth
]

    \coordinate (O) at (0,0,0);
    \coordinate (K1_end) at (5,0,0);
    \coordinate (K2_end) at (0,5,0);
    \coordinate (K3_end) at (0,0,5);

    \draw[->, thick, gray!60] (O) -- (K1_end) node[anchor=north] {$K_1$ (Hub)};
    \draw[->, thick, gray!60] (O) -- (K2_end) node[anchor=west] {$K_2$ (Intermediate)};
    \draw[->, thick, gray!60] (O) -- (K3_end) node[anchor=east] {$K_3$ (Spoke)};

    
    \draw[ultra thick, blue] (1,0,0) -- (3,0,0);
    \node[blue, font=\scriptsize, anchor=north] at (2,0,0) {Permissible $K_1$};

    \draw[ultra thick, blue] (0,1,0) -- (0,3,0);
    \node[blue, font=\scriptsize, anchor=east] at (0,2,0) {Permissible $K_2$};

    \draw[ultra thick, blue] (0,0,1) -- (0,0,3);
    \node[blue, font=\scriptsize, anchor=east] at (0,0,2) {Permissible $K_3$};

    
    \draw[dotted, blue!40] (1,0,0) -- (1,3,0); 
    \draw[dotted, blue!40] (3,0,0) -- (3,3,0); 
    \draw[dotted, blue!40] (0,1,0) -- (3,1,0); 
    \draw[dotted, blue!40] (0,3,0) -- (3,3,0); 

    \filldraw[fill=red!20, draw=red!80, opacity=0.8] 
        (1.2, 1.2, 0) -- (2.8, 1.4, 0) -- (2.6, 2.8, 0) -- (1.4, 2.6, 0) -- cycle;
    \node[red!80!black, font=\footnotesize] at (3.5, 1.2, 0) {$R \subsetneq K_1 \times K_2$};

    
    \draw[dotted, blue!40] (0,1,0) -- (0,1,3); 
    \draw[dotted, blue!40] (0,3,0) -- (0,3,3); 
    \draw[dotted, blue!40] (0,0,1) -- (0,3,1); 
    \draw[dotted, blue!40] (0,0,3) -- (0,3,3); 

    \filldraw[fill=orange!20, draw=orange!80, opacity=0.8] 
        (0, 1.4, 1.2) -- (0, 2.8, 1.5) -- (0, 2.6, 2.7) -- (0, 1.2, 2.5) -- cycle;
    \node[orange!80!black, font=\footnotesize, anchor=south] at (0, 2, 3.2) {$S \subsetneq K_2 \times K_3$};

    
    \draw[dotted, blue!40] (1,0,0) -- (1,0,3); 
    \draw[dotted, blue!40] (3,0,0) -- (3,0,3); 
    \draw[dotted, blue!40] (0,0,1) -- (3,0,1); 
    \draw[dotted, blue!40] (0,0,3) -- (3,0,3); 

    \filldraw[fill=blue!15, draw=blue!80, opacity=0.6] 
        (1.3, 0, 1.3) -- (2.7, 0, 1.6) -- (2.5, 0, 2.6) -- (1.5, 0, 2.4) -- cycle;
    
    \node[blue!80!black, font=\bfseries, anchor=north] at (2.2, 0, 2.0) {$S \circ R \subseteq K_1 \times K_3$};

    
    \coordinate (pxy) at (2.0, 2.2, 0);   
    \coordinate (pyz) at (0, 2.2, 1.8);   
    \coordinate (pxz) at (2.0, 0, 1.8);   
    \coordinate (p3d) at (2.0, 2.2, 1.8); 

    \fill[red] (pxy) circle (1.5pt);
    \fill[red] (pyz) circle (1.5pt);
    \fill[blue] (pxz) circle (2pt);

    \draw[dashed, gray] (pxy) -- (p3d);
    \draw[dashed, gray] (pyz) -- (p3d);
    \draw[dashed, blue] (p3d) -- (pxz);
    
    \node[font=\scriptsize, anchor=west, align=left] at (p3d) {$\exists y \in K_2$ s.t. $(x, y)\in R$ and $(y, z)\in S$};

\end{tikzpicture}
\caption{3D visualization of vertical composition. The relations $R$ (red) and $S$ (orange) are irregular polygons strictly contained within the product of the permissible sets (indicated by the dotted blue frames). The relation $S\circ R$ (blue) is a polygonal subset of the $K_1$--$K_3$ plane.}
\label{fig:vertical-composition-3d-subsets}
\end{figure}
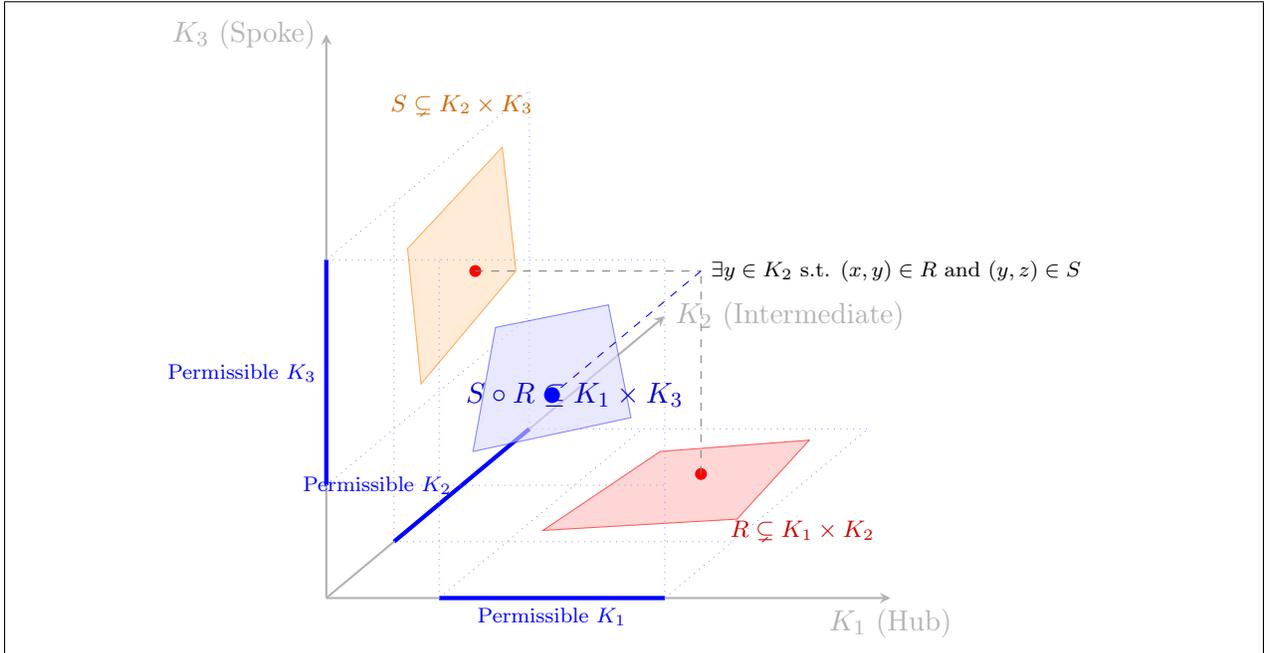

\begin{lemma}[Closedness of Vertical Composition]
\label{lem:vertical-closed}
If $R$ and $S$ are vertical morphisms (closed relations), then $S \circ R$ is
also a vertical morphism (a closed relation).
\end{lemma}

\begin{proof}
Given $R \subseteq K_1 \times K_2$ and $S \subseteq K_2 \times K_3$, both
closed. We need $S \circ R$ to be closed in $K_1 \times K_3$.
Consider the subset $E = \{((x, y), (y', z)) \in (R \times S) \mid y = y'\}$ of
the product $R \times S$. This is closed.
The objects $K_1, K_2, K_3$ are all compact. Thus $R$ and $S$ are compact. $R \times S$ is
compact, and $E$ is compact (as a closed subset of a compact space).
The map $\pi: E \to K_1 \times K_3$ defined by $\pi((x, y), (y, z)) = (x, z)$
is continuous. The image $S \circ R = \pi(E)$ is compact (continuous image of
compact), hence closed in the Hausdorff space $K_1 \times K_3$.
\end{proof}

\begin{example}[Multi-Layer Asset Allocation]
\label{ex:asset-alloc-comp}
Let $K_1$ be the space of portfolios consisting of individual securities (stocks, bonds, cash), and let $K_2$ be the space of portfolios of asset classes (e.g., US Equity, International Equity, Fixed Income). Let $g: K_1 \to K_2$ be the aggregation map that sums the weights of securities within each asset class.

Consider the following two alignment relations:
\begin{itemize}
    \item $R : K_1 \dashrightarrow K_2$: The ``Implementation'' alignment. A portfolio $x \in K_1$ aligns with a recommended tactical asset allocation $y \in K_2$ if the portfolio implements the allocation within a given tolerance $\epsilon$:
    $R = \{ (x, y) \mid \| g(x) - y \| \le \epsilon \} \subseteq K_1 \times K_2$.
    \item $S : K_2 \dashrightarrow K_2$: The ``Policy'' alignment. A tactical asset allocation $y \in K_2$ aligns with a long-term strategic baseline $z \in K_2$ if the tactical shifts are within a permitted band $\delta$:
    $S = \{ (y, z) \mid \| y - z \| \le \delta \} \subseteq K_2 \times K_2$.
\end{itemize}

The vertical composition $S \circ R : K_1 \dashrightarrow K_2$ relates the actual security holdings $x$ directly to the strategic baseline $z$. The condition $(x, z) \in S \circ R \subseteq K_1 \times K_2$ means that the portfolio $x$ is a valid implementation of \emph{some} tactical allocation $y$ that is itself compliant with the strategic baseline $z$.
\end{example}

\clearpage  \section{The Double Category $\HSSimp$}
\label{sec:double}

\subsection{Definition of the Double Category}

\begin{definition}[Double Category $\HSSimp$]
The \textbf{double category of permissible (Hub-and-Spoke) portfolio spaces}
$\HSSimp$ consists of:
\par
\begin{itemize}
\item \textbf{Objects}: permissible portfolio spaces $K$ (closed subsets of
simplices $\Delta^n$)
\item \textbf{Horizontal morphisms}: Continuous re-implementations $f: K_1 \to
K_2$
\item \textbf{Vertical morphisms}: Closed alignment relations $R \subseteq K_1
\times K_3$
\item \textbf{Horizontal composition}: Standard function composition $g \circ f$.
\item \textbf{Vertical composition}: Standard relational composition $S \circ R$.
\item \textbf{Horizontal identities}: $\text{id}_K: K \to K$ (identity map).
\item \textbf{Vertical identities}: $\Delta_K = \{(x, x) \mid x \in K\}$
(diagonal relation).
\item \textbf{2-cells}: For a square
\begin{center}
\begin{tikzcd}
K_1 \arrow[r, "f"] \arrow[d, dashed, "R"'] & K_2 
\arrow[d, dashed, "S"] \\ 
K_3 \arrow[r, "g"'] & K_4
\end{tikzcd}
\end{center}
a 2-cell exists if and only if the inclusion $\Graph(g) \circ R \subseteq S
\circ \Graph(f)$ holds. This makes $\HSSimp$ a \textbf{thin} double category.
\end{itemize}
\end{definition}

\begin{remark}
Note that this definition is strongly reminiscent of the simplest example of a double category: the double category of sets, functions and relations. Its objects are sets, horizontal morphisms are set functions, vertical morphisms are relations between sets, and 2-cells exist when the above inclusion holds. However, our definition has a much more restricted class of objects, and is also topological, leading to a double category with less structure.
\end{remark}

\begin{remark}[Graph closedness]
\label{rem:graph-closed}
For any continuous map $f: K_1 \to K_2$ between compact Hausdorff spaces, the graph $\Graph(f) = \{(x, f(x)) : x \in K_1\}$ is a closed subset of $K_1 \times K_2$. This follows because the map $x \mapsto (x, f(x))$ is continuous and $K_1$ is compact, so the image is compact, hence closed in the Hausdorff space $K_1 \times K_2$. Therefore $\Graph(f)$ is a valid vertical morphism in $\HSSimp$.
\end{remark}

\begin{remark}
The thinness assumption reflects the idea that alignment checking yields a binary yes/no answer---a simplification that limits expressivity but is precisely what makes the coherence laws automatic. Thinness is dropped in Section~\ref{sec:2-cells}, which introduces a richer framework incorporating evidence of alignment.
\end{remark}

\begin{theorem}[Well-definedness of $\HSSimp$]
\label{thm:coherence}
$\HSSimp$ is a well-defined double category.
\end{theorem}

  \begin{proof}
  We verify the double category axioms:
  \begin{itemize}
  \item \emph{Horizontal category:} Composition of continuous maps
    is continuous; identity maps are continuous.
  \item \emph{Vertical category:} Composition of closed relations
    preserves closedness (Lemma~\ref{lem:vertical-closed});
    diagonal relations $\Delta_K$ are closed in $K \times K$.
  \item \emph{Associativity and units:} Horizontal associativity
    and units are inherited from function composition; vertical
    associativity and units from relational composition.
  \item \emph{2-cells and interchange:}
    The 2-cell condition
    $\Graph(g) \circ R \subseteq S \circ \Graph(f)$
    is a property (not structure), so $\HSSimp$ is thin.
    The interchange law holds because relational composition
    is monotone with respect to inclusion.  For horizontal
    pasting: given composable 2-cells
    $\Graph(g) \circ R \subseteq S \circ \Graph(f)$ (left square)
    and $\Graph(g') \circ S \subseteq T \circ \Graph(f')$ (right square),
    monotonicity gives
    $\Graph(g') \circ \Graph(g) \circ R
    \subseteq \Graph(g') \circ S \circ \Graph(f)
    \subseteq T \circ \Graph(f') \circ \Graph(f)$,
    i.e., $\Graph(g' \circ g) \circ R \subseteq T \circ \Graph(f' \circ f)$.
    The vertical pasting argument is analogous, using
    $\Graph(g) \circ R \subseteq S \circ \Graph(f)$ and
    $\Graph(h) \circ S' \subseteq T' \circ \Graph(g)$
    to obtain
    $\Graph(h) \circ (S' \circ R) \subseteq (T' \circ S) \circ \Graph(f)$
    via the same monotonicity.
  \end{itemize}
  \end{proof}

\begin{remark}[2-Cell Interpretation]
The inclusion $\Graph(g) \circ R \subseteq S \circ \Graph(f)$ has a concrete reading:
\par
\begin{itemize}
\item $\Graph(g) \circ R = \{(x, z) \mid \exists w \in K_3: (x, w) \in R, g(w) = z\}$
(Go ``down'' via $R$ to $w$, then ``right'' via $g$ to $z=g(w)$)
\item $S \circ \Graph(f) = \{(x, z) \mid (f(x), z) \in S\}$
(Go ``right'' via $f$ to $y=f(x)$, then ``down'' via $S$ to $z$)
\end{itemize}
The 2-cell condition means: ``For \textit{any} portfolio $x$, if the 'down-then-right'
path is possible, then the 'right-then-down' path must also be possible and
yield the same alignment.''
\end{remark}

\begin{remark}[Financial Intuition for the Inclusion Direction]
The direction of the 2-cell inclusion, $\Graph(g) \circ R \subseteq S \circ \Graph(f)$, is chosen to ensure that the framework acts as a \textbf{safety guarantee} rather than a rigid identity.

In financial operations, $R$ and $S$ typically represent compliance or alignment \emph{bounds} (e.g., tracking error $\le 10$\,bps). The inclusion asserts that if a hub portfolio is compliant ($R$), any portfolio generated by the transformation process ($g \circ R$) must land within the permissible downstream bounds ($S \circ f$).

If the inclusion were reversed ($S \circ \Graph(f) \subseteq \Graph(g) \circ R$), it would imply that \textit{every} compliant spoke portfolio must have originated from a compliant hub. This is often false: a re-implementation might ``accidentally'' satisfy a downstream sector limit even if the parent hub violated it. The inclusion as defined preserves investment intent without over-constricting the valid outcomes of the re-implementation.
\end{remark}

\begin{example}[Example Diagram]
\label{ex:diagram}
\par
Consider the following diagram

\begin{center}
\begin{tikzcd}[column sep=huge, row sep=huge]
K_{\text{Funds}} \arrow[r, "f_{\text{ETF}}"]
\arrow[d, dashed, "R_{\text{Track}}"'] & K_{\text{ETFs}} \arrow[d, dashed, 
"S_{\text{Factor}}"] \\ 
K_{\text{Bench}} \arrow[r, "g_{\text{Factor}}"'] &
K_{\text{Factors}}
\end{tikzcd}
\end{center}
\par
\noindent where
\begin{itemize}
\item $K_{\text{Funds}}$: permissible active mutual fund portfolios.
\item $K_{\text{ETFs}}$: permissible passive ETF baskets.
\item $K_{\text{Bench}}$: Benchmark portfolios.
\item $K_{\text{Factors}}$: Factor exposure vectors.
\item $f_{\text{ETF}}$: Re-implementation (e.g., optimal replacement).
\item $g_{\text{Factor}}$: Re-implementation (e.g., factor decomposition).
\item $R_{\text{Track}}$: Alignment relation (e.g., tracking error is sufficiently low).
\item $S_{\text{Factor}}$: Alignment relation (e.g., factor loadings are sufficiently close).
\end{itemize}
A 2-cell here asserts compatibility: if a fund is aligned with a benchmark (down) and
that benchmark maps to a bundle of factor exposures (right), then that fund must map to
an ETF (right) that is aligned with that same bundle of factor exposures (down). Intuitively, this will be the case if the tolerance in factor loading differences is ``looser'' than the tracking error tolerance.
\end{example}

\subsection{Alignment Constraints on Re-implementation}
\label{subsec:reimplalign}

An important special case of 2-cell existence is the following case:

\begin{center}
\begin{tikzcd}
K_1 \arrow[r, "id"] \arrow[d, dashed, "id"'] & K_1 
\arrow[d, dashed, "R"] \\ 
K_1 \arrow[r, "f"'] & K_2
\end{tikzcd}
\end{center}

In this case, a 2-cell exists if and only if the re-implementation morphism $f$ respects the alignment relation $R$ between portfolios in its domain and codomain.

A typical portfolio management problem is: given an alignment relation $R$, construct a re-implementation morphism $f$ that respects $R$, and in addition is ``optimal'' in some sense. For further discussion, see Section~\ref{subsec:opt-align}.

\begin{remark}[Relation to optimization and to menus]
    The condition that a re-implementation $f$ respects an alignment $R$ (i.e., the existence of the 2-cell) can be viewed in two complementary ways, which we develop later:
    \begin{itemize}
        \item \textbf{Constructive View (Section~\ref{subsec:opt-align}):} Instead of checking if a pre-existing map $f$ satisfies $R$, we often \textit{construct} $f$ by selecting the optimal spoke portfolio from the set of all aligned possibilities.
        \item \textbf{Algebraic View (Section~\ref{sec:DOTS}):} In the DOTS framework, we treat the alignment $R$ not as a constraint on a single map, but as an operator that generates the entire ``menu'' of permissible spoke portfolios (denoted $K \odot R$). A valid horizontal morphism $f$ is then simply a continuous selection from this menu.
    \end{itemize}
\end{remark}

\clearpage  \section{Pullback and Pushforward Operations}
\label{sec:operations}

\subsection{Pullback: Reindexing Alignments}

\begin{definition}[Pullback]
\label{def:pullback}
Given $f: K_1 \to K_2$ and $S \subseteq K_2 \times K_3$, the
\textbf{pullback} is:
\[
f^*S = \{(x, z) \in K_1 \times K_3 \mid (f(x), z) \in S\}
\]
\end{definition}

\begin{proposition}[Pullback Preserves Closedness]
\label{prop:pullback-closed}
If $f$ is a horizontal morphism and $S$ is a vertical morphism, $f^*S$ is a
vertical morphism.
\end{proposition}

\begin{proof}
The map $f \times \text{id}_{K_3}: K_1 \times K_3 \to K_2 \times K_3$ is
continuous (product of continuous maps). $S$ is closed in $K_2 \times K_3$ (by
definition). $f^*S = (f \times \text{id}_{K_3})^{-1}(S)$ is the continuous
preimage of a closed set, so it is closed in $K_1 \times K_3$.
\end{proof}

\begin{remark}[Portfolio Interpretation: Pre-alignment]
Pullback answers: \textbf{``Which hub portfolios in $K_1$ re-implement to
portfolios in $K_2$ that satisfy $S$?''}
\end{remark}

\subsection{Pushforward: Forward Propagation}

\begin{definition}[Pushforward]
\label{def:pushforward}
Given $f: K_1 \to K_2$ and $R \subseteq K_1 \times K_3$, the
\textbf{pushforward} is:
\[
f_!R = \{(y, z) \in K_2 \times K_3 \mid \exists x \in K_1: f(x) = y
\text{ and } (x, z) \in R\}
\]
\end{definition}

\begin{remark}[Portfolio Interpretation: Forward Alignment Propagation]
Pushforward answers: \textbf{``Which spoke portfolios in $K_2$ are
re-implementations of hub portfolios in $K_1$ that satisfy $R$?''}
\end{remark}

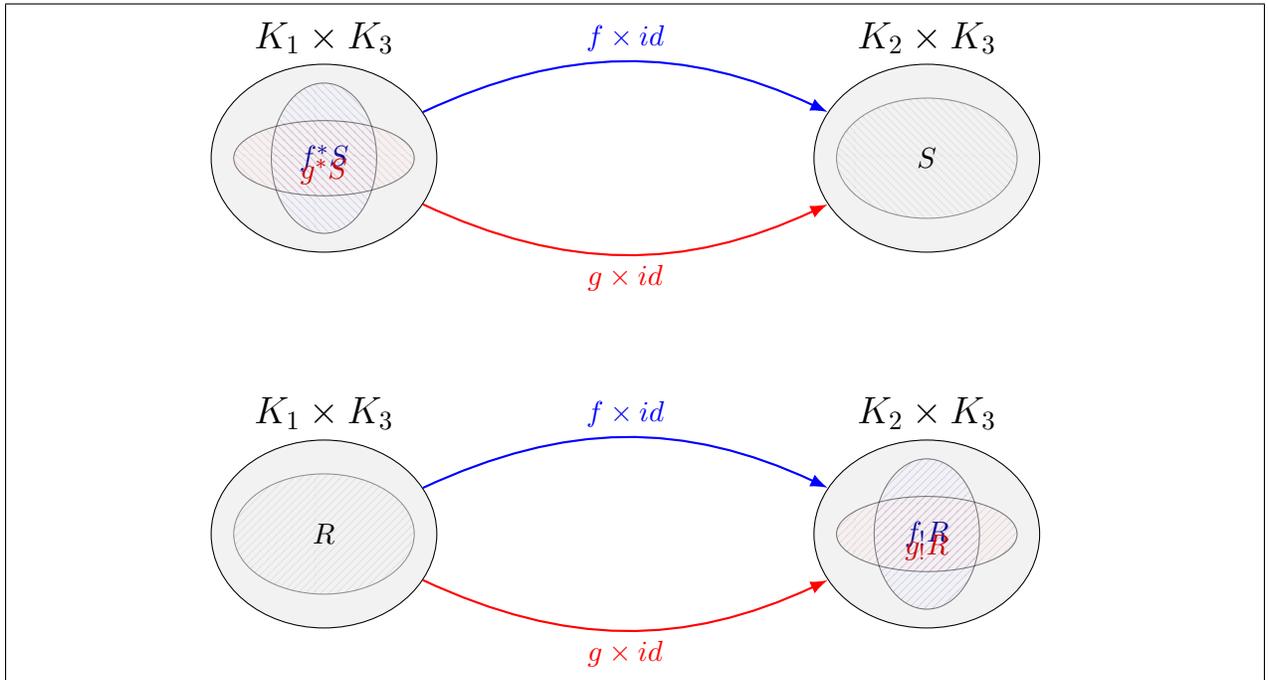
\begin{figure}[h]
    \centering
    \begin{tikzpicture}[
        obj/.style={draw, ellipse, minimum width=3cm, minimum height=2.5cm, fill=gray!10},
        lbl/.style={font=\Large}
    ]
        \begin{scope}
            \node (K1) [obj, label={[lbl]above:$K_1 \times K_3$}] {};
            \node (K2) [obj, label={[lbl]above:$K_2 \times K_3$}, right=5cm of K1] {};

            \draw [->, thick, blue, -Latex] (K1) to [bend left=25]
                  node [midway, above] {$f \times id$} (K2);
            \draw [->, thick, red, -Latex] (K1) to [bend right=25]
                  node [midway, below] {$g \times id$} (K2);

            \draw [fill=gray, opacity=0.4, pattern=north west lines, pattern color=gray!80]
                  (K2.center) ellipse (1.2cm and 0.8cm);
            \node at (K2.center) {$S$};

            \draw [fill=blue, opacity=0.5, pattern=north west lines, pattern color=blue!50]
                  (K1.center) ellipse (0.7cm and 1.0cm);
            \node at (K1.center) [text=blue!60!black] {$f^*S$};

            \draw [fill=red, opacity=0.5, pattern=north west lines, pattern color=red!50]
                  (K1.center) ellipse (1.2cm and 0.5cm);
            \node at (K1.center) [yshift=-5pt, text=red!80!black] {$g^*S$};
        \end{scope}

        \begin{scope}[yshift=-5cm]
            \node (K1R) [obj, label={[lbl]above:$K_1 \times K_3$}] {};
            \node (K2R) [obj, label={[lbl]above:$K_2 \times K_3$}, right=5cm of K1R] {};

            \draw [->, thick, blue, -Latex] (K1R) to [bend left=25]
                  node [midway, above] {$f \times id$} (K2R);
            \draw [->, thick, red, -Latex] (K1R) to [bend right=25]
                  node [midway, below] {$g \times id$} (K2R);

            \draw [fill=gray, opacity=0.4, pattern=north east lines, pattern color=gray!80]
                  (K1R.center) ellipse (1.2cm and 0.8cm);
            \node at (K1R.center) {$R$};

            \draw [fill=blue, opacity=0.5, pattern=north east lines, pattern color=blue!50]
                  (K2R.center) ellipse (0.7cm and 1.0cm);
            \node at (K2R.center) [text=blue!60!black] {$f_!R$};

            \draw [fill=red, opacity=0.5, pattern=north east lines, pattern color=red!50]
                  (K2R.center) ellipse (1.2cm and 0.5cm);
            \node at (K2R.center) [yshift=-5pt, text=red!80!black] {$g_!R$};
        \end{scope}
    \end{tikzpicture}
    \caption{Conceptual diagram comparing Pullback (top) and Pushforward (bottom) for two different horizontal morphisms, $f$ and $g$. The operations project the alignment relations to different (and differently shaped) regions in the hub or spoke space.}
    \label{fig:pull-push-offset}
\end{figure}

In general, the pushforward of a closed set under a continuous map is not closed---this is precisely where our topological requirements earn their keep.

\begin{example}[Pushforward Fails Closedness]
\label{ex:pushforward-not-closed}
\par
This example shows why we cannot use non-closed objects.
Let $K_1 = [0, 1)$ (an open interval, \emph{not} a valid object in $\HSSimp$).
Let $f: K_1 \to \Delta^1$ be the inclusion $f(x) = (x, 1-x)$. This is continuous.
Let $K_3 = [0,1] \subset \Delta^1$ and $R = \{(x, x) \mid x \in K_1\}$ (closed in $K_1 \times K_3$).
Then:
\[
f_!R = \{( (y_0, y_1), z) \mid y_0 = z, (y_0, y_1) = (x, 1-x), x \in [0, 1) \}
\]
\[
f_!R = \{ ((x, 1-x), x) \mid x \in [0, 1) \}
\]
This set is \emph{not closed} in $\Delta^1 \times K_3$. The limit point
$((1, 0), 1)$ (as $x \to 1$) is missing.

\textbf{Portfolio problem}: A ``phantom portfolio'' $((1, 0), 1)$ appears to
satisfy the propagated alignment, but has no actual hub portfolio
satisfying the original alignment. This failure is precisely because $K_1$ was
not closed (and therefore not compact).
\end{example}

The consequences are not confined to this example. Phantom portfolios break the audit trail (no permissible hub produced them), destroy the adjunction between hub-side and spoke-side verification (Theorem~\ref{thm:adjunction}), create numerical instability (a rounding error can flip a portfolio from valid to invalid), and prevent logical composition (the Beck--Chevalley and Frobenius laws cease to hold, so the order of compliance checks changes the result). Section~\ref{sec:failure} demonstrates these failures in detail.

\clearpage  \section{The Central Role of Properness}
\label{sec:proper}

\subsection{Properness and Pushforward}

\begin{definition}[Proper Map]
A continuous map $f: X \to Y$ between topological spaces is \textbf{proper} if
for every compact $K \subseteq Y$, the preimage $f^{-1}(K)$ is compact in $X$.
\end{definition}

\begin{remark}
A continuous map $f : K_1 \to K_2$ from a compact space $K_1$ to any topological space $K_2$ has the property that the preimage of any compact set in $K_2$ is compact in $K_1$. If $K_2$ is Hausdorff, then compact subsets are closed, which is used below to deduce closedness of pushforwards.
\end{remark}

\begin{proposition}[Horizontal Morphisms are Proper]
\label{prop:reimpl-are-proper}
Every horizontal morphism $f: K_1 \to K_2$ in $\HSSimp$ is a proper map.
\end{proposition}
\begin{proof}
By definition, $f$ is a continuous map. The object $K_1$ is compact. The object $K_2$ is a subset of
$\Delta^m$, which is Hausdorff. A continuous map from a compact space ($K_1$)
to a Hausdorff space ($K_2$) is automatically a proper map.
\end{proof}

Properness means: \textbf{no boundary escape; limits behave well.}
Our framework's design (using compact objects $K$) ensures this.

Automatic properness is precisely what is needed to ensure the pushforward
operation is well-defined. The key ``theorem'' is trivial, but crucial.

\begin{theorem}[Proper Pushforward Preserves Closedness]
\label{thm:proper-closed}
Let $f: K_1 \to K_2$ be a horizontal morphism and $R \subseteq K_1 \times K_3$ be
a vertical morphism. Then $f_!R$ is a vertical morphism (i.e., it is
closed in $K_2 \times K_3$).
\end{theorem}

\begin{proof}
Since $K_1$ and $K_3$ are compact and $R$ is closed, $R$ is compact. The map $f \times \mathrm{id}_{K_3}: K_1 \times K_3 \to K_2 \times K_3$ is continuous, so
\[
f_!R = (f \times \mathrm{id}_{K_3})(R)
\]
is compact. In the Hausdorff space $K_2 \times K_3$, compact sets are closed.
\end{proof}

\begin{remark}[Why This Matters in Practice]
\label{rmk:properness-operational}
This result deserves emphasis beyond its brevity.  In portfolio construction, constraint specifications are built from weak inequalities ($\le$, $\ge$) and equalities --- half-spaces, hyperplanes, norm balls --- whose intersections with the simplex are always closed.  The theorem says: if your inputs (constraint specifications) are closed, your outputs (downstream portfolio sets) are automatically closed.  No phantom portfolios can appear.  No special care is needed at the boundary.  The compactness of the simplex does the work silently.  This is not a mathematical convenience but a \emph{design requirement}: the entire coherence theory of Sections~\ref{sec:Beck--Chevalley}--\ref{sec:frobenius} rests on it, and Section~\ref{sec:failure} shows precisely what breaks when it is violated.
\end{remark}

\subsection{Functoriality and Adjunction}

\begin{proposition}
\label{prop:functoriality}
Pushforward is covariant functorial: $(g \circ f)_! = g_! \circ f_!$.
Pullback is contravariant functorial: $(g \circ f)^* = f^* \circ g^*$.
\end{proposition}

\begin{proof}
These follow from the definitions.
For pullback:
$(x, z) \in (g \circ f)^*S \iff ((g \circ f)(x), z) \in S$
$\iff (g(f(x)), z) \in S \iff (f(x), z) \in g^*S$
$\iff (x, z) \in f^*(g^*S) \iff (x, z) \in (f^* \circ g^*)S$.

For pushforward:
$(y, z) \in (g \circ f)_!R \iff \exists x \in K_1: (g \circ f)(x) = y, (x, z) \in R$
$\iff \exists x \in K_1: g(f(x)) = y, (x, z) \in R$
$\iff \exists x \in K_1, w \in K_2: w = f(x), g(w) = y, (x, z) \in R$
$\iff \exists w \in K_2: g(w) = y \text{ and } (\exists x \in K_1: w = f(x), (x, z) \in R)$
$\iff \exists w \in K_2: g(w) = y, (w, z) \in f_!R$
$\iff (y, z) \in g_!(f_!R)$.
\end{proof}

\begin{theorem}[Adjunction $f_! \dashv f^*$]
\label{thm:adjunction}
For \textit{any} $f: K_1 \to K_2$, $R \subseteq K_1 \times K_3$, and $S \subseteq K_2
\times K_3$:
\[
R \subseteq f^*S \iff f_!R \subseteq S
\]
\end{theorem}

\begin{proof}
\textbf{($\Rightarrow$)} Suppose $R \subseteq f^*S$. Let $(y, z) \in f_!R$.
Then there exists $x \in K_1$ with $f(x) = y$ and $(x, z) \in R$. Since $R
\subseteq f^*S$, we have $(x, z) \in f^*S$, which means $(f(x), z) \in S$.
Therefore $(y, z) = (f(x), z) \in S$. Thus $f_!R \subseteq S$.

\textbf{($\Leftarrow$)} Suppose $f_!R \subseteq S$. Let $(x, z) \in R$.
Since $x \in K_1 = \dom(f)$, we have $f(x) \in K_2$ and $(f(x), z) \in f_!R$ by
definition of pushforward.
Since $f_!R \subseteq S$, we have $(f(x), z) \in S$.
By definition of pullback, this means $(x, z) \in f^*S$. Thus $R \subseteq f^*S$.
\end{proof}

\begin{remark}
In order-theoretic language, this adjunction constitutes a \emph{Galois connection} between the posets of alignment relations on $K_1$ and  $K_2$. In particular, this guarantees that $f^*$ preserves intersections (meets) and $f_!$ preserves unions (joins).
\end{remark}

\begin{remark}
The fixed points of the closure operator $c=f^*\circ f_!$ are alignment relations on the hub space that are ``spoke-compatible'', i.e. do not rely on any information that is lost during re-implementation via $f$. The fixed points of the interior operator $i=f_!\circ f^*$ are alignment relations on the spoke space that are ``fully reachable'', i.e. every spoke portfolio allowed by $S$ is the image of a hub portfolio under $f$.
\end{remark}

\subsection{Invalid Maps and Objects}
\label{sec:non-proper-examples}

\begin{example}[Invalid Object: Ratio Constraints]
\par
Consider a constraint requiring that asset~$i$ be held at least $c$ times the weight of asset~$j$: the set $\{x \in \Delta^n \mid x_i / x_j \ge c\}$ is \textbf{not closed}, because sequences with $x_j \to 0^+$ and $x_i / x_j \ge c$ can converge to boundary points where $x_j = 0$ and $x_i = 0$, which are not in the set (or are, depending on convention---precisely the ambiguity that makes it a bad definition). The correct formulation is the \emph{linear} constraint $\{x \in \Delta^n \mid x_i \ge c \cdot x_j\}$, which \emph{is} closed (an intersection of a half-space with $\Delta^n$) and is therefore a valid object in $\HSSimp$.
\textbf{Portfolio problem}: The open version admits sequences of portfolios whose limit violates the constraint---the phantom-portfolio pathology again. The lesson is simple: ratio constraints should always be reformulated as linear inequalities before being encoded as permissible spaces.
\end{example}

\clearpage  \section{Beck--Chevalley Condition}
\label{sec:Beck--Chevalley}

The Beck--Chevalley condition is, in a sense, the theorem the entire framework exists to prove. Everything else---the closed objects, the properness, the adjunction---is scaffolding for this result. The question it answers is simple: if I check compliance at the hub level, do I need to check again at the spoke level? If the answer is no, we have path independence, and the audit trail is sound.

\paragraph*{Guide to this section.}  This section presents two results.  The \textbf{lax} Beck--Chevalley inclusion (Section~\ref{sub:lax-bc}) holds for \emph{every} commuting square in $\HSSimp$, with no extra hypotheses.  It is a one-way containment --- upstream checking is \emph{conservative}: it may reject portfolios that would pass downstream, but never approves one that would fail.  This is the operational workhorse; it is what makes hub-level compliance checking safe.  The \textbf{strict} equality (Section~\ref{sub:strict-bc}) holds only when the square is \emph{pointwise cartesian} --- a condition that requires each fiber of the top map to surject onto the corresponding fiber of the bottom map.  This is a restrictive condition, and in multi-stage optimization pipelines it may fail: the optimizer at the hub level may not explore the full space of spokes.  The strict result is a bonus, not the main event.

In this section we state and prove a Beck--Chevalley (BC) law that uses only
objects and morphisms \emph{already present} in the double category $\HSSimp$.
We do \emph{not} form new fiber-product objects. Instead, we assume we are given
a commuting square of existing horizontal morphisms and impose a pointwise
cartesian condition on that square.

\subsection{Lax Beck--Chevalley}
\label{sub:lax-bc}

The full Beck--Chevalley law holds only under the restrictive pointwise cartesianness condition discussed above. Fortunately, the lax half is completely general, and it is this half that matters in practice. We call it the ``Audit Safety Inequality'':

\begin{proposition}[Lax Beck--Chevalley Inclusion]
\label{prop:lax-BC}
Let the following be a commuting square of horizontal morphisms (continuous maps between compact permissible spaces) in $\HSSimp$:
\[
\begin{tikzcd}
K_A \arrow[r, "g"] \arrow[d, "f'"'] & K_B \arrow[d, "f"] \\
K_C \arrow[r, "h"'] & K_D
\end{tikzcd}
\]
satisfying $f \circ g = h \circ f'$. Let $Z$ be a further permissible space and let $R \subseteq K_B \times Z$ be a closed vertical morphism (alignment relation). Then:
\[
f'_! (g^* R) \subseteq h^* (f_! R).
\]
\end{proposition}

\begin{proof}
Let $(y, z) \in f'_! (g^* R)$. By the definition of pushforward, there exists $x \in K_A$ such that:
\begin{enumerate}
   \item $f'(x) = y$, and
   \item $(x, z) \in g^* R$.
\end{enumerate}
By the definition of pullback, condition (2) implies:
\[
(g(x), z) \in R.
\]
Let $w = g(x) \in K_B$. Then $(w, z) \in R$.
We now map $w$ to $K_D$ via $f$. By the commutativity of the square ($f \circ g = h \circ f'$):
\[
f(w) = f(g(x)) = h(f'(x)) = h(y).
\]
Since $(w, z) \in R$ and $f(w) = h(y)$, it follows by the definition of pushforward that:
\[
(h(y), z) \in f_! R.
\]
Finally, by the definition of pullback along $h$, this implies:
\[
(y, z) \in h^* (f_! R).
\]
Thus, $f'_! (g^* R) \subseteq h^* (f_! R)$.
\end{proof}

\subsubsection*{Financial Motivation: The ``Safety'' of hub-side Verification}

For example, consider a commuting square of portfolio re-implementations where $g: K_{Hub} \to K_{Agg}$ aggregates a detailed portfolio into a summary view (e.g., sectors), and $f': K_{Hub} \to K_{Spoke}$ re-implements the portfolio into a spoke universe (e.g., ETFs). Let $R$ be an alignment constraint defined on the aggregate space $K_{Agg}$.

The lax Beck--Chevalley condition, given by the inclusion
\[
    f'_{!}(g^{*}R) \subseteq h^{*}(f_{!}R),
\]
admits a portfolio management interpretation regarding the timing of alignment verification:

\begin{itemize}
    \item \textbf{LHS (Early Verification):} The set $f'_{!}(g^{*}R)$ represents the universe of spoke portfolios generated by filtering the \textit{Hub} portfolios first. We identify all $x \in K_{Hub}$ that satisfy the sector limit $R$ (via $g^*$), and then re-implement them ($f'_{!}$).
    
    \item \textbf{RHS (Late Verification):} The set $h^{*}(f_{!}R)$ represents the universe of spoke portfolios that appear compliant after the fact. These are portfolios in $K_{Spoke}$ whose aggregate characteristics align with \textit{some} valid sector profile in $K_{Agg}$.
    
    \item \textbf{The Inclusion as Safety:} The inclusion asserts that \textit{hub-side verification is conservative}. Any portfolio generated from an aligned hub will satisfy the downstream alignment check.
\end{itemize}

\subsection{Pointwise Cartesianness}

The reverse inclusion only holds under a rather restrictive condition. Suppose we have a commuting square
in $\HSSimp$ as above.

\begin{definition}[Pointwise cartesian square]
\label{def:pointwise-cartesian}
The commuting square is \emph{pointwise cartesian} if for every
pair $(y,z)\in K_{B}\times K_{C}$ with $f(y)=h(z)$ there \emph{exists}
$x\in K_{A}$ such that $g(x)=y$ and $f'(x)=z$.
\end{definition}

\begin{remark}
Definition~\ref{def:pointwise-cartesian} refers only to the four given objects
and arrows. It does not introduce any new object such as a fiber product, which would not exist in $\HSSimp$. It is equivalent to the existence, for every pair $(y, z)$ with $f(y) = h(z)$, of a hub $x$ such that $g(x) = y$ and $f'(x) = z$. This is a covering property for the set-theoretic pullback, and ensures that the square commutes up to surjectivity on points.

It is worth being frank about what this condition demands. Pointwise cartesianness is not something one typically checks directly; rather, it is something one hopes to verify by recognizing the square as belonging to a class for which it holds automatically (linear aggregation chains, for instance). In practice, many multi-stage pipelines involve optimization steps whose outputs depend on solver state, numerical tolerances, or ordering effects---and for these, the condition may simply fail. The lax Beck--Chevalley inclusion remains valid regardless, which is why the paper emphasizes it separately; but the full equality, which is what path independence really requires, should be treated as a property to be verified, not assumed.
\end{remark}

\subsubsection*{When does the condition hold?}

\paragraph{Example 1: Linear Aggregation Chains}
Stock $\to$ Sector $\to$ Asset Class mappings are total, continuous maps on closed simplices.
Every consistent pair of sector and asset class portfolios arises from a stock portfolio.
\par
\begin{itemize}
 \item $g$: aggregate stocks into sectors
 \item $f'$: map sector weights to asset class exposures
\end{itemize}
Therefore, the square is pointwise cartesian and Beck--Chevalley holds (path independence).

\paragraph{Example 2: Proper Optimization-Based Re-implementations}
Consider a commutative square of vertical morphisms, i.e. alignment relations:
\[
\begin{tikzcd}
K_{n} \arrow[r, "S", dashed] \arrow[d, "R'"', dashed] & K_{p} \arrow[d, "R", dashed] \\ 
K_{m} \arrow[r, "S'"', dashed] & K_{q}
\end{tikzcd}
\]
Suppose that in the context of Section~\ref{subsec:opt-align} we have constructed horizontal morphisms
\[
\begin{tikzcd}
K_{n} \arrow[r, "g"] \arrow[d, "f'"'] & K_{p} \arrow[d, "f"] \\ 
K_{m} \arrow[r, "h"'] & K_{q}
\end{tikzcd}
\]
where the horizontal morphisms are re-implementations that respect the corresponding
alignment constraints and are the unique \emph{optimal} re-implementations among those that do (see Section~\ref{subsec:opt-align} for the setup). \emph{If this square commutes}, then the restricted square
\[
\begin{tikzcd}
K_{n} \arrow[r, "g"] \arrow[d, "f'"'] & \Image{K_{n}} \arrow[d, "f"] \\ 
\Image{K_{n}} \arrow[r, "h"'] & K_{q}
\end{tikzcd}
\]
is pointwise cartesian. (Note that by properness, the images are closed and hence well-defined objects of $\HSSimp$.)

Unfortunately, we will see in Section~\ref{subsec:coherence} that the square of re-implementations need not, in general, commute---the result of successive optimizations may be path-dependent. We must impose additional assumptions to ensure that the square commutes, and hence that the restricted square is pointwise cartesian.

\begin{figure}[h]
    \centering
    \begin{tikzpicture}[node distance=3.5cm, auto,
        obj/.style={draw, ellipse, minimum width=2.2cm, minimum height=3.2cm, fill=gray!10},
        pt/.style={circle, fill, inner sep=1.5pt}]

        \node (Kn) [obj, label=above:$K_A$] {};
        \node (Kp) [obj, label=above:$K_B$, below right=2cm and 4cm of Kn] {};
        \node (Km) [obj, label=below:$K_C$, below left=2cm and 4cm of Kn] {};
        \node (Kq) [obj, label=below:$K_D$, below right=2cm and 4cm of Km] {};

        \node (Hole) [draw, ellipse, minimum width=0.8cm, minimum height=1.2cm, fill=white] at (Kn.center) {};

        \node (y) [pt, label=right:$y$] at ([xshift=-0.4cm]Kp.center) {};
        \node (z) [pt, label=left:$z$] at ([xshift=0.4cm]Km.center) {};
        \node (w) [pt, label=above:$w$] at ([yshift=0.6cm]Kq.center) {};

        \draw [->, thick] (Kn) -- node[above, sloped, pos=0.4] {$g$} (Kp);
        \draw [->, thick] (Kn) -- node[above, sloped, pos=0.4] {$f'$} (Km);
        \draw [->, thick] (Kp) -- node[below, sloped, pos=0.4] {$f$} (Kq);
        \draw [->, thick] (Km) -- node[below, sloped, pos=0.4] {$h$} (Kq);

        \draw [->, dashed, shorten >=1pt] (y) -- (w);
        \draw [->, dashed, shorten >=1pt] (z) -- (w);
        \node [text width=3cm, align=center] at (Kq.center) {Consistent Pair: \\ $f(y) = h(z) = w$};

        \draw [dashed, opacity=0.7] (y) -- node[below, sloped, pos=0.3] {$g^{-1}(y)$} (Hole.east);
        \draw [dashed, opacity=0.7] (z) -- node[below, sloped, pos=0.3] {$(f')^{-1}(z)$} (Hole.west);

        \node [red, text width=3cm, align=center, font=\bfseries] at (Kn.south) {Failure --- \\ No $x \in K_A$ exists such that $g(x)=y$ and $f'(x)=z$.};
    \end{tikzpicture}
    \caption{Illustration of a commuting square that fails the pointwise cartesian condition (Definition~\ref{def:pointwise-cartesian}). A consistent pair $(y, z)$ exists (since $f(y) = h(z)$), but the required hub portfolio $x$ would lie in a ``hole'' (a non-permissible region) of the object $K_A$. }
    \label{fig:cartesian-fail-revised}
\end{figure}
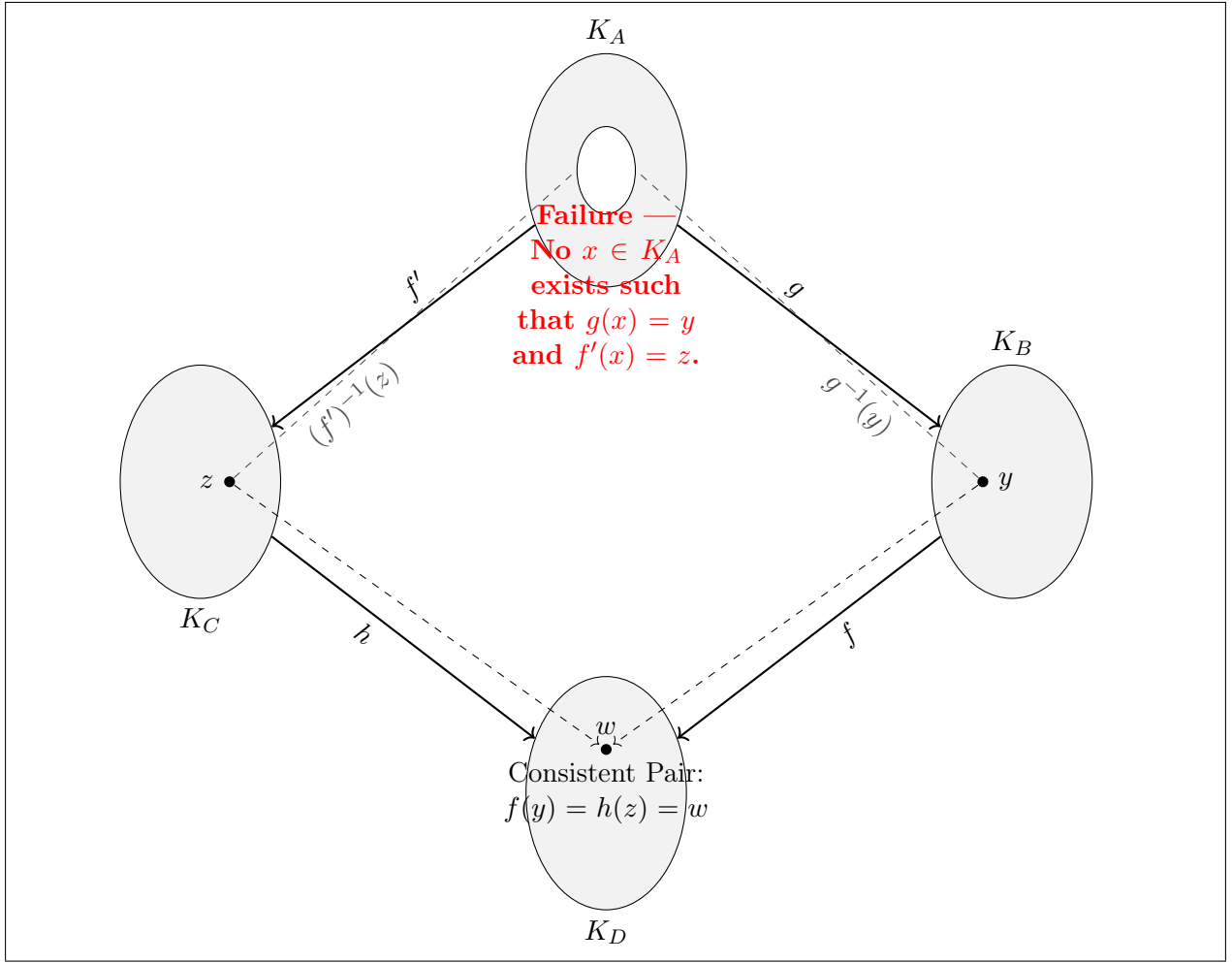

\subsection{Strict Beck--Chevalley}
\label{sub:strict-bc}

Recall the operators $f_{!}$ (pushforward) and $f^{*}$ (pullback) on closed
relations from Section~\ref{sec:operations}. Properness of horizontal maps
(Section~\ref{sec:proper}) ensures closedness is preserved.

\begin{theorem}[Beck--Chevalley for existing pointwise cartesian squares]
\label{thm:Beck--Chevalley}
Let
\[
\begin{tikzcd}
K_{A} \arrow[r, "g"] \arrow[d, "f'"'] & K_{B} \arrow[d, "f"] \\ 
K_{C} \arrow[r, "h"'] & K_{D}
\end{tikzcd}
\]
be a commuting square of existing objects and horizontal morphisms in $\HSSimp$
that is pointwise cartesian. Let $Z$ be a further permissible space. Then for every closed relation
$R\subseteq K_{B}\times Z$,
\[
f'_{!}\big(g^{*}R\big) \;=\; h^{*}\big(f_{!}R\big).
\]
\end{theorem}

\begin{proof}
We show both inclusions.

\emph{($\subseteq$)} Let $(y,z)\in f'_{!}(g^{*}R)$. Then there exists
$x\in K_{A}$ with $f'(x)=y$ and $(x,z)\in g^{*}R$, i.e., $(g(x),z)\in R$.
Set $y':=g(x)\in K_{B}$. Commutativity gives
$f(y')=f(g(x))=h(f'(x))=h(y)$. Since $(y',z)\in R$, we have
$(f(y'),z)\in f_{!}R$, hence $(h(y),z)\in f_{!}R$. By definition of pullback,
this is equivalent to $(y,z)\in h^{*}(f_{!}R)$.

\emph{($\supseteq$)} Let $(y,z)\in h^{*}(f_{!}R)$. Then
$(h(y),z)\in f_{!}R$, so there exists $y'\in K_{B}$ with
$f(y')=h(y)$ and $(y',z)\in R$. By pointwise cartesianness, there exists
$x\in K_{A}$ with $g(x)=y'$ and $f'(x)=y$. From $(y',z)\in R$ and $g(x)=y'$ we
get $(x,z)\in g^{*}R$, and from $f'(x)=y$ we conclude
$(y,z)\in f'_{!}(g^{*}R)$.
\end{proof}

The following restates Theorem~\ref{thm:Beck--Chevalley} with an alternative proof via the mates identity, which clarifies the categorical structure.

\begin{theorem}[Beck--Chevalley via Mates Identity]
\label{cor:pc-implies-cart}
If the square is pointwise cartesian in the sense of
Definition~\ref{def:pointwise-cartesian}, then for every permissible space $Z$ and every alignment relation $R\subseteq K_{B}\times Z$ we have
\[
f'_{!}(g^{*}R) \;=\; h^{*}(f_{!}R).
\]
\end{theorem}

  \begin{proof}
  We prove the equivalent \emph{mates identity}: as relations from $K_C$ to $K_B$,
  \begin{equation}\label{eq:mates}
  \Graph(g) \circ \Graph(f')^{\dagger} \;=\; \Graph(f)^{\dagger} \circ \Graph(h),
  \end{equation}
  where $\Graph(k)^{\dagger} = \{(k(x),x)\}$ denotes the reverse graph (conjoint).
  Theorem~\ref{thm:Beck--Chevalley} is then recovered from~\eqref{eq:mates}
  by left-composing with an arbitrary alignment relation $S\subseteq K_B\times Z$.

  \emph{($\subseteq$)}
  Let $(y,y')\in \Graph(g)\circ\Graph(f')^{\dagger}$.
  By definition, there exists $x\in K_A$ with
  $(y,x)\in\Graph(f')^{\dagger}$ and $(x,y')\in\Graph(g)$,
  i.e., $f'(x)=y$ and $g(x)=y'$.
  By commutativity of the square, $f(y')=f(g(x))=h(f'(x))=h(y)$,
  so $(y,y')\in\Graph(f)^{\dagger}\circ\Graph(h)$.

  \emph{($\supseteq$)}
  Let $(y,y')\in \Graph(f)^{\dagger}\circ\Graph(h)$.
  Then $(y,h(y))\in\Graph(h)$ and $(h(y),y')\in\Graph(f)^{\dagger}$,
  so $f(y')=h(y)$.
  By pointwise cartesianness
  (Definition~\ref{def:pointwise-cartesian}),
  there exists $x\in K_A$ with $g(x)=y'$ and $f'(x)=y$.
  Hence $(y,y')\in \Graph(g)\circ\Graph(f')^{\dagger}$.

  \medskip
  \noindent\emph{Recovery of Theorem~\ref{thm:Beck--Chevalley}.}
  Since pushforward and pullback are expressible as relational composition
  with (reverse) graphs---$f'_{!}(T) = T\circ\Graph(f')^{\dagger}$
  and $h^{*}(T) = T\circ\Graph(h)$---the
  Beck--Chevalley equality $f'_{!}(g^{*}S) = h^{*}(f_{!}S)$ for
  $S\subseteq K_B\times Z$ follows by left-composing~\eqref{eq:mates}
  with $S$:
  \[
  f'_{!}(g^{*}S)
  = S\circ\Graph(g)\circ\Graph(f')^{\dagger}
  \;=\; S\circ\Graph(f)^{\dagger}\circ\Graph(h)
  = h^{*}(f_{!}S).  \qedhere
  \]
  \end{proof}

\subsubsection*{Portfolio Interpretation}

The Beck--Chevalley equality
\[
f'_{!}\circ g^{*} \;=\; h^{*}\circ f_{!}
\]
formalizes \emph{path independence} for alignment propagation across an existing
two-stage pipeline. Properness guarantees closedness at each step,
and pointwise cartesianness guarantees that the two paths capture exactly the
same portfolios.

\clearpage  \section{Frobenius Reciprocity}
\label{sec:frobenius}

Frobenius reciprocity, in the classical representation-theoretic setting, says that induction and restriction are adjoint---that filtering before or after a change of representation gives the same answer. In our setting the content is more concrete: if you intersect an alignment relation with a pulled-back constraint and then push forward, you get the same result as pushing forward first and intersecting afterward. The practical force of this is that a portfolio manager can filter at the hub level or at the spoke level and arrive at the same compliant universe. This is the kind of result that, once stated, feels inevitable; the point is that it \emph{is} inevitable, given the framework, and does not require a separate proof strategy.

\subsection{Statement}

\begin{theorem}
\label{thm:frobenius}
Let $f: K_1 \to K_2$ be a horizontal morphism in $\HSSimp$. For closed
alignment relations $R \subseteq K_1 \times K_3$ and $S \subseteq K_2 \times
K_3$:
\[
f_!(R \cap f^*S) = f_!R \cap S
\]
\end{theorem}

\begin{proof}
(Note: $f$ is total on its domain $K_1$ by definition.)

\textbf{Compute LHS}: $f_!(R \cap f^*S)$
\begin{align*}
R \cap f^*S &= \{(x, z) \in K_1 \times K_3 \mid (x, z) \in R
\text{ and } (f(x), z) \in S\}
\\
f_!(R \cap f^*S) &= \{(y, z) \in K_2 \times K_3 \mid \exists x \in K_1:
f(x) = y, (x, z) \in R, (f(x), z) \in S\}
\end{align*}

\textbf{Compute RHS}: $f_!R \cap S$
\begin{align*}
f_!R &= \{(y, z) \in K_2 \times K_3 \mid \exists x \in K_1: f(x) = y, (x,
z) \in R\}
\\
f_!R \cap S &= \{(y, z) \mid (\exists x \in K_1: f(x) = y, (x, z) \in R) \text{ and
} (y, z) \in S\}
\end{align*}

\textbf{Show LHS $\subseteq$ RHS}:
Suppose $(y, z) \in f_!(R \cap f^*S)$. Then there exists $x \in K_1$ with $f(x) =
y$, $(x, z) \in R$, and $(f(x), z) \in S$.
Therefore $(y, z) \in S$ and $(y, z) \in f_!R$. Thus $(y, z) \in f_!R \cap S$.

\textbf{Show RHS $\subseteq$ LHS}:
Suppose $(y, z) \in f_!R \cap S$. Then:
\par
\begin{itemize}
\item There exists $x \in K_1$ with $f(x) = y$ and $(x, z) \in R$
\item $(y, z) \in S$
\end{itemize}
Since $f(x) = y$, we have $(f(x), z) = (y, z) \in S$.
Therefore $(x, z) \in R$ and $(f(x), z) \in S$, which means $(x, z) \in R \cap
f^*S$.
Since $f(x) = y$, we have $(y, z) \in f_!(R \cap f^*S)$.

\textbf{Closedness}: Both sides are closed by Theorem~\ref{thm:proper-closed}
(since $f$ is proper) and closure under intersection.
\end{proof}

\begin{figure}
\centering
\begin{tikzpicture}[scale=1, >=stealth]
    
    \node[draw, rectangle, minimum width=3.5cm, minimum height=3.5cm] (hub1) at (0,4) {};
    \node[anchor=north] at (hub1.north) [yshift=-2pt] {\textbf{Hub} $K_1$};
    
    \fill[blue!20] (-0.4, 3.8) circle (0.8); 
    \node[blue, font=\small] at (-1.0, 3.8) {$R$};
    
    \fill[red!20, opacity=0.7] (0.4, 3.8) circle (0.8); 
    \node[red, font=\small] at (1.0, 3.8) {$f^*S$};
    
    \begin{scope}
        \clip (-0.4, 3.8) circle (0.8);
        \fill[purple!60] (0.4, 3.8) circle (0.8);
    \end{scope}
    \node[white, font=\tiny] at (0, 3.8) {$\cap$};
    
    \node[draw, rectangle, minimum width=3.5cm, minimum height=3.5cm] (spoke1) at (7,4) {};
    \node[anchor=north] at (spoke1.north) [yshift=-2pt] {\textbf{Spoke} $K_2$};

    \begin{scope}
        \clip (6.6, 3.8) circle (0.8);
        \fill[purple!60] (7.4, 3.8) circle (0.8);
    \end{scope}
    \node[purple!40!black, font=\small, right] (Res1) at (7.5, 3.8) {Result 1};

    \draw[->, thick] (hub1.east) -- (spoke1.west) node[midway, above] {$f_!$ (Push)};

    
    \node[draw, rectangle, minimum width=3.5cm, minimum height=3.5cm] (hub2) at (0,0) {};
    \node[anchor=north] at (hub2.north) [yshift=-2pt] {\textbf{Hub} $K_1$};

    \fill[blue!20] (0, -0.2) circle (0.8);
    \node[blue, font=\small] at (0, -0.2) {$R$};
    
    \node[draw, rectangle, minimum width=3.5cm, minimum height=3.5cm] (spoke2) at (7,0) {};
    \node[anchor=north] at (spoke2.north) [yshift=-2pt] {\textbf{Spoke} $K_2$};

    \fill[blue!20] (6.6, -0.2) circle (0.8);
    \node[blue, font=\small] at (5.9, -0.2) {$f_!R$};
    
    \fill[red!20, opacity=0.7] (7.4, -0.2) circle (0.8);
    \node[red, font=\small] at (8.1, -0.2) {$S$};
    
    \begin{scope}
        \clip (6.6, -0.2) circle (0.8);
        \fill[purple!60] (7.4, -0.2) circle (0.8);
    \end{scope}
    \node[purple!40!black, font=\small, right] (Res2) at (7.8, -1.2) {Result 2};

    \draw[->, thick] (hub2.east) -- (spoke2.west) node[midway, above] {$f_!$ (Push)};

    
    \draw[dashed, gray] (8.8, 3.8) -- (9.5, 3.8); 
    \draw[dashed, gray] (8.8, -0.2) -- (9.5, -0.2); 
    
    \draw[<->, thick, green] (9.5, 3.8) -- (9.5, -0.2) 
        node[midway, right, align=left] {Sets are\\Identical};

\end{tikzpicture}
\caption{Visualization of Frobenius reciprocity}
\label{fig:frobenius}
\end{figure}
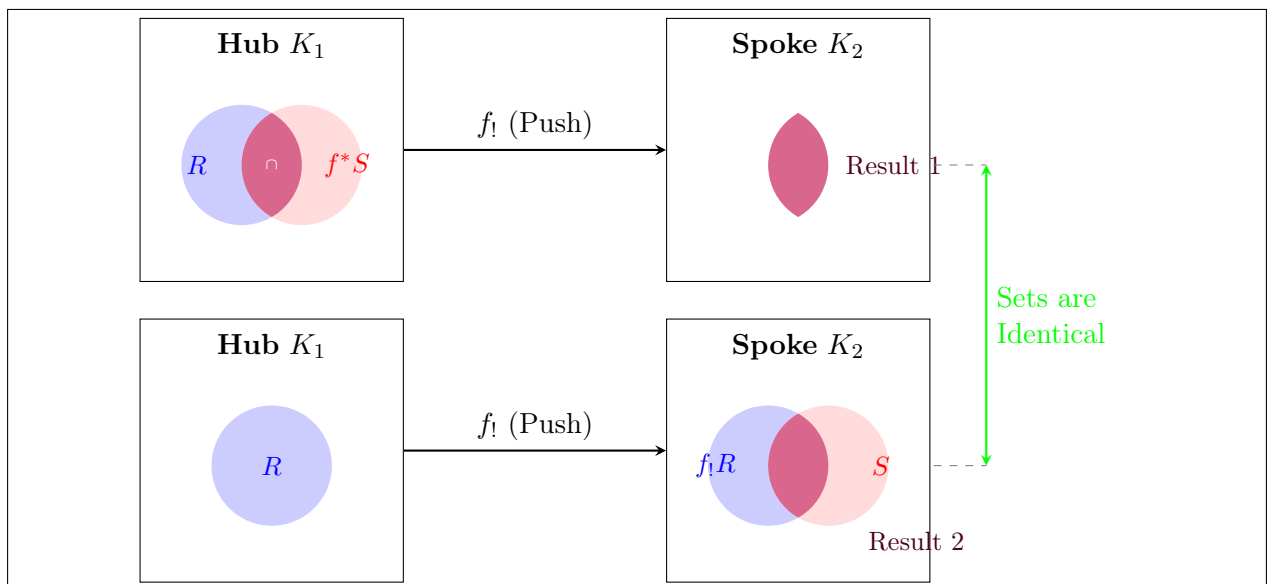

\subsection{Portfolio Interpretation}

Frobenius reciprocity guarantees that filtering portfolios by both a
hub-side alignment $R$ and a spoke-side alignment $S$ produces the same
set of permissible spoke portfolios whether the filtering is done before or after
re-implementation.

\begin{center}
\begin{tikzcd}[column sep=huge, row sep=huge]
K_1 \arrow[r, "f"] \arrow[d, dashed, "R \cap f^{*}S"'] & 
K_2 \arrow[d, dashed, "f_{!}R \cap S"] \\ 
K_3 \arrow[r, "\mathrm{id}_{K_3}"'] & K_3
\end{tikzcd}
\end{center}
(Here $R \subseteq K_1 \times K_3$ and $S \subseteq K_2 \times K_3$ for some $K_3$, and $f: K_1 \to K_2$)

\begin{example}[Customization with Multiple Constraints]
\label{ex:esg-benchmark}
\qquad\par
\begin{itemize}
\item hub: $K_1$ (Permissible Mutual Fund portfolios)
\item spoke: $K_2$ (Permissible ETF baskets)
\item $f: K_1 \to K_2$: Optimal re-implementation
\item $R = $ ``hub portfolio is diversified'': $\{(x, \checkmark) \mid
\text{HHI}(x) \le H_0\}\subseteq K_1 \times \{\checkmark\}$
\item $S = $ ``spoke ETF basket is low-fee'': $\{(y, \checkmark) \mid \text{Fee}(y)
\le 0.1\%\}\subseteq K_2 \times \{\checkmark\}$
\end{itemize}

(This uses Remark~\ref{rem:compliance} to interpret the constraints/eligibility/compliance conditions as alignment relations $R$ and $S$.)

Which ETF baskets in $K_2$ are implementations of diversified hub portfolios, \emph{and} low-fee? We can construct the desired portfolios in two ways:

\begin{enumerate}
\item Find hub portfolios ($x \in K_1$) that are
diversified ($R$) and then re-implement to low-fee ETFs ($f^*S$). Implement them via $f$.

\item Implement all diversified hub portfolios ($f_!R$).
Filter the resulting ETF baskets for low fees ($S$).
\end{enumerate}

Frobenius guarantees identical answers.

\end{example}

\begin{remark}
The elegance of Frobenius reciprocity in this setting---the proof is essentially two lines of set algebra---may obscure a genuine operational difficulty. The theorem assumes that the constraints $R$ and $S$ are both expressible as closed alignment relations. In practice, many constraints (turnover limits, for instance, or constraints that depend on the current state of the portfolio rather than on the target allocation alone) are not naturally of this form. Frobenius tells you that \emph{if} you can express your constraints as closed relations, the order of filtering is irrelevant; it does not tell you how to express them so.
\end{remark}

\clearpage  \section{Why Closedness Is Not Optional}
\label{sec:failure}

The closedness requirement on objects and vertical morphisms in $\HSSimp$ is a design principle, not a mathematical convenience.  Every coherence result in this paper --- properness (Theorem~\ref{thm:proper-closed}), adjunction (Theorem~\ref{thm:adjunction}), Beck--Chevalley (Proposition~\ref{prop:lax-BC}), Frobenius (Theorem~\ref{thm:frobenius}) --- depends on it.  This section demonstrates that the requirement cannot be relaxed, and that the natural attempt to repair the damage (taking topological closures of pushforwards) introduces new pathologies.  It then shows that the same principle extends to the probabilistic setting of Part~III: closedness is essential not only for the deterministic framework but for every extension developed in this paper.

\subsection{Loss of Properness, and an Unsuccessful Solution}

If $K_1$ is not closed in $\Delta^n$, it is not compact. A continuous map $f:
K_1 \to K_2$ from a non-compact space $K_1$ is \emph{not} guaranteed to be
proper.

As shown in Example~\ref{ex:pushforward-not-closed}, a non-proper map $f$
(like the inclusion $f: [0, 1) \to [0, 1]$) fails to preserve closedness under
pushforward. $f_!R$ will not be closed.

If we relax the requirement for objects to be closed, the pushforward $f_!R$ may fail to be closed (as shown in Example~\ref{ex:pushforward-not-closed}). A naive fix might be to force closedness by defining a modified pushforward $\overline{f_!}R = \cl(f_!R)$. This does preserve composition, adjunction and lax Beck--Chevalley; but it breaks strict Beck--Chevalley and, more importantly, Frobenius.

\begin{theorem}[Stability of the Closure Fix]
Let $f: X \to Y$ and $g: Y \to Z$ be continuous maps between topological spaces (not necessarily compact).
Let vertical morphisms be \textbf{closed relations}.
Define the ``patched'' pushforward as $\overline{f_!}R := \operatorname{cl}(f(R))$.
This operation preserves the following structural laws:

\begin{enumerate}
    \item \textbf{Compositionality (Functoriality):}
    \[ \overline{(g \circ f)_!} R = \overline{g_!} (\overline{f_!} R) \]
    \begin{proof}
    For any continuous map $g$, we have $g(\operatorname{cl}(A)) \subseteq \operatorname{cl}(g(A))$: indeed, $g^{-1}(\operatorname{cl}(g(A)))$ is closed (by continuity) and contains $A$ (since $g(A) \subseteq \operatorname{cl}(g(A))$), hence contains $\operatorname{cl}(A)$.
    Consequently, $g(A) \subseteq g(\operatorname{cl}(A)) \subseteq \operatorname{cl}(g(A))$.
    Taking closures: $\operatorname{cl}(g(A)) \subseteq \operatorname{cl}(g(\operatorname{cl}(A))) \subseteq \operatorname{cl}(\operatorname{cl}(g(A))) = \operatorname{cl}(g(A))$, so $\operatorname{cl}(g(\operatorname{cl}(A))) = \operatorname{cl}(g(A))$.
    Applying this to the relation $R$:
    \begin{align*}
        \text{RHS} &= \operatorname{cl}( g( \operatorname{cl}( f(R) ) ) ) \\
        &= \operatorname{cl}( g( f(R) ) ) \quad \text{(by topological identity)} \\
        &= \operatorname{cl}( (g \circ f)(R) ) \\
        &= \overline{(g \circ f)_!} R \quad \text{(LHS)}
    \end{align*}
    \end{proof}

    \item \textbf{Adjunction ($f_! \dashv f^*$):}
    \[ \overline{f_!}R \subseteq S \iff R \subseteq f^*S \]
    \begin{proof}
    The definition of the patched pushforward is $\operatorname{cl}(f(R))$.
    The condition $\operatorname{cl}(f(R)) \subseteq S$ holds if and only if $f(R) \subseteq S$, because $S$ is a \textbf{closed set} (by definition of vertical morphisms).
    The condition $f(R) \subseteq S$ is the standard set-theoretic adjunction, equivalent to $R \subseteq f^{-1}(S) = f^*S$.
    Thus, the adjunction is preserved exactly.
    \end{proof}

    \item \textbf{Lax Beck--Chevalley (Safety):}
    \[ \overline{f'_!} (g^* R) \subseteq h^* (\overline{f_!} R) \]
    \begin{proof}
    Let $A = f'(g^* R)$ be the standard (non-closed) pushforward.
    From the standard set-theoretic logic, we know $A \subseteq h^*(f_! R)$.
    Since $f_! R \subseteq \operatorname{cl}(f_! R) = \overline{f_!} R$, monotonicity implies:
    \[ A \subseteq h^*(\overline{f_!} R) \]
    The set $h^*(\overline{f_!} R)$ is the preimage of a closed set under a continuous map $h$, so it is closed.
    If a set $A$ is contained in a closed set $C$, then its closure is also contained in $C$:
    \[ \operatorname{cl}(A) \subseteq C \]
    Substituting terms back:
    \[ \overline{f'_!} (g^* R) \subseteq h^* (\overline{f_!} R) \]
    \end{proof}
\end{enumerate}
\end{theorem}

  \begin{theorem}[Closure Fix Breaks Path Independence (Strict Beck--Chevalley)]
  \label{thm:closure-breaks-BC}
  Let $\overline{f_!}R := \cl(f_!R)$. There exists a commuting square and
  a closed relation $S$ for which:
  \[
  \overline{f'_!}(g^* S) \;\subsetneq\; h^*(\overline{f_!}\, S)
  \]
  so the strict Beck--Chevalley equality fails for the patched pushforward.
  \end{theorem}

  \begin{proof}[Counterexample]
  Consider the commuting square
  \[
  \begin{tikzcd}
  K_A \arrow[r, "g"] \arrow[d, "f'"'] & K_B \arrow[d, "f"] \\
  K_C \arrow[r, "h"'] & K_D
  \end{tikzcd}
  \]
  with the following data:
  \begin{itemize}
  \item $K_A = [0, 1)$ (non-closed, hence non-compact).
  \item $K_B = K_C = K_D = [0, 1]$.
  \item $g = f' = \iota \colon [0,1) \hookrightarrow [0,1]$
    (the inclusion), and $f = h = \operatorname{id}_{[0,1]}$.
  \item $Z = [0, 1]$, and $S = \{(1, 1)\} \subseteq K_B \times Z$ (a boundary constraint pairing the endpoint with
    itself).
  \end{itemize}
  Commutativity holds: $f \circ g = \iota = h \circ f'$.
  Note that $f' = \iota$ is \emph{not} proper: properness fails because $K_A = [0,1)$ is non-compact. (A continuous map from a compact space to a Hausdorff space is automatically proper, but $K_A$ is not compact. In particular, $\iota^{-1}([0,1]) = [0,1)$ is non-compact.)

  \medskip
  \noindent\textbf{LHS (Early Audit):}
  First pull back $S$ to the hub, then push forward.
  \[
  g^* S = \{(x, z) \in K_A \times Z : (g(x), z) \in S\}
        = \{(x, z) \in [0,1) \times [0,1] : x = 1,\; z = 1\}
        = \emptyset
  \]
  since $1 \notin [0, 1)$.  Therefore
  $\overline{f'_!}(\emptyset) = \operatorname{cl}(\emptyset) = \emptyset$.

  \medskip
  \noindent\textbf{RHS (Late Audit):}
  First push forward, then pull back.
  \[
  f_! S = \{(f(w), z) : (w, z) \in S\}
        = \{(\operatorname{id}(1), 1)\} = \{(1, 1)\}
  \]
  This is already closed, so
  $\overline{f_!}\, S = \{(1, 1)\}$.  Then:
  \[
  h^*(\overline{f_!}\, S)
  = \{(y, z) \in K_C \times Z : (h(y), z) \in \{(1,1)\}\}
  = \{(1, 1)\}
  \]

  \medskip
  \noindent\textbf{Conclusion:}
  $\overline{f'_!}(g^* S) = \emptyset \subsetneq \{(1,1)\}
  = h^*(\overline{f_!}\, S)$.
  The lax inclusion $\text{LHS} \subseteq \text{RHS}$ holds
  (consistent with the previous theorem), but strict equality fails.
  The RHS accepts a ``phantom portfolio'' at the boundary $x = 1$ that
  has no pre-image in the hub space $K_A = [0, 1)$; the LHS correctly
  identifies that no such hub exists.
  \end{proof}
  
\begin{theorem}[Closure Fix Breaks Frobenius Reciprocity]
\label{thm:closure-fails}
Let $\overline{f_!}R := \cl(f_!R)$ denote the ``patched'' pushforward for non-proper maps.
This operation violates the Frobenius reciprocity law. Specifically, there exist a map $f$ and relations $R, S$ such that:
\[
\overline{f_!}(R \cap f^*S) \neq \overline{f_!}R \cap S
\]
\end{theorem}

\begin{proof}[Counterexample]
Let $K_1 = [0, 1)$ (non-closed hub) and $K_2 = [0, 1]$ (closed spoke).
Let $f: K_1 \to K_2$ be the inclusion map $f(x) = x$.
Let $R$ be the diagonal relation on $K_1$ (``all valid hubs'') and $S = \{(1, 1)\}$ be a boundary constraint on $K_2$.

\textbf{LHS (Filter then Push):}
First, compute the pullback $f^*S = \{(x, z) \in K_1 \times K_2 \mid (f(x), z) \in S\} = \{(x, z) \mid x \in [0,1),\, f(x) = 1,\, z = 1\}$.
Since $1 \notin K_1$, the pullback is empty: $f^*S = \emptyset$.
Consequently, the intersection $R \cap f^*S$ is empty, and the pushforward is empty.
\[ \text{LHS} = \emptyset \]

\textbf{RHS (Push then Filter):}
Compute the standard pushforward $f_!R = \{(x,x) \mid x \in [0, 1)\}$.
The closure is $\overline{f_!R} = \{(x,x) \mid x \in [0, 1]\}$, which \emph{includes} the point $(1,1)$.
Intersecting with $S$:
\[ \text{RHS} = \overline{f_!R} \cap S = \{(x,x) \mid x \in [0,1]\} \cap \{(1,1)\} = \{(1,1)\} \]

\textbf{Conclusion:} $\emptyset \neq \{(1,1)\}$.
The RHS accepts a ``phantom portfolio'' at the boundary that implies a valid hub existed, while the LHS correctly identifies that no such hub exists.
\end{proof}

\subsubsection*{Portfolio Interpretation of Failures}

\textbf{Phantom portfolios}: Closure adds limit points to $f_!R$ that have no
actual pre-image in $R$. These are spoke portfolios that \emph{appear} to
satisfy propagated alignments but don't actually result from \textit{any} hub
portfolio satisfying the original alignment.

\textbf{Conclusion}: Closedness (hence compactness) of objects is the single
assumption that makes all horizontal morphisms proper, pushforward
well-defined, and the compositional coherence laws (Adjunction, BC, Frobenius)
available.

\subsection{The Insufficiency of the Closure Fix for Probabilistic Models}
\label{subsec:closure-insufficiency}

As we will see later, in production settings where re-implementation is automated via the use of imperfect optimizers, alignment must be understood in a probabilistic setting. While the ``closure fix'' $\overline{f_!}R = \operatorname{cl}(f_!R)$ can technically restore basic categorical composition in the deterministic setting, it is mathematically insufficient for the probabilistic frameworks developed in Sections \ref{sec:Polish} and \ref{sec:wasserstein}. The following discussion assumes that background.

Even if we require that horizontal morphisms are \emph{tight Feller kernels} (preventing mass from escaping to infinity, as per Definition~\ref{def:tight}), we cannot relax the requirement for vertical morphisms to be \emph{closed relations}, because the two conditions play different roles in the theory. While tightness manages the behaviour of the system at infinity, closedness is essential for managing behaviour at the \emph{boundaries}. If alignment relations $R$ are not closed (e.g., strict inequalities $x > 0$), the link between probability and geometry established in Section~\ref{sec:Polish} breaks down, leading to the following failure modes:

\paragraph*{Failure of the Portmanteau Theorem (Limit Instability)}
The mathematical engine of the probabilistic framework is the \textbf{Portmanteau Theorem} (Theorem~\ref{thm:portmanteau}). It establishes that weak convergence of probability measures (solver stability) implies the upper semi-continuity of mass on \textbf{closed sets}:
\[
\mu_n \Rightarrow \mu \implies \limsup_{n \to \infty} \mu_n(F) \le \mu(F) \quad \text{for all \textbf{closed} sets } F.
\]
\begin{itemize}
    \item \textbf{The Failure:} If the alignment relation $R$ defines a non-closed constraint set $S$ (e.g., ``Tracking Error $< 15$bps''), the function $x \mapsto P(x, S)$ is no longer guaranteed to be upper semi-continuous.
    \item \textbf{The Consequence:} We lose the ability to trust the limit. A sequence of hub portfolios $x_n \to x$ could all be ``safe'' (high probability of landing in $S$), but the limit portfolio $x$ could be ``unsafe'' (probability drops discontinuously at the boundary). This destroys the stability required for the sensitivity analysis discussed in Section~\ref{sub:wasserstein-compute}.
\end{itemize}

\paragraph*{Instability of the $\epsilon$-Support}
The definition of the \textbf{$\epsilon$-Highest Density Region} (Definition~\ref{def:epsilon-supp}) relies on level sets of the density function to define a closed support.
\begin{itemize}
    \item \textbf{The Failure:} Verification of a non-closed relation $S$ against these closed supports becomes erratic.
    \item \textbf{The Consequence:} Consider a density that ``touches'' the boundary of a strict inequality constraint. A numerically negligible shift in the distribution could flip the check from \textit{True} to \textit{False} because the boundary points are excluded from $S$. Closed relations ensure that ``touching the boundary'' is a valid, stable state.
\end{itemize}

\paragraph*{Ill-Posed Wasserstein Projections}
In the Transport model (Section~\ref{sec:wasserstein}), safety is defined by the distance to the valid set $\mathcal{P}_S$.
\[
W_1(P(x), \mathcal{P}_S) = \inf_{\nu \in \mathcal{P}_S} W_1(P(x), \nu).
\]
\begin{itemize}
    \item \textbf{The Failure:} If $S$ is not closed, the set of valid measures $\mathcal{P}_S$ is not closed in the weak topology. The infimum might not be attained by any valid distribution.
    \item \textbf{The Consequence:} The ``cure'' portfolio doesn't strictly exist (one can get arbitrarily close but never reach it). This breaks the computational algorithms in Section~\ref{sub:wasserstein-compute}, which rely on projecting samples onto a valid set. A solver trying to project onto $x > 0$ will oscillate or fail numerically, whereas projecting onto $x \ge 0$ is a standard, stable operation.
\end{itemize}

It can also be shown that if alignment relations are not required to be closed, then adjunction, lax Beck--Chevalley and Frobenius fail, for all three approaches: Safety Radius, HDR and Transport-Based.

\paragraph*{The Design Principle and Its Consequences for Part III}
These failures are not edge cases.  Every constraint specification in a real portfolio system is built from weak inequalities ($\le$, $\ge$) and equalities --- and these produce closed sets.  Strict inequalities ($<$, $>$) arise only from specification errors or ambiguous documentation.  The design principle is therefore: \emph{model constraints as closed sets, and the framework works; model them as open sets, and it doesn't.}

This principle has a deeper consequence.  When horizontal morphisms become \emph{stochastic} (Part~III), the deterministic pushforward $f_!R$ is replaced by a probabilistic operation that ``fattens'' $R$ according to a risk budget $\epsilon$.  The natural response is to keep the fattened set open --- after all, it represents an \emph{approximate} constraint.  But the analysis above shows that this is precisely the wrong move.  Open constraint sets destroy the Portmanteau-based stability arguments that underpin the probabilistic framework (Section~\ref{sec:Polish}), make Wasserstein projections ill-posed (Section~\ref{sec:wasserstein}), and break the HDR support calculus.  Part~III therefore maintains the closedness requirement throughout, even for probabilistic operations: the graded double category $\HSPr$ (Definition~\ref{def:HSP-r}) uses \emph{closed} relations as its vertical morphisms, and the three approaches to probabilistic compliance (Safety Radius, HDR, Wasserstein) are all defined so that the resulting sets remain closed.  The ``closure fix'' attempted in Section~\ref{subsec:closure-insufficiency} is not needed, because the problem it was designed to solve --- open pushforwards --- never arises when the inputs are correctly specified.

\clearpage  \section{Constructing Optimal Re-implementations}
\label{sec:optimal-reimpl}

In previous sections, horizontal morphisms $f: K_1 \to K_2$ were treated as \emph{given} continuous maps between feasible portfolio spaces. We now demonstrate a method for \emph{constructing} such morphisms, including their domains $K_1$ and codomains $K_2$, based on practical optimization goals. This grounds the abstract framework in concrete portfolio management tasks, such as finding a low-fee ETF portfolio that optimally replicates a given mutual fund portfolio. This process naturally yields the compact objects and continuous (proper) morphisms required by $\HSSimp$.

\subsection{Setup: Permissible Spaces, Quantitative Attributes, and Objectives}

Assume we start with two ambient portfolio spaces $\Delta^n(\A)$ and $\Delta^m(\mathcal{B})$ (e.g., mutual funds and ETFs). We may initially only know the baseline \textbf{permissible sets} $P_1 \subseteq \Delta^n$ and $P_2 \subseteq \Delta^m$, which must be closed subsets (and hence compact objects in $\HSSimp$). These represent the maximum possible spaces satisfying basic constraints, but the actual domain and codomain of our constructed morphism might be smaller.

We require the following components:
\par
\begin{itemize}
   \item A ``space of quantitative attributes'' $(V, \|\cdot\|_V)$, a finite-dimensional Banach space, used to compare portfolios from different universes.
   \item Continuous \textbf{attribute maps} defined on the permissible sets: $g_\A: P_1 \to V$ and $g_\mathcal{B}: P_2 \to V$. These map portfolios to their relevant characteristics (e.g., asset class exposures, factor loadings).
   \item A continuous \textbf{objective function} $u: V \to \mathbb{R}$ defined on the space of quantitative attributes. This function quantifies the desirability of a portfolio based on its metric representation (e.g., utility derived from risk/return metrics, negative fees associated with factor metrics). We assume higher values of $u$ are better.
\end{itemize}

(Note that the term ``metric'' in the above definition has nothing to do with the mathematical notion of metric space. It refers to a quantitative measurement of a financial characteristic.)

\begin{definition}[Portfolio space of quantitative attributes]
A \textbf{portfolio space of quantitative attributes} is a finite-dimensional Banach space $(V,
\|\cdot\|_V)$.

A \textbf{portfolio quantitative attributes map} is a continuous function $g: P \to V$ defined on a permissible portfolio space $P \subseteq \Delta^k$.
\end{definition}

\begin{example}[Spaces of quantitative attributes]
\qquad\par
\begin{itemize}
\item \textbf{Asset Class Exposure}: $V = \mathbb{R}^k$. $g_\A: P_1 \to V$ and $g_\mathcal{B}: P_2 \to V$ summarize asset class exposures.
\item \textbf{Factor Exposure}: $V = \mathbb{R}^k$ (risk factors). $g(x) = Fx$ is the portfolio's factor loading.
\item \textbf{Risk Contribution}: $V = \mathbb{R}^n$ where $g(x)_i$ is the risk contribution of asset $i$.
\end{itemize}
\end{example}

\begin{example}[Mutual Funds vs. ETFs]
\label{ex:fund-etf-metric}
\par
Let $P_1 \subseteq \Delta^n(\A)$ be the permissible space of mutual fund portfolios and
$P_2 \subseteq \Delta^m(\mathcal{B})$ the permissible space of ETF portfolios. Both aim to provide exposure to $k$ asset classes. $g_\A: P_1 \to \mathbb{R}^k$ and $g_\mathcal{B}: P_2 \to \mathbb{R}^k$ map portfolios to their asset class exposures in $V=\mathbb{R}^k$.
\end{example}

\begin{example}[Objective Functions on Quantitative Attributes]
\qquad\par
\begin{itemize}
\item \textbf{Metric-Based Utility}: If $v \in V$ represents expected return and variance metrics, $u(v) = v_{\text{return}} - \lambda v_{\text{variance}}$.
\item \textbf{Negative Cost associated with quantitative attributes}: If $v$ represents factor exposures, $u(v) = -\text{Cost}(v)$, where cost might relate to turnover or taxes implied by achieving those exposures.
\end{itemize}

Note: If the primary objective (like fees) depends directly on the spoke portfolio $y \in P_2$ rather than its metric $g_\mathcal{B}(y)$, the setup can be adapted using a composite objective in the optimization below, but here we focus on $u: V \to \mathbb{R}$.
\end{example}

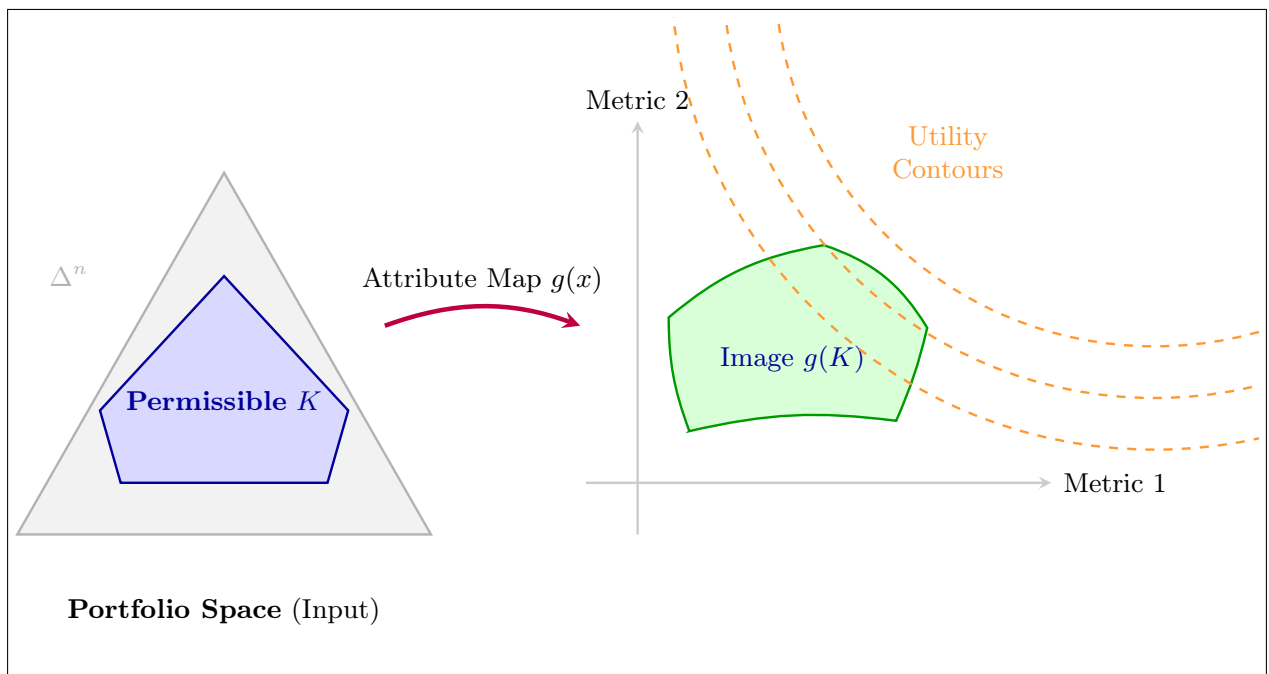
\begin{figure}
    \centering
\resizebox{\linewidth}{!}{%
    \begin{tikzpicture}[
        scale=1.2,
        >=stealth,
        font=\footnotesize, 
        simplex/.style={fill=gray!10, draw=gray!60, thick},
        permissible/.style={fill=blue!15, draw=blue!60!black, thick},
        metric_space/.style={draw=gray!40, ->, thick},
        map_arrow/.style={->, thick, shorten >= 2pt, shorten <= 2pt, color=purple},
        contour/.style={draw=orange!80, thick},
        optimal_pt/.style={circle, fill=red, inner sep=1.5pt}
    ]

        
        \coordinate (A) at (0,0);       
        \coordinate (B) at (4,0);       
        \coordinate (C) at (2,3.5);     
        
        \filldraw[simplex] (A) -- (B) -- (C) -- cycle;
        \node[gray!60] at (0.5, 2.5) {$\Delta^n$};

        \coordinate (K1) at (1, 0.5);
        \coordinate (K2) at (3, 0.5);
        \coordinate (K3) at (3.2, 1.2);
        \coordinate (K4) at (2, 2.5); 
        \coordinate (K5) at (0.8, 1.2);
        
        \filldraw[permissible] (K1) -- (K2) -- (K3) -- (K4) -- (K5) -- cycle;
        \node[blue!60!black] at (2, 1.3) {\textbf{Permissible} $K$};
        
        \node[below] at (2,-0.5) {\textbf{Portfolio Space} (Input)};

        
        \begin{scope}[xshift=6cm, yshift=0.5cm]
            \clip (-1, -1) rectangle (6.0, 4.5);
        
            \draw[metric_space] (-0.5,0) -- (4,0) node[right] {Metric 1};
            \draw[metric_space] (0,-0.5) -- (0,3.5) node[above] {Metric 2};
            
            \coordinate (GK1) at (0.5, 0.5);
            \coordinate (GK2) at (2.5, 0.6);
            \coordinate (GK3) at (2.8, 1.5);
            \coordinate (GK4) at (1.8, 2.3); 
            \coordinate (GK5) at (0.3, 1.6);
            
            \draw[fill=green!15, draw=green!60!black, thick] 
                (GK1) to[bend left=10] 
                (GK2) to[bend right=5] 
                (GK3) to[bend right=20] 
                (GK4) to[bend right=15] 
                (GK5) to[bend right=10] 
                cycle;
                
            \node[blue!60!black] at (1.5, 1.2) {Image $g(K)$};
            
            \coordinate (UtilCenter) at (5.0, 5.0);
            
            \draw[contour, dashed] (UtilCenter) circle (4.18cm);
            
            \draw[contour, dashed] (UtilCenter) circle (3.68cm);
            
            \draw[contour, dashed] (UtilCenter) circle (4.68cm);
            
            \node[orange!80, align=center, fill=white, inner sep=1pt] at (3.0, 3.2) {Utility\\Contours};
            
            \node[below] at (1.5,-1) {\textbf{Space of quantitative attributes} $V$};
            
        \end{scope}

        
        \draw[map_arrow, line width=1.5pt] (3.5, 2) to[bend left=20] node[midway, above, text=black] {Attribute Map $g(x)$} (5.5, 2);

    \end{tikzpicture}%
}
    \caption{A portfolio space, a space of quantitative attributes for it, and an objective function (utility function).}
    \label{fig:space-of-quantitative-attributes}
\end{figure}

\subsection{Optimal Re-implementation via Constrained Optimization}

The goal is, for a given hub portfolio $x \in P_1$, to find a spoke portfolio $y \in P_2$ that best achieves two potentially conflicting aims:
\begin{enumerate}
   \item Match the hub metric: $g_\mathcal{B}(y)$ should be ``close'' to $g_\A(x)$.
   \item Optimize the objective: $u(g_\mathcal{B}(y))$ should be maximized.
\end{enumerate}
We combine these into a single optimization problem. For a fixed $x \in P_1$, we seek $y \in P_2$ that solves:
\[
\min_{y \in P_2} \left( \|g_\A(x) - g_\mathcal{B}(y)\|_V^p - \lambda u(g_\mathcal{B}(y)) \right)
\]
Here, $p \ge 1$ (typically $p=2$ for squared Euclidean distance) and $\lambda \ge 0$ is a weighting parameter balancing metric matching (minimizing distance) against objective maximization (maximizing $u$, hence minimizing $-u$).

Let $F(x,y)$ be the objective function depending on parameter $\lambda$:
\[
F(x,y) = \|g_{\mathcal{A}}(x) - g_{\mathcal{B}}(y)\|_V^p - \lambda u(g_{\mathcal{B}}(y))
\]

The set of optimal solutions for a given $x \in P_1$ is:
\[
y^*(x) = \operatorname*{arg\,min}_{y \in P_2} F(x, y)
\]

We use Berge's Maximum Theorem \cite{Berge} (applied to minimization) to establish the conditions under which this optimization problem constructs a valid horizontal morphism $f: K_1 \to K_2$ in $\HSSimp$.

\begin{theorem}[Berge's Theorem for Minimizers (informal)]
Let $X, Y$ be topological spaces. Let $F: X \times Y \to \mathbb{R}$ be a continuous function. Let $C \subseteq Y$ be a non-empty compact set. Define the correspondence $\Phi: X \rightrightarrows C$ by $\Phi(x) = \arg \min_{y \in C} F(x, y)$. Then:
\begin{enumerate}
   \item The value function $v(x) = \min_{y \in C} F(x, y)$ is continuous on $X$.
   \item The argmin correspondence $\Phi(x)$ is upper hemicontinuous and compact-valued on $X$.
   \item If, for each $x \in X$, the function $y \mapsto F(x, y)$ has a unique minimum over $C$, then $\Phi(x)$ is single-valued and the resulting function $f(x) = \Phi(x)$ is continuous on $X$.
\end{enumerate}
\end{theorem}

\begin{theorem}
\label{thm:berge-optimal}
Assume:
\begin{itemize}
   \item The permissible sets $P_1 \subseteq \Delta^n$ and $P_2 \subseteq \Delta^m$ are closed (hence compact).
   \item The attribute maps $g_\A: P_1 \to V$, $g_\mathcal{B}: P_2 \to V$ are continuous.
   \item The objective function $u: V \to \mathbb{R}$ is continuous.
   \item The combined objective $F(x, y) = \|g_\A(x) - g_\mathcal{B}(y)\|_V^p - \lambda u(g_\mathcal{B}(y))$ (for fixed $p \ge 1, \lambda \ge 0$) has a \textbf{unique minimum} $y^*(x)$ over $P_2$ for each $x \in P_1$.
\end{itemize}

Then:
\begin{enumerate}
   \item The map $f: P_1 \to P_2$ defined by $f(x) = y^*(x) = \arg \min_{y \in P_2} F(x, y)$ is continuous.
   \item Let $K_1 = P_1$. $K_1$ is a valid object in $\HSSimp$.
   \item Let $K_2 = f(K_1) \subseteq P_2$. $K_2$ is compact (hence closed) and is a valid object in $\HSSimp$.
   \item The map $f: K_1 \to K_2$ constitutes a \textbf{valid horizontal morphism in $\HSSimp$}.
\end{enumerate}
\end{theorem}

Sufficient conditions for uniqueness include: $P_2$ is convex, $g_\mathcal{B}$ is affine, the function $v \mapsto \|v\|^p$ is strictly convex (which requires both $p>1$ and that $\|\cdot\|_V$ is a strictly convex norm, e.g., the Euclidean norm or any $\ell^q$ norm with $1 < q < \infty$), and the function $v \mapsto -u(v)$ is strictly convex (i.e., $u$ is strictly concave).

\begin{proof}
(1) The function $F(x, y)$ is continuous on $P_1 \times P_2$ because $g_\A, g_\mathcal{B}, u$, the norm, and exponentiation ($(\cdot)^p$) are continuous. The constraint set $P_2$ is compact (as a closed subset of $\Delta^m$) and independent of $x$. Since $F(x, y)$ has a unique minimum for each $x$ by assumption, part 3 of Berge's Theorem applies directly (with $X=P_1$, $C=P_2$), showing that the function $f(x) = \arg \min_{y \in P_2} F(x, y)$ is continuous.
(2) $K_1=P_1$ is closed by assumption, hence a valid object.
(3) $K_2 = f(K_1)$ is the continuous image of the compact set $K_1$. Therefore, $K_2$ is compact. Compact subsets of Hausdorff spaces (like $\Delta^m$) are closed. Thus $K_2$ is a valid object.
(4) $K_1$ and $K_2$ are valid objects, and $f: K_1 \to K_2$ is a continuous map. This matches the definition of a horizontal morphism in $\HSSimp$.
\end{proof}

We restate this more simply:

\begin{theorem}[Construction of $\mathbb{HS}$ Morphisms via Optimization]
  Assume $P_1 \subseteq \Delta^n$ and $P_2 \subseteq \Delta^m$ are closed
  (hence compact), $g_\A : P_1 \to V$, $g_\mathcal{B} : P_2 \to V$ are continuous,
  and $u : V \to \mathbb{R}$ is continuous. If $g_\mathcal{B}$ is affine and injective, $P_2$ is convex,
  $1 < p < \infty$, $\|\cdot\|_V$ is a strictly convex norm (e.g., the Euclidean norm or any $\ell^q$ norm with $1 < q < \infty$; this ensures $v \mapsto \|v\|_V^p$ is strictly convex),
  and $v \mapsto -u(v)$ is strictly convex, then for each $x \in P_1$,
  the minimization problem has a unique solution $y^*(x)$, and the map
  $f(x) = y^*(x)$ is continuous.
\end{theorem}

\begin{remark}[Role of Injectivity]
The injectivity hypothesis on $g_\mathcal{B}$ ensures that the objective $F(x,y) = \|g_\A(x) - g_\mathcal{B}(y)\|_V^p - u(g_\mathcal{B}(y))$ is \emph{strictly} convex in $y$, not merely convex. If $g_\mathcal{B}$ is not injective, distinct portfolios $y \neq y'$ with $g_\mathcal{B}(y) = g_\mathcal{B}(y')$ would yield identical objective values, destroying uniqueness. In practice, uniqueness can also be recovered by adding a regulariser (e.g., a small $\|y\|^2$ term) or by restricting $P_2$ so that $g_\mathcal{B}|_{P_2}$ is injective.
\end{remark}

\begin{remark}[Significance]
This theorem shows that, under standard conditions ensuring a unique optimal solution (for example, strict convexity/concavity), the process of finding an `optimal' re-implementation that balances metric matching and objective optimization naturally constructs a continuous map $f$ between compact feasible spaces $K_1, K_2$. The hub space $K_1$ is simply the initial permissible space $P_1$, and the spoke space $K_2$ is the set of achievable optimal portfolios. This guarantees the resulting re-implementation $f$ is automatically proper and fits coherently into the $\HSSimp$ double categorical structure.
\end{remark}

\begin{figure}
\centering
\begin{tikzpicture}[scale=1, >=stealth]
    \node[circle, fill, inner sep=1.5pt, label=left:$x$] (x) at (0,2) {};
    \node at (0, 3) {Hub $K_1$};

    \draw[thick] (4,0) -- (9,0) -- (6.5,4) -- cycle;
    \node at (8, 3.5) {Spoke $K_2$};

    \draw[fill=blue!15, draw=blue!60!black, thick] 
        (5.3, 1.2) -- (6.2, 1.0) -- (6.8, 2.0) -- (5.8, 2.4) -- (5.3, 1.8) -- cycle;
    \node[text=blue!60!black, font=\footnotesize, align=center] at (5.8, 1.5) {Permissible\\Set $F_R(x)$};

    \coordinate (Ideal) at (7.8, 2.2);
    \node[circle, fill=gray, inner sep=1pt] at (Ideal) {};
    \node[text=gray, font=\footnotesize, right] at (Ideal) {Ideal metric (infeasible)};

    \draw[orange!80, thin] (Ideal) circle (1.02);
    \draw[orange!80, thin] (Ideal) circle (1.4);
    \draw[orange!80, thin] (Ideal) circle (1.8);
    
    \node[circle, fill=red, inner sep=1.5pt, label=above:$y^*$] (opt) at (6.8, 2.0) {};
    \node[text=red, font=\footnotesize, align=center] at (6.8, 0.4) {Optimal Spoke\\(Boundary Solution)};

    \draw[dashed, red, opacity=0.5] (Ideal) -- (opt);

    \draw[->, dashed] (x) -- (5.3, 1.8) node[midway, above] {$R$};
\end{tikzpicture}
\caption{Finding the optimal spoke from a permissible set of spokes}
\label{fig:optimal-spoke}
\end{figure}
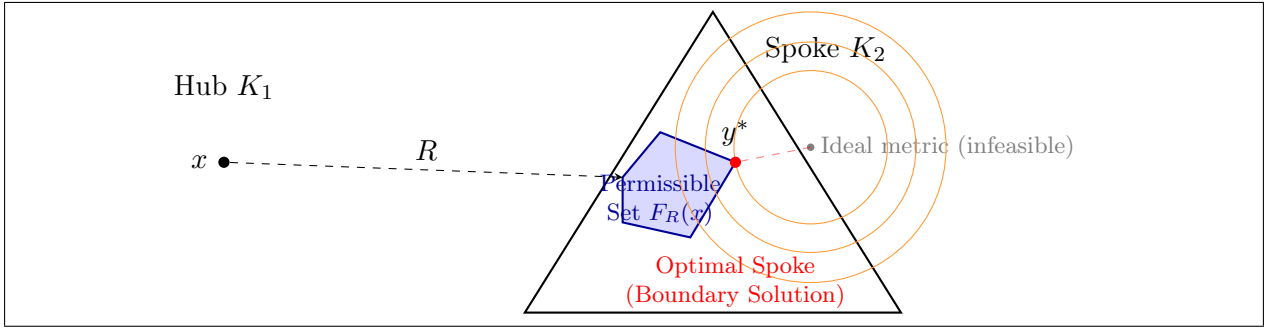

\subsubsection*{Non-uniqueness}

In cases where the objective is flat (e.g., non-unique solutions to linear programs) or absent (pure constraint satisfaction), we must rely on Michael's Selection Theorem to guarantee the existence of a continuous re-implementation.

\begin{theorem}[Michael's Selection Theorem \cite{Michael}]
\label{thm:michael-selection}
Let $X$ be a paracompact space (e.g., a metric space) and $Y$ be a Banach space. Let $\Phi: X \rightrightarrows Y$ be a lower hemicontinuous correspondence with non-empty, closed, and convex values. Then there exists a continuous function $f: X \to Y$ such that $f(x) \in \Phi(x)$ for all $x \in X$.
\end{theorem}

\begin{corollary}[Existence of Continuous Re-implementation]
Let permissible spaces $K_1$ and $K_2$ be convex, closed subsets of simplices. Let $R \subseteq K_1 \times K_2$ be a closed alignment relation such that for every $x \in K_1$, the feasible set $F_R(x)$ is convex and non-empty. If the correspondence $x \mapsto F_R(x)$ is lower hemicontinuous, then there exists a continuous horizontal morphism $f: K_1 \to K_2$ in $\HSSimp$ such that $\Graph(f) \subseteq R$.  (Since $K_2$ is a closed convex subset of $\mathbb{R}^m$, Michael's theorem guarantees a continuous selection into $\mathbb{R}^m$ with values in the convex images of the correspondence, so the selection lands in $K_2$.)
\end{corollary}

\begin{remark}
The convexity requirement is naturally satisfied by linear constraints (e.g., budget sets, sector limits). The lower hemicontinuity condition intuitively requires that the ``menu'' of options does not abruptly collapse in the interior of the hub space. When measurability rather than continuity suffices, the Kuratowski--Ryll-Nardzewski measurable selection theorem \cite{Kuratowski} provides an alternative under weaker hypotheses. This theorem ensures that the ``Design'' phase (defining the menu of allowable spoke portfolios, cf. Section~\ref{sec:DOTS}) is mathematically compatible with the ``Execution'' phase (selecting a continuous map $f$).
\end{remark}

\subsection{Optimal Re-implementation via Alignment Relation Constraint}
\label{subsec:opt-align}

An alternative scenario arises when the primary goal is to optimize an objective function $u: K_2 \to \mathbb{R}$ on the spoke space, subject to a pre-defined alignment relation $R \subseteq K_1 \times K_2$ acting as a constraint, as in Section~\ref{subsec:reimplalign}. Here, $K_1 \subseteq \Delta^n$ and $K_2 \subseteq \Delta^m$ are the given compact (closed) permissible portfolio spaces, and $R$ is a closed relation (a vertical morphism).

For a given hub portfolio $x \in K_1$, the set of feasible spoke portfolios is determined by the relation $R$:
\[
F_R(x) = \{ y \in K_2 \mid (x, y) \in R \} \subseteq K_2
\]
This set $F_R(x)$ represents all portfolios in $K_2$ that are considered aligned with $x$ according to $R$. The optimization problem is then to select the best portfolio from this set:
\[
\max_{y \in F_R(x)} u(y)
\]
Let $\Phi_R(x)$ be the set of optimal solutions:
\[
\Phi_R(x) = \arg \max_{y \in F_R(x)} u(y) \subseteq K_2
\]
If optimal solutions are unique, $\Phi_R$ defines a function $\phi: K_1\rightarrow K_2$.

We again use Berge's Maximum Theorem to determine when this process yields a valid horizontal morphism in $\HSSimp$.

\begin{theorem}[Construction via Alignment Constraint]
\label{thm:berge-relation-optimal}
Assume:
\begin{enumerate}
\item $K_1 \subseteq \Delta^n$ and $K_2 \subseteq \Delta^m$ are permissible portfolio spaces, i.e.\ closed (hence compact) subsets.
\item $R \subseteq K_1 \times K_2$ is a vertical morphism, i.e.\ a closed relation.
\item $u : K_2 \to \mathbb{R}$ is a continuous objective function.
\item For each $x \in K_1$, define the feasible correspondence
\[
F_R(x) := \{ y \in K_2 \mid (x,y) \in R \}.
\]
Assume that $F_R(x)$ is nonempty for all $x \in K_1$.
\item The correspondence $F_R : K_1 \rightrightarrows K_2$ is \emph{upper hemicontinuous} with compact values. (Note: by Remark~\ref{rem:UHC} below, this condition is \emph{automatic} when $R$ is a closed relation between compact Hausdorff spaces---which is guaranteed by assumptions (1) and (2). It is listed here for logical clarity, but in practice need not be verified separately.)
\item For every $x \in K_1$, the maximization problem
\[
\max_{y \in F_R(x)} u(y)
\]
has a unique solution.
\end{enumerate}
Then:
\begin{enumerate}
\item The induced map
\[
f : K_1 \to K_2, \qquad f(x) := \arg\max_{y \in F_R(x)} u(y),
\]
is continuous.
\item Let
\[
K_1' := \operatorname{dom}(R) = \{ x \in K_1 \mid F_R(x) \neq \varnothing \}.
\]
Then $K_1'$ is closed and hence a valid object of $\HSSimp$.
\item Let $K_2' := f(K_1')$. Then $K_2'$ is compact (hence closed) and a valid object of $\HSSimp$.
\item The map $f : K_1' \to K_2'$ is a horizontal morphism in $\HSSimp$ satisfying $\Graph(f) \subseteq R$.
\end{enumerate}

\end{theorem}

\begin{remark}[Sufficient conditions for upper hemicontinuity]
\label{rem:UHC}
A convenient sufficient condition for the correspondence
\[
F_R : K_1 \rightrightarrows K_2, \qquad
F_R(x) := \{y \in K_2 \mid (x,y) \in R\},
\]
to be upper hemicontinuous with compact values is that the alignment relation
\[
R \subseteq K_1 \times K_2
\]
be closed, with $K_1$ and $K_2$ compact Hausdorff spaces.

Indeed, if $R$ is closed in $K_1 \times K_2$, then for each $x \in K_1$ the fiber
\[
F_R(x) = \{y \in K_2 \mid (x,y) \in R\}
\]
is compact. Moreover, the closedness of $R$ implies that the correspondence
$F_R$ has a closed graph. For correspondences with compact values into a
Hausdorff space, the closed--graph property is equivalent to upper
hemicontinuity.

In practical portfolio-construction settings, relations $R$ are typically
defined by systems of inequalities of the form
\[
g_j(x,y) \le 0, \qquad j = 1,\dots,m,
\]
where each $g_j : K_1 \times K_2 \to \mathbb{R}$ is continuous. Such relations
are automatically closed, and hence induce upper hemicontinuous feasible-set
correspondences. Consequently, the continuity of the optimizer constructed in
Theorem~\ref{thm:berge-relation-optimal} follows from standard topological facts, without requiring any
interiority or convexity assumptions beyond those ensuring uniqueness of the
optimum.
\end{remark}

\begin{proof}
(1) For a fixed $x$, $F_R(x)$ is the projection onto $K_2$ of the closed set $R \cap (\{x\} \times K_2)$. Since $R$ is closed in $K_1 \times K_2$ and $\{x\} \times K_2$ is closed, the intersection is closed. As $K_1 \times K_2$ is compact, $R$ is compact, and the intersection is compact. The continuous projection of a compact set is compact, so $F_R(x)$ is compact. The upper hemicontinuity of $F_R$ follows from $R$ being a closed relation between compact Hausdorff spaces (a standard result, related to the Closed Graph Theorem for correspondences).
(2) $K'_1 = \operatorname{dom}(R) = \pi_1(R)$ is the image of the compact set $R$ under the continuous projection $\pi_1: K_1 \times K_2 \to K_1$, hence compact and therefore closed in $K_1 \subseteq \Delta^n$.
(3) $K'_2 = f(K'_1)$ is the continuous image of the compact set $K'_1$, hence compact and therefore closed in $K_2 \subseteq \Delta^m$.
(4) $K'_1$ and $K'_2$ are valid objects, and $f: K'_1 \to K'_2$ is a continuous map, matching the definition of a horizontal morphism.
\end{proof}

\begin{remark}[Significance]
This theorem shows that optimization subject to alignment relation constraints also naturally constructs valid $\HSSimp$ objects and morphisms, provided the relation $R$ and objective $u$ satisfy standard conditions (ensuring the domain $K'_1$ is closed and solutions are unique). This reinforces the suitability of the $\HSSimp$ framework for modeling practical portfolio construction tasks where constraints are given by relations rather than metric proximity. The necessary condition on $R$ for $K'_1$ to be closed (related to lower hemicontinuity of $F_R$) often holds in practical cases, for example, if $R$ itself is defined by inequalities involving continuous functions.
\end{remark}

\begin{remark}[Categorical interpretation of the feasible set]
Recall from Section~\ref{subsec:reimplalign} that a valid horizontal morphism $f$ must satisfy the inclusion $Graph(f) \subseteq R$. In the context of optimization, the feasible set $F_R(x)$ is precisely the fiber of the relation $R$ at $x$:
\[
    F_R(x) = \{ y \in K_2 \mid (x,y) \in R \}.
\]
Therefore, the optimization problem is structurally equivalent to finding a continuous selection $f(x) \in F_R(x)$. If such a continuous selection exists, the resulting map $f$ automatically satisfies the 2-cell condition of Section~\ref{subsec:reimplalign}. The union of these feasible sets across all $x$ constitutes the ``action'' defined in Section~\ref{sec:DOTS}, $K_1 \odot R = \bigcup_{x \in K_1} F_R(x)$.
\end{remark}

\subsection{Coherence of Squares Constructed from Alignment Constraints}
\label{subsec:coherence}

We now analyse the coherence of a square of horizontal morphisms where, instead of being constructed from metric-based optimization (Theorem~\ref{thm:berge-optimal}), they are constructed via \textbf{alignment relation constraints} as defined in Theorem~\ref{thm:berge-relation-optimal}.

Assume we are given four objects ($K_1, K_2, K_3, K_4$) and a ``square'' of four vertical morphisms (alignment relations) that serve as constraints:
\begin{itemize}
   \item $R_f \subseteq K_1 \times K_2$
   \item $R_g \subseteq K_1 \times K_3$
   \item $R_{g'} \subseteq K_2 \times K_4$
   \item $R_{f'} \subseteq K_3 \times K_4$
\end{itemize}
We also assume given continuous objective functions $u_2: K_2 \to \mathbb{R}$, $u_3: K_3 \to \mathbb{R}$, and $u_4: K_4 \to \mathbb{R}$.

The horizontal morphisms $f, g, f', g'$ are then constructed (per Theorem~\ref{thm:berge-relation-optimal}) as the unique optimizers over the feasible sets defined by these relations. For example, for a given $x \in dom(f)$, $f(x)$ is the unique portfolio $y$ that maximizes $u_2(y)$ subject to the constraint that $(x,y) \in R_f$.

\begin{figure}[h]
\centering
\begin{tikzcd}[row sep=large, column sep=large]
K_1 
    \arrow[r, blue, "f"] 
    \arrow[d, blue, "g"'] 
    \arrow[r, dashed, "R_f"', below, yshift=-5pt] 
    \arrow[d, dashed, "R_g", right, xshift=5pt] 
& K_2 
    \arrow[d, blue, "g'"'] 
    \arrow[d, dashed, "R_{g'}", right, xshift=5pt] 
\\
K_3 
    \arrow[r, blue, "f'"] 
    \arrow[r, dashed, "R_{f'}"', below, yshift=-5pt] 
& K_4
\end{tikzcd}
\caption{A square of horizontal morphisms $f, g, f', g'$ (\textcolor{blue}{solid blue}) constructed from a square of vertical alignment constraints $R_f, R_g, R_{f'}, R_{g'}$ (dashed).}
\end{figure}
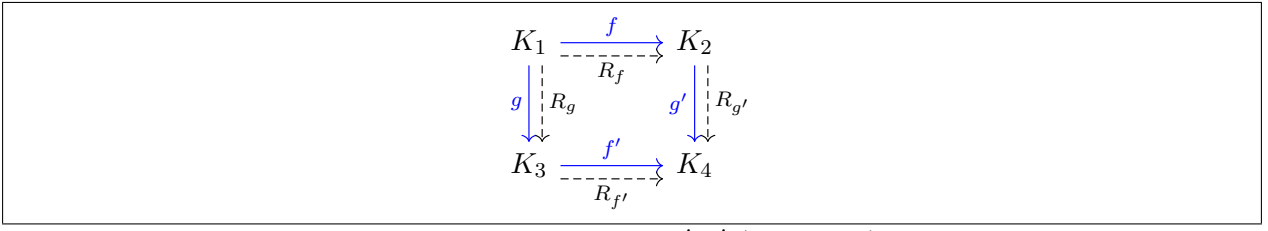

\subsubsection*{Analysis of Commutativity}

In this framework, the square \textbf{does not commute} in general. The construction is inherently path-dependent.
Let's trace the two paths for a given $x \in K_1$:

\begin{itemize}
   \item \textbf{Path 1 ($f' \circ g$):}
   \begin{enumerate}
       \item $z = g(x) = \arg\max_{z' \in F_{R_g}(x)} u_3(z')$. This $z$ is the ``best'' portfolio in $K_3$ (according to $u_3$) that is aligned with $x$ (via $R_g$).
       \item $w = f'(z) = \arg\max_{w' \in F_{R_{f'}}(z)} u_4(w')$. This $w$ is the ``best'' portfolio in $K_4$ (according to $u_4$) that is aligned with $z$ (via $R_{f'}$).
   \end{enumerate}
   \item \textbf{Path 2 ($g' \circ f$):}
   \begin{enumerate}
       \item $y = f(x) = \arg\max_{y' \in F_{R_f}(x)} u_2(y')$. This $y$ is the ``best'' portfolio in $K_2$ (according to $u_2$) that is aligned with $x$ (via $R_f$).
       \item $w' = g'(y) = \arg\max_{w'' \in F_{R_{g'}}(y)} u_4(w'')$. This $w'$ is the ``best'' portfolio in $K_4$ (according to $u_4$) that is aligned with $y$ (via $R_{g'}$).
   \end{enumerate}
\end{itemize}

Commutativity ($w = w'$) fails because the intermediate ``greedy'' optimizations are different. The portfolio $z$ is chosen using $R_g$ and $u_3$, while $y$ is chosen using $R_f$ and $u_2$. Since $y$ and $z$ are generally different, the inputs to the second-stage optimizations are different, leading to different final portfolios ($w \neq w'$).

We can obtain commutativity by imposing some rather special conditions on the optimization problems defining the square.

\begin{theorem}[Commutativity via Bellman Consistency]
\label{thm:bellman_commutativity}
Consider alignment relations $R_f \subseteq K_1 \times K_2$, $R_g \subseteq K_1 \times K_3$, $R_{g'} \subseteq K_2 \times K_4$, $R_{f'} \subseteq K_3 \times K_4$ forming a \emph{commuting square of relations}: $R_{f'} \circ R_g = R_{g'} \circ R_f$ (i.e., the set of final portfolios reachable via both paths is the same).

Let $u_4: K_4 \to \mathbb{R}$ be a continuous, strictly concave objective function. Define the intermediate objective functions $u_2: K_2 \to \mathbb{R}$ and $u_3: K_3 \to \mathbb{R}$ by the Bellman value condition:
\begin{align*}
u_2(y) &= \max \{ u_4(w) \mid (y, w) \in R_{g'} \} \\
u_3(z) &= \max \{ u_4(w) \mid (z, w) \in R_{f'} \}
\end{align*}
Construct the horizontal morphisms $f, g, f', g'$ as the unique optimizers subject to alignment constraints:
\begin{align*}
f(x) &= \arg \max_{y \in K_2, (x,y) \in R_f} u_2(y), & g'(y) &= \arg \max_{w \in K_4, (y,w) \in R_{g'}} u_4(w), \\
g(x) &= \arg \max_{z \in K_3, (x,z) \in R_g} u_3(z), & f'(z) &= \arg \max_{w \in K_4, (z,w) \in R_{f'}} u_4(w).
\end{align*}
Assume that for every $x \in K_1$, the set of reachable final portfolios $W_x = \{ w \in K_4 \mid (x, w) \in R_{f'} \circ R_g \}$ is non-empty and \emph{convex}. Then the square of re-implementations commutes:
\[ f' \circ g = g' \circ f \]
\end{theorem}

\begin{proof}
Let $x \in K_1$ be an arbitrary hub portfolio. We determine the image of $x$ under both paths.

\noindent \textbf{Step 1: Uniqueness of the global optimum.}
Define $W_x \subseteq K_4$ as the set of all final portfolios reachable from $x$. Since the relations commute ($R_{f'} \circ R_g = R_{g'} \circ R_f$), the set of reachable portfolios is identical for both paths. Since $K_i$ are compact and relations are closed, $W_x$ is compact. By hypothesis, $W_x$ is convex. Since $u_4$ is strictly concave and $W_x$ is convex and compact, there exists a unique global maximizer $w^* \in W_x$:
\[ w^* = \arg \max_{w \in W_x} u_4(w). \]

\noindent \textbf{Step 2: Analysis of Path 1 ($x \xrightarrow{g} z \xrightarrow{f'} w$).}
The morphism $g(x)$ selects the portfolio $z^* \in K_3$ that maximizes $u_3(z)$ subject to $(x, z) \in R_g$. Substituting the definition of $u_3$:
\[ z^* = \arg \max_{z: (x,z) \in R_g} \left( \max_{w: (z,w) \in R_{f'}} u_4(w) \right). \]
This nested maximization is equivalent to maximizing $u_4(w)$ over the composite relation $R_{f'} \circ R_g$. To see that the witnesses coincide: let $w^*$ be the unique maximiser of $u_4$ over $W_x$ (unique by strict concavity on the convex compact set $W_x$). Since $w^* \in W_x = \{w : \exists z,\, (x,z) \in R_g \text{ and } (z,w) \in R_{f'}\}$, there exists $z_0$ with $(x, z_0) \in R_g$ and $(z_0, w^*) \in R_{f'}$. Then $u_3(z_0) = \max_{w:(z_0,w) \in R_{f'}} u_4(w) \ge u_4(w^*)$. But $u_3(z^*) \ge u_3(z_0)$ by optimality of $z^*$, and $u_3(z^*) = u_4(w')$ for some $w'$ with $(z^*, w') \in R_{f'}$, so $w' \in W_x$ and $u_4(w') \ge u_4(w^*)$. By uniqueness, $w' = w^*$, whence $(z^*, w^*) \in R_{f'}$.
The second morphism $f'(z^*)$ then computes:
\[ f'(z^*) = \arg \max_{w: (z^*, w) \in R_{f'}} u_4(w). \]
Since $w^*$ is reachable from $z^*$ and $w^*$ is the global maximum on $W_x$, it follows that $f'(z^*) = w^*$. Thus, $(f' \circ g)(x) = w^*$.

\noindent \textbf{Step 3: Analysis of Path 2 ($x \xrightarrow{f} y \xrightarrow{g'} w$).}
By symmetry, $f(x)$ selects $y^* \in K_2$ maximizing $u_2(y)$. By the Bellman definition of $u_2$, $y^*$ is selected precisely because it allows access to the maximal value of $u_4$, which is attained at $w^*$. The second morphism $g'(y^*)$ then selects $w^*$. Thus, $(g' \circ f)(x) = w^*$.

\noindent \textbf{Conclusion.}
Both paths map every $x \in K_1$ to the same unique global optimizer $w^*$. The square commutes.
\end{proof}

\begin{remark}[Convexity of $W_x$]
The convexity hypothesis on $W_x$ is a mild condition when the alignment relations $R_g, R_{f'}$ are defined by affine constraints (linear equalities and inequalities on portfolio weights), since the image of a convex set under composition of affine-constrained relations remains convex.  It may fail when the relations encode non-convex constraints (e.g., cardinality bounds or discrete sector assignments).
\end{remark}

\begin{remark}[Set-valued intermediate optimizers]
When the intermediate value functions $u_2$ and $u_3$ are not strictly concave, the argmax correspondences $f$ and $g$ may be set-valued. The commutativity conclusion---that both paths reach the same terminal portfolio $w^*$---remains valid regardless: it depends only on uniqueness of the final optimum, not of the intermediate ones. If single-valued horizontal morphisms are required, sufficient conditions include strict concavity of the intermediate objectives, which holds when $u_4$ is strictly concave and the constraint correspondences $z \mapsto \{w : (z,w) \in R_{f'}\}$ are affine.
\end{remark}

Unfortunately, this theorem is rarely useful in practice, since it involves replacing the original objective functions with the intermediate objective functions defined by the Bellman value condition.

\subsection{Summary}
\label{subsec:opt-summary}

It is convenient to summarize the results of this section as follows.

\begin{theorem}[Coherent construction of re-implementation pipelines]
\label{thm:pipeline}
Let $K_1,K_2,K_3,K_4$ be permissible portfolio spaces, i.e.\ closed (hence compact)
subsets of simplices. 

Let
\[
R_f \subseteq K_1 \times K_2, \quad
R_g \subseteq K_1 \times K_3, \quad
R_{g'} \subseteq K_2 \times K_4, \quad
R_{f'} \subseteq K_3 \times K_4
\]
be closed alignment relations forming a commuting square:
\[
R_{f'} \circ R_g = R_{g'} \circ R_f .
\]

Assume continuous objective functions
\[
u_2 : K_2 \to \mathbb{R}, \quad
u_3 : K_3 \to \mathbb{R}, \quad
u_4 : K_4 \to \mathbb{R}.
\]

\paragraph{(A) Existence and continuity.}
Assume that for each $h \in \{f,g,f',g'\}$ and each admissible input portfolio,
the feasible correspondence defined by the relation $R_h$ is nonempty, compact,
upper hemicontinuous, and admits a unique optimizer with respect to the relevant
objective function. 

Then the induced maps
\[
f : K_1 \to K_2, \quad
g : K_1 \to K_3, \quad
g' : K_2 \to K_4, \quad
f' : K_3 \to K_4
\]
are continuous horizontal morphisms in $\HSSimp$ whose graphs lie in the
corresponding relations.

\paragraph{(B) General (path-dependent) behaviour.}
In general, the resulting two-stage re-implementations
\[
f' \circ g \quad \text{and} \quad g' \circ f
\]
need not coincide. 
Thus, even when the alignment relations commute, pipelines constructed via
greedy stagewise optimization are generically path-dependent.

\paragraph{(C) Bellman-consistent (path-independent) case.}
Suppose in addition that the final objective $u_4 : K_4 \to \mathbb{R}$ is
strictly concave, and define intermediate value functions by
\[
u_2(y) := \max \{ u_4(w) \mid (y,w) \in R_{g'} \}, \qquad
u_3(z) := \max \{ u_4(w) \mid (z,w) \in R_{f'} \}.
\]

Assume the set of reachable final portfolios $W_x$ is convex for each $x$, so that strict concavity of $u_4$ on $W_x$ guarantees these maxima are attained uniquely (Theorem~\ref{thm:bellman_commutativity}).
Construct all four maps $f,g,g',f'$ as unique optimizers relative to these value
functions.

Then the resulting square of horizontal morphisms commutes:
\[
f' \circ g = g' \circ f .
\]

\paragraph{(D) Interpretation.}
Part (B) reflects genuine economic path dependence of greedy local optimizations.
Part (C) identifies a dynamic-programming regime in which re-implementation
pipelines become coherent: each stage optimizes future value, and the final
portfolio is globally optimal among all portfolios reachable via the composite
alignment relation.
\end{theorem}

\clearpage  \section{Propagation of Portfolio Changes and New Requirements}
\label{sec:propagation}

\subsection{Propagation of Portfolio Changes}

A fundamental task in portfolio management is determining how changes in a ``hub''
portfolio should flow through to its ``spokes.'' The $\HSSimp$ double category
provides rigorous support for this process.

The horizontal morphism $f: K_1 \to K_2$ acts as the propagation mechanism.
A portfolio change $x \leadsto x'$ at the hub $K_1$ is coherently
transformed into the change $y \leadsto y'$ at the spoke $K_2$,
where $y = f(x)$ and $y' = f(x')$.

Alignment is modeled by a closed relation $R \subseteq K_1 \times K_3$.
\par
\begin{itemize}
\item \textbf{Role of Properness}: Because $f$ is automatically proper
(Prop.~\ref{prop:reimpl-are-proper}), the pushforward $f_{!}R$ is
a closed relation (Thm.\ \ref{thm:proper-closed}).
This ensures that the re-implemented portfolio $y'$ cannot be a
``phantom portfolio''---a limit point that appears to satisfy the
propagated alignment but has no compliant pre-image.

\item \textbf{Role of Adjunction}: The equivalence
$R \subseteq f^{*}S \iff f_{!}R \subseteq S$ (Thm.~\ref{thm:adjunction})
allows alignment to be checked in either the hub space $K_1$ or
the spoke space $K_2$.

\item \textbf{Compositional Laws}: Beck--Chevalley (Thm.
\ref{thm:Beck--Chevalley}) and Frobenius (Thm.~\ref{thm:frobenius})
ensure path independence and that filtering commutes with
re-implementation.
\end{itemize}

The coherence of $\HSSimp$ means that a change in the \textbf{Hub} portfolio
($x \leadsto x'$) that respects the alignment relation $R$ will produce
a change in the \textbf{Spoke} portfolio ($y \leadsto y'$) that respects the propagated closed alignment relation $f_{!}R$.

\begin{remark}[Diagnostic versus prescriptive use of pullback]
The pullback $f^{*}S$ admits two distinct modes of use---a distinction that is mathematical in origin but operational in consequence. In \emph{diagnostic} mode, we ask: ``given a proposed change to the hub, will the spoke's alignment break?'' Here we compute $f^{*}S$ and check whether the new hub portfolio lies in the appropriate fiber; the pullback serves as a feasibility test. In \emph{prescriptive} mode, we reverse the question: ``how should we redesign the hub so that a desired spoke alignment can be achieved?'' Here $f^{*}S$ characterizes the set of hub portfolios compatible with the spoke constraint, and we optimize within it.

A subtler case---one the basic framework handles but does not emphasize---arises when spoke constraints are systematic: regulatory position limits, say, that apply uniformly to all client accounts. Here the spokes collectively constrain the hub. The intersection $\bigcap_i f_i^{*}S_i$ over all spoke constraints defines the feasible region for the hub, and the direction of prescription reverses.
\end{remark}

\subsection{Propagation of New Alignment Requirements}

New alignment or compliance requirements can be imposed on the hub and must be propagated downstream to the spokes---possibly through several stages. In practice, portfolio changes are more often triggered by new alignment requirements than by changes in investment conviction.

To illustrate the practical necessity of the topological and categorical results established in previous sections, consider a multi-stage investment pipeline. In this scenario, a centralized ``Hub'' strategy is transformed first into an intermediate ``Model'' portfolio, and subsequently into hundreds of customized ``Client'' accounts.

\subsubsection*{The Scenario: New CIO Asset Allocation Guidance}

A Chief Investment Officer (CIO) has determined that portfolios must now conform to centralized asset allocation guidance. This investment framework is formalized as follows: the asset allocation guidance corresponds to a portfolio $z\in K_{\text{Asset Classes}}$, and ``acceptable conformance'' of actual portfolios to guidance corresponds to an alignment relation \[R_{\text{CIO}}: K_{\text{Hub}}\dashrightarrow K_{\text{Asset Classes}}\] 

The ongoing portfolio management problem is as follows: whenever the CIO asset allocation guidance changes, the hub portfolio may need to be changed in order to conform to the new guidance. Changes then need to be propagated to the spoke portfolios, and this may be a multi-stage process.

That is, the propagation of this guidance may involve a composite re-implementation process:
\[
K_{\text{Hub}} \xrightarrow{f} K_{\text{Model}} \xrightarrow{g} K_{\text{Client}}
\]
New constraints may be introduced at each stage of this pipeline:
\begin{enumerate}
    \item \textbf{Stage 1 (Hub $\to$ Model):} The map $f$ re-implements the Hub strategy into a Model portfolio. At this stage, a \textit{Sector Constraint} $S_{\text{Sec}}$ is applied (e.g., strictly aligning with maximum sector weights).
    \item \textbf{Stage 2 (Model $\to$ Client):} The map $g$ re-implements the Model into final Client accounts. At this stage, \textit{Tax and Fee Constraints} $T_{\text{Tax}}$ are applied (e.g., limiting turnover for tax efficiency).
\end{enumerate}

The question is whether the final client portfolio $z \in K_{\text{Client}}$---which has been subjected to sector constraints $S_{\text{Sec}}$ and tax constraints $T_{\text{Tax}}$---actually satisfies the original CIO guidance $R_{\text{CIO}}$. The double category $\HSSimp$ ensures that it does, and that the audit trail remains unbroken across this multi-stage construction:

\begin{itemize}
    \item \textbf{Conformance via Properness (Theorem~\ref{thm:proper-closed}):}
    Because the re-implementation maps $f$ and $g$ are proper, the pushforward of the CIO's guidance moves securely through the chain. The set of valid client portfolios derived from the CIO's guidance, $(g \circ f)_! R_{\text{CIO}}$, is guaranteed to be a \textit{closed} set. This ensures that the intersection with local constraints at each stage is well-defined; no ``phantom'' client portfolios can be generated that technically satisfy the solver tolerances but violate the fundamental CIO directive.

    \item \textbf{Audit Conservatism via Lax Beck--Chevalley (Proposition~\ref{prop:lax-BC}):}
    While strict path independence requires specific conditions (pointwise cartesianness), the \textit{lax} Beck--Chevalley condition holds in full generality:
    \[
    f'_! (g^* R_{\text{CIO}}) \subseteq h^* (f_! R_{\text{CIO}})
    \]
    In operational terms, this inclusion acts as a \emph{safety guarantee for hub-side verification}. It asserts that if we verify conformance \textit{early} (filtering the hub portfolios against the CIO guidance before re-implementation), the resulting client portfolios are mathematically guaranteed to satisfy the downstream conformance checks. While some conforming portfolios might theoretically be found via ``late verification'' that cannot be traced back to the hub (unless the square is pointwise cartesian), the lax condition ensures that the upstream audit is strictly conservative and never produces false positives.

    \item \textbf{Coherent Filtering via Frobenius (Theorem~\ref{thm:frobenius}):}
    The Frobenius law ensures that sector compliance can be checked at either the hub or the model level without affecting the outcome. Applying Frobenius to $f: K_{\text{Hub}} \to K_{\text{Model}}$ (with $R_{\text{CIO}}$ on $K_{\text{Hub}}$ and $S_{\text{Sec}}$ on $K_{\text{Model}}$):
    \[
    f_!(R_{\text{CIO}} \cap f^* S_{\text{Sec}}) = f_!(R_{\text{CIO}}) \cap S_{\text{Sec}}
    \]
    The LHS filters at the hub (pulling the sector constraint back via $f^*$ before re-implementation); the RHS filters at the model (intersecting after re-implementation). Pushing both sides forward by $g_!$ and intersecting with $T_{\text{Tax}}$:
    \[
    g_!\bigl(f_!(R_{\text{CIO}} \cap f^* S_{\text{Sec}})\bigr) \cap T_{\text{Tax}} = g_!\bigl(f_!(R_{\text{CIO}}) \cap S_{\text{Sec}}\bigr) \cap T_{\text{Tax}}
    \]
    In business terms, the set of permissible client portfolios is identical regardless of whether we incorporate sector compliance at the hub level (via the pullback $f^* S_{\text{Sec}}$, before Stage~1) or at the model level (via $S_{\text{Sec}}$, between Stages~1 and~2). The CIO's guidance is not ``diluted'' or lost as it passes through the subsequent layers of constraints.
\end{itemize}

\clearpage  \section{Computational Considerations}
\label{sec:computation}

\subsection{Representing Objects and Morphisms}

For computational work in $\HSSimp$:
\par
\begin{itemize}
\item \textbf{Objects $K$}: Represented as polyhedral sets
$K = \{x \in \Delta^n \mid Ax \leq b\}$. These are intersections of closed
half-spaces, so they are always closed (and thus valid objects).
\item \textbf{Vertical Morphisms $R$}: Represented as $R = \{(x, y) \in K_1
\times K_2 \mid Cx + Dy \leq e\}$. These are also closed.
\item \textbf{Horizontal Morphisms $f$}:
   \par
   \begin{itemize}
       \item If $f(x) = Fx + g$ is affine, it is continuous.
       \item If $f(x)$ is the solution to the optimization problem in
       Section~\ref{sec:optimal-reimpl}, it is continuous (by Thm.
       \ref{thm:berge-optimal}).
   \end{itemize}
\end{itemize}

\subsection{Computing Pullback}

Given $f: K_1 \to K_2$ and $S = \{(y, z) \in K_2 \times K_3 \mid Q(y, z) \leq 0\}$:
\[
f^*S = \{(x, z) \in K_1 \times K_3 \mid Q(f(x), z) \leq 0\}
\]
If $f$ is affine ($f(x) = Fx+g$) and $S$ is polyhedral ($Py+Qz \le q$):
\[
f^*S = \{(x, z) \in K_1 \times K_3 \mid P(Fx + g) + Qz \leq q\}
\]
This is a simple matrix-vector operation.

\subsection{Computing Pushforward}

Given $f: K_1 \to K_2$ and $R = \{(x, z) \in K_1 \times K_3 \mid C(x, z) \leq 0\}$:
\[
f_!R = \{(y, z) \in K_2 \times K_3 \mid \exists x \in K_1: y = f(x), C(x, z) \leq 0\}
\]
This is a projection operation. If $f$ and $C$ are linear, this is
Fourier-Motzkin elimination \cite{Ziegler}, which has high computational complexity.

In general, however, computing the pushforward is even more difficult because it requires projecting a nonlinear set
\[
\{(x,y,z) \mid y=f(x),\ (x,z)\in R\}
\]
onto the $(y,z)$ coordinates. This is a nonlinear projection or quantifier-elimination problem, which
is intractable in general. 

Fortunately, in a high dimensional production setting, we rarely need to visualize or delineate the boundary of $f_!R$ precisely. We only need to answer: ``Does this specific spoke portfolio $y$ satisfy $(y,z)\in f_!R$\,? In other words, does there exist $x\in K_1$ such that $f(x)=y$ and $(x, z)\in R$\,?'' Even for $N=500$, if $f$ is defined by convex optimization and $R$ is convex, this verification step is a standard convex feasibility problem, which is tractable for modern solvers.

There are, however, several special cases where delineating $f_!R$ is tractable in practice.

\subsubsection*{1. Smooth and Locally Invertible Maps}
If $f$ is smooth and locally invertible (Jacobian non-singular), one can invert the relation locally:
\[
f_!R = \{(f(x), z) \mid (x,z)\in R\}.
\]
In such cases, the pushforward can be computed directly by evaluating $f$ pointwise on samples of
$R$, making it practical for low-dimensional applications.

\subsubsection*{2. Piecewise Linear or Affine Maps}
When $f$ is piecewise affine (common in optimization-based mappings that switch regimes) and $R$
is convex or polyhedral, $f_!R$ is a finite union of convex projections:
\[
f_!R = \bigcup_i (A_i R_i),
\]
where each region $R_i$ corresponds to one affine piece of $f$. Each projection is computationally
tractable using linear programming or convex elimination.

\subsubsection*{3. Low-Dimensional Nonlinear Maps}
If the domain of $f$ is low-dimensional (e.g., two to five assets) and $R$ is defined by linear or
quadratic constraints, the pushforward can be approximated numerically by discretizing $x$,
evaluating $y=f(x)$, and taking the convex hull of the resulting $(y,z)$ pairs. This method provides
accurate permissible-region estimates for practical portfolio problems.

\subsubsection*{4. Linear or Convex Codomain Constraints}
If one wishes to verify inclusion $f_!R \subseteq C$ and $C$ is convex (for instance, linear
alignment constraints), the check reduces to
\[
\forall (x,z) \in R, \quad f(x) \in C_z, \quad C_z = \{y \mid (y,z) \in C\}.
\]
When $f$ is continuous and $C_z$ is convex, this inclusion can be verified via convex relaxation or
interval enclosures of $f$, such as McCormick envelopes \cite{McCormick} or semidefinite relaxations.

\subsubsection*{5. Linearization via Tangent Cones}
When $R$ describes local alignment conditions, the pushforward can be approximated by the action
of the Jacobian of $f$ on tangent cones:
\[
T_{(f(x),z)}(f_!R) = Df(x)\, T_{(x,z)}R.
\]
This linearized pushforward is analytically tractable and corresponds to standard sensitivity
analysis or first-order propagation in differentiable portfolio transformations.

\paragraph{Financial Interpretation.} In portfolio settings, $R$ typically encodes linear alignment
constraints (e.g., exposure matching), while $f$ represents a nonlinear optimization-based
re-implementation (e.g., fee or tax optimization). When $f$ is continuous and proper, pushforwards
can be reliably delineated through sampling or convexification. However, certain realistic cases, such as discrete constraints arising from lot size requirements, may not be handled well.

\subsection{Computing Optimal Re-implementations}

The optimization problem defined in Section~\ref{sec:optimal-reimpl} is often
a standard convex optimization problem:
\[
\max_{y \in \Delta^m} U_\mathcal{B}(y) \quad \text{subject to} \quad \|g_\A(x)
- g_\mathcal{B}(y)\|_V \leq \epsilon
\]
(and $y \in \Delta^m$, which is polyhedral). This is often a Linear Program
(LP), Quadratic Program (QP), or Second-Order Cone Program (SOCP), all of
which are efficiently solvable \cite{BoydVandenberghe, Bertsekas}.

\subsection{A Concrete Worked Example}

To make the pullback and pushforward operations concrete, consider a simple system:

\begin{itemize}
    \item \textbf{Ambient Spaces:} Let the hub be $\Delta^2 = \{(x_0, x_1, x_2) \in \mathbb{R}^3 \mid \sum x_i = 1, x_i \ge 0\}$ (3 assets) and the spoke be $\Delta^1 = \{(y_0, y_1) \in \mathbb{R}^2 \mid y_0+y_1 = 1, y_i \ge 0\}$ (2 assets).

    \item \textbf{Object $K_1$:} A permissible space with a cap on asset 0:
    \[ K_1 = \{x \in \Delta^2 \mid x_0 \le 0.5\} \]

    \item \textbf{Object $K_2$ and $K_3$:} Let $K_2 = K_3 = \Delta^1$.

    \item \textbf{Horizontal Morphism $f: K_1 \to K_2$:} An aggregation map that combines assets 0 and 1:
    \[ f(x_0, x_1, x_2) = (x_0 + x_1, x_2) = (y_0, y_1) \]
    This is continuous and maps $K_1$ to $K_2$ (since $x_0+x_1+x_2 = 1 \implies y_0+y_1=1$).
\end{itemize}

\paragraph{Computing the Pullback $f^*S$ (Pre-alignment)}

Let's define a spoke-space alignment $S \subseteq K_2 \times K_3$ that limits the aggregated asset $y_0$ relative to some benchmark $z_0$:
\[ S = \{(y, z) \in \Delta^1 \times \Delta^1 \mid y_0 \le 2z_0 \} \]
By Definition~\ref{def:pullback}, the pullback $f^*S$ finds all hub portfolios $x \in K_1$ that satisfy this alignment after re-implementation:
\begin{align*}
    f^*S &= \{(x, z) \in K_1 \times K_3 \mid (f(x), z) \in S\} \\
         &= \{(x, z) \in K_1 \times \Delta^1 \mid ( (x_0+x_1, x_2), z) \in S\} \\
         &= \{(x, z) \in K_1 \times \Delta^1 \mid (x_0+x_1) \le 2z_0 \}
\end{align*}
The resulting computation is simple, yielding a new closed relation on the hub space.

\paragraph{Computing the Pushforward $f_!R$ (Forward Propagation)}

Now, define a hub-space alignment $R \subseteq K_1 \times K_3$ that constrains asset 0 against the benchmark $z_0$:
\[ R = \{(x, z) \in K_1 \times \Delta^1 \mid x_0 \ge z_0 \} \]
By Definition~\ref{def:pushforward}, the pushforward $f_!R$ finds all spoke-space pairs $(y, z)$ that can be ``hit'' by an aligned hub portfolio $x \in R$:
\[ f_!R = \{(y, z) \in K_2 \times K_3 \mid \exists x \in K_1 : f(x) = y \text{ and } (x, z) \in R\} \]
To compute this, we must eliminate $x$ subject to the constraints:
\begin{enumerate}
    \item $x \in K_1 \implies x_0 \le 0.5$ (and $x \in \Delta^2$)
    \item $f(x) = y \implies x_0+x_1 = y_0$ and $x_2 = y_1$
    \item $(x, z) \in R \implies x_0 \ge z_0$
    \item $x_i \ge 0 \implies x_0 \ge 0, x_1 \ge 0, x_2 \ge 0$
\end{enumerate}
From (2) and (4), $x_1 = y_0 - x_0 \ge 0 \implies x_0 \le y_0$.
We need to find if there \emph{exists} an $x_0$ that satisfies all conditions on $x_0$:
\[ x_0 \ge 0 \quad \text{and} \quad x_0 \le 0.5 \quad \text{and} \quad x_0 \ge z_0 \quad \text{and} \quad x_0 \le y_0 \]
Such an $x_0$ exists if and only if the ``floor'' is less than or equal to the ``ceiling'':
\[ \max(0, z_0) \le \min(0.5, y_0) \]
Thus, the pushforward (the result of the projection) is:
\[ f_!R = \{(y, z) \in \Delta^1 \times \Delta^1 \mid \max(0, z_0) \le \min(0.5, y_0) \} \]
As guaranteed by Theorem~\ref{thm:proper-closed}, this is a closed set, but its definition is not completely trivial and emerges from the quantifier elimination.

\clearpage

\part{Axiomatization and Extensions}

\clearpage

\section*{Guide to Models}
\addcontentsline{toc}{section}{Guide to Models}

\begin{table}[h]
\centering
\renewcommand{\arraystretch}{1.4}
\begin{tabular}{p{0.22\textwidth} p{0.48\textwidth} p{0.22\textwidth}}
\toprule
\textbf{Model Variant} & \textbf{Primary Use Case} & \textbf{Mathematical Cost} \\
\midrule
\textbf{$\HSSimp$ (Reference)} \newline \textit{Part I} 
& \textbf{Architecture \& Safety.} The baseline model for defining objects, guaranteeing properness, and establishing audit safety (Lax BC). 
& \textbf{Low} \newline (Point-set Topology) \\
\midrule
\textbf{$\mathbb{S}\mathbf{pan}$ (Evidence)} \newline \textit{Sec. \ref{sec:2-cells}} 
& \textbf{Quantitative Audit.} Replaces binary ``aligned'' checks with quantitative evidence (e.g., tracking error magnitude, fee budgets). 
& \textbf{Medium} \newline (Bicategories) \\
\midrule
\textbf{$\HSSimp/B$ (Personas)} \newline \textit{Sec. \ref{sec:personas}} 
& \textbf{Mass Customization.} Modeling \textit{families} of portfolios indexed by age or risk parameters (e.g., Glidepaths, Target Date Funds). 
& \textbf{Medium} \newline (Slice Categories) \\
\midrule
\textbf{DOTS (Menus)} \newline \textit{Sec. \ref{sec:DOTS}} 
& \textbf{Strategy Design.} Generating ``menus'' of permissible spokes for human selection or OCIO workflows (Actions $K \odot R$). 
& \textbf{High} \newline (Operads) \\
\midrule
\textbf{HSP-r (Stochastic)} \newline \textit{Sec. \ref{sec:Polish}} 
& \textbf{Execution \& Liquidity.} Handling randomness in solvers (rounding, tie-breaking) and messy, non-compact state spaces. 
& \textbf{High} \newline (Measure Theory) \\
\bottomrule
\end{tabular}
\caption{A guide to the different models presented in this paper. The Reference Model ($\HSSimp$) provides the architectural foundation, while the extensions in Parts II and III adapt the framework to specific industrial constraints.}
\label{tab:model_guide}
\end{table}

\clearpage

\section{An Axiomatic Framework}
\label{sec:abstract}

The double category $\HSSimp$ constructed in Part~I is a \emph{thin equipment} (framed bicategory) in the sense of Shulman~\cite{Shulman} and Grandis--Par\'e~\cite{GrandisPare}: companions are graphs of maps, conjoints are their converses, and the adjunction $f_! \dashv f^*$ follows from the unit/counit inequalities.  This is standard categorical machinery, and the axioms are stated below for completeness.

The interest of the axiomatic treatment lies elsewhere --- in identifying what \emph{fails} when the framework is extended to the probabilistic setting.  In $\HSPr$ (Part~III), horizontal composition of stochastic kernels introduces error accumulation: the composite pushforward $(g \circ f)_!$ is contained in, but generally strictly smaller than, the staged pushforward $g_!(f_!)$.  The companion construction, which requires strict functoriality, cannot be carried out.  The pushforward and pullback operations become \emph{primitive predicate transformers} rather than compositions with a distinguished vertical arrow.  Section~\ref{subsec:lax-axioms} develops this observation into a ``Lax Logical Framework'' that identifies exactly which axioms survive (adjunction, sub-functoriality, lax Beck--Chevalley, lax Frobenius) and which do not (companions, conjoints, strict functoriality).  This is the section's contribution: not the axioms themselves, but the precise diagnosis of where they break.

\paragraph*{Guide to this section.}  Readers familiar with equipments may skim Sections~\ref{sec:axioms-subsec}--\ref{sec:equipment-frobenius} and proceed directly to Section~\ref{subsec:lax-axioms}.  Readers whose primary interest is the concrete models may skip this section entirely; the variant semantics in Sections~\ref{sec:2-cells}--\ref{sec:DOTS} and Part~III do not depend on the abstract axioms.

\subsection{Axioms}
\label{sec:axioms-subsec}
\noindent\textbf{Notation.} We write $\meet$ for the binary meet in the vertical poset; in the concrete model $\HSSimp$, this is set-theoretic intersection $\cap$.

We assume a \emph{thin} double category (a framed bicategory / proarrow equipment)
\[
\mathbb{D} = (\mathrm{Obj}, \mathbb{D}_h, \mathbb{D}_v)
\]
consisting of (with the financial interpretation in $\HSSimp$ noted parenthetically):
\begin{itemize}
\item Objects $A,B,C,\dots$ \emph{(permissible portfolio spaces)}.
\item Horizontal arrows $f:A\to B$ forming a category $\mathbb{D}_h$ (compose strictly; identities $1_A$) \emph{(re-implementations)}.
\item Vertical arrows $R:A\dashrightarrow B$ forming, for each $(A,B)$, a poset $\mathbb{D}_v(A,B)$
with vertical composition $S\vcomp R$ and identities $\id_A:A\dashrightarrow A$ \emph{(alignment relations; the poset order is inclusion of constraints)}.
\end{itemize}
Thinness means all $2$-cells are inclusions in these posets: if a 2-cell $\alpha$ exists, it is unique.
 
\[
    \begin{tikzcd}
    A \arrow[r, "f"] \arrow[d, "R"', dashed] & B \arrow[d, "S", dashed] \\
    C \arrow[r, "g"'] & D
    \arrow[from=1-1, to=2-2, phantom, "\Downarrow \alpha"]
    \end{tikzcd}
\]
 
\begin{definition}[Vertical daggered allegory-like structure]\label{def:vertical}
Assume:
\begin{enumerate}
\item[(V1)] Each hom $\mathbb{D}_v(A,B)$ is a poset $(\le)$ with binary meets $R\meet R'$ (in $\HSSimp$, $\meet$ is set-theoretic intersection $\cap$).
\item[(V2)] Vertical composition is monotone and distributes over meets when one factor is a companion or conjoint: for any vertical morphisms $R, R'$ and any companion $\langle f\rangle$ or conjoint $[f]$,
\[
\langle f\rangle\vcomp (R\meet R') = (\langle f\rangle\vcomp R)\meet(\langle f\rangle\vcomp R'), \qquad
(R\meet R')\vcomp [f] = (R\vcomp [f])\meet(R'\vcomp [f]),
\]
and similarly with the roles of companion and conjoint exchanged.  Meet-distribution \emph{fails} for the composition of two arbitrary vertical morphisms (see Remark~\ref{rem:V2-scope}).
\item[(V3)] A dagger (converse) $(-)^\dagger:\mathbb{D}_v(A,B)\to\mathbb{D}_v(B,A)$ with
$(R^\dagger)^\dagger=R$ and $(S\vcomp R)^\dagger=R^\dagger\vcomp S^\dagger$.
\item[(V4)] The modular law (Dedekind's formula): for composable $S$ and $R$,
\[
(S \vcomp R) \meet T \;\le\; S \vcomp (R \meet (S^\dagger \vcomp T)).
\]
\end{enumerate}
In $\HSSimp$, axiom V4 is a standard identity of the calculus of relations: if $(y,z) \in (S \circ R) \cap T$, the witness $x$ with $(y,x) \in R$ and $(x,z) \in S$ also satisfies $(y,z) \in T$; since $(y,z) \in T$ and $(x,z) \in S$ (equivalently $(z,x) \in S^\dagger$), we get $(y,x) \in S^\dagger \circ T$, hence $(y,x) \in R \cap (S^\dagger \circ T)$, placing $(y,z)$ in $S \circ (R \cap (S^\dagger \circ T))$.  The modular law holds in every model of this paper whose vertical morphisms are relations (including HSP-r and the three probabilistic compliance frameworks), but it is only \emph{used} in the deterministic equipment setting where companions exist (Theorem~\ref{thm:Frob}).  The probabilistic models lack companions and use only the forward Frobenius inclusion, which does not require V4.
\end{definition}

On allegories as an abstraction of categories of relations, see \cite{FreydScedrov}; for the logical foundation of adjoint quantifiers ($f_! \dashv f^*$), see \cite{Lawvere}. For restriction categories, see \cite{CockettLack}.

\begin{remark}[Why V2 is restricted to companions and conjoints]
\label{rem:V2-scope}
In $\HSSimp$, meet-distribution holds when one factor is a function graph: $\Graph(f) \circ (R \cap R') = (\Graph(f) \circ R) \cap (\Graph(f) \circ R')$ and $(R \cap R') \circ \Graph(f) = (R \circ \Graph(f)) \cap (R' \circ \Graph(f))$, because a function graph assigns a unique image to each element.  The unrestricted version (distribution for arbitrary vertical composites) \emph{fails} for the same reason it fails in $\mathbf{Rel}$: the standard counterexample ($R = \{(a,b)\}$, $R' = \{(a,c)\}$, $T = \{(b,d),(c,d)\}$, giving $T \circ (R \cap R') = \emptyset \subsetneq \{(a,d)\} = (T \circ R) \cap (T \circ R')$) applies to finite compact spaces with closed relations.
\end{remark}
 
\begin{definition}[Equipment: companions and conjoints]\label{def:equipment}
For each horizontal $f:A\to B$ there exist verticals
\[
\langle f\rangle: A\into B \qquad \text{and} \qquad [f]: B\into A
\]
(the \emph{companion} and \emph{conjoint}) such that:
\begin{enumerate}
\item[(E1)] $\langle f\rangle \adj [f]$ in the vertical posets:
\[
\id_A \le [f]\vcomp\langle f\rangle
\qquad\text{and}\qquad
\langle f\rangle\vcomp [f] \le \id_B.
\]
\item[(E2)] Functoriality: $\langle 1_A\rangle=\id_A$, $[1_A]=\id_A$, and standard mates
compatibilities for horizontal composition (via cartesian comparison $2$-cells).
\item[(E3)] (Required for Frobenius) Dagger-compatibility: $\langle f\rangle^\dagger=[f]$ and $[f]^\dagger=\langle f\rangle$. This axiom is not needed for the adjunction (Proposition~\ref{prop:adjunction}) or the Beck--Chevalley condition (Theorem~\ref{thm:BC}), but is essential for the Frobenius reciprocity law (Theorem~\ref{thm:Frob}), which requires the modular law (Dedekind formula) for vertical composition.
\end{enumerate}
\end{definition}

For example, in $\HSSimp$ the companion of a horizontal morphism $f$ is the graph $\Graph(f)$ and the conjoint is the reverse graph $\Graph(f)^\dagger$, related by E3.

\begin{lemma}[Mates comparison]
\label{lemma:mates-comp}
For any commuting square
\[
\begin{tikzcd}[column sep=large, row sep=large]
A \arrow[r, "g"] \arrow[d, "f'"'] & B \arrow[d, "f"] \\
C \arrow[r, "h"'] & D
\end{tikzcd}
\]
(so $f \circ g = h \circ f'$), and any $R: B \into Z$, we have $f'_!(g^* R) \le h^*(f_! R)$.
\end{lemma}

\begin{proof}
By the adjunction $h_! \dashv h^*$ (which follows from E1), the inequality $f'_!(g^* R) \le h^*(f_! R)$ is equivalent to $h_!(f'_!(g^* R)) \le f_! R$.
By functoriality (E2), $h_! \circ f'_! = (h \circ f')_!$ and $f_! \circ g_! = (f \circ g)_!$.
Since $f \circ g = h \circ f'$, we have $h_! \circ f'_! = f_! \circ g_!$, so the LHS becomes $f_!(g_!(g^* R))$.
By the counit inequality (E1), $g_!(g^* R) \le R$, and by monotonicity of $f_!$:
\[
f_!(g_!(g^* R)) \le f_!(R).
\]
The inequality holds.
\end{proof}

\subsection{Pushforward and pullback along horizontals}
\label{sec:equipment-pushpull}
Fix an object $Z$. For $f:A\to B$ define monotone maps between vertical hom-posets:
\[
f_!: \mathbb{D}_v(A,Z)\to \mathbb{D}_v(B,Z),
\quad
f_!(R):=\langle f\rangle\vcomp R,
\qquad
f^{\!*}: \mathbb{D}_v(B,Z)\to \mathbb{D}_v(A,Z),
\quad
f^{\!*}(S):=[f]\vcomp S.
\]
 
\begin{center}
\begin{tikzcd}[column sep=large, row sep=large]
A \arrow[r, "f"] \arrow[d, dashed, "R"'] & B \arrow[d, dashed, "f_!R"] \\
Z & Z
\end{tikzcd}
\qquad
\begin{tikzcd}[column sep=large, row sep=large]
A \arrow[r, "f"] \arrow[d, dashed, "f^{\!*}S"'] & B \arrow[d, dashed, "S"] \\
Z & Z
\end{tikzcd}
\end{center}
 
\begin{proposition}[Adjunction $f_!\adj f^{\!*}$]\label{prop:adjunction}
For all $R:A\into Z$ and $S:B\into Z$,
\[
R \le f^{\!*}S \quad\Longleftrightarrow\quad f_!R \le S.
\]
\emph{Proof.}
$(\Rightarrow)$\; If $R\le [f]\vcomp S$, then by monotonicity of vertical composition $\langle f\rangle\vcomp R \le \langle f\rangle\vcomp [f]\vcomp S \le \id_B\vcomp S=S$, where the second step uses the counit $\langle f\rangle\vcomp [f] \le \id_B$ from (E1).
$(\Leftarrow)$\; If $\langle f\rangle\vcomp R\le S$, then by monotonicity
$[f]\vcomp\langle f\rangle\vcomp R \le [f]\vcomp S$, and the unit $\id_A \le [f]\vcomp\langle f\rangle$ from (E1) gives $R\le [f]\vcomp\langle f\rangle\vcomp R \le [f]\vcomp S$. \qed
\end{proposition}
 
\begin{proposition}[Functoriality]\label{prop:abstract-functoriality}
For composable horizontals $A\xrightarrow{f}B\xrightarrow{g}C$,
\[
(g\circ f)_! = g_!\vcomp f_!,
\qquad
(g\circ f)^{\!*} = f^{\!*}\vcomp g^{\!*}.
\]
\emph{Proof.} By (E2) and associativity of vertical composition. \qed
\end{proposition}
 
\subsection{Beck--Chevalley}

\begin{proposition}[Lax Beck--Chevalley]
Let $f, g, f', h$ be horizontal morphisms in the equipment $\mathbb{D}$ forming a commuting square as in Lemma~\ref{lemma:mates-comp} (so $f \circ g = h \circ f'$). For any vertical morphism $R: B \dashrightarrow Z$, we have the inequality:
\[
f'_! (g^* R) \leq h^* (f_! R).
\]
This is the content of Lemma~\ref{lemma:mates-comp}.
\end{proposition}
 
\begin{definition}\label{def:cartesian}
A commuting square of horizontals
\[
\begin{tikzcd}[column sep=large, row sep=large]
A \arrow[r, "g"] \arrow[d, "f'"'] & B \arrow[d, "f"] \\
C \arrow[r, "h"'] & D
\end{tikzcd}
\]
is \emph{cartesian} if $f'_!(g^*R) = h^*(f_!R)$ for every $R:B\into Z$.
\end{definition}

\begin{theorem}[Strict Beck--Chevalley]\label{thm:BC}
For a cartesian horizontal square (cf. \cite{CarboniWalters, Wood}) and every $R:B\into Z$,
\[
f'_! (g^* R) = h^* (f_! R).
\]
\end{theorem}
 
\subsection{Frobenius reciprocity}
\label{sec:equipment-frobenius}

\begin{theorem}[Frobenius]\label{thm:Frob}
For any horizontal $f:A\to B$, any $R:A\into Z$ and $S:B\into Z$,
\[
f_!\bigl(R \meet f^{\!*}S\bigr) \;=\; f_!R \meet S.
\]
\emph{Proof.} Work in $\mathbb{D}_v(B,Z)$. Using (V2) and (E1),
\[
\begin{aligned}
f_!\bigl(R \meet f^{\!*}S\bigr)
&= \langle f\rangle \vcomp \bigl(R \meet ([f]\vcomp S)\bigr) \\
&= (\langle f\rangle \vcomp R) \meet \bigl(\langle f\rangle \vcomp [f] \vcomp S\bigr) \\
&\le (\langle f\rangle \vcomp R) \meet (\id_B \vcomp S) = f_!R \meet S.
\end{aligned}
\]
For the reverse inequality, we use axiom V4 (the modular law):
\[
(\langle f\rangle \vcomp R) \meet S
\;\le\;
\langle f\rangle \vcomp \bigl(R \meet \langle f\rangle^\dagger \vcomp S\bigr).
\]
By E3, $\langle f\rangle^\dagger = [f]$, so the right-hand side is
$\langle f\rangle \vcomp (R \meet [f] \vcomp S) = f_!(R \meet f^{\!*}S)$.
Since $f_!R = \langle f\rangle \vcomp R$, the left-hand side is $f_!R \meet S$.
Thus $f_!R \meet S \le f_!(R \meet f^{\!*}S)$, and equality holds. \qed
\end{theorem}
 
\begin{remark}
These abstract arguments use:
(i) vertical hom-posets with binary meets (V1);
(ii) meet-preservation of vertical composition (V2, for companions only);
(iii) companions/conjoints with unit and counit inequalities (E1);
(iv) the modular law (V4) together with dagger-compatibility (E3) for the Frobenius reverse direction.
\end{remark}

\subsection{Lax Semantics and the Failure of Companions}
\label{subsec:lax-axioms}

The axioms of a standard equipment (Definition~\ref{def:equipment}) imply strict functoriality for the pushforward and pullback operations: $(g \circ f)_! = g_! \circ f_!$. However, as seen in the probabilistic extensions (Section~\ref{sec:Polish}, specifically Theorems \ref{thm:pushforward-comp} and \ref{thm:pullback-comp}), practical risk management often obeys only \emph{lax} functoriality, where composition leads to an accumulation of error bounds or risk budgets.

We therefore outline a generalized framework where the rigorous equality of the equipment is relaxed to an inequality.

  \begin{definition}[Lax Logical Framework]
  \label{def:lax-logical}
  A \textbf{Lax Logical Framework} consists of a double category structure where:
  \begin{enumerate}
      \item The horizontal category $\mathbb{D}_h$ and vertical posets
        $\mathbb{D}_v$ are defined as in Definition~\ref{def:vertical}.
      \item For every horizontal morphism $f: A \to B$, there exist
        monotonic mappings
        $f_!: \mathbb{D}_v(A, \cdot) \to \mathbb{D}_v(B, \cdot)$ and
        $f^*: \mathbb{D}_v(B, \cdot) \to \mathbb{D}_v(A, \cdot)$.
      \item \textbf{Adjunction:} $f_! \dashv f^*$ holds as a Galois
        connection: $f_! R \le S \iff R \le f^* S$.
      \item \textbf{Unitality:} $(\mathrm{id}_A)_! R = R$ and
        $(\mathrm{id}_A)^* S = S$ for all $R, S$.
      \item \textbf{Lax Composition (Sub-functoriality):}
      \[
      (g \circ f)_! R \;\le\; g_!(f_! R) \quad \text{and} \quad
      f^*(g^* S) \;\le\; (g \circ f)^* S
      \]
      Note: these two inequalities are mates under the adjunction in
      item~(3) and are therefore equivalent.

      \emph{Pushforward $\Rightarrow$ pullback.}  Suppose $(g \circ f)_! R \le g_!(f_! R)$ for all $R$.  We wish to show $f^*(g^* S) \le (g \circ f)^* S$ for all $S$.  From $(g \circ f)_! R \le g_!(f_! R)$, apply $(g \circ f)^*$ to both sides and use the unit $R \le (g \circ f)^*((g \circ f)_! R)$ to obtain $R \le (g \circ f)^*(g_!(f_! R))$.  By the adjunction $g_! \dashv g^*$, $g_!(f_! R) \le S$ iff $f_! R \le g^* S$, and similarly $f_! R \le g^* S$ iff $R \le f^*(g^* S)$.  Chaining: whenever $(g \circ f)_! R \le S$, we get $R \le f^*(g^* S)$; but $(g \circ f)_! R \le S$ iff $R \le (g \circ f)^* S$ by item~(3), so $f^*(g^* S) \le (g \circ f)^* S$.

      \emph{Pullback $\Rightarrow$ pushforward.}  Suppose $f^*(g^* S) \le (g \circ f)^* S$ for all $S$.  We wish to show $(g \circ f)_! R \le g_!(f_! R)$.  Set $S = g_!(f_! R)$.  Then $(g \circ f)^* S \ge f^*(g^* S) \ge f^*(f_! R) \ge R$, where the last step uses the unit of $f_! \dashv f^*$.  By item~(3), $R \le (g \circ f)^* S$ iff $(g \circ f)_! R \le S = g_!(f_! R)$.

      The dual framework
      (with reversed inequalities, modeling \emph{super-functorial}
      error accumulation) can be obtained by swapping $f_!$ and $f^*$.
  \end{enumerate}
  \end{definition}

\begin{remark}[The unavailability of Companions]
In a strict equipment, the pushforward is defined constructively via the companion: $f_! R = \langle f \rangle \vcomp R$. This definition forces strict functoriality because the assignment $f \mapsto \langle f \rangle$ is a functor from the horizontal category to the vertical bicategory.

In a Lax Logical Framework, \textbf{pushforwards and pullbacks cannot be defined via companions and conjoints}.
Proof: Suppose $f_! R = \langle f \rangle \vcomp R$ for some vertical arrow $\langle f \rangle$. Then:
\[
(g \circ f)_! \Delta_A = \langle g \circ f \rangle \quad \text{and} \quad g_!(f_! \Delta_A) = \langle g \rangle \vcomp \langle f \rangle.
\]
If the physics of the system dictates a strict inequality (e.g., $(g \circ f)_! < g_! \circ f_!$ due to convexity accumulation), then $\langle g \circ f \rangle \neq \langle g \rangle \vcomp \langle f \rangle$. This means the mapping from horizontal maps to vertical relations is not a functor, destroying the structural coherence required for a standard equipment.

Consequently, in lax settings such as HSP-r, the operations $f_!$ and $f^*$ must be treated as \textbf{primitive predicate transformers} rather than as compositions with distinguishing morphisms. We lose the ability to represent the function $f$ as a ``perfect'' relation, because in a lax system, no relation is perfectly composable.
\end{remark}

  \begin{remark}[Lax Beck--Chevalley: A Guiding Principle]
  \label{rem:lax-BC}
  In a Lax Logical Framework, the strict Beck--Chevalley equality
  $f'_!(g^* S) = h^*(f_! S)$ of Theorem~\ref{thm:Beck--Chevalley}
  generally fails.  However, the \emph{safe} (conservative) direction
  \begin{equation}\label{eq:lax-BC}
  f'_!(g^* S) \;\le\; h^*(f_! S)
  \end{equation}
  can survive under suitable hypotheses.  Informally, if the laxness
  is \emph{sub-functorial in the pushforward direction}
  ($(g \circ f)_! R \le g_!(f_! R)$, i.e., the composite pushforward
  under-approximates the staged pushforward), and if the
  pullback--pushforward adjunction is preserved at each stage,
  then~\eqref{eq:lax-BC} holds by the same mates argument as in
  the strict case, with equalities replaced by inequalities throughout.

  Concrete instances of this principle, with precise hypotheses and proofs,
  appear in the probabilistic setting:
  \begin{itemize}
  \item \textbf{Safety Radius:}  Theorem~\ref{thm:pullback-comp}
    establishes $P^{*,\delta}(Q^{*,\epsilon} S) \subseteq
    (Q \circ P)^{*,\delta+\epsilon} S$, which is the pullback
    analogue of~\eqref{eq:lax-BC} with additive risk budgets.
  \item \textbf{HDR:}  The Lax Beck--Chevalley for
    $\epsilon$-pushforward (Theorem~\ref{thm:pushforward-comp})
    gives $(Q \circ P)_!^{\delta+\epsilon} R \subseteq
    Q_!^\epsilon(P_!^\delta R)$ under the Convolution Stability
    Condition.
  \end{itemize}
  We do not formulate a general abstract theorem here, as the
  precise hypotheses depend on the choice of $\epsilon$-support
  semantics (Safety Radius vs.\ HDR vs.\ Wasserstein).
  \end{remark}

\subsubsection*{Adjunction as the Primary Safety Guarantee}

While the failure of the companion construction prevents us from treating horizontal morphisms as ``perfect'' relations, the adjunction between pushforward and pullback survives. In a lax framework, this adjunction becomes the primary definition of system safety.

\begin{proposition}[Adjunction in a Lax Framework]
Even in a Lax Logical Framework where horizontal composition is sub-functorial (i.e., $(g \circ f)_! \le g_! \circ f_!$), the operations for a \emph{single} morphism $f$ form a Galois connection:
\[
f_! R \le S \iff R \le f^* S
\]
\end{proposition}

\begin{remark}
This persistence of adjunction is critical. As established in the probabilistic model (Theorem~\ref{thm:safety-adjunction}), this equivalence provides the logical justification for ``Stress Testing''. It guarantees that checking if a hub portfolio $x$ lies in the preimage $f^* S$ is sufficient to ensure that the resulting spoke portfolio $f(x)$ lands in $S$, even if $f$ is noisy or imperfect.
\end{remark}

\subsubsection*{Lax Frobenius: The Cost of Leakage}

In a strict equipment (like $\HSSimp$), the Frobenius reciprocity law is an equality: $f_!(R \meet f^*S) = f_!R \meet S$. This implies that filtering a portfolio for compliance before re-implementation is identical to filtering it after.

In a Lax Logical Framework, ``leakage'' of risk or mass (as seen in the HSP-r model) breaks this equality. It is replaced by a one-way inclusion that dictates the safe order of operations.

\begin{proposition}[Lax Frobenius / Conservative Filtering]
In a Lax Logical Framework, filtering upstream is strictly more conservative than filtering downstream:
\[
f_!(R \meet f^* S) \;\le\; f_!R \meet S
\]
\end{proposition}

\begin{proof}
Let $(y, z) \in f_!(R \meet f^*S)$. Then there exists $x$ with $(x, z) \in R \meet f^*S$ and $f(x) = y$. Since $(x, z) \in R$, we have $(f(x), z) = (y, z) \in f_!R$. Since $(x, z) \in f^*S$, we have $(f(x), z) = (y, z) \in S$. Therefore $(y, z) \in f_!R \meet S$.

In a lax setting (e.g., HSP-r), the inclusion is strict because the probabilistic pushforward $f_!^\epsilon$ ``leaks'' mass: a hub portfolio satisfying both $R$ and $f^*S$ with individual risk budgets may fail the joint budget under the composite pushforward, while a spoke portfolio could ``luckily'' satisfy $S$ after the fact.
\end{proof}

\emph{Intuition.}
The term $f_!(R \meet f^* S)$ represents the set of spoke portfolios derived from hub portfolios that are \emph{guaranteed} to satisfy $S$ (pre-verified). The term $f_!R \meet S$ represents the set of spoke portfolios that happen to satisfy $S$ after the fact.

In a lax system (e.g., probabilistic), there exist ``lucky'' outcomes where a hub portfolio with a 10\% failure risk generates a compliant spoke (accepted by the RHS), but would have been rejected by the strict upstream check (LHS). To maintain a valid audit trail, the framework must enforce the LHS (upstream filtering) to prevent this ``tail risk leakage''.

\clearpage  \section{Evidence of Alignment}
\label{sec:2-cells}

\subsection{Sketch: Extension to Non-Thin 2-Cells}

The axiomatic framework presented above assumes a \textbf{thin double category}, where the vertical arrows in any given hom-set $\mathbb{D}_v(A, B)$ form a poset. This is reflected in the main model $\HSSimp$, where 2-cells are simply inclusions $Graph(g) \circ R \subseteq S \circ Graph(f)$. Such a 2-cell is unique if it exists, signifying only \emph{whether} a re-implementation preserves an alignment.

This, however, is a simplification of practical financial reality. Alignment may be binary; verification of alignment almost never is. It is measured by continuous metrics---tracking error, factor loading differences, fee budgets---and a portfolio is not just ``aligned'' but ``aligned with 35 basis points of tracking error.''

In this particular case, one could replace the binary alignment relation $R$ with a family of relations $R_\epsilon$, indexed by a tracking error tolerance $\epsilon$. However, such \emph{ad hoc} solutions become unsatisfactory in the case of more complex alignment relations, where it is desirable to organize verification data more systematically.

A natural extension, therefore, is to accommodate \textbf{non-thin} (or ``non-locally-posetal'') 2-cells. This would allow the framework to model situations where there are \emph{multiple, distinct ways} for a re-implementation to be compatible with an alignment.

A non-thin 2-cell $\alpha$ could represent a \emph{specific, quantitative measure} or \emph{method} of this compatibility. For example, consider two different optimal re-implementation strategies, $\alpha_1$ (minimizing fees) and $\alpha_2$ (minimizing tracking error). Both might satisfy the same alignment inclusion, but they represent different, non-equivalent processes. In a non-thin framework, $\alpha_1$ and $\alpha_2$ could be distinct 2-cells filling the same square.

To accommodate this, the axiomatic framework would be modified:

\begin{enumerate}
    \item \textbf{Drop Thinness:} The assumption that the double category is thin would be removed.

    \item \textbf{2-Cell Hom-Sets:} For any square of 1-morphisms ($f, g, R, S$), the 2-cells $\alpha$ filling that square would form a \emph{set} (or a new poset, in a ``locally posetal'' model) rather than being a single boolean fact.

    \begin{center}
    \begin{tikzcd}
    A \arrow[r, "f"] \arrow[d, "R"', dashed] & B \arrow[d, "S", dashed] \\
    C \arrow[r, "g"'] & D
    \arrow[from=1-1, to=2-2, phantom, "\Downarrow \alpha"]
    \end{tikzcd}
    \end{center}

    \item \textbf{Composition and Interchange:} This change would require introducing explicit axioms for the horizontal and vertical composition of these 2-cells. Crucially, the \textbf{interchange law}, which is automatic in the thin HS model, would become a non-trivial axiom that must be explicitly satisfied.
\end{enumerate}

\subsection{``Evidence'' 2-Cells via the Double Category of Spans}

We now sketch a concrete instance of the non-thin framework just described, which is the most natural candidate for practical applications.

We replace the category of relations with a category of \textbf{spans}. We assume an ambient category $\mathcal{E}$ with pullbacks (e.g., $\mathbf{Top}$ or $\mathbf{CompHaus}$, consistent with the proper maps in the main text).

\subsubsection*{Vertical Morphisms as Spans of Evidence}

In this richer framework, a vertical morphism (alignment) from permissible space $A$ to $B$ is not a subset of $A \times B$, but a \textbf{span} in $\mathcal{E}$:
\[
A \xleftarrow{p} E \xrightarrow{q} B
\]
\textbf{Interpretation:}
\begin{itemize}
   \item $E$ is the \textbf{Space of Evidence} (or the space of alignment certificates).
   \item A point $e \in E$ is a specific witness that portfolio $p(e) \in A$ is aligned with portfolio $q(e) \in B$.
   \item The physical alignment relation from the thin model is recovered by taking the joint image $\langle p, q \rangle(E) \subseteq A \times B$.
\end{itemize}

\begin{example}[Tracking Error with Evidence]
Let $A$ be a space of funds and $B$ a benchmark. Let
\[
E = \{(a, b, \delta) \in A \times B \times \mathbb{R}_{\ge 0} \mid \|g_A(a) - g_B(b)\| \le \delta\}.
\]
The projection maps $p(a,b,\delta)=a$ and $q(a,b,\delta)=b$ verify the portfolios, while the third component $\delta$ retains the \emph{degree} of alignment. This allows us to distinguish between ``barely aligned'' and ``perfectly aligned.''
\end{example}

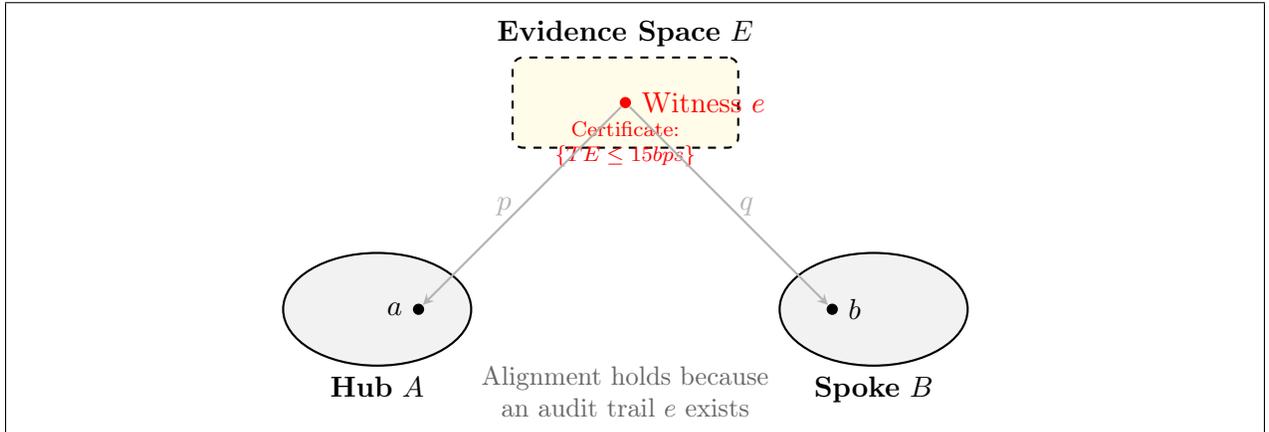
\begin{figure}[h]
\centering
\begin{tikzpicture}[
    scale=1.1,
    >=stealth,
    space/.style={draw, ellipse, minimum width=2.5cm, minimum height=1.5cm, fill=gray!10, thick},
    evidence/.style={draw, rectangle, rounded corners, minimum width=3cm, minimum height=1.2cm, fill=yellow!10, dashed, thick},
    point/.style={circle, fill, inner sep=1.5pt},
    map/.style={->, thick, gray!60}
]

    \node[space, label=below:\textbf{Hub} $A$] (A) at (0,0) {};
    \node[space, label=below:\textbf{Spoke} $B$] (B) at (6,0) {};

    \node[evidence, label=above:\textbf{Evidence Space} $E$] (E) at (3, 2.5) {};

    \node[point, label=left:$a$] (pa) at (0.5, 0) {};
    \node[point, label=right:$b$] (pb) at (5.5, 0) {};

    \node[point, red, label=right:\textcolor{red}{Witness $e$}] (pe) at (3, 2.5) {};
    \node[font=\scriptsize, align=center, text=red] at (3, 2.0) {Certificate:\\$\{TE \le 15bps\}$};

    \draw[map] (pe) -- (pa) node[midway, left] {$p$};
    \draw[map] (pe) -- (pb) node[midway, right] {$q$};

    \node[align=center, font=\small, gray!80!black] at (3, -1) {Alignment holds because\\an audit trail $e$ exists};

\end{tikzpicture}
\caption{A Vertical Morphism in the span model. Unlike a simple relation which merely asserts $R(a, b)$, a span requires a concrete object $e$ (the evidence) that projects down to identifying the specific hub and spoke portfolios it links.}
\label{fig:span-bridge}
\end{figure}

\subsubsection*{Composition via Pullback}

Vertical composition models the transitivity of alignment. If we have evidence $E$ that $A$ aligns with $B$, and evidence $F$ that $B$ aligns with $C$:
\[
A \xleftarrow{p_1} E \xrightarrow{q_1} B \quad \text{and} \quad B \xleftarrow{p_2} F \xrightarrow{q_2} C
\]
The composite vertical morphism is the \textbf{pullback} $E \times_B F$:
\[
E \times_B F = \{ (e, f) \in E \times F \mid q_1(e) = p_2(f) \}
\]
\textbf{Interpretation:} A witness for the composite alignment $A \dashrightarrow C$ consists of a pair of witnesses $(e, f)$ that match at the intermediate portfolio in $B$. This formalizes the ``no phantom portfolios'' requirement: we must successfully identify an actual portfolio in $B$ that bridges the gap.

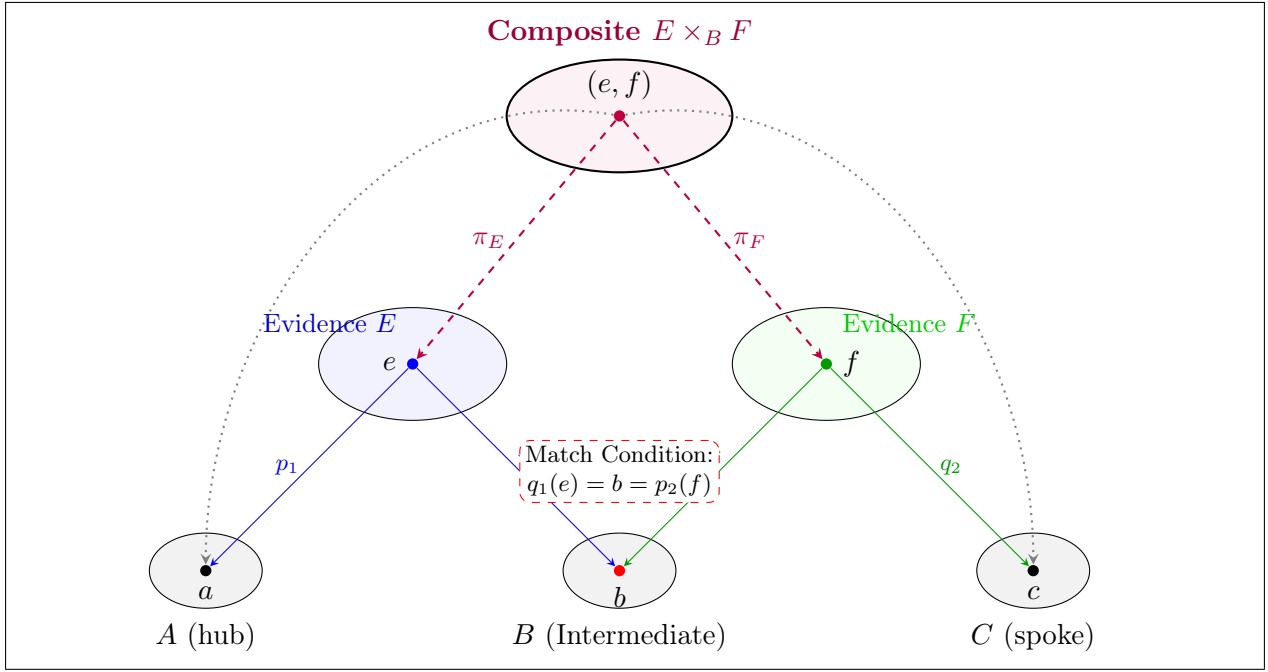
\begin{figure}[h]
\centering
\begin{tikzpicture}[scale=1.1, >=stealth]

    
    \node (A) at (0,0) [draw, ellipse, minimum width=1.5cm, minimum height=1cm, fill=gray!10] {};
    \node at (0, -0.8) {$A$ (hub)};
    
    \node (B) at (5,0) [draw, ellipse, minimum width=1.5cm, minimum height=1cm, fill=gray!10] {};
    \node at (5, -0.8) {$B$ (Intermediate)};
    
    \node (C) at (10,0) [draw, ellipse, minimum width=1.5cm, minimum height=1cm, fill=gray!10] {};
    \node at (10, -0.8) {$C$ (spoke)};

    \node (E) at (2.5, 2.5) [draw, ellipse, minimum width=2.5cm, minimum height=1.5cm, fill=blue!5] {};
    \node at (1.5, 3.0) [blue!80!black, font=\small] {Evidence $E$};

    \node (F) at (7.5, 2.5) [draw, ellipse, minimum width=2.5cm, minimum height=1.5cm, fill=green!5] {};
    \node at (8.5, 3.0) [green!80!black, font=\small] {Evidence $F$};

    \node (P) at (5, 5.5) [draw, ellipse, minimum width=3cm, minimum height=1.5cm, fill=purple!5, thick] {};
    \node at (5, 6.5) [purple!80!black, font=\bfseries] {Composite $E \times_B F$};


    \node (a_pt) at (0,0) [circle, fill=black, inner sep=1.5pt, label=below:$a$] {};
    
    \node (b_pt) at (5,0) [circle, fill=red, inner sep=1.5pt, label=below:$b$] {};
    
    \node (c_pt) at (10,0) [circle, fill=black, inner sep=1.5pt, label=below:$c$] {};

    \node (e_pt) at (2.5, 2.5) [circle, fill=blue, inner sep=1.5pt, label=left:$e$] {};
    
    \node (f_pt) at (7.5, 2.5) [circle, fill=green!60!black, inner sep=1.5pt, label=right:$f$] {};

    \node (ef_pt) at (5, 5.5) [circle, fill=purple, inner sep=1.5pt, label=above:{$(e,f)$}] {};


    \draw[->, dashed, purple, thick] (ef_pt) -- (e_pt) node[midway, left, font=\small] {$\pi_E$};
    \draw[->, dashed, purple, thick] (ef_pt) -- (f_pt) node[midway, right, font=\small] {$\pi_F$};

    \draw[->, blue] (e_pt) -- (a_pt) node[midway, left, font=\footnotesize] {$p_1$};
    \draw[->, blue] (e_pt) -- (b_pt) node[midway, above, font=\footnotesize, sloped] {$q_1$};

    \draw[->, green!60!black] (f_pt) -- (b_pt) node[midway, above, font=\footnotesize, sloped] {$p_2$};
    \draw[->, green!60!black] (f_pt) -- (c_pt) node[midway, right, font=\footnotesize] {$q_2$};

    
    \node[draw=red, dashed, inner sep=2pt, rounded corners, fill=white, font=\footnotesize, align=center] at (5, 1.2) 
        {Match Condition:\\$q_1(e) = b = p_2(f)$};

    \draw[->, dotted, thick, gray] (ef_pt) to[out=170, in=90] (a_pt);
    \draw[->, dotted, thick, gray] (ef_pt) to[out=10, in=90] (c_pt);

\end{tikzpicture}
\caption{Detail of element-wise mappings in the pullback. A composite witness $(e,f)$ can exist in the top space only because the individual witnesses $e$ and $f$ map to the same intermediate portfolio $b$. The projection maps $\pi_E$ and $\pi_F$ recover the original evidence from the composite.}
\label{fig:span-composition-elements}
\end{figure}

\subsubsection*{The Double Category $\mathbb{S}\mathbf{pan}(\mathcal{E})$}

We construct the double category $\mathbb{S}\mathbf{pan}(\mathcal{E})$ as follows:

\begin{enumerate}
   \item \textbf{Objects:} Permissible portfolio spaces $A, B, \dots$ (objects of $\mathcal{E}$).
   \item \textbf{Horizontal Morphisms:} Continuous maps $f: A \to B$ (the re-implementations). Note that the horizontal category $\mathbb{D}_h$ is simply $\mathcal{E}$ itself.
   \item \textbf{Vertical Morphisms:} Spans $S = (A \xleftarrow{p_S} E_S \xrightarrow{q_S} C)$.
   \item \textbf{2-Cells (Evidence Transformations):} A 2-cell $\alpha$
   is a morphism of spans in the standard sense.
\end{enumerate}

Note that pullbacks in $\mathcal{E}$ are required to define composition of spans.

\textbf{Motivation:} In the thin framework, a 2-cell asserted that ``if $x$ aligns with $z$ via $S$, then $f(x)$ aligns with $g(z)$ via $T$.'' In this framework, the 2-cell $\alpha$ acts as an \textbf{evidence transformer}. It takes a certificate of alignment for the hub portfolios and systematically transforms it into a certificate of alignment for the spoke portfolios.

This is essential for \textbf{audit trails}. It is not enough to know that the ETF basket \emph{is} compliant; we require a map $\alpha$ that derives the ETF's alignment report directly from the Mutual Fund's alignment report.

\subsubsection*{Horizontal Morphisms as a Full Subcategory}

To integrate standard re-implementations into the vertical structure, we use the \textbf{companion} construction. Any horizontal map $f: A \to B$ induces a vertical span (its graph):
\[
\text{Graph}(f) := (A \xleftarrow{\text{id}} A \xrightarrow{f} B)
\]
This allows us to view horizontal re-implementations as a special case of alignment (``perfect deterministic alignment''). The fact that horizontal morphisms embed fully into the span structure ensures that we can treat ``re-implementation'' and ``alignment'' uniformly when needed, while retaining the directional distinction for the user.

\subsubsection*{Technical Trade-off: Weak Associativity}

The move from relations to spans introduces a layer of algebraic complexity. While relations form an associative algebra (sets satisfy $(A \cap B) \cap C = A \cap (B \cap C)$), spans form a \textbf{bicategory}. Pullbacks are associative only \emph{up to isomorphism}:
\[
(E \times_B F) \times_C G \cong E \times_B (F \times_C G)
\]
For practical portfolio systems, this requires handling coherence isomorphisms (bookkeeping). The thin $\mathbb{HS}$ model in the main text suppresses this complexity to focus on topological properness, but a fully automated ``Evidence Engine'' would likely require the span-based bicategorical structure to manage the data structures of alignment reports.

\subsubsection*{Worked Example: Evidence Spans in a CIO Pipeline}

To make the construction concrete, consider the three-stage pipeline from Section~\ref{sec:propagation}: a CIO model portfolio $x \in K_{\mathrm{CIO}}$ is re-implemented into a sector model $y = f(x) \in K_{\mathrm{Sec}}$, which is then re-implemented into a client portfolio $z = g(y) \in K_{\mathrm{Client}}$.

In the relational (thin) setting, the alignment constraint between $K_{\mathrm{CIO}}$ and $K_{\mathrm{Sec}}$ is a closed relation $R \subseteq K_{\mathrm{CIO}} \times K_{\mathrm{Sec}}$ encoding, say, ``tracking error $\le 50$bps.''  The 2-cell is binary: either $(x, y) \in R$ or it is not.

In the span setting, we replace $R$ with a span
\[
K_{\mathrm{CIO}} \xleftarrow{\;\pi_1\;} E_R \xrightarrow{\;\pi_2\;} K_{\mathrm{Sec}}
\]
where $E_R$ is the \emph{evidence space}: a set of triples $(x, y, c)$ where $c$ is a \textbf{compliance certificate} --- a data structure recording the tracking error $\|g_A(x) - g_B(y)\|$, the fee differential, the date of computation, and the solver parameters used.  The projections $\pi_1, \pi_2$ recover the hub and spoke portfolios.

Suppose a second alignment span $S$ represents a fee-budget constraint between $K_{\mathrm{Sec}}$ and $K_{\mathrm{Client}}$, with evidence space $E_S$ recording fee computations and budget attestations.  The composite span $S \circ R$ is formed by pullback:
\[
E_{S \circ R} \;=\; E_R \times_{K_{\mathrm{Sec}}} E_S \;=\; \{(x, y, c_1, z, c_2) : \pi_2(x, y, c_1) = \pi_1'(y, z, c_2)\}
\]
An element of $E_{S \circ R}$ is a \emph{composite audit trail}: a tracking-error certificate $c_1$ linking CIO to sector, paired with a fee-budget certificate $c_2$ linking sector to client, joined at the common intermediate portfolio $y$.  This is precisely the data structure an Evidence Ledger (Appendix~\ref{app:architecture}) would store.  The thin framework records only whether $(x, z)$ satisfies the composite constraint; the span framework records \emph{which certificates were used at each stage and through which intermediate portfolio they compose}.

\subsection{Beck--Chevalley and Frobenius}

In the baseline $\mathbb{HS}$ framework, the Beck--Chevalley and Frobenius laws were expressed as equalities of sets (relations). In the $\mathbb{S}\mathbf{pan}(\mathcal{E})$ framework, these laws become \textbf{natural isomorphisms} between spans. This shift is significant for software implementation: it means we can structurally convert audit trails between different stages of the pipeline without information loss.

The following results are standard.

\begin{theorem}[Lax Beck--Chevalley Condition for General Squares]
Given a commuting square of horizontal morphisms (continuous maps) in $\mathcal{E}$:
\[
\begin{tikzcd}
    A \arrow[r, "g"] \arrow[d, "f'"'] & B \arrow[d, "f"] \\
    C \arrow[r, "h"'] & D
\end{tikzcd}
\]
such that $f \circ g = h \circ f'$.

For any vertical morphism (alignment span) $S: B \dashrightarrow Z$, there exists a canonical 2-cell (morphism of spans) $\alpha$:
\[
\alpha : f'_! (g^* S) \Longrightarrow h^* (f_! S)
\]
\end{theorem}

Now consider a commuting square of horizontal re-implementations that is a \emph{pullback} in $\mathcal{E}$ (the analog of the ``pointwise cartesian'' condition from Definition~\ref{def:pointwise-cartesian}). In the span framework, the Beck--Chevalley condition asserts an isomorphism of evidence spaces for any vertical alignment $S$ on $B$:
\[
f'_! \, g^* S \cong h^* \, f_! S
\]
\textbf{Operational Meaning:} Suppose $S$ represents a benchmark constraint on $B$.
\begin{itemize}
   \item The LHS ($f'_! \, g^* S$) represents: ``Pull the benchmark constraint back to $A$ (pre-alignment), find evidence of alignment there, and push that evidence forward to $C$.''
   \item The RHS ($h^* \, f_! S$) represents: ``Push the benchmark constraint forward to $D$ (post-alignment), and pull back to check if $C$'s re-implementation satisfies it.''
\end{itemize}
The isomorphism guarantees that an alignment certificate generated by the ``early verification'' process (LHS) can be \emph{losslessly} translated into an alignment certificate for the ``late verification'' process (RHS). The audit trail is preserved across the lattice of re-implementations.

\paragraph{Frobenius Reciprocity (Preservation of Filtered Evidence).}
For a horizontal morphism $f: A \to B$ and vertical alignments $R$ (on $A$) and $S$ (on $B$), the Frobenius law describes the interaction between filtering and re-implementation:
\[
f_!(R \times_A f^*S) \cong f_!R \times_B S
\]
\textbf{Operational Meaning:}
\begin{itemize}
   \item \textbf{LHS (Filter then Re-implement):} We find portfolios in $A$ that satisfy alignment $R$ \emph{and} whose re-implementations satisfy $S$. We bundle this combined evidence and push it to $B$.
   \item \textbf{RHS (Re-implement then Filter):} We take evidence that $A$ satisfies $R$, push it to $B$, and pair it with independent evidence that the result satisfies $S$.
\end{itemize}
In the relational model, this was a statement about the set of permissible portfolios remaining constant. In the span model, the isomorphism explicitly constructs the composite witness. It tells us that we do not need to re-run the alignment check $S$ from scratch after re-implementation; we can essentially ``carry over'' the witness from the hub analysis.

\clearpage  \section{Persona–Indexed Hub–and–Spoke}
\label{sec:personas}

This section develops a variant of the $\HSSimp$ calculus in which every ``object'' is a \emph{family of permissible portfolio spaces indexed by a base of personas}. Formally, we work in the slice over a fixed compact Hausdorff space \(B\) of personas (age bands, client segments, life–cycle states, etc.). Objects are \emph{display maps} \(p:K\to B\) whose fibers \(K_b\subseteq \Delta^{n(b)}\) are closed permissible portfolio spaces. Horizontal morphisms are continuous, fiberwise re–implementations over \(B\); vertical morphisms are closed, persona–indexed alignment relations. All coherence results (closed pushforwards, Adjunction, Beck–Chevalley, Frobenius) lift \emph{fiberwise} and hence globally.

The intuition is that \(B\) is a space of \emph{investor personas}. For example:
\begin{itemize}
    \item \(B\) is a range of ages, e.g. the closed interval $[25, 95]$. This can be used to model \emph{target date glide paths} (see \cite{Singh}), capturing the requirement that re-implementations and alignment specifications be compatible with the way that the portfolio in a glide path evolves as the investor ages.
    \item \(B\) is a product space incorporating both the age and educational attainment of the investor.
    \item \(B = [0,1]\) is an interpolation parameter between two existing hub portfolios. Given \(K_0, K_1 \subseteq \Delta^n\), define \(K_t = (1-t)K_0 + tK_1\) (Minkowski combination; closed in \(\Delta^n\) since both summands are compact). The family \(K = \bigsqcup_{t \in [0,1]} K_t\) is a persona-indexed object, and any \(t \in (0,1)\) defines a ``virtual hub'' interpolating between the two endpoint strategies---allowing a client to hold a portfolio ``between'' two existing products without designing a new hub from scratch.
\end{itemize}

Throughout we assume \(B\) compact Hausdorff unless stated otherwise. This matches key use–cases (closed age interval, finite demographic factors, aggregate wealth, etc.), and makes properness automatic again. In practice \(B\) will usually be a simple space.

While the main applications below involve connected \(B\), the finite (hence discrete) case commonly occurs---a series of model portfolios indexed by risk level, say. Mathematically trivial; organizationally, not.

\subsection{Persona–indexed Permissible Spaces}
\label{sub:display}

\begin{definition}[Persona–indexed object]

Fix a compact Hausdorff base \(B\). A \emph{persona–indexed permissible portfolio space} is a continuous map \(p:K\to B\) such that \(K\) is a closed subset of \(B\times \Delta^{n}\) for some \(n\), and for each \(b\in B\) the fiber

\[
K_b \;=\; \{x\in \Delta^{n} \mid (b,x)\in K\}
\]

is a closed subset of \(\Delta^{n}\). We often write \(K=\bigsqcup_{b\in B} K_b\) suggestively.

\end{definition}

\noindent\emph{Interpretation.} For each persona \(b\) (e.g., age \(b\)), \(K_b\) is the permissible region under baseline constraints \emph{for that persona}. Continuity in \(b\) lets constraints (caps, budgets, liquidity) vary smoothly along demographic axes.

\begin{lemma}[Compactness]

If \(B\) is compact Hausdorff and \(K\subseteq B\times \Delta^{n}\) is closed, then \(K\) is compact and the projection \(p:K\to B\) is proper.

\end{lemma}

 Let \(p_i:K_i\to B\) and \(q_j:L_j\to B\) be objects over \(B\).

\begin{definition}[Horizontal morphisms over \(B\)]

A \emph{horizontal morphism} \(f:K_1\to K_2\) \emph{over \(B\)} is a continuous map with \(p_2\circ f=p_1\). Equivalently, a family \(\{f_b:K_{1,b}\to K_{2,b}\}_{b\in B}\) of continuous re–implementations varying continuously in \(b\).

\end{definition}

\begin{definition}[Vertical morphisms over \(B\)]

A \emph{vertical morphism} (alignment relation indexed by \(B\)) from \(K_1\) to \(L_1\) is a closed subset

\[
R \;\subseteq\; K_1 \times_{B} L_1 \;=\; \{(x,z)\in K_1\times L_1 \mid p_1(x)=q_1(z)\},
\]

\noindent i.e., a closed, fiberwise relation \(R_b\subseteq K_{1,b}\times L_{1,b}\).

\end{definition}

\noindent

Horizontal composition and vertical (relational) composition are computed over \(B\). We write

\[
\Graph_B(f)\;=\;\{(x,f(x))\in K_1\times_B K_2\}
\]

\noindent for the graph of a horizontal map over \(B\).

\begin{proposition}[Automatic properness over a compact base]

Let \(B\) be compact Hausdorff. For any horizontal morphism \(f:K_1\to K_2\) over \(B\) between objects \(p_i:K_i\to B\) as in \S\ref{sub:display}, the map \(f\) is proper. Moreover, for any object \(r:M\to B\), the product map

\[
f\times_B \mathrm{id}_M:\; K_1\times_B M \;\longrightarrow\; K_2\times_B M
\]

is proper and hence closed.

\end{proposition}

\begin{definition}[Pullback and pushforward over \(B\)]

Given \(f:K_1\to K_2\) over \(B\) and a closed relation \(S\subseteq K_2\times_B M\), define

\[
f^{*}S\;=\;\{(x,w)\in K_1\times_B M \mid (f(x),w)\in S\}.
\]

Given a closed relation \(R\subseteq K_1\times_B M\), define its pushforward

\[
f_{!}R\;=\;(f\times_B \mathrm{id}_M)(R)\;\subseteq\; K_2\times_B M.
\]

\end{definition}

\begin{theorem}[Closedness of pushforward over \(B\)]

If \(B\) is compact Hausdorff, \(f:K_1\to K_2\) is horizontal, and \(R\subseteq K_1\times_B M\) is closed, then \(f_{!}R\) is closed in \(K_2\times_B M\).

\end{theorem}

\noindent\emph{Proof sketch.} By \S\ref{sub:display}, \(K_i\) and \(M\) are compact; thus the fiber product \(K_1\times_B M\) is compact and \(R\) is compact. The map \(f\times_B \mathrm{id}_M\) is continuous, hence it maps the compact \(R\) to a compact (therefore closed) subset of \(K_2\times_B M\).

\subsection{The Indexed Double Category \(\HSSimp/B\) and its Laws}

\begin{definition}[Double category \(\HSSimp/B\)]

Objects: display maps \(p:K\to B\) as in \S\ref{sub:display}. Horizontal morphisms: continuous maps over \(B\). Vertical morphisms: closed relations inside fiber products over \(B\). A \(2\)-cell for a square

\[
\begin{tikzcd}[row sep=small, column sep=large]
K_1 \ar[r,"f"] \ar[d, dashed, "R"'] & K_2 \ar[d, dashed, "S"]\\
L_1 \ar[r,"g"'] & L_2
\end{tikzcd}
\qquad\text{(all over \(B\))}
\]

\noindent exists iff the inclusion $\Graph_B(g)\circ R \;\subseteq\; S\circ \Graph_B(f)
$ holds inside the appropriate fiber products. 

\end{definition}

Note that like $\HSSimp$, \(\HSSimp/B\) is thin.

\begin{proposition}[Adjunction \(f_{!}\dashv f^{*}\) in \(\HSSimp/B\)]

For any horizontal \(f:K_1\to K_2\) and closed relations \(R\subseteq K_1\times_B M\), \(S\subseteq K_2\times_B M\),

\[
R\;\subseteq\; f^{*}S \quad\Longleftrightarrow\quad f_{!}R\;\subseteq\; S.
\]

\end{proposition}

\begin{definition}[Pointwise–cartesian square over \(B\)]

A commuting square of horizontals over \(B\),

\[
\begin{tikzcd}[row sep=small, column sep=large]
K \ar[r,"g"] \ar[d,"f'"'] & L \ar[d,"f"]\\
K' \ar[r,"h"'] & L',
\end{tikzcd}
\qquad f\circ g=h\circ f',
\]

\noindent is \emph{pointwise–cartesian} if for every \(b\in B\) and every pair \((y,z)\in L_b\times K'_b\) with \(f_b(y)=h_b(z)\) there exists \(x\in K_b\) with \(g_b(x)=y\) and \(f'_b(x)=z\).

\end{definition}

\begin{theorem}[Beck–Chevalley in \(\HSSimp/B\)]

For a pointwise–cartesian square of horizontals over \(B\) and any closed relation \(R\subseteq L\times_B M\),

\[
f'_{!}\, g^{*} R \;=\; h^{*}\, f_{!} R.
\]

\end{theorem}

\begin{theorem}[Frobenius reciprocity in \(\HSSimp/B\)]

For a horizontal \(f:K_1\to K_2\) and closed relations \(R\subseteq K_1\times_B M\), \(S\subseteq K_2\times_B M\),

\[
f_{!}\bigl(R\cap f^{*}S\bigr)\;=\; f_{!}R \;\cap\; S.
\]

\end{theorem}

\noindent\emph{Comments.} Proofs are word–for–word fiberwise adaptations of HS: compactness (from compact \(B\) and closedness of objects) makes properness automatic; pushforwards preserve closedness; the elementwise set–chases are identical in each fiber and hence globally.

\subsection{Examples of Bases \(B\) and Modeling Patterns}

\paragraph{(1) target–date glidepath.} \(B=[b_{\min},b_{\max}]\subset \mathbb{R}\) (age). Let \(K\subseteq B\times \Delta^n\) encode age–dependent caps (e.g., equity cap decreasing continuously with age). A horizontal \(f:K\to K'\) over \(B\) may implement an age–specific “active\(\to\)ETF” mapping, with optimization parameters (tracking tolerance, fee budgets) continuous in \(b\). Vertical relations (tracking, factor budgets) live fiberwise in \(K_b\times K'_b\).

  \paragraph{(2) Age \(\times\) educational attainment.}
  \(B=[b_{\min},b_{\max}]\times E\) with \(E\) a finite discrete set
  (e.g., \(\{\)HS, College, Graduate\(\}\)). Finite \(E\) is compact;
  thus \(B\) is compact. Piecewise–constant features in \(e\in E\) and
  continuous in age satisfy the hypotheses; all laws apply componentwise
  over~\(E\).

  \paragraph{(3) Age \(\times\) wealth.}
  \(B=[b_{\min},b_{\max}]\times [w_{\min},w_{\max}]\subset\mathbb{R}^{2}\)
  (age \(\times\) wealth level).
  In standard lifecycle portfolio theory
  (Merton \cite{Merton}, Samuelson \cite{Samuelson}),
  the optimal asset allocation is a function of \emph{both} the
  investor's time horizon (age) and current wealth---through mechanisms
  such as decreasing relative risk aversion, subsistence floors, or
  housing considerations.
  Let \(K\subseteq B\times\Delta^{n}\) encode joint age--wealth
  constraints: for instance, the equity cap may decrease with age
  (de-risking) while the allocation to alternatives increases with
  wealth (access to illiquid vehicles requires higher minimums).
  A horizontal morphism \(f:K\to K'\) over \(B\) implements a
  re-implementation that is continuous in both \(b\) (age) and \(w\)
  (wealth); vertical alignment relations are fiberwise in
  \(K_{(b,w)}\times K'_{(b,w)}\).
  The product of two compact intervals is compact Hausdorff, so all
  results of Section~\ref{sec:personas} apply.

  \paragraph{(4) Non–contractible bases.}
  Non–contractibility does not impede the calculus; compactness does
  the work. However, the practical use cases are more speculative.

  \begin{itemize}\itemsep3pt

  \item \(B=S^{1}\) \textbf{(spending seasonality):}
    Annual spending patterns are approximately periodic---tuition
    payments in autumn, tax payments in spring, holiday expenditures in
    winter---so the ``time of year'' parameter lives naturally on the
    circle \(S^{1}\cong\mathbb{R}/\mathbb{Z}\).
    A family \(K\subseteq S^{1}\times\Delta^{n}\) encodes seasonal
    cash and duration targets (e.g., higher cash allocation in the
    months preceding large outflows), with constraints that wrap
    continuously around the year-end boundary rather than exhibiting a
    discontinuity at December/January.
    All coherence laws hold fiberwise over~\(S^{1}\).

  \item \(B=S^{1}\) \textbf{(investment-cycle regimes):}
    Tactical asset allocation often depends on where the economy sits
    in a recurring regime cycle.
    For example, in the growth--inflation clock
    (cf.\ Greer \cite{Greer}),
    the economy rotates through four quadrants---recovery, expansion,
    stagflation, contraction---and the recommended factor tilts
    (value vs.\ growth, duration, commodities) vary continuously around
    the cycle.
    Modeling the regime parameter as \(b\in S^{1}\) captures the
    periodicity: the ``post-contraction'' regime borders
    ``pre-recovery'' without a boundary discontinuity.
    Families \(K\subseteq S^{1}\times\Delta^{n}\) then encode
    regime-dependent permissible allocations, and re-implementations
    (e.g., from a strategic to a tactical sleeve) are continuous in the
    regime parameter.

  \item \(B\) \textbf{a finite connected graph}
    (a 1--dimensional CW complex)\textbf{:}
    Vertices represent discrete \emph{states}---career stages
    (student, employed, unemployed, retired), product-lifecycle phases
    (launch, growth, maturity, decline), or plan statuses
    (accumulation, payout, rollover)---and edges represent
    \emph{permissible transitions} between them.
    Loops are allowed: an employee can cycle between ``employed'' and
    ``unemployed'' multiple times, and a product can re-enter a growth
    phase after repositioning.
    Under the standard CW topology, any finite graph is compact
    Hausdorff.
    A display map \(p:K\to B\) assigns to each state \(v\) a fiber
    \(K_{v}\subseteq\Delta^{n}\) of permissible allocations, and to
    each edge a continuous family of intermediate allocations
    interpolating between the endpoint fibers (modeling a transition
    glide path).
    Coherence results (adjunction, Beck--Chevalley, Frobenius) persist
    because they depend only on compactness of~\(B\), not on its
    homotopy type.

  \end{itemize}
  
\begin{definition}[Reindexing]

Given \(u:B'\to B\) and \(p:K\to B\), define the pullback

\[
u^{*}K \;=\; \{(b',x)\in B'\times \Delta^{n}\mid (u(b'),x)\in K\} \;\xrightarrow{\;p'\;}\; B',
\]

\noindent an object over \(B'\). Horizontals and verticals pull back along \(u\) via fiber products. This makes \(B\mapsto \HSSimp/B\) a pseudofunctor \( (\mathbf{Top}_{\mathrm{cpt,Haus}})^{\mathrm{op}}\to \mathbf{DblCat}\) (cf. fibrations/display maps; \emph{compare} \cite{Shulman, Street}).

\end{definition}

\noindent\emph{Portfolio reading.} Aggregating persona attributes or refining them induces functorial reindexing of constraints, maps, and alignments---useful when migrating between segmentation schemes.

\subsection{Cross-Vintage Constraints}

A key reason to consider persona-indexed portfolio spaces is that it is often insufficient to consider constraints on individual portfolios. It is often the case that a ``series'' of portfolios is required to satisfy \emph{cross-portfolio} constraints. The $\HSSimp/B$ framework allows us to specify this formally.

For example, for suites of target date funds indexed by retirement year $r$, it is often desirable to enforce cross-vintage constraints. Consider a constraint on the steepness of the glide path: as calendar time advances by $\Delta$, the allocation for the $r$-vintage at age $b$ differs from the allocation for the $(r+\Delta)$-vintage at age $b+\Delta$ by no more than $\epsilon$.

Standard slice categories restrict vertical morphisms to the fiber product $K \times_B L$, forcing aligned portfolios to map to the exact same base point $b$. This prevents us from modeling constraints that link different personas, such as these cross-vintage constraints. To address this, we generalize the vertical morphisms in $\HSSimp/B$.

\begin{definition}[Generalized Vertical Morphism over $B$]
Let $p: K \to B$ and $q: L \to B$ be objects (display maps). A vertical morphism is a pair $(r, R)$ where:
\begin{enumerate}
    \item $r \subseteq B \times B$ is a closed relation on the base space (e.g., the ``time shift’’ relation $r = \{(b, b+\Delta) \mid b \in [0, T-\Delta]\}$).
    \item $R \subseteq K \times L$ is a closed relation on the total spaces.
    \item \textbf{Compatibility:} The relation $R$ lies over $r$. That is, for every $(x, y) \in R$, the pair of base points $(p(x), q(y))$ must belong to $r$.
\end{enumerate}
\end{definition}

  \begin{example}[Formulation of Glide Path Steepness Constraint]
  Let $x^{(t)}_b$ denote the allocation for vintage $t$ at age $b$.
  A smoothness constraint requires
  $|x^{(t)}_b - x^{(t)}_{b+\Delta}| \le \epsilon$.
  This is a relation $R$ covering the shift relation
  $r = \{(b, b+\Delta)\}$ on the age axis $B$.
  Because $R$ and $r$ are closed, vertical composition of generalized
  morphisms preserves closedness (by the same argument as
  Lemma~\ref{lem:vertical-closed}, applied fiberwise).

  The standard coherence laws (Beck--Chevalley, Frobenius) for the slice
  $\HSSimp/B$ hold when vertical morphisms are restricted to the fiber
  product $K \times_B L$.  For generalized
  morphisms which allow
  cross-fiber relations, the coherence laws require additional hypotheses
  on the interaction between the base relation $r$ and the display maps.
  In the glide path setting, these conditions are satisfied when the
  shift map $b \mapsto b + \Delta$ is a proper self-map of~$B$.
  \end{example}

\subsection{Summary}

By replacing single spaces \(K\subseteq \Delta^n\) with persona–indexed display maps \(p:K\to B\) and working in the slice \(\HSSimp/B\), we obtain a calculus in which \emph{every $\HSSimp$ notion is fibered over personas}. With \(B\) compact Hausdorff (a natural condition for many segmentations), properness is again automatic, pushforwards preserve closedness, and Adjunction, Beck–Chevalley, and Frobenius hold verbatim, persona by persona.

This supports glidepaths, demographic conditioning, and even non–contractible persona bases (circles/graphs) without changing the proofs---only the ambient bookkeeping. It also supports practical considerations in glide path design, in particular cross-portfolio (i.e. cross-vintage) constraints.

\clearpage

\section{DOTS for Portfolio Systems}
\label{sec:DOTS}

This variant framework is based on a suggestion from David Jaz Myers. The idea is simple: instead of committing to a single spoke per hub---which is what the $\HSSimp$ framework does when it asks for a deterministic re-implementation $f$---define a \emph{menu} of permissible spokes and let the portfolio manager choose from it. The choice may involve qualitative judgment, client preferences, or constraints that resist formalization; what matters is that the menu itself is well-defined, that it composes sensibly, and that compliance can be verified at the menu level rather than at the level of individual selections. This is the DOTS approach (see \cite{LibkindMyers}, \cite{MyersDOTS}). This section provides only a sketch of the theory.

The primitive operation replaces the deterministic re-implementation with a set-valued one. Given an alignment relation $R$ specifying a constraint on a desired re-implementation (as in Section~\ref{subsec:reimplalign}), we act on a ``state''---a closed subset representing the current portfolio universe---and obtain a new closed subset representing the menu of permissible outcomes:
\[
K\odot R:=\{\,y: \exists x\in K,\ (x,y)\in R\,\}.
\]
As we will see this unifies screens, budgets, and optimizers as algebraic operations; it supports parallel assembly via a symmetric monoidal product; and it elevates strategy blueprints to operadic operations. When unique optimizers exist, DOTS collapses to $\HSSimp$ and all $\HSSimp$ guarantees apply. Otherwise, DOTS retains and exposes path dependence so that it can be managed, audited, and---where possible---reduced.

\subsection{The DOTS Equipment}
\subsubsection*{The thin double category}
Define $\DOTS$ with the following data:
\begin{itemize}[leftmargin=*]
  \item \textbf{Objects:} simplices $\Delta^n$, $n\ge 0$.
  \item \textbf{Horizontal arrows:} continuous maps $f\colon \Delta^n\to\Delta^m$.
  \item \textbf{Vertical arrows:} closed relations $R\subseteq\Delta^n\times\Delta^m$.
  \item \textbf{2-cells:} inclusions $\Graph(g)\circ R\subseteq S\circ \Graph(f)$.
\end{itemize}
Composition is by ordinary function composition horizontally and relational composition vertically. Thinness (only inclusions) keeps algebraic reasoning close to elementary set manipulations \cite{Shulman,GrandisPare} while still supporting a full calculus of pushforwards, pullbacks, and Frobenius reciprocity.

We could also define our objects to be closed permissible subsets $K\subseteq \Delta^n$ as before, but for simplicity, we use the above definition for now.

\subsubsection*{States as loose morphisms and the action}
A \emph{state} on $\Delta^n$ is a closed subset $K\subseteq\Delta^n$, viewed as a loose morphism $\Delta^0\Rightarrow\Delta^n$. Given a closed state $K\subseteq\Delta^n$ and a closed relation $R\subseteq\Delta^n\times\Delta^m$, the action is the pushforward
\[
K\odot R := \{\, y\in\Delta^m : \exists x\in K,\ (x,y)\in R \,\}.
\]

\begin{proposition}[Action laws]\label{prop:action}
For closed $K$ and closed relations $R,S$ the following hold:
\begin{enumerate}[label=(\alph*),leftmargin=*]
  \item \textbf{Closedness:} $K\odot R$ is closed.
  \item \textbf{Unitality:} $K\odot \Delta_{\Delta^n}=K$, where $\Delta_{\Delta^n}$ is the identity relation on $\Delta^n$.
  \item \textbf{Associativity:} $(K\odot R)\odot S=K\odot (S\circ R)$, where $S \circ R$ denotes the standard relational composition (right-to-left: first apply $R$, then $S$).
  \item \textbf{Isotonicity:} $K\subseteq K'$ implies $K\odot R\subseteq K'\odot R$.
  \item \textbf{Projectors:} If $R=\{\,(x,x):x\in E\,\}$ with $E$ closed, then $K\odot R=K\cap E$ and $(K\odot R)\odot R=K\odot R$.
\end{enumerate}
\end{proposition}

\begin{proof}
Closedness:
Let $S:=(K\times\Delta^m)\cap R$. Both $K\times\Delta^m$ and $R$ are closed in $\Delta^n\times\Delta^m$, which is compact; hence $S$ is compact. The projection $\pi_y\colon S\to\Delta^m$ is continuous, so $\pi_y(S)=K\odot R$ is compact and closed.

Associativity:
We compute directly:
\[
(K\odot R)\odot S = \{ z : \exists y\in\Delta^m, \exists x\in K, (x,y)\in R, (y,z)\in S\} 
\]
\[
= \{ z : \exists x\in K, (x,z)\in S\circ R\} = K\odot (S\circ R).
\]

Projectors:
If $R=\{(x,x): x\in E\}$ with $E$ closed, then
\[
K\odot R=\{ x\in \Delta^n : \exists x'\in K, (x',x)\in R\}=K\cap E.
\]
Idempotence follows since $(K\cap E)\cap E=K\cap E$.
\end{proof}

\begin{remark}[Finance mapping]
(a) \emph{Closedness} means menus are stable under limits in the portfolio space: there are no phantom spokes created by limit points. (b) \emph{Unitality} expresses that ``doing nothing changes nothing.'' (c) \emph{Associativity} says we may act by $R$ and then by $S$, or equivalently by the composite constraint $S\circ R$. (d) \emph{Isotonicity} states that broadening the hub can only expand the menu. (e) \emph{Projectors} model idempotent screens such as ESG, liquidity, or admissibility filters.
\end{remark}

\begin{remark}[Generalizing horizontal morphisms]
    Although the action $K \odot R$ is defined using vertical composition, it serves the functional role of a generalized horizontal transition. In the main HS framework (Section~\ref{subsec:reimplalign}), we sought a single deterministic re-implementation $f: K \to L$ respecting $R$. Here, the object $K \odot R$ represents the \textit{space of all possible re-implementations}---effectively a ``cloud'' or ``menu'' of spoke portfolios permissible under the constraint $R$. 
    
    Thus, the DOTS action algebraizes the ``feasible region'' concept from Section~\ref{subsec:opt-align}, allowing us to manipulate sets of potential spokes before a final selection is made.
\end{remark}

\begin{remark}[Deterministic slice and bridge back to HS]
If $R=\Graph(f)$ for a continuous $f\colon\Delta^n\to\Delta^m$, then $K\odot R=f(K)$, which is closed because $K$ is compact. Thus, whenever an action reduces to a graph (e.g., a unique optimizer), we are back in the original $\HSSimp$ theory and inherit all its guarantees.
\end{remark}

\begin{remark}[Non-determinism and selection]
In the general case where $K\odot R$ is not a graph, implementing menu selection means making a \emph{continuous} selection of a spoke portfolio in the menu for each hub portfolio in $K$, e.g. via some optimization process. Theorem~\ref{thm:michael-selection} guarantees that this is possible in practical settings.
\end{remark}

\begin{remark}[DOTS vs. HSP-r]
It is instructive to compare the DOTS approach with the stochastic HSP-r framework (Section~\ref{sec:Polish}).
\begin{itemize}
\item \textbf{HSP-r (Distributions):} Models \textit{uncontrolled} variance (e.g., solver noise, liquidity slippage). The system produces a distribution, and we require that \textit{almost all} outcomes are compliant.
\item \textbf{DOTS (Menus):} Models \textit{controlled} choice (e.g., human manager discretion, regime selection). The system produces a set of valid options, and we require that the \textit{selection capability} is preserved.
\end{itemize}
Ideally, a pipeline uses DOTS for the design phase---narrowing the menu of permissible strategies---and HSP-r for the execution phase, modeling the noisy implementation of the selected strategy. Design is about choice; execution is about noise. The two require different mathematics.

The gap between the two frameworks is, however, not as clean as this distinction suggests. In practice, the ``design'' phase involves solver runs whose outputs are themselves noisy, and the ``execution'' phase involves human overrides that are better modeled as menu selections than as stochastic perturbations. The frameworks are complementary, but neither alone captures the full pipeline---a point the operadic formalism acknowledges implicitly by treating composition abstractly, without specifying whether its components are deterministic, set-valued, or stochastic.
\end{remark}

\begin{remark}[Multi-hub alignment]
The operadic composition above handles \emph{building a spoke from multiple inputs}: the core--satellite wiring diagram, for instance, takes two hub spaces and produces one output. DOTS handles this naturally. A distinct---and harder---situation arises when a single spoke must remain aligned with multiple \emph{independent} hubs simultaneously: a portfolio that must track both a regional equity benchmark and a global fixed-income benchmark managed by separate teams. This is a multi-objective alignment problem. The spoke's constraint set is the intersection of pullbacks from each hub, and feasibility of the intersection is not guaranteed. The current framework does not claim to solve this problem in general, though it does provide the tools to formulate it precisely: the intersection $\bigcap_j f_j^{*}S_j$ is well-defined and closed (each pullback is closed, and a finite intersection of closed sets is closed), so the question reduces to nonemptiness---a feasibility problem the framework poses but does not resolve.
\end{remark}

\begin{figure}
\centering
\begin{tikzpicture}[scale=1.1, >=stealth]
    \draw[thick] (0,0) -- (3.5,0) -- (1.75, 3) -- cycle;
    \node at (1.75, 3.3) {Ambient $\Delta^n$};
    
    \draw[fill=blue!15, draw=blue!60!black, thick] 
        (0.8, 0.5) -- (2.7, 0.5) -- (2.1, 1.8) -- (1.4, 1.8) -- cycle;
    \node[text=blue!60!black, font=\footnotesize] at (1.75, 2.0) {Permissible $K_{Hub}$};

    \draw[fill=blue!15, draw=blue!60!black] (1.2, 0.8) rectangle (1.9, 1.3);
    \node[text=blue!60!black, font=\footnotesize] at (1.55, 1.05) {$K$};
    
    \node[circle, fill=black, inner sep=1.5pt, label=below:$x$] (x) at (1.7, 0.9) {};

    \draw[thick] (6,0) -- (9.5,0) -- (7.75, 3) -- cycle;
    \node at (7.75, 3.3) {Ambient $\Delta^m$};

    \draw[fill=blue!15, draw=blue!60!black, thick] 
        (6.8, 0.5) -- (8.7, 0.5) -- (8.3, 1.8) -- (7.2, 1.8) -- cycle;
    \node[text=blue!60!black, font=\footnotesize] at (7.75, 2.0) {Permissible $K_{Spoke}$};

    \draw[fill=green!15, draw=green!60!black, thick] 
        (7.1, 0.7) -- (8.4, 0.7) -- (8.1, 1.5) -- (7.4, 1.5) -- cycle;
    
    \node[text=green!60!black, font=\bfseries, align=left] at (8.0, 1.2) {Menu\\$K \odot R$};

    \draw[dashed, draw=black, fill=white, opacity=0.6] (7.4, 0.9) circle (0.2);
    
    \node[font=\tiny, pin={[pin edge={black, thin}, align=left]below:$F_R(x)$}] at (7.4, 0.65) {};

    \draw[->, dashed, thick] (x) to[bend left=15] (7.2, 0.9);
    \node[align=center, font=\footnotesize, fill=white] at (4.75, 1.5) {Action $\odot R$\\ (One-to-Many)};

\end{tikzpicture}
\caption{Action of an alignment relation on a portfolio, generating a menu of spokes. Pre-existing permissible subsets are also shown in blue.}
\end{figure}
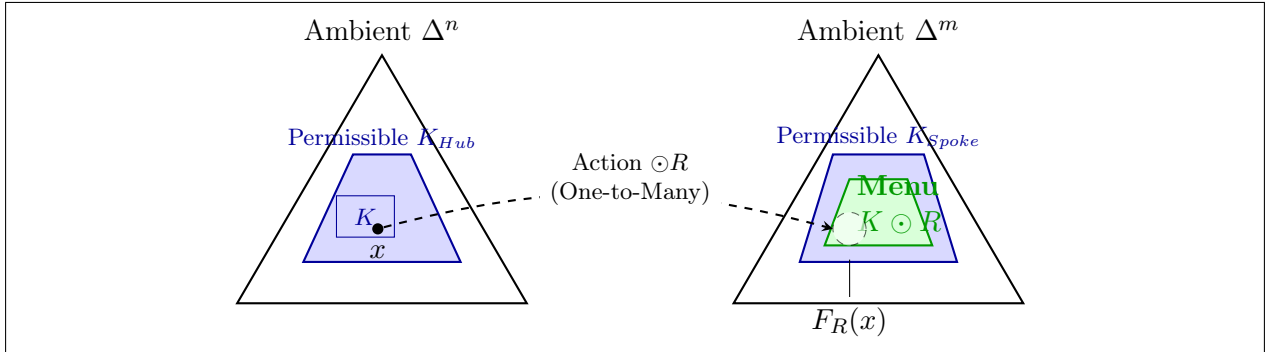

\subsection{Encoding Frequent Portfolio Relations}
\label{subsec:frequent-relations}
We list canonical closed relations used repeatedly in practice and indicate how to encode them. Throughout, the emphasis is on making each modeling choice explicit so that its mathematical properties (closedness, convexity) are clear to downstream consumers.

\subsubsection*{Tracking and exposure relations}
Let $g_A\colon\Delta^n\to\mathbb{R}^k$ and $g_B\colon\Delta^m\to\mathbb{R}^k$ be continuous exposure maps (asset classes, factors, or risk contributions). Define
\begin{align*}
R_{\mathrm{track},\epsilon}&:=\big\{(x,y): \lVert g_A(x)-g_B(y)\rVert_2\le \epsilon\big\},\\
R_{\mathrm{fac}\le b}&:=\big\{(y,z): (Fy)_j-z_j\le b_j\ \text{ for all } j\big\}.
\end{align*}
Closedness follows from continuity of $g_A,g_B$ and the closedness of Euclidean balls.

\subsubsection*{Fee, liquidity, and capacity}
With continuous fee, ADV, and position-size maps we set
\begin{align*}
R_{\mathrm{fee}\le \tau}&:=\{(y,y): \Fee(y)\le \tau\},\\
R_{\mathrm{liq}}&:=\{(y,y): \mathrm{ADV}(y)\ge \mathrm{ADV}_{\min},\ y_i\le c_i\ \forall i\}.
\end{align*}
These are projectors (subidentities), hence idempotent and order-insensitive with respect to themselves.

\subsubsection*{Tax and turnover}
Assuming continuous models for realized gain and turnover over a step, define
\begin{align*}
R_{\mathrm{tax}\le \theta}&:=\{(x,y): \mathrm{RealizedGain}(x,y)\le \theta\},\\
R_{\mathrm{turn}\le \kappa}&:=\{(x,y): \lVert y-x\rVert_1\le \kappa\}.
\end{align*}
The $\ell_1$ ball in the simplex is closed; realized-gain constraints are closed if the tax-lot model is finite with continuous aggregation.

\subsubsection*{Optimization-generated correspondences}
Let $F\colon\Delta^n\times\Delta^m\to\mathbb{R}$ be continuous with compact feasible sets in $y$ for each fixed $x$. For a closed feasible set $E$ (e.g., a liquidity face) define the correspondence
\[
\Gamma_F:=\big\{(x,y): y\in\operatorname*{arg\,min}_{y'\in \Delta^m\cap E} F(x,y')\big\}.
\]
Then $\Gamma_F$ is closed. If $F(x,\cdot)$ is strictly convex on $E$ for all $x$, the correspondence is a graph $\Graph(f)$ and the action determinizes.

\subsection{Operadic Templates}
In the DOTS paradigm, operadic templates encode multi-input wiring patterns. Each template consists of an input profile, a family of internal relations, and an output interface; evaluating the template on input states returns a menu. This formalism describes a simple framework for modularity and reusability in portfolio construction.

\subsubsection*{Core--satellite}
To illustrate the practical utility of operadic composition, we define a reusable wiring diagram $\Phi$ for a \textbf{Core-Satellite} strategy. In the operadic view (see \cite{Leinster}), $\Phi$ is an operation with arity 2: it accepts two permissible hub spaces (the inputs) and produces a single permissible total portfolio space (the output).

Let the inputs be:
\begin{itemize}
    \item $K_{\mathrm{core}} \subseteq \Delta^n$: A high-capacity, low-fee permissible space (e.g., an Index Fund hub).
    \item $K_{\mathrm{sat}} \subseteq \Delta^m$: A high-alpha, concentrated permissible space (e.g., a Hedge Fund hub).
\end{itemize}

The wiring diagram $\Phi$ encapsulates the logic of combining these independent sleeves into a unified strategy. It is defined by a weight parameter $w \in [0,1]$ (the allocation to Core) and a global constraint relation $R_{\mathrm{global}}$. The composition is given by:
\[
\Phi(K_{\mathrm{core}}, K_{\mathrm{sat}}) \;=\; 
\left( K_{\mathrm{core}} \otimes K_{\mathrm{sat}} \right) \odot R_{\mathrm{mix}, w} \odot R_{\mathrm{global}}
\]
Here, $R_{\mathrm{mix}, w}$ is the structural relation that creates the composite portfolio $y = w \cdot x_c + (1-w) \cdot x_s$, and $R_{\mathrm{global}}$ represents portfolio-level constraints that cannot be checked at the sleeve level, such as a total fee cap or a portfolio-level volatility cap.

\begin{figure}[h]
\centering
\begin{tikzpicture}[
    scale=1.1, 
    >=stealth, 
    hub/.style={
        draw=blue!60!black, 
        top color=blue!5, 
        bottom color=blue!15, 
        rounded corners=2pt, 
        minimum width=2.5cm, 
        minimum height=1.2cm,
        align=center,
        drop shadow
    },
    process/.style={
        circle,
        draw=gray!80,
        fill=white,
        thick,
        minimum size=0.8cm
    },
    constraint/.style={
        draw=red!60!black, 
        fill=red!5, 
        dashed, 
        rounded corners,
        inner sep=5pt
    },
    wire/.style={->, thick, draw=gray!60}
]

    \node[hub] (Core) at (0, 2) {\textbf{Input 1}\\$K_{\mathrm{core}}$\\(Index)};
    \node[hub] (Sat) at (0, 0) {\textbf{Input 2}\\$K_{\mathrm{sat}}$\\(Alpha)};

    
    \node[process] (Mix) at (4, 1) {\Large $\Sigma$};
    \node[below=0.1cm of Mix, text=gray] {Mix $w$};
    
    \draw[wire] (Core.east) -- (Mix) node[midway, above] {Weight $w$};
    \draw[wire] (Sat.east) -- (Mix) node[midway, below] {Weight $1-w$};

    \node[hub, right=2.0cm of Mix, fill=green!10, draw=green!60!black] (Result) {\textbf{Output}\\$K_{\mathrm{total}}$};
    
    \draw[wire] (Mix) -- (Result);

    \draw[red, very thick] (5.5, 0.7) -- (5.5, 1.3);
    
    \node[text=red!60!black, font=\scriptsize, align=center] at (5.5, 1.6) {Fee Cap\\Filter};
    
    \begin{scope}[on background layer]
        \node[constraint, fit=(Mix) (Result) (Core) (Sat), inner sep=15pt, label={[red!60!black]above:\textbf{Wiring Diagram} $\Phi$}] (Template) {};
    \end{scope}

\end{tikzpicture}
\caption{Visualizing Operadic Composition. The simple wiring diagram $\Phi$ (red dashed box) acts as a template. It takes two permissible spaces as inputs ($K_{\mathrm{core}}, K_{\mathrm{sat}}$), mixes them according to a structural relation, applies a global constraint (Fee Cap), and outputs the resulting menu of total portfolios.}
\label{fig:operad-wiring}
\end{figure}
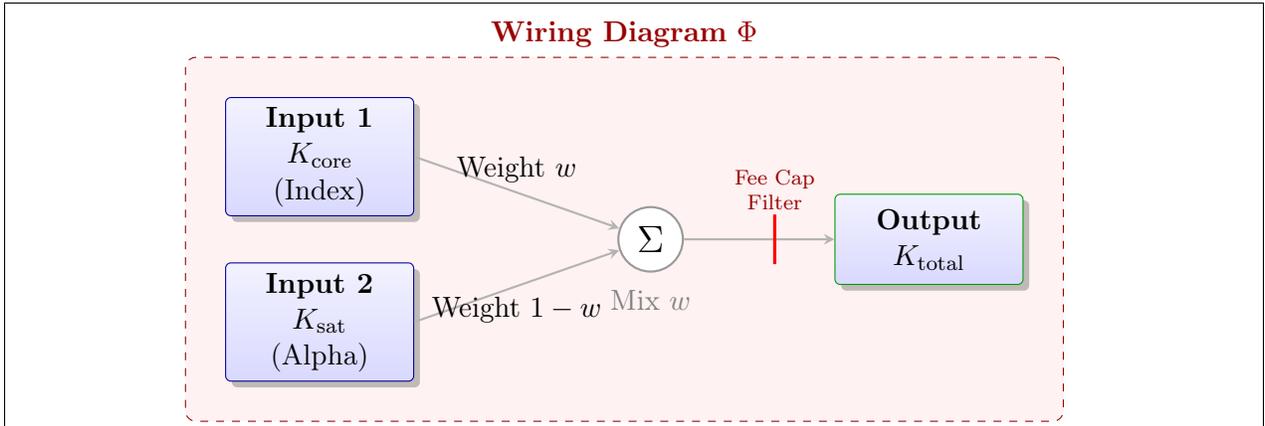

We treat $\Phi$ as a ``black box'' circuit that wires the two inputs together (Figure~\ref{fig:operad-wiring}). The power of the framework is that $K_{\mathrm{core}}$ and $K_{\mathrm{sat}}$ can be swapped out for any other valid objects (e.g., replacing the Index Fund with an ETF basket) without redesigning the compliance logic encoded in $\Phi$.

\subsubsection*{Blend-with-caps}
\textbf{Inputs.} $K_A\subseteq\Delta^{n_A}$, $K_B\subseteq\Delta^{n_B}$. A blend relation $R_{\mathrm{blend}}$ encodes budgets and a mapping to the spoke universe; fee and capacity projectors enforce operational limits. The pattern exposes where order matters when optimizers are non-unique.

\subsubsection*{Glidepath transitions}
\textbf{Inputs.} Time-indexed hubs $K_t$. The template wires $K_t$ to $K_{t+1}$ via turnover and tax constraints, with tracking to a spoke-date benchmark. The result is a menu of valid next-step spokes.

\subsection{DOTS Beck--Chevalley and Frobenius}

Consider the commuting square
\[
\begin{tikzcd}
\Delta^n \arrow[r, "g"] \arrow[d, "f'"'] & \Delta^p \arrow[d, "f"] \\ 
\Delta^m \arrow[r, "h"'] & \Delta^q
\end{tikzcd}
\]
As before, we say that it is \emph{pointwise cartesian} if for all $y\in\Delta^p$, $z\in\Delta^m$ with $f(y)=h(z)$ there exists $x\in\Delta^n$ such that $g(x)=y$ and $f'(x)=z$ (informally: no new objects need to be created to make the square a pullback).

Define pullback and pushforward on relations by
\[
g^{*}(R):=\{(x,z): (g(x),z)\in R\},\qquad f_{!}(R):=\{(y,z): \exists x,\ f(x)=y,\ (x,z)\in R\}.
\]

\begin{theorem}[DOTS Beck--Chevalley]\label{thm:DOTS-BC}
Consider a commuting square of continuous maps between permissible spaces (objects):
\[
\begin{tikzcd}
\Delta^n \arrow[r, "g"] \arrow[d, "f'"'] & \Delta^p \arrow[d, "f"] \\
\Delta^m \arrow[r, "h"'] & \Delta^q
\end{tikzcd}
\]
If the square is pointwise cartesian and all maps are proper (automatic on compact domains), then for any closed alignment relation $R\subseteq\Delta^p\times Z$ (where $Z$ is a compact Hausdorff parameter space), the following equality of relations holds:
\[
f'_{!}(g^{*} R) \;=\; h^{*}(f_{!} R).
\]
\end{theorem}

\begin{proof}
We verify the equality of the relations as subsets of $\Delta^m \times Z$.

\noindent \textbf{($\subseteq$)} Let $(y,z)\in f'_{!}(g^{*}R)$. By definition of pushforward, there exists $x\in \Delta^n$ such that $f'(x)=y$ and $(x,z)\in g^{*}R$. The latter implies $(g(x),z)\in R$. Let $w=g(x)$. By the commutativity of the square, $f(w)=f(g(x))=h(f'(x))=h(y)$. Combining $f(w)=h(y)$ with $(w,z)\in R$ implies that $(h(y),z) \in f_{!}R$. By definition of pullback, this yields $(y,z)\in h^{*}(f_{!}R)$.

\noindent \textbf{($\supseteq$)} Let $(y,z)\in h^{*}(f_{!}R)$. Then $(h(y),z)\in f_{!}R$, meaning there exists $w\in \Delta^p$ such that $f(w)=h(y)$ and $(w,z)\in R$. Since the square is pointwise cartesian, there exists $x\in \Delta^n$ such that $g(x)=w$ and $f'(x)=y$. Substituting $w=g(x)$ into the relation condition yields $(g(x),z)\in R$, which means $(x,z)\in g^*R$. Since $f'(x)=y$, we conclude that $(y,z)\in f'_{!}(g^*R)$.
\end{proof}

\begin{theorem}[DOTS Frobenius]\label{thm:DOTS-Frob}
Let $f\colon\Delta^n\to\Delta^m$ be a continuous map (re-implementation). For any closed relation $R\subseteq\Delta^n\times Z$ (hub constraint) and $S\subseteq\Delta^m\times Z$ (spoke constraint) with $Z$ compact,
\[
f_{!}\big(R\cap f^{*}S\big)\;=\; f_{!}R \;\cap\; S.
\]
\end{theorem}

\begin{proof}
\noindent \textbf{($\subseteq$)} Let $(y,z)\in f_{!}(R\cap f^{*}S)$. There exists $x\in \Delta^n$ such that $f(x)=y$ and $(x,z)\in R\cap f^{*}S$. This implies $(x,z)\in R$ and $(x,z)\in f^*S$. By definition of pullback, $(f(x),z) \in S$, so $(y,z)\in S$. Since $f(x)=y$ and $(x,z)\in R$, we also have $(y,z)\in f_{!}R$. Thus $(y,z)\in f_{!}R \cap S$.

\noindent \textbf{($\supseteq$)} Let $(y,z)\in f_{!}R \cap S$. Then $(y,z)\in S$ and there exists $x\in \Delta^n$ with $f(x)=y$ and $(x,z)\in R$. Since $f(x)=y$, the condition $(y,z)\in S$ implies $(f(x),z)\in S$, which means $(x,z)\in f^*S$. Therefore, $(x,z)\in R \cap f^*S$. Combining this with $f(x)=y$, we obtain $(y,z)\in f_{!}(R \cap f^*S)$.
\end{proof}

\subsubsection*{When BC and Frobenius fail}
BC can fail if the square is not pointwise cartesian (missing preimages) or if some map loses properness (objects are not closed). Frobenius fails if images are forcibly closed \emph{after} the fact (post hoc closure), because topological closure does not distribute over intersections or composition. Financially this manifests as phantom spokes and order-sensitive surprises that are hard to debug.

\subsection{Case Studies: End-to-End Pipelines}
Each case specifies (i) inputs, (ii) relations and their order, (iii) points of determinization, (iv) guaranteed equalities (BC, Frobenius), and (v) computational notes. The aim is to make explicit which steps commute and which encode genuine choices.

\subsubsection*{Multi-manager blend with fee and factor budgets}
\textbf{Inputs.} $K_A\subseteq\Delta^{n_A}$, $K_B\subseteq\Delta^{n_B}$ (closed).\\
\textbf{Relations.} $R_{\mathrm{blend}}\subseteq(\Delta^{n_A}\oplus\Delta^{n_B})\times\Delta^{m}$, fee-cap projector $R_{\mathrm{fee}\le\tau}$, factor projector $R_{\mathrm{fac}\le b}$. Here $\oplus$ denotes the direct sum of simplices (juxtaposition of two independent weight vectors into a single higher-dimensional space) and $\otimes$ denotes the product of permissible spaces $K_A \times K_B$. The menu is
\[
K^{\star}=(K_A\otimes K_B)\odot R_{\mathrm{blend}}\odot R_{\mathrm{fee}\le\tau}\odot R_{\mathrm{fac}\le b}.
\]
\textbf{Determinization.} If blending solves a strictly convex program, replace $R_{\mathrm{blend}}$ by $\Graph(f_{\mathrm{blend}})$; then, by Frobenius, fee and factor projectors commute with $f_{\mathrm{blend}}$.
\\\textbf{Computation.} $R_{\mathrm{fac}\le b}$ and $R_{\mathrm{fee}\le \tau}$ are typically polyhedral; $R_{\mathrm{blend}}$ is often conic (SOCP). Use elimination to derive $K^{\star}$ or sample-and-hull for outer approximations, retaining certificates of inclusion.

\subsubsection*{Tax-aware ETF transition with turnover gate}
\textbf{Inputs.} Legacy state $K_{\mathrm{legacy}}\subseteq\Delta^{n}$.\\
\textbf{Relations.} $R_{\mathrm{turn}\le \kappa}$, $R_{\mathrm{tax}\le \theta}$, tracking $R_{\mathrm{track},\epsilon}$. The pipeline is
\[
K^{\star}=K_{\mathrm{legacy}}\odot R_{\mathrm{turn}\le \kappa}\odot R_{\mathrm{tax}\le \theta}\odot R_{\mathrm{track},\epsilon}.
\]
\textbf{Order.} If $R_{\mathrm{track},\epsilon}$ is functional (unique optimizer), Frobenius allows moving projectors across it. Otherwise, the order encodes a real choice: ``turnover-first'' versus ``tax-first'' yields different menus.\\
\textbf{Practice.} Apply the liquidity projector before tracking to avoid infeasible optimizers; log the exact order of operations for audit.

\subsubsection*{Factor neutralization overlay}
\textbf{Inputs.} $K_{\mathrm{base}}\subseteq\Delta^{n}$.\\
\textbf{Relations.} Overlay relation $R_{\mathrm{hedge}}$ targeting factor exposures and max gross notional; fee projector. The output is
\[
K^{\star}=K_{\mathrm{base}}\odot R_{\mathrm{hedge}}\odot R_{\mathrm{fee}\le \tau}.
\]
\textbf{Determinization.} With strictly convex penalties on residual exposure and notional, $R_{\mathrm{hedge}}=\Graph(f_{\mathrm{hedge}})$. By Frobenius, fee filtering may occur before or after $f_{\mathrm{hedge}}$ without changing the menu.

\subsubsection*{Benchmark-aware customization}
\textbf{Inputs.} Client hubs $K_{\mathrm{client}}\subseteq\Delta^{n}$; benchmark map $b\colon\Delta^{n}\to\Delta^{m}$.\\
\textbf{Relations.} Tracking band around $b(x)$; ESG projector; liquidity projector. The menu reads
\[
K^{\star}=K_{\mathrm{client}}\odot \Graph(b)\odot R_{\mathrm{track},\epsilon}\odot R_{\mathrm{ESG}}\odot R_{\mathrm{liq}}.
\]
\textbf{BC.} If a later factor map commutes pointwise cartesianly with $b$, BC allows reordering checks across maps without changing the menu.

\subsection{A Small Worked Example}
\label{subsec:worked-example}
\subsubsection*{Setup}
Hubs $\Delta^2$ with constraint $x_1\le 0.6$ (a closed face). Spokes $\Delta^2$ with fee map $\Fee(y)=10y_1+5y_2+0y_3$ (bps). Exposure maps $g_\A,g_\mathcal{B}$ both sum weights by two sleeves; tracking tolerance $\epsilon=0.05$ under identity covariance.

\subsubsection*{Relations}
$R_{\mathrm{track},\epsilon}=\{(x,y): \lVert g_A(x)-g_B(y)\rVert_2\le 0.05\}$,
$R_{\mathrm{fee}\le 6}=\{(y,y): \Fee(y)\le 6\}$.

\subsubsection*{Menus}
Compute $K\odot R_{\mathrm{track},\epsilon}$ by sampling hubs on a $21\times 21$ lattice inside the face and enumerating spoke lattice points at step $0.01$, filtering by the tracking ball. Intersect with the fee projector (idempotent) to obtain the final menu. The result visibly shrinks when the fee cap is applied first if optimizers are not unique. Make the optimizer unique by adding $\alpha\,\lVert y\rVert_2^2$ with $\alpha>0$, thereby determinizing the action.

\clearpage

\part{Probabilistic Theory}

\clearpage

\begin{center}
\begin{tikzpicture}[
    node distance=0.8cm and 0.5cm,
    >=stealth,
    font=\small,
    core/.style={
        rectangle,
        draw=blue!60!black,
        top color=blue!5,
        bottom color=blue!15,
        rounded corners=5pt,
        minimum width=13cm,
        minimum height=1.5cm,
        align=center,
        drop shadow,
        font=\bfseries
    },
    approach_box/.style={
        rectangle,
        draw=gray!50,
        fill=white,
        rounded corners=5pt,
        minimum width=13cm,
        minimum height=4.2cm,
        align=center,
        thick
    },
    header/.style={
        font=\bfseries,
        text=black!90,
        anchor=north west
    },
    math_label/.style={
        font=\scriptsize\bfseries,
        text=blue!60!black,
        fill=blue!5,
        inner sep=3pt,
        rounded corners=2pt,
        anchor=north east
    },
    desc_text/.style={
        align=left, 
        text width=7.5cm, 
        font=\footnotesize,
        anchor=west
    },
    viz_anchor/.style={
        anchor=north west 
    }
]

    \node[core] (Kernel) {
        The Stochastic Foundation\\
        \normalfont Horizontal Morphisms are \textbf{Tight Feller Kernels} $P: X \rightsquigarrow Y$\\
        \textit{\footnotesize ``Bounded inputs produce bounded probability mass.''}
    };


    \node[approach_box, below=1.0cm of Kernel] (Radius) {};
    
    \node[header] at ($(Radius.north west)+(0.3,-0.2)$) {1. Safety Radius};
    \node[math_label] at ($(Radius.north east)+(-0.2,-0.2)$) {Robust Safety (Geometry)};

    \node[viz_anchor] at ($(Radius.north west)+(0.5,-1.0)$) {
        \begin{tikzpicture}[scale=0.45]
            \draw[thick, green!60!black, fill=green!5] (-2,-1.5) rectangle (2,1.5);
            
            \node[green!60!black, anchor=south west, inner sep=1pt] at (2, -1.5) {$S$};
            
            \fill[black] (0,0) circle (3pt);
            \draw[thick, red!80, fill=red!10, fill opacity=0.5] (0,0) circle (1.2cm);
            \draw[->, red!80] (0,0) -- (1.2,0) node[midway, above, font=\tiny] {$r_\epsilon$};
        \end{tikzpicture}
    };

    \node[desc_text] at ($(Radius.north west)+(4.5,-2.1)$) {
        \textbf{Concept:} Metric Concentration\\
        \textbf{Math:} $d(f(x), S^c) \ge r_\epsilon$\\
        \textbf{Pros:} Fast ($O(1)$ check)\\
        \textbf{Cons:} Conservative\\
        \textbf{Use:} Real-time replication
    };

    \node[approach_box, below=0.6cm of Radius] (HDR) {};
    
    \node[header] at ($(HDR.north west)+(0.3,-0.2)$) {2. Highest Density Region};
    \node[math_label] at ($(HDR.north east)+(-0.2,-0.2)$) {Robust Safety (Probability)};

    \node[viz_anchor] at ($(HDR.north west)+(0.5,-1.0)$) {
        \begin{tikzpicture}[scale=0.45]
            \draw[thick, green!60!black, fill=green!5] (-2,-1.5) rectangle (2,1.5);
            
            \node[green!60!black, anchor=south west, inner sep=1pt] at (2, -1.5) {$S$};

            \draw[thick, red!80, fill=red!10, fill opacity=0.5]
                plot [smooth cycle] coordinates {(-1, -0.5) (0.5, -1) (1.2, 0) (0.5, 1) (-0.8, 0.8)};
            \node[red!80, font=\tiny] at (0,0) {$\operatorname{supp}_\epsilon$};
        \end{tikzpicture}
    };

    \node[desc_text] at ($(HDR.north west)+(4.5,-2.1)$) {
        \textbf{Concept:} Chance Constraint\\
        \textbf{Math:} $\operatorname{supp}_\epsilon(P) \subseteq S$\\
        \textbf{Pros:} Exact (Tightest set)\\
        \textbf{Cons:} Slow ($O(e^d)$)\\
        \textbf{Use:} Complex/Non-convex
    };

    \node[approach_box, below=0.6cm of HDR] (Wass) {};
    
    \node[header] at ($(Wass.north west)+(0.3,-0.2)$) {3. Transport (Cure)};
    \node[math_label] at ($(Wass.north east)+(-0.2,-0.2)$) {Economic Safety (Cost)};

    \node[viz_anchor] at ($(Wass.north west)+(0.5,-1.0)$) {
        \begin{tikzpicture}[scale=0.45]
            \draw[thick, green!60!black, fill=green!5] (1.5,-1.5) rectangle (3.5,1.5);
            
            \node[green!60!black, anchor=south west, inner sep=1pt] at (3.5, -1.5) {$S$};
            
            \draw[thick, red!80, fill=red!10, fill opacity=0.5] (-2, -0.5) circle (0.8);
            \node[red!80, font=\tiny] at (-2, -0.5) {$\mu$};
            
            \draw[->, thick, dashed, blue] (-1.2, -0.5) -- (1.5, -0.5);
            
            \node[blue, font=\tiny, align=center] at (-0.5, -0.5) {cost\\$\le \epsilon$};
        \end{tikzpicture}
    };

    \node[desc_text] at ($(Wass.north west)+(4.5,-2.1)$) {
        \textbf{Concept:} Economic Cost\\
        \textbf{Math:} $W_1(P, \mathcal{P}_S) \le \epsilon$\\
        \textbf{Pros:} Handles ``Near Misses''\\
        \textbf{Cons:} Requires LP/QP\\
        \textbf{Use:} Trading Budgets
    };

    \node[below=0.6cm of Wass, draw=gray, dashed, rounded corners, inner sep=8pt, fill=gray!5, text width=12.5cm, align=center] {
        \textbf{Key Insight:} Geometric approaches ($1, 2$) verify \textit{possibility}; Transport ($3$) calculates \textit{remediation}.
    };

    
    \coordinate (BusXCoord) at ($(Radius.west) + (-1.0, 0)$);
    
    \coordinate (BusTop) at (BusXCoord |- Kernel.west);
    \coordinate (BusBottom) at (BusXCoord |- Wass.west);
    
    \draw[very thick, gray!80] (BusTop) -- (BusBottom);
    
    \draw[very thick, gray!80] (Kernel.west) -- (BusTop);
    \fill[gray!80] (Kernel.west) circle (2pt);
    
    \foreach \mybox in {Radius, HDR, Wass} {
        \draw[->, very thick, gray!80] 
            (BusXCoord |- \mybox.west) -- (\mybox.west);
        \fill[gray!80] (BusXCoord |- \mybox.west) circle (2pt);
    }

\end{tikzpicture}
\end{center}

\clearpage

\section{Probabilistic Hub-and-Spoke Theory}
\label{sec:Polish}

\subsection{Motivation}

This variant framework, based on a suggestion of Greg Zitelli's, extends the basic topological theory to cover practical situations in which optimization is approximate, nondeterministic, or subject to ``measure zero'' exceptions. The underlying categorical structure mimics Part~I, but the mathematical foundations shift significantly---from point-set topology to infinite-dimensional analysis \cite{AliprantisBorder}, set-valued analysis \cite{AubinFrankowska, RockafellarWets}, and measure theory \cite{Bogachev, Billingsley}. The probabilistic theory is, as a consequence, considerably more technical than the basic theory. This is not merely a matter of exposition; the complications are genuine.

\subsection*{From Audit to Risk Management}

In the idealized $\HSSimp$ framework, we assumed portfolio re-implementation is deterministic. In the graded extension developed below (HSP-r, Definition~\ref{def:HSP-r}), we will see that \textit{any} possible outcome, no matter how unlikely, could in principle be tracked for audit purposes. However, production environments often require a \emph{probabilistic} perspective:

\begin{enumerate}
\item \textbf{Stochastic Solvers:} Real-world constraints may be non-convex, requiring heuristic solvers that produce distributions of outcomes.
\item \textbf{Tail Risk vs. Core Behavior:} A strategy should not be invalidated by a ``one-in-a-billion'' solver glitch or an extreme tail event. We require a formalism that ignores events below a specified probability threshold (the risk budget).
\item \textbf{Value-at-Risk (VaR):} Compliance rules are often probabilistic (e.g., ``The loss shall not exceed $L$ with 99\% confidence'').
\end{enumerate}

To model these phenomena, we extend our category to one in which horizontal morphisms are nondeterministic, and 2-cells carry a ``risk budget'' parameter $\epsilon$.

Informally, the picture is this: a deterministic re-implementation sends each hub portfolio to exactly one spoke portfolio. A stochastic re-implementation sends each hub portfolio to a \emph{probability distribution} over spoke portfolios. The question ``is the spoke aligned with the hub?'' becomes ``is the spoke aligned with the hub \emph{with high enough probability}?''---and the answer depends on a tolerance parameter $\epsilon$ that quantifies how much tail risk is acceptable. The mathematical machinery below (Polish spaces, tight Feller kernels, Chapman--Kolmogorov composition) exists to make this question precise and to ensure that the answer composes correctly across multiple stages. The reader who finds the machinery heavy should keep in mind that its purpose is modest: to define what ``almost always compliant'' means in a way that is stable under composition.

\subsection{The Graded Double Category HSP-r}
\label{sub:HSP-r}

We define a double category based on \textbf{tight Feller kernels} equipped with an $\epsilon$-support semantics. Unlike the strict audit framework, which tracks all possibilities, \textbf{HSP-r} (Hub-and-Spoke Risk) tracks the ``effective'' behaviour of the system up to a defined tolerance.

  \begin{definition}[HSP-r Data]
  \label{def:HSP-r}
  The $[0,1)$-\textbf{graded double category} \textbf{HSP-r}
  (Hub-and-Spoke Risk) consists of:
  \par
  \begin{itemize}
  \item \textbf{Objects:} Polish spaces $X$ (separable, completely metrizable
    topological spaces, equipped with their Borel $\sigma$-algebras).

  \item \textbf{Horizontal morphisms:} \textbf{Tight Feller kernels}
    $P: X \rightsquigarrow Y$ (Section~\ref{sub:tightness}).

  \item \textbf{Vertical morphisms:} \textbf{Closed relations}
    $R \subseteq X \times Z$.

  \item \textbf{Horizontal composition:} Chapman--Kolmogorov composition.
    Given kernels $P: X \rightsquigarrow Y$ and $Q: Y \rightsquigarrow Z$,
    their composite $(Q \circ P): X \rightsquigarrow Z$ is defined by
    \[
      (Q \circ P)(x, C) \;=\; \int_Y Q(y, C)\, P(x, dy)
    \]
    for any measurable $C \subseteq Z$.

  \item \textbf{Vertical composition:} Standard relational composition
    $S \circ R = \{(x, z) : \exists\, y,\; (x, y) \in R \text{ and }
    (y, z) \in S\}$.

  \item \textbf{Horizontal identities:} The Dirac (deterministic) kernel
    $\delta_{\operatorname{id}}: X \rightsquigarrow X$ defined by
    $\delta_{\operatorname{id}}(x, C) = \mathbf{1}_C(x)$.

  \item \textbf{Vertical identities:} The diagonal relation
    $\Delta_X = \{(x, x) \mid x \in X\}$.

  \item \textbf{2-Cells ($\epsilon$-graded):}
    Given a kernel $P: X \rightsquigarrow Y$ and closed relations
    $R \subseteq X \times Z$, $S \subseteq Y \times Z$ over a common
    base $Z$, a \textbf{2-cell at risk level $\epsilon \in [0,1)$} filling
    the square
    \begin{center}
    \begin{tikzcd}
    X \arrow[r, "P", squiggly] \arrow[d, dashed, "R"'] & Y
    \arrow[d, dashed, "S"] \\
    Z \arrow[r, equal] & Z
    \end{tikzcd}
    \end{center}
    exists if and only if $R \subseteq P^{*,\epsilon} S$, i.e.,
    for every $(x, z) \in R$:
    \begin{equation}\label{eq:HSP-r-2cell}
      P(x,\, S_z) \;\geq\; 1 - \epsilon
    \end{equation}
    where $S_z = \{y \in Y : (y, z) \in S\}$ denotes the $z$-slice of $S$.
    A 2-cell at level $\epsilon$ is automatically a 2-cell at any level
    $\epsilon' > \epsilon$.

    A general square with distinct bottom kernel $Q: Z \rightsquigarrow W$,
    \begin{center}
    \begin{tikzcd}
    X \arrow[r, "P", squiggly] \arrow[d, dashed, "R"'] & Y
    \arrow[d, dashed, "S"] \\
    Z \arrow[r, "Q"', squiggly] & W
    \end{tikzcd}
    \end{center}
    admits a 2-cell at risk level $\epsilon$ if for every $(x, z) \in R$:
    \begin{equation}\label{eq:HSP-r-2cell-full}
      \int_W P(x,\, S_w)\; Q(z, dw) \;\geq\; 1 - \epsilon
    \end{equation}
    \emph{Interpretation:} drawing outcomes $y \sim P(x, \cdot)$
    and $w \sim Q(z, \cdot)$ \emph{independently} from the respective marginal kernels, the pair $(y, w)$ lies
    in $S$ with probability at least $1 - \epsilon$ under the product measure $P(x, \cdot) \otimes Q(z, \cdot)$.  (The independence assumption is built into the product measure; in practice, this requires that the solver noise for different pipeline stages is uncorrelated.)  When
    $Q = \delta_g$ is a Dirac kernel,~\eqref{eq:HSP-r-2cell-full}
    reduces to $P(x, S_{g(z)}) \geq 1 - \epsilon$; when additionally
    $P = \delta_f$, it reduces to $(f(x), g(z)) \in S$, recovering the
    deterministic $\HSSimp$ condition.
  \end{itemize}

  \medskip
  \noindent\textbf{Composition of 2-cells.}
  Risk budgets accumulate \emph{additively} under composition:
  \begin{enumerate}
  \item \textbf{Horizontal composition} (stages of a pipeline):
    if the left square has a 2-cell at level $\delta$ and the right square
    at level $\epsilon$, the composed square has a 2-cell at level
    $\delta + \epsilon$ (Theorem~\ref{thm:pullback-comp}).  Concretely,
    the risk of a two-stage pipeline is bounded by the sum of the
    per-stage risks.

  \item \textbf{Vertical composition} (stacking alignment constraints):
    if the top square has a 2-cell at level $\delta$ and the bottom
    at level $\epsilon$, the composed square (with vertically composed
    relations $T \circ R$ and $U \circ S$) has a 2-cell at level
    $\delta + \epsilon$.

    \emph{Proof sketch.}
    Consider two vertically composable squares sharing a
    middle kernel $Q: Z \rightsquigarrow W$.  The top square has
    kernel $P: X \rightsquigarrow Y$, relations
    $R \subseteq X \times Z$, $S \subseteq Y \times W$,
    and 2-cell at level $\delta$.
    The bottom square has kernel $Q': V \rightsquigarrow V'$,
    relations $T \subseteq Z \times V$, $U \subseteq W \times V'$,
    and 2-cell at level $\epsilon$.
    The composed square has kernel $P$ (top), kernel $Q'$
    (bottom), relations $T \circ R$ and $U \circ S$.

    Let $(x, v) \in T \circ R$, witnessed by $z$ with
    $(x, z) \in R$ and $(z, v) \in T$.
    Introduce the auxiliary triple product measure
    $P(x, \cdot) \otimes Q(z, \cdot) \otimes Q'(v, \cdot)$
    on $Y \times W \times V'$.
    The top 2-cell gives
    $\Pr_{P \otimes Q}[(y, w) \in S] \ge 1 - \delta$,
    and the bottom 2-cell gives
    $\Pr_{Q \otimes Q'}[(w, v') \in U] \ge 1 - \epsilon$.
    By the union bound on the triple product:
    \[
    \Pr_{P \otimes Q \otimes Q'}
    \bigl[(y, w) \in S \text{ and } (w, v') \in U\bigr]
    \;\ge\; 1 - (\delta + \epsilon).
    \]
    On this event, $y \in S_w$ and $w \in U_{v'}$,
    hence $y \in (U \circ S)_{v'}$.
    Marginalising out $w$:
    $\Pr_{P \otimes Q'}[(y, v') \in U \circ S]
    \ge 1 - (\delta + \epsilon)$,
    which is the composed 2-cell condition.

  \item \textbf{Identity 2-cells:}
    For any relation $R \subseteq X \times Z$, the identity square
    (with Dirac kernel $\delta_{\operatorname{id}}$ on top and $R$ on
    both sides) admits a 2-cell at risk level $\epsilon = 0$, since
    $\delta_{\operatorname{id}}(x, R_z) = \mathbf{1}_{R_z}(x) = 1$
    for every $(x,z) \in R$.
  \end{enumerate}
  \end{definition}

This definition embeds the deterministic theory $\HSSimp$ as the ``zero-risk'' limit ($\epsilon=0$). Horizontal morphisms in $\HSSimp$ become Dirac kernels in \textbf{HSP-r}, i.e. deterministic kernels.

\subsection{Tight Feller Kernels: Modeling Imperfect Optimization}
\label{sub:tightness}

We restrict horizontal morphisms to \textbf{tight Feller kernels}. This ensures that the system is stable under perturbations and that the ``bulk'' of the probability mass remains within a compact region, allowing $\epsilon$-supports to be well-defined compact sets.

\begin{definition}[Feller Property]
A Markov kernel $P(x, dy)$ is \textbf{Feller} if the map $x \mapsto \int \phi(y) P(x, dy)$ is continuous for every bounded continuous function $\phi$.
\end{definition}
\textit{Interpretation:} Small changes in the input (hub) lead to small changes in the distribution of outputs (spoke). This is the analog of continuity in the deterministic case.

\begin{definition}[Tightness]
A kernel $P$ is \textbf{tight} if for every compact $K \subseteq X$ and $\delta > 0$, there exists a compact $L \subseteq Y$ such that $P(x, L) \ge 1 - \delta$ for all $x \in K$.
\end{definition}
\textit{Interpretation:} This guarantees that valid bounded inputs produce valid bounded outputs with high probability. Mass does not ``escape to infinity.'' This is the analog of properness in the deterministic case.

\begin{example}[Failure of Tightness: Unbounded Leverage]
Consider an optimization problem where short selling is permitted without gross exposure limits. If the solver detects an arbitrage opportunity, the solution sequence $\{y^{(k)}\}$ may diverge ($w_A \to +\infty, w_B \to -\infty$).
The probability mass describing the solution \emph{escapes}. For any compact subset $L$, $\lim_{k \to \infty} P(x, L) = 0$. This violates tightness, confirming that unbounded solvers are not valid morphisms.
\end{example}

\begin{remark}
The Feller and tightness conditions are natural in the sense that they generalize exactly the continuity and properness on which the deterministic theory depends. Whether they are \emph{checkable} for a given production solver is a different question. Most commercial optimizers do not expose enough of their internal structure to verify the Feller property analytically; one can test it empirically (perturb the input, observe the output distribution) but this is not a proof. The honest position is that the HSP-r framework provides the correct \emph{target} specification for a well-behaved stochastic pipeline, while verification that a particular implementation meets the specification remains, in most cases, an empirical matter.
\end{remark}

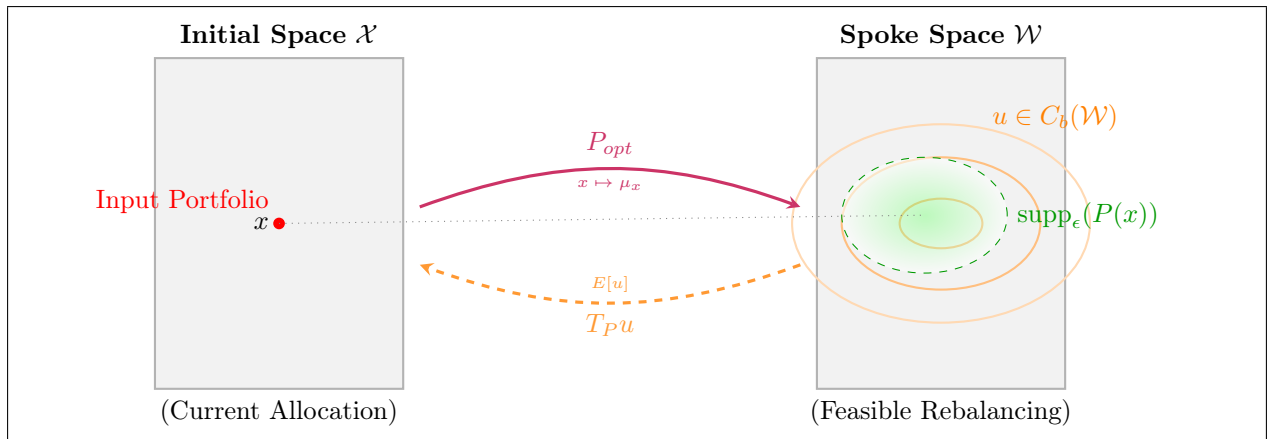
\begin{figure}[h]
    \centering
    \begin{tikzpicture}[scale=1.1, >=stealth, font=\small]

        \draw[thick, fill=gray!10, draw=gray!60] (-4, -2) rectangle (-1, 2);
        \node[above] at (-2.5, 2) {\textbf{Initial Space} $\mathcal{X}$};
        \node[below] at (-2.5, -2) {(Current Allocation)};

        \coordinate (x) at (-2.5, 0);
        \fill[red] (x) circle (2pt) node[left, black] {$x$};
        \node[red, above left] at (x) {Input Portfolio};

        \draw[thick, fill=gray!10, draw=gray!60] (4, -2) rectangle (7, 2);
        \node[above] at (5.5, 2) {\textbf{Spoke Space} $\mathcal{W}$};
        \node[below] at (5.5, -2) {(Feasible Rebalancing)};

        \draw[thick, orange!30] (5.5, 0) ellipse (1.8cm and 1.2cm);
        \draw[thick, orange!50] (5.5, 0) ellipse (1.2cm and 0.8cm);
        \draw[thick, orange!80] (5.5, 0) ellipse (0.5cm and 0.3cm);
        
        \node[orange, above right] at (6.0, 1.0) {$u \in C_b(\mathcal{W})$};

        \shade[inner color=green!40, outer color=white, opacity=0.6] (5.3, 0.1) ellipse (1.0cm and 0.7cm);
        \draw[green!60!black, dashed] (5.3, 0.1) ellipse (1.0cm and 0.7cm);
        \node[green!60!black, right] at (6.3, 0.1) {$\supp_\epsilon(P(x))$};

        \draw[->, very thick, purple!80] (-0.8, 0.2) to[out=20, in=160] 
            node[midway, above] {$P_{opt}$} 
            node[midway, below, font=\tiny] {$x \mapsto \mu_x$}
            (3.8, 0.2);

        \draw[->, very thick, dashed, orange!80] (3.8, -0.5) to[out=200, in=-20] 
            node[midway, below] {$T_{P} u$} 
            node[midway, above, font=\tiny] {$E[u]$} 
            (-0.8, -0.5);

        \draw[dotted, gray] (x) -- (5.3, 0.1);
    \end{tikzpicture}
    \caption{Stochastic Optimization Consistency in HSP-r. The ``smear'' represents the $\epsilon$-support: the region containing the vast majority of the probability mass. Tail events outside this region are discarded by the risk budget.}
    \label{fig:portfolio_feller_kernel}
\end{figure}

\begin{example}[Production optimizers: \textbf{Gurobi}]
    \label{ex:gurobi-feller}
    Standard MIP solvers are discontinuous (not Feller). To use them in HSP-r, we apply \textit{randomized smoothing} by adding noise $\sigma \xi$ to the objective (where $\sigma > 0$ is the noise magnitude and $\xi$ is a standard random variable, e.g.\ Gaussian).
    \[
        P(x, dy) = \text{Law}\left( \operatorname{argmin}_{y \in K} \left[ C(x,y) + \sigma \xi \right] \right)
    \]
    While this creates a kernel with potentially infinite topological support (if $\xi$ is Gaussian), the \textbf{$\epsilon$-support} concentrates tightly around the optimal solution lattice. The HSP-r framework allows us to ignore the Gaussian tails that extend to the boundary of the feasible region, treating the solver as effectively ``sharp'' up to a confidence interval.
\end{example}

\begin{example}[Production optimizers: \textbf{Axioma}]
    Similarly, a smoothed Axioma kernel $P_\sigma$ is tight because the long-only budget constraint is compact. While the strict topological support is the entire simplex (due to Gaussian smoothing of alpha signals), the $\epsilon$-support correctly identifies the ``likely'' portfolios.
\end{example}

\subsection{Probabilistic Operations: Two Approaches}
\label{sub:stoch_ops}

To define the operations of pullback and pushforward in a probabilistic setting, we must identify a geometric object that captures the ``bulk'' of the probability mass. We explore two approaches: \emph{concentration of measure} and \emph{highest density region}. Before developing the results, we introduce and compare these two approaches.

The first approach (suggested by Greg Zitelli) models the safety region as a metric ball centered on the solver's expected output. This approach relies on foundational concentration inequalities. This approach works well for production solvers (e.g. Gurobi, Axioma) in cases where there is a global optimum and the error is typically sub-Gaussian.

\begin{definition}[Centered Kernel and Safety Radius]
Let $P: X \rightsquigarrow Y$ be a tight Feller kernel between metric spaces $(X, d_X)$ and $(Y, d_Y)$. We assume $P$ is a perturbation of a deterministic center map $f: X \to Y$ (e.g., the noise-free solution).
For a risk budget $\epsilon > 0$, the \textbf{Safety Radius} $r_\epsilon$ is a scalar satisfying:
\[
P(x, B(f(x), r_\epsilon)) \ge 1 - \epsilon \quad \text{for all } x \in X
\]
where $B(c, r) = \{ y \in Y \mid d_Y(c, y) \le r \}$ denotes the closed metric ball.
\end{definition}

\begin{remark}[Derivation via Concentration Inequalities]
If the solver error is bounded or sub-Gaussian, $r_\epsilon$ can be derived analytically. For example, using Hoeffding's inequality for bounded perturbations \cite{Hoeffding}, or Talagrand's concentration inequality for product spaces \cite{Talagrand}, one typically obtains a radius scaling as $r_\epsilon \propto \sqrt{\ln(1/\epsilon)}$. This replaces the complex density integration of the HDR approach with a single global parameter.
\end{remark}

The second approach tracks the exact density contours.

\begin{definition}[$\epsilon$-Highest Density Region]
\label{def:epsilon-supp}
Let $Y$ be a Polish space and $P(x, \cdot)$ a probability measure on $Y$ with a continuous density function $\rho_x$.
The \textbf{$\epsilon$-Highest Density Region ($\epsilon$-HDR)}, denoted $\operatorname{supp}_{\epsilon}(P(x))$, is the super-level set of the density containing probability mass $1-\epsilon$:
\[
\operatorname{supp}_{\epsilon}(P(x)) := \{ y \in Y \mid \rho_x(y) \ge \lambda_\epsilon \}
\]
where $\lambda_\epsilon = \sup \{ \lambda \ge 0 \mid P(x, \{y \mid \rho_x(y) \ge \lambda\}) \ge 1-\epsilon \}$.
\end{definition}

\begin{remark}
If the density $\rho_x$ has no flat regions (i.e., the measure of level sets $\{y \mid \rho_x(y) = c\}$ is zero for all $c$), then $\lambda_\epsilon$ is the unique threshold satisfying the mass condition with equality.
\end{remark}

\begin{remark}[Stability of Production Solvers]
For the \textbf{Gaussian-smoothed Gurobi} example (Example~\ref{ex:gurobi-feller}), the addition of Gaussian noise $\sigma \xi$ typically induces a quasi-concave (unimodal) density around the optimal lattice points.
This ensures that the $\epsilon$-HDR is a finite union of compact convex sets (or a single convex set in the relaxed case) that evolves continuously with the input $x$.
Similarly, for \textbf{Axioma}, signal smoothing ensures the alpha landscape is sufficiently regular that the HDR is stable.
Unlike the intersection definition, the HDR remains robust for both discrete and continuous error models.
\end{remark}

\paragraph{Computational Tractability.}
In the HDR approach, calculating the density-based $\epsilon$-support of a distribution defined by a solver (e.g., a convex optimization problem with random perturbations) requires numerical integration or Monte Carlo sampling to delineate the boundary of the high-density region. This is $O(e^d)$ in dimension $d$.
In contrast, the Safety Radius approach requires only the calculation of the Euclidean distance $d(f(x), S_z^c)$. For polyhedral constraints (common in finance), this is a standard convex projection problem, solvable in polynomial time.

\paragraph{Architecture Impact.}
For the Safety Radius approach, in the system architecture described in Appendix~\ref{app:architecture}, the \emph{Verification Worker} can now be stateless and purely geometric. It does not need to simulate thousands of solver runs to estimate a failure probability, as it does for the HDR approach. Instead, it queries the \emph{Solver Interface} for the deterministic solution $f(x)$ and the pre-calibrated Lipschitz constant/variance parameters to determine $r_\epsilon$. It then performs a geometric containment check against the constraints stored in the \emph{Morphism Registry}.

\begin{table}[htbp]
    \centering
\begin{tabular}{p{0.2\textwidth} p{0.35\textwidth} p{0.35\textwidth}}
\toprule
\textbf{Feature} & \textbf{Safety Radius} & \textbf{Highest Density Region} \\
\midrule
\textbf{Geometry} & Symmetric Metric Balls (Spheres) & Irregular Density Contours \\
\textbf{Verification Cost} & \textbf{Low} (Polynomial time distance check) & \textbf{High} (Numerical integration or Sampling) \\
\textbf{Precision} & \textbf{Conservative} (Approximation) & \textbf{Exact} (Smallest possible set) \\
\textbf{Composition} & \textbf{Linear} (Triangle Inequality) & \textbf{Complex} (Convolution Stability) \\
\textbf{Ideal Use Case} & Real-time production software & Theoretical analysis or low-dimensional models \\
\bottomrule
\end{tabular}    
\caption{Comparison of the Safety Radius and Highest Density Region (HDR) approaches}
    \label{tab:comparison}
\end{table}

The Safety Radius approach is generally preferred in a portfolio replication context. Here are some examples where portfolio construction can generate irregular probability distributions (multimodal or ``kidney-shaped''), forcing the use of the HDR approach.

\paragraph*{Minimum Position Sizes (The ``Split Peak'')}

This is the most common cause of \emph{multimodal} distributions in institutional systems.

\begin{itemize}
    \item \emph{The Constraint:} A rule requiring that if an asset is held, it must constitute at least a certain percentage of the portfolio (e.g., ``Hold $0\%$ or $\ge 2\%$''). This is common to avoid ``dust'' positions that are costly to trade.
    \item \emph{The Scenario:} An optimizer wants to allocate an ideal weight of $1.0\%$ to a stock.
    \item \emph{The Resulting Distribution:}
    \begin{itemize}
        \item Because $1.0\%$ is illegal, a randomized solver (or a solver with input noise) will ``dither.''
        \item In 50\% of runs, it rounds \emph{down} to $0\%$.
        \item In 50\% of runs, it forces \emph{up} to $2\%$.
    \end{itemize}
    \item \emph{Why Metric Concentration Fails:}
    \begin{itemize}
        \item The ``average'' or ``center'' of this distribution is $1\%$.
        \item A metric ball centered at $1\%$ suggests that $1\%$ is a valid, safe holding.
        \item In reality, $1\%$ is explicitly \emph{forbidden}. The Highest Density Region (HDR) correctly identifies two separate islands of safety ($0$ and $2+$), while the metric ball includes the ``valley of death'' between them.
    \end{itemize}
\end{itemize}

\paragraph*{Tax-Aware Transitioning (The ``Kidney'' or ``Banana'')}

This scenario creates \emph{non-convex, curved} distributions, often referred to as banana-shaped or kidney-shaped.

\begin{itemize}
    \item \emph{The Constraint:} A combination of a maximum tracking error (an ellipsoid constraint) and a maximum realized capital gain tax budget (a linear or piecewise-linear constraint).
    \item \emph{The Scenario:} You are transitioning a legacy portfolio to a new target. You want to minimize distance to the target, but you have a strict limit on how much tax you can generate.
    \item \emph{The Resulting Distribution:}
    \begin{itemize}
        \item The solver tries to sell assets to buy the new target.
        \item It hits the tax cap. To get closer to the target without selling high-cost basis stock, it starts buying ``proxy'' assets (highly correlated substitutes) that don't trigger taxes.
        \item As correlation shifts, the efficient frontier of solutions curves around the tax obstacle.
    \end{itemize}
    \item \emph{Why Metric Concentration Fails:}
    \begin{itemize}
        \item The set of valid solutions curves around the ``expensive'' tax event.
        \item A metric ball centered on the mean would cut \textit{through} the tax constraint violation region.
        \item The HDR correctly wraps around the obstacle in a kidney shape, excluding the high-tax center.
    \end{itemize}
\end{itemize}

\paragraph*{Regime Switching Signals (The ``Bifurcation'')}

This generates \emph{highly separated multimodal} distributions.

\begin{itemize}
    \item \emph{The Constraint:} A quantitative signal that acts as a binary switch. For example, ``If Volatility $< 15$, leverage is $150\%$; if Volatility $> 15$, leverage is $100\%$.''
    \item \emph{The Scenario:} Current market volatility is measuring exactly at 15.0, but there is measurement noise ($\pm 0.1$).
    \item \emph{The Resulting Distribution:}
    \begin{itemize}
        \item Small perturbations in the input data cause the solver to flip aggressively between two distinct portfolios (High Leverage vs. Low Leverage).
        \item The output distribution has two dense clouds of solutions that are far apart in Euclidean space.
    \end{itemize}
    \item \emph{Why Metric Concentration Fails:}
    \begin{itemize}
        \item The mean of these two clouds is a portfolio with 125\% leverage.
        \item However, the strategy \textit{never} allows 125\% leverage.
        \item A metric ball would declare the 125\% leverage portfolio ``safe,'' creating a false positive for audit.
    \end{itemize}
\end{itemize}

\paragraph*{Long/Short Gross Exposure Limits (The ``Hollow Shell'')}

This creates distributions where the ``center'' is empty.

\begin{itemize}
    \item \emph{The Constraint:} A fund must maintain exactly $0\%$ net exposure (dollar neutral) but maximize exposure to a specific factor, subject to a Gross Exposure limit of $200\%$.
    \item \emph{The Scenario:} The solver pushes positions to the absolute limit of the Gross Exposure constraint to maximize alpha.
    \item \emph{The Resulting Distribution:}
    \begin{itemize}
        \item The solutions essentially live on the ``skin'' (boundary) of a sphere defined by the Gross Exposure limit.
        \item Random noise pushes the solution along the surface of this sphere, not toward the center.
    \end{itemize}
    \item \emph{Why Metric Concentration Fails:}
    \begin{itemize}
        \item The expected value (mean) of points distributed on a sphere's surface is inside the sphere (lower gross exposure).
        \item While the mean is technically ``valid'' (it satisfies the limit), it is distinct from the actual behaviour of the solver, which is always at the boundary. The metric ball overestimates the stability of the solution by assuming mass exists in the interior, where there is none.
    \end{itemize}
\end{itemize}

\section{Probabilistic Operations: Safety Radius}
\label{sec:safety-radius}

Before developing either approach in detail, we state the central trade-off.  The Safety Radius approach (this section) replaces the solver's output distribution with a metric ball $B(f(x), r_\epsilon)$ centered on the expected output.  This is geometrically clean: balls compose under Minkowski sums, risk budgets accumulate linearly, and the Frobenius law holds as a tight geometric equality.  The cost is a spherical assumption --- the error ball is symmetric, so non-spherical distributions (bimodal solvers, banana-shaped tax constraints) are either missed or conservatively over-approximated.

The Highest Density Region (HDR) approach (Section~\ref{sec:HDR-section}) replaces the ball with the smallest high-probability region of the actual distribution.  This adapts to the distribution's shape, but at a price: composition requires a Convolution Stability Condition (Proposition~\ref{prop:convolution-stability}), and the Frobenius law degrades from an equality to a one-way inclusion.  The difference between the two sides represents \emph{tail risk leakage} --- portfolios that pass a chance constraint despite carrying residual tail risk.  For multi-stage pipelines, this asymmetry is the reason to prefer robust constraints (the Safety Radius or Geometric Frobenius) at intermediate stages and to reserve the HDR for terminal diagnostics.  Section~\ref{subsec:compositional-safety} develops this comparison in detail after both approaches have been presented.

\subsection{Safety Radius Approach}
\label{sub:concentration}

\begin{definition}[Probabilistic Pullback]
Given a kernel $P$ centered at $f$ with safety radius $r_\epsilon$, and a relation $S \subseteq Y \times Z$, the \textbf{$\epsilon$-pullback} identifies hubs that are ``deeply'' compliant:
\[
P^{*,\epsilon} S = \{ (x, z) \in X \times Z \mid d_Y(f(x), S_z^c) \ge r_\epsilon \}
\]
where $S_z = \{y \mid (y,z) \in S\}$ is the slice of permissible spokes, and $d_Y(y, A) = \inf_{a \in A} d_Y(y, a)$ is the distance to a set.
\end{definition}

\noindent \textbf{Interpretation:} A hub portfolio $x$ is deemed safe if its expected re-implementation $f(x)$ lies deep enough inside the permissible spoke region $S_z$ such that the distance to the boundary ($S_z^c$) exceeds the safety radius. Geometrically, this ensures $B(f(x), r_\epsilon) \subseteq S_z$.

\begin{definition}[Probabilistic Pushforward]
Given a kernel $P$ centered at $f$ with safety radius $r_\epsilon$, and a relation $R \subseteq X \times Z$, the \textbf{$\epsilon$-pushforward} propagates the safety buffer:
\[
P_!^\epsilon R = \overline{ \bigcup_{(x, z) \in R} \left( B(f(x), r_\epsilon) \times \{z\} \right) } \subseteq Y \times Z
\]
\end{definition}

\noindent \textbf{Interpretation:} $P_!^\epsilon R$ is simply the Minkowski sum \cite{Schneider} of the deterministic image of $R$ and the error ball of radius $r_\epsilon$. It answers: ``Where might the spoke portfolio land?'' by inflating the deterministic projection $f(R)$ by the worst-case deviation $r_\epsilon$.

\begin{remark}
Pullback and pushforward along an identity morphism are no longer trivial. When the center map is the identity $f = \mathrm{id}_X$, the operations correspond to metric erosion and dilation by the safety radius $r_\epsilon$:

\begin{itemize}
    \item \textbf{Pullback (Erosion/Inner Parallel Set):} 
    The $\epsilon$-pullback identifies points deep enough inside the permissible region $S$ that they remain valid even after a worst-case perturbation of size $r_\epsilon$:
    \[
    \mathrm{id}^{*,\epsilon} S = \{ (x, z) \in X \times Z \mid d(x, S_z^c) \ge r_\epsilon \} = \{ (x,z) \mid B(x, r_\epsilon) \subseteq S_z \}
    \]
    
    \item \textbf{Pushforward (Dilation/Minkowski Sum):} 
    The $\epsilon$-pushforward expands the valid relation $R$ by the safety ball, representing all locations reachable by the "jitter" $r_\epsilon$:
    \[
    \mathrm{id}_!^\epsilon R = R \oplus B(0, r_\epsilon) = \overline{ \bigcup_{(x, z) \in R} \left( B(x, r_\epsilon) \times \{z\} \right) }
    \]
\end{itemize}
\end{remark}

\subsection{Composition and Accumulation of Risk}
\label{sub:accumulation}

To analyse multi-stage pipelines (e.g., Hub $\to$ Model $\to$ Spoke), we must determine how the safety radius accumulates across compositions. The choice of accumulation function determines the logical relationship between the \emph{step-by-step} audit and the \emph{end-to-end} audit.

\subsubsection*{Defining the Composite Safety Radius}

Let $P: X \rightsquigarrow Y$ and $Q: Y \rightsquigarrow Z$ be tight Feller kernels centered at $f_P$ and $f_Q$ with safety radii $r_P$ and $r_Q$ respectively. Assume the center map $f_Q$ is $L$-Lipschitz continuous.

The \textbf{Sequential Risk Budget} represents the worst-case stack-up of errors:
\[ r_{seq} \;=\; L r_P + r_Q \]

However, in many high-dimensional portfolio applications, errors are uncorrelated (e.g., orthogonal Gaussian noise from randomized smoothing). In such cases, the Linear budget is overly conservative. We therefore define a \textbf{Composite Safety Radius} consistent with the independence of the error terms.

\begin{definition}[Quadratic Composition]
We define the quadratic composition radius $r_{quad}$ for the kernel $Q \circ P$ based on independence of errors:
    \[ r_{quad} = \sqrt{(L r_P)^2 + r_Q^2} \]
Since $\sqrt{a^2+b^2} \le a+b$ for non-negative reals (square both sides), we have $r_{quad} \le r_{seq}$.
\end{definition}

\subsubsection*{Compositionality: Linear Composition (Worst Case Error Accumulation)}

\begin{theorem}[Composition of Safety Radius Operations, Deterministic Case]
Let $f: X \to Y$ and $g: Y \to Z$ be deterministic maps. Let $r_f, r_g > 0$ be the safety radii associated with $f$ and $g$ respectively. Assume $g$ is $L$-Lipschitz continuous. The composition of probabilistic operations follows the linearity of the risk budget, resulting in the following conservative inclusions:

\begin{enumerate}
    \item \textbf{Pushforward Composition (Conservative Risk Envelope):} 
    The detailed sequential propagation is contained within the scalar approximation:
    \[ g_!^{r_g} \left( f_!^{r_f} R \right) \;\subseteq\; (g \circ f)_!^{L r_f + r_g} R \]
    
    \item \textbf{Pullback Composition (Conservative Verification):}
    The scalar approximation imposes a stricter verification standard than the sequential check:
    \[ (g \circ f)^{*, L r_f + r_g} S \;\subseteq\; f^{*, r_f} \left( g^{*, r_g} S \right) \]
\end{enumerate}
\end{theorem}

\begin{proof}
\textbf{1. Pushforward:}
Let $z \in g_!^{r_g} ( f_!^{r_f} R )$. By definition, there exists an intermediate point $y \in f_!^{r_f} R$ such that $d_Z(z, g(y)) \le r_g$.
Since $y \in f_!^{r_f} R$, there exists a hub $x$ with $(x, \cdot) \in R$ such that $d_Y(y, f(x)) \le r_f$.
We evaluate the distance of $z$ from the composite center $g(f(x))$:
\begin{align*}
d_Z(z, g(f(x))) &\le d_Z(z, g(y)) + d_Z(g(y), g(f(x))) \\
&\le r_g + L \cdot d_Y(y, f(x)) \quad \text{(by Lipschitz property)} \\
&\le r_g + L r_f
\end{align*}
Thus, $z$ lies in the ball $B_Z(g(f(x)), L r_f + r_g)$, which implies $z \in (g \circ f)_!^{L r_f + r_g} R$.

\textbf{2. Pullback:}
Let $(x, w) \in (g \circ f)^{*, L r_f + r_g} S$. This implies that the large ball $B_Z(g(f(x)), L r_f + r_g)$ is entirely contained in the safe set $S_w$.
We must show $(x, w) \in f^{*, r_f} ( g^{*, r_g} S )$. This requires that for all $y \in B_Y(f(x), r_f)$, the ball $B_Z(g(y), r_g)$ is contained in $S_w$.
Let $y$ be any point such that $d_Y(y, f(x)) \le r_f$.
Let $z$ be any point such that $d_Z(z, g(y)) \le r_g$.
By the triangle inequality and Lipschitz condition, we established above that:
\[
d_Z(z, g(f(x))) \le r_g + L r_f.
\]
Since the composite ball is contained in $S_w$, $z \in S_w$.
Thus, for every likely intermediate $y$, every resulting $z$ is safe. The sequential condition is satisfied.
\end{proof}

\begin{theorem}[Composition of Safety Radius Operations, Probabilistic Case]
Let $P: X \rightsquigarrow Y$ and $Q: Y \rightsquigarrow Z$ be tight Feller kernels.
Assume $P$ is centered at $f_P$ with safety radius $r_P$, and $Q$ is centered at $f_Q$ with safety radius $r_Q$.
Assume the center map $f_Q$ is $L$-Lipschitz continuous.
We define the composite kernel $Q \circ P$ to be centered at $f_Q \circ f_P$ with composite safety radius $r_{seq} = L r_P + r_Q$.

The probabilistic operations satisfy the following conservativeness inclusions:

\begin{enumerate}
    \item \textbf{Pushforward Composition (Risk Envelope):}
    The sequential safety region is contained within the composite safety envelope:
    \[
    Q_!^{r_Q} \left( P_!^{r_P} R \right) \;\subseteq\; (Q \circ P)_!^{r_{seq}} R
    \]

    \item \textbf{Pullback Composition (Verification Safety):}
    Composite verification is a conservative proxy for sequential verification:
    \[
    (Q \circ P)^{*, r_{seq}} S \;\subseteq\; P^{*, r_P} \left( Q^{*, r_Q} S \right)
    \]
\end{enumerate}
\end{theorem}

\begin{proof}
\textbf{1. Pushforward Inclusion:}
Recall that the probabilistic pushforward $P_!^{r_P} R$ is the union of balls $B_Y(f_P(x), r_P)$.
Let $z$ be an element of the sequential pushforward $Q_!^{r_Q} (P_!^{r_P} R)$.
This means $z$ belongs to a ball $B_Z(f_Q(y), r_Q)$ for some $y$ which itself belongs to a ball $B_Y(f_P(x), r_P)$.
Using the Lipschitz property of $f_Q$, the image of the ball $B_Y(f_P(x), r_P)$ is contained within $B_Z(f_Q(f_P(x)), L r_P)$.
The Minkowski sum of this image with the error budget $r_Q$ is contained within the ball of radius $L r_P + r_Q$.
Therefore, $z \in B_Z(f_Q(f_P(x)), L r_P + r_Q)$, validating the inclusion.

\textbf{2. Pullback Inclusion:}
Let $(x, w) \in (Q \circ P)^{*, r_{seq}} S$. This implies that the entire ball $B_Z(f_Q(f_P(x)), L r_P + r_Q)$ is contained in the safe slice $S_w$.
We need to show that $x$ is valid sequentially.
Sequential validity requires that for any likely intermediate outcome $y \in B_Y(f_P(x), r_P)$, the resulting outcome $z$ from $Q$ is safe.
The set of all such reachable $z$ is contained in the union of $r_Q$-balls centered on the image of the $r_P$-ball:
\[ \mathcal{Z} = \bigcup_{y \in B_Y(f_P(x), r_P)} B_Z(f_Q(y), r_Q). \]
As shown in the deterministic case, $\mathcal{Z} \subseteq B_Z(f_Q(f_P(x)), L r_P + r_Q)$.
Since the larger composite ball is safe (contained in $S_w$), the subset $\mathcal{Z}$ is also safe.
Therefore, the sequential check passes.
\end{proof}

\subsubsection*{Compositionality: Quadratic Composition (Independent Errors)}

Adopting the Quadratic definition allows us to derive tighter bounds for the composite system. This yields the structurally important property that \emph{sequential verification is a conservative proxy for end-to-end verification}.

These results require an additional structural condition on the downstream center map $f_Q$, ensuring that it does not collapse directions in the output space.

\begin{definition}[Metric Openness]
\label{def:metric-openness}
A map $f_Q: Y \to Z$ between metric spaces is \textbf{$\ell$-co-Lipschitz} (or \textbf{metrically open with constant $\ell > 0$}) if for all $y \in Y$ and $r > 0$:
\[
f_Q(B_Y(y, r)) \supseteq B_Z(f_Q(y), \ell r).
\]
The co-Lipschitz constant satisfies $\ell \le L$ (the Lipschitz constant), with equality when $f_Q$ is a scalar multiple of an isometry.
\end{definition}

\noindent We define the \textbf{Coverage Radius}:
\[ r_{cov} \;=\; \ell r_P + r_Q \]
which measures the guaranteed reach of the sequential risk tube. Since $\ell \le L$, we have $r_{cov} \le r_{seq}$.

\begin{remark}[Financial Interpretation of Metric Openness]
\label{rem:metric-openness}
Where the Lipschitz constant $L$ bounds the \emph{maximum} amplification of errors (the gross leverage), the co-Lipschitz constant $\ell$ bounds the \emph{minimum}: a unit perturbation in model space produces at least an $\ell$-unit displacement in spoke space. For a linear re-implementation $f_Q(y) = Ay + b$, we have $L = \sigma_{\max}(A)$ and $\ell = \sigma_{\min}(A)$ (the largest and smallest singular values of $A$).

The condition $\ell > 0$ means the solver is \emph{non-degenerate}---it does not collapse any direction in the output space. This holds for full replication ($\ell = L = 1$), leveraged strategies ($\ell = L = |\lambda|$), and any optimizer whose Jacobian has full rank. It fails when the solver maps a high-dimensional model space into a lower-dimensional spoke space with rank deficiency---for instance, a strategy that ignores certain asset classes entirely, making some output directions unreachable from model perturbations. In such degenerate cases, the sequential audit tube does not cover the full composite ball, and the results below require modification (e.g., restricting the composite ball to the image of $f_Q$).
\end{remark}

\begin{theorem}[Lax Functoriality: Risk Concentration]
Assume the center map $f_Q$ is $\ell$-co-Lipschitz (Definition~\ref{def:metric-openness}), and assume $r_{quad} \le r_{cov}$, where $r_{quad} = \sqrt{(L r_P)^2 + r_Q^2}$ is the quadratic composition radius and $r_{cov} = \ell r_P + r_Q$ is the coverage radius.  (This holds when the condition number $L/\ell$ of $f_Q$ is moderate; it can fail for ill-conditioned maps.) The probabilistic pushforward satisfies the following inclusion:
\[
(Q \circ P)_!^{r_{quad}} R \;\subseteq\; Q_!^{r_Q} \left( P_!^{r_P} R \right)
\]
\end{theorem}

\begin{proof}
We must show that every point in the composite risk envelope (a ball of radius $r_{quad}$) is covered by the sequential risk tube.  Throughout the proof, we use the standing hypothesis $r_{quad} \le r_{cov} := \ell\, r_P + r_Q$.

Let $(z, w) \in (Q \circ P)_!^{r_{quad}} R$. By the definition of the pushforward, there exists $x$ with $(x, w) \in R$ such that $d_Z(z, f_Q(f_P(x))) \le r_{quad}$.

We exhibit a witness for membership in $Q_!^{r_Q}(P_!^{r_P} R)$ by finding $y' \in B_Y(f_P(x), r_P)$ with $d_Z(z, f_Q(y')) \le r_Q$.

Consider the Minkowski-sum tube $\mathcal{T} = \bigcup_{y \in B_Y(f_P(x), r_P)} B_Z(f_Q(y), r_Q)$. This is the $r_Q$-neighborhood of the image $I = f_Q(B_Y(f_P(x), r_P))$. By the co-Lipschitz property of $f_Q$ (Definition~\ref{def:metric-openness}), $I \supseteq B_Z(f_Q(f_P(x)), \ell r_P)$. In a normed space, the distance from a point to a closed ball satisfies $d(z, B(a, r)) = \max(0,\, \|z - a\| - r)$; this follows from the triangle inequality: $\|z - a\| \le \|z - y\| + r$ for any $y \in B(a,r)$, and the bound is attained by $y = a + r(z-a)/\|z-a\|$. Our portfolio weight spaces are convex subsets of $\mathbb{R}^n$ with the induced norm, so the formula applies. Therefore:
\[
d_Z(z, I) \;\le\; d_Z(z, B_Z(f_Q(f_P(x)), \ell r_P)) \;=\; \max\bigl(0,\, d_Z(z, f_Q(f_P(x))) - \ell r_P\bigr) \;\le\; \max(0,\, r_{cov} - \ell r_P) \;=\; r_Q.
\]
Hence there exists $w_0 \in I$ with $d_Z(z, w_0) \le r_Q$. Since $w_0 \in I = f_Q(B_Y(f_P(x), r_P))$, there exists $y' \in B_Y(f_P(x), r_P)$ with $f_Q(y') = w_0$. Then $d_Z(z, f_Q(y')) \le r_Q$, and since $(x, w) \in R$ with $d_Y(y', f_P(x)) \le r_P$, we have $y' \in P_!^{r_P} R_w$, placing $z$ in $Q_!^{r_Q}(P_!^{r_P} R)$.
\end{proof}

\begin{theorem}[Conservative Sequential Verification]
Assume the same co-Lipschitz hypothesis as above ($r_{quad} \le r_{cov} = \ell r_P + r_Q$). Passing the sequential audit implies passing the composite audit:
\[
P^{*, r_P} \left( Q^{*, r_Q} S \right) \;\subseteq\; (Q \circ P)^{*, r_{quad}} S
\]
\end{theorem}

\begin{proof}
Let $(x, w)$ be a pair accepted by the sequential verification (LHS). The sequential condition requires that the Worst-Case Tube
\[
\mathcal{T} \;:=\; \bigcup_{y \in B_Y(f_P(x), r_P)} B_Z(f_Q(y), r_Q) \;\subseteq\; S_w.
\]
We must show $B_Z(f_Q(f_P(x)), r_{quad}) \subseteq S_w$.  The tube $\mathcal{T}$ is the $r_Q$-neighbourhood of the image set $I = f_Q(B_Y(f_P(x), r_P))$.  By the co-Lipschitz property of $f_Q$ (Definition~\ref{def:metric-openness}), $I \supseteq B_Z(f_Q(f_P(x)), \ell\, r_P)$.  For any $z$ with $d_Z(z, f_Q(f_P(x))) \le r_{quad} \le \ell\, r_P + r_Q = r_{cov}$, the triangle inequality gives $d_Z(z, I) \le \max(0, r_{quad} - \ell\, r_P) = r_Q$, so $z \in \mathcal{T}$.  Therefore $B_Z(f_Q(f_P(x)), r_{quad}) \subseteq \mathcal{T} \subseteq S_w$.
\end{proof}

\paragraph*{Guidance on Definition Selection}
The choice between Linear and Quadratic composition depends on the nature of the solver errors:

\begin{itemize}
    \item \textbf{Use Linear Composition ($L r_P + r_Q$) when:}
    The errors are systematic or correlated (e.g., both stages rely on the same biased risk model), or when strict regulatory guarantees are required against adversarial conditions. In this regime, the inclusions become equalities, and the audit is exact but costly.

    \item \textbf{Use Quadratic Composition ($\sqrt{\dots}$) when:}
    The errors are driven by distinct randomized processes (e.g., Monte Carlo noise, independent solver seeds). This regime enables \emph{Audit Safety}: verifying each step individually is a sufficient condition for global compliance, allowing the system to decouple intermediate checks without compromising end-to-end safety.
\end{itemize}

\paragraph*{Final remarks on Safety Radius approach}

The Safety Radius approach pays an upfront ``complexity tax'' (calculating the Lipschitz constant $L$) to buy computational simplicity later (checking $d(x,y) \le r$). However, in many practical situations---in particular, customization and replication---we can simply assume $L=1$, as the following remark explains. 

\begin{remark}[Interpretation of $L$ as Financial Leverage]
\label{rem:lipschitz-leverage}
In the context of the Safety Radius approach, the Lipschitz constant $L$ of a re-implementation map $f: K_{Hub} \to K_{Spoke}$ admits a direct financial interpretation as the \textbf{Sensitivity} or \textbf{Gross Leverage} of the transformation.

Recall the definition of Lipschitz continuity:
\[
\| f(x) - f(x') \| \le L \cdot \| x - x' \|
\]
Let $\delta = x - x'$ represent a perturbation or ``error'' in the hub portfolio (e.g., input noise). Let $\Delta = f(x) - f(x')$ represent the resulting deviation in the spoke portfolio. The inequality states:
\[
\frac{\| \Delta \|}{\| \delta \|} \le L
\]
Thus, $L$ bounds the ratio of output deviation to input deviation. We distinguish three regimes:

\begin{enumerate}
    \item \textbf{De-risking / Dampening ($L < 1$):} 
    Consider a map that moves assets into Cash or performs heavy diversification (shrinking the variance). The output is \textit{more stable} than the input. The safety radius contracts ($r_{spoke} < r_{hub}$), meaning the process naturally absorbs upstream errors.
    
    \item \textbf{Unitary / Pass-through ($L \approx 1$):} 
    Consider a simple re-custody or long-only replication of an index. A $1\%$ error in the input weights results in exactly a $1\%$ error in the output weights. The safety radius is preserved.
    
    \item \textbf{Leveraging / Amplification ($L > 1$):} 
    Consider a strategy that applies a leverage factor $\lambda > 1$ to the hub signals (e.g., a $130/30$ fund derived from a standard model).
    \[
    f(x) \approx \lambda x \implies \| \lambda x - \lambda x' \| = |\lambda| \| x - x' \|
    \]
    Here, $L = \lambda$. A $1\%$ error in the hub signal is magnified into a $\lambda\%$ error in the spoke. The safety radius must expand ($r_{spoke} = \lambda \cdot r_{hub}$) to maintain the same confidence level.
\end{enumerate}
\end{remark}

\subsection{Coherence Laws for the Safety Radius Approach}

\begin{theorem}[Safety Guarantee (one-way adjunction)]
\label{thm:metric-adjunction}
Let $P$ be a centered kernel with safety radius $r_\epsilon$. For any closed relations $R$ and $S$, we have the one-way implication:
\[
P_!^\epsilon R \subseteq S \implies R \subseteq P^{*,\epsilon} S
\]
\end{theorem}

\begin{proof}
Assume the inclusion $P_!^\epsilon R \subseteq S$. Let $(x, z)$ be an arbitrary element of $R$.
By the definition of the probabilistic pushforward:
\[
B(f(x), r_\epsilon) \times \{z\} \subseteq P_!^\epsilon R
\]
Given the hypothesis $P_!^\epsilon R \subseteq S$, it follows that:
\[
B(f(x), r_\epsilon) \times \{z\} \subseteq S
\]
This implies that for every $y \in B(f(x), r_\epsilon)$, the pair $(y, z)$ is in $S$. In terms of slices, this means the entire ball is contained in the permissible slice:
\[
B(f(x), r_\epsilon) \subseteq S_z
\]
The condition that a closed ball of radius $r_\epsilon$ is contained in a set is equivalent to the center being at distance at least $r_\epsilon$ from the complement of the set. Thus:
\[
d_Y(f(x), S_z^c) \ge r_\epsilon
\]
By definition of the probabilistic pullback, this implies $(x, z) \in P^{*,\epsilon} S$.

It remains to handle limit points arising from the closure in the definition of $P_!^\epsilon R$.  If $(y,z)$ is a limit point of $\bigcup_{(x,z') \in R} B(f(x), r_\epsilon) \times \{z'\}$ (necessarily with $z' = z$ eventually) rather than an element of the union itself, then there exists a sequence $(y_n, z) \to (y, z)$ with $y_n \in B(f(x_n), r_\epsilon)$ for some $(x_n, z) \in R$.  Since $P_!^\epsilon R \subseteq S$, each $y_n \in S_z$.  Because $S$ is closed, the slice $S_z$ is closed, so $y \in S_z$.  Thus the one-way implication holds for limit points as well.
\end{proof}

\begin{remark}[Why the converse fails]
Unlike the deterministic adjunction $f_! \dashv f^*$ of Theorem~\ref{thm:adjunction}, the stochastic version provides only a one-way implication.  The converse fails because the closure in the definition of $P_!^\epsilon$ can introduce limit points not witnessed by any element of $R$: a pair $(y,z)$ may belong to $P_!^\epsilon R$ (via a limiting sequence) without any hub $x$ satisfying both $(x,z) \in R$ and $y \in B(f(x), r_\epsilon)$.
\end{remark}

\begin{theorem}[Beck--Chevalley / Early-Audit Safety]
\label{thm:metric-bc}
Consider a commuting square of portfolio transformations:
\[
\begin{tikzcd}
X \arrow[r, "P"] \arrow[d, "g"'] & Z \arrow[d, "h"] \\
Y \arrow[r, "Q"'] & W
\end{tikzcd}
\]
Assume $g$ and $h$ are deterministic aggregation maps, and $P, Q$ are stochastic kernels centered at $f_P, f_Q$ with safety radii $r_\delta, r_\epsilon$.
Assume the square commutes on centers ($h \circ f_P = f_Q \circ g$) and that $h$ is $L$-Lipschitz continuous.
If the downstream risk budget satisfies $r_\epsilon \ge L \cdot r_\delta$, then:
\[
P_!^\delta (g^* R) \subseteq h^* (Q_!^\epsilon R)
\]
\end{theorem}

\begin{proof}
Let $(y, w) \in P_!^\delta (g^* R)$. By definition, there exists $x \in X$ such that:
1. $(x, w) \in g^* R \implies (g(x), w) \in R$.
2. $y \in B(f_P(x), r_\delta)$.

We must show $(y, w) \in h^* (Q_!^\epsilon R)$, which requires $(h(y), w) \in Q_!^\epsilon R$.
This, in turn, requires the existence of an element $z \in Y$ such that $(z, w) \in R$ and $h(y) \in B(f_Q(z), r_\epsilon)$.

Let us choose $z = g(x)$. From (1), we know $(z, w) \in R$.
It remains to verify the distance condition. We know $d_Y(y, f_P(x)) \le r_\delta$.
Applying the map $h$:
\[
d_W(h(y), h(f_P(x))) \le L \cdot d_Y(y, f_P(x)) \le L \cdot r_\delta
\]
Using the commutativity of centers, $h(f_P(x)) = f_Q(g(x)) = f_Q(z)$. Thus:
\[
d_W(h(y), f_Q(z)) \le L \cdot r_\delta
\]
Given the stability assumption $r_\epsilon \ge L \cdot r_\delta$, we have:
\[
d_W(h(y), f_Q(z)) \le r_\epsilon
\]
Thus $h(y) \in B(f_Q(z), r_\epsilon)$, which implies $(h(y), w) \in Q_!^\epsilon R$.
\end{proof}

\begin{theorem}[Frobenius / Conservative Filtering]
\label{thm:metric-frobenius}
Let $P$ be a centered kernel with safety radius $r_\epsilon$. For closed relations $R$ and $S$:
\[
P_!^\epsilon (R \cap P^{*,\epsilon} S) \subseteq P_!^\epsilon R \cap S
\]
\end{theorem}

\begin{proof}
Let $(y, z)$ be an element of the LHS. This implies the existence of $x$ such that:
1. $(x, z) \in R$
2. $(x, z) \in P^{*,\epsilon} S \implies B(f(x), r_\epsilon) \subseteq S_z$
3. $y \in B(f(x), r_\epsilon)$

From (1) and (3), it is immediate that $(y, z) \in P_!^\epsilon R$.
From (2) and (3), since $y$ is in the ball and the ball is contained in $S_z$, it follows that $y \in S_z$, and thus $(y, z) \in S$.
Therefore, $(y, z)$ is in the intersection.
\end{proof}

\clearpage  \section{Probabilistic Operations: Highest Density Regions}
\label{sec:HDR-section}

The Safety Radius approach replaces a distribution with a ball.  The HDR approach keeps the distribution and asks: where does most of the mass actually land?  The \textbf{$\epsilon$-Highest Density Region} of a distribution $P(x, \cdot)$ is the smallest closed set containing probability mass at least $1 - \epsilon$.  For a distribution with a continuous density, this is a superlevel set $\{y : \rho_x(y) \ge \lambda_\epsilon\}$ --- the region where the density is highest.  The resulting ``safety region'' adapts to the distribution's actual shape: elongated for anisotropic solver noise, bimodal for discrete-constraint solvers, banana-shaped for tax-constrained transitions.  The cost of this adaptability is compositional: unlike metric balls, HDRs do not compose by simple addition, and the Frobenius law holds only as an inclusion.

\subsection{Definitions}
\label{sub:HDR}

\begin{definition}[Probabilistic Pullback]
Given $P: X \rightsquigarrow Y$ and $S \subseteq Y \times Z$, the \textbf{$\epsilon$-pullback} identifies hubs that succeed with high probability:
\[
P^{*,\epsilon} S = \{ (x, z) \in X \times Z \mid P(x, S_z) \ge 1 - \epsilon \}
\]
\end{definition}

\begin{definition}[Probabilistic Pushforward]
Given a tight Feller kernel $P: X \rightsquigarrow Y$ and a relation $R \subseteq X \times Z$, the \textbf{$\epsilon$-pushforward} propagates the safety region:
\[
P_!^\epsilon R = \overline{ \bigcup_{(x, z) \in R} \left( \operatorname{supp}_\epsilon(P(x)) \times \{z\} \right) } \subseteq Y \times Z
\]
\end{definition}

\begin{remark}[Operational Meaning]
$P_!^\epsilon R$ defines the \textbf{Canonical Safety Region} for the spokes.
It answers: ``If the hub $x$ is valid, where will the resulting spoke portfolio land with high confidence?''
By tracking the HDR, we exclude tail events (outliers) that would otherwise expand the reachable set to the entire space, rendering the audit meaningless.
\end{remark}

\begin{remark}
Again, pullback and pushforward along an identity morphism are not trivial. When the kernel $P$ represents noise centered around the identity, the operations correspond to chance constraints and robust regions:

\begin{itemize}
    \item \textbf{Pullback (Chance Constraint):} 
    The $\epsilon$-pullback selects inputs that satisfy the constraint $S$ with probability at least $1-\epsilon$:
    \[
    P^{*,\epsilon} S = \{ (x, z) \in X \times Z \mid P(x, S_z) \ge 1 - \epsilon \}
    \]
    \item \textbf{Pushforward (Robust Safety Region):} 
    The $\epsilon$-pushforward is the closure of the union of all $\epsilon$-supports, representing the set of all "likely" outcomes generated by the valid inputs:
    \[
    P_!^\epsilon R = \overline{ \bigcup_{(x, z) \in R} \left( \operatorname{supp}_\epsilon(P(x)) \times \{z\} \right) }
    \]
\end{itemize}
\end{remark}

The $\epsilon$-pullback and $\epsilon$-pushforward operations compose coherently, so that multi-stage pipelines (e.g., Hub $\to$ Model $\to$ Client) can be analysed step by step without losing track of accumulated risk.

\begin{theorem}[Composition of $\epsilon$-Pullback]
\label{thm:pullback-comp}
Let $P: X \rightsquigarrow Y$ and $Q: Y \rightsquigarrow Z$ be tight Feller kernels.
For any closed relation $S \subseteq Z \times W$ and risk budgets $\delta, \epsilon \ge 0$ such that $\delta + \epsilon < 1$:
\[
P^{*,\delta} \left( Q^{*,\epsilon} S \right) \;\subseteq\; (Q \circ P)^{*,\delta+\epsilon} S
\]
\end{theorem}

\begin{proof}
Let $(x, w)$ be an arbitrary element of the nested pullback $P^{*,\delta} (Q^{*,\epsilon} S)$.
We must show that the composite kernel maps $x$ to the target set $S_w$ with probability at least $1 - (\delta + \epsilon)$.

\textbf{Step 1: Analyse the Inner Pullback.}
Define the set of ``successful'' intermediate states in $Y$ as:
\[
A := \{ y \in Y \mid (y, w) \in Q^{*,\epsilon} S \} = \{ y \in Y \mid Q(y, S_w) \ge 1 - \epsilon \}
\]
For any $y \in A$, the kernel $Q$ maps $y$ into the target slice $S_w$ with high probability.

\textbf{Step 2: Analyse the Outer Pullback.}
The condition $(x, w) \in P^{*,\delta} (Q^{*,\epsilon} S)$ implies that $x$ maps into the set $A$ with probability at least $1-\delta$:
\[
P(x, A) \ge 1 - \delta
\]

\textbf{Step 3: Evaluate the Composite Kernel.}
The probability that the composite kernel $(Q \circ P)$ maps $x$ into $S_w$ is given by the Chapman--Kolmogorov integral:
\[
(Q \circ P)(x, S_w) = \int_Y Q(y, S_w) \, P(x, dy)
\]
We split this integral over the success set $A$ and its complement $A^c$:
\[
= \int_A Q(y, S_w) \, P(x, dy) + \int_{A^c} Q(y, S_w) \, P(x, dy)
\]
Since probabilities are non-negative, we can lower-bound the second term by 0.
For the first term, we know that for all $y \in A$, $Q(y, S_w) \ge 1 - \epsilon$. Thus:
\[
\int_A Q(y, S_w) \, P(x, dy) \ge \int_A (1 - \epsilon) \, P(x, dy) = (1 - \epsilon) P(x, A)
\]
Using the bound $P(x, A) \ge 1 - \delta$, we have:
\[
(Q \circ P)(x, S_w) \ge (1 - \epsilon)(1 - \delta) = 1 - \epsilon - \delta + \epsilon\delta
\]
Since $\epsilon, \delta \ge 0$, the term $\epsilon\delta$ is non-negative. Dropping it preserves the inequality:
\[
(Q \circ P)(x, S_w) \ge 1 - (\epsilon + \delta)
\]
By the definition of the probabilistic pullback, this implies $(x, w) \in (Q \circ P)^{*,\delta+\epsilon} S$.
\end{proof}

\begin{theorem}[Composition of $\epsilon$-Pushforward]
\label{thm:pushforward-comp}
Let $P: X \rightsquigarrow Y$ and $Q: Y \rightsquigarrow Z$ be tight Feller kernels.
Assume the Highest Density Regions (HDRs) satisfy the \textbf{Convolution Stability Condition}:
\[
\operatorname{supp}_{\delta+\epsilon}((Q \circ P)(x)) \;\subseteq\; \overline{\bigcup_{y \in \operatorname{supp}_\delta(P(x))} \operatorname{supp}_\epsilon(Q(y))}
\]
Then for any closed relation $R \subseteq X \times W$:
\[
(Q \circ P)_!^{\delta + \epsilon} R \;\subseteq\; Q_!^\epsilon \left( P_!^\delta R \right)
\]
\end{theorem}

\begin{proof}
Let $(z, w) \in (Q \circ P)_!^{\delta + \epsilon} R$. By the definition of the pushforward as the closure of the union of supports, there exists a sequence $\{(z_k, w_k)\}_{k=1}^\infty$ converging to $(z, w)$, where for each $k$, there exists a hub $x_k$ such that $(x_k, w_k) \in R$ and $z_k \in \operatorname{supp}_{\delta+\epsilon}((Q \circ P)(x_k))$.

We apply the Convolution Stability Condition to each term in the sequence. For every $k$:
\[
z_k \in \operatorname{supp}_{\delta+\epsilon}((Q \circ P)(x_k)) \subseteq \overline{\bigcup_{y \in \operatorname{supp}_\delta(P(x_k))} \operatorname{supp}_\epsilon(Q(y))}.
\]
This implies that $(z_k, w_k)$ lies in the closure of the set
\[
\mathcal{S} := \bigcup_{(x, w') \in R} \bigcup_{y \in \operatorname{supp}_\delta(P(x))} \left( \operatorname{supp}_\epsilon(Q(y)) \times \{w'\} \right).
\]
We observe that the set on the RHS of the theorem statement, $Q_!^\epsilon (P_!^\delta R)$, is by definition the closure of $\mathcal{S}$. Since $(z_k, w_k) \in \overline{\mathcal{S}}$ for all $k$, and the sequence converges to $(z, w)$, it follows that $(z, w) \in \overline{\mathcal{S}} = Q_!^\epsilon (P_!^\delta R)$.
\end{proof}

\begin{remark}
It is instructive to reconcile the ``simple'' set-theoretic inclusion of the HDR approach with the ``complicated'' arithmetic of the Safety Radius approach. The complication (the Lipschitz constant $L$) is not absent in the HDR formulation; it is implicit in the geometry of the image set.
\end{remark}

\paragraph*{The Convolution Stability Condition}

We have just seen that the composition of pushforwards relies on the behaviour of the \emph{Highest Density Regions} (HDRs) under convolution. To ensure that the risk budget adds linearly (i.e., the risk of the composite process is bounded by $\delta + \epsilon$), we required the Convolution Stability Condition. This condition asserts that the ``likely'' outcomes of the end-to-end process $(Q \circ P)$ are contained within the set of outcomes generated by chaining the ``likely'' intermediate steps. It prevents the probability mass from drifting into the ``tails'' of the distribution during composition, a phenomenon that would invalidate the audit trail.

While this condition does not hold for arbitrary probability measures (particularly multi-modal ones), the following proposition identifies sufficient conditions.

\begin{proposition}[Sufficient Conditions for Convolution Stability]
\label{prop:convolution-stability}
The Convolution Stability Condition holds in each of the following cases:
\begin{enumerate}
    \item \textbf{Log-Concave Additive Noise.} If $P(x, \cdot)$ has a log-concave density, and $Q(y, \cdot)$ acts as additive noise---i.e., the output has the form $y + \xi$ where $\xi$ is independent of $y$ with log-concave density---then the composite kernel $(Q \circ P)(x, \cdot)$ is a convolution of log-concave densities, which is log-concave (by the Pr\'ekopa-Leindler inequality, see \cite{Prekopa}, \cite{Leindler}, \cite{Gardner}). For log-concave measures, the $\epsilon$-HDR is a convex superlevel set, and the Minkowski sum of convex superlevel sets contains the superlevel set of the convolution at the combined threshold. The mass bound $\Pr(\xi \in A_\delta + B_\epsilon) \ge 1 - (\delta + \epsilon)$ always holds (by a union-bound argument).  However, the stronger HDR containment $\mathrm{HDR}_{1-(\delta+\epsilon)} \subseteq A_\delta + B_\epsilon$ requires additional concentration: it is verified for Gaussian and sub-Gaussian noise (where $\Pr(\xi \in B) \ge 1 - \epsilon$ holds for balls of explicit radius), but does \emph{not} follow from log-concavity alone.  For general log-concave noise, only the mass bound is established; HDR containment should be checked for the specific noise model in applications (see Appendix~\ref{app:convolution-stability-proof} for details).

    \item \textbf{Compact Support.}
    If the kernels $P(x, \cdot)$ and $Q(y, \cdot)$ have compact support, the condition holds by the sub-additivity of the support diameter: $\operatorname{supp}(Q \circ P) \subseteq \overline{\bigcup_y \operatorname{supp}(Q(y, \cdot))}$ where the union is taken over $y \in \operatorname{supp}(P(x, \cdot))$, and restricting to $\epsilon$-HDRs only tightens the inclusion.
\end{enumerate}
\end{proposition}

\begin{proof}
See Appendix~\ref{app:convolution-stability-proof}.
\end{proof}

In production, these conditions cover the principal use cases:

\begin{enumerate}
    \item \textbf{Gaussian-Smoothed Solvers:}
    When solvers are smoothed via additive Gaussian noise (e.g., applying Gaussian input perturbation to Gurobi \cite{Gurobi} as in Example~\ref{ex:gurobi-feller}), the resulting probability densities are log-concave, so case (1) of Proposition~\ref{prop:convolution-stability} applies.

    \item \textbf{Signal Regularization:}
    In systems like Axioma \cite{Axioma}, where the optimization is formulated as continuous QP over regularized inputs, the mapping $x \mapsto P(x, \cdot)$ avoids the bifurcation and discontinuity characteristic of raw Mixed-Integer Programming. The stability of the density function $\rho_x$ ensures that the $\epsilon$-HDR evolves continuously, preventing the ``mass escape'' that would violate stability.

    \item \textbf{Bounded Errors:}
    If the error kernels $P$ and $Q$ have compact support (e.g., uniform noise within a small tracking error band), case (2) of Proposition~\ref{prop:convolution-stability} applies.
\end{enumerate}

\begin{remark}[Failure Mode]
The condition fails if the solver is \emph{strictly multi-modal} (e.g., a raw MIP solver flipping between two distinct sectors based on numerical noise). In such cases, the convolution could place the highest density of the composite distribution in a region that corresponds to the ``valley'' between the modes of the individual steps. The failure of Convolution Stability correctly flags this solver behaviour as \emph{audit-unsafe}.
\end{remark}

In the probabilistic setting, we trade strict categorical equivalence for one-directional safety guarantees.
Instead of proving that upstream and downstream checks are identical ($=$), we prove that upstream checks are \emph{conservative} ($\implies$).
This ensures that the system may reject some valid portfolios (type I error) but will never approve an invalid one (type II error), satisfying the primary requirement of a risk management system.

\begin{theorem}[The Safety Guarantee (One-Way Adjunction)]
\label{thm:safety-adjunction}
Let $P: X \rightsquigarrow Y$ be a tight Feller kernel.
The probabilistic pushforward and pullback satisfy a one-way implication:
\[
P_!^\epsilon R \subseteq S \implies R \subseteq P^{*,\epsilon} S
\]
\end{theorem}

\begin{proof}
Assume the inclusion of relations $P_!^\epsilon R \subseteq S$.
We must show that for any pair $(x,z) \in R$, the probability of compliance is sufficiently high.

Let $(x,z)$ be an arbitrary element of $R$.
By the definition of the probabilistic pushforward:
\[
P_!^\epsilon R = \overline{ \bigcup_{(x', z') \in R} \left( \operatorname{supp}_\epsilon(P(x')) \times \{z'\} \right) }
\]
Since the union is a subset of its closure, we have:
\[
\operatorname{supp}_\epsilon(P(x)) \times \{z\} \subseteq P_!^\epsilon R \subseteq S
\]
This implies that for every $y \in \operatorname{supp}_\epsilon(P(x))$, the pair $(y,z)$ is in $S$.
Consequently, the entire $\epsilon$-HDR is contained within the $z$-slice of $S$:
\[
\operatorname{supp}_\epsilon(P(x)) \subseteq S_z := \{y \in Y \mid (y,z) \in S\}
\]
  We now evaluate the probability mass. By Definition~\ref{def:epsilon-supp},
  the $\epsilon$-HDR is constructed to contain mass \emph{at least}
  $1-\epsilon$:
  \[
  P(x, \operatorname{supp}_\epsilon(P(x))) \;\ge\; 1-\epsilon
  \]
  (Equality holds when the density has no flat regions;
  see the remark following Definition~\ref{def:epsilon-supp}.)
  By the monotonicity of the probability measure, since
  $\operatorname{supp}_\epsilon(P(x)) \subseteq S_z$:
  \[
  P(x, S_z) \ge P(x, \operatorname{supp}_\epsilon(P(x))) \ge 1-\epsilon
  \]
This inequality $P(x, S_z) \ge 1-\epsilon$ is precisely the condition for $(x,z) \in P^{*,\epsilon} S$.
Thus, $R \subseteq P^{*,\epsilon} S$.
\end{proof}

\begin{remark}[Audit Implication]
This theorem guarantees that if we constrain the \textbf{Safety Region} of the spoke to lie within the permissible set $S$, we effectively guarantee that the \textbf{probability of compliance} is at least $1-\epsilon$.
Note that the reverse implication does not hold for HDRs: a set $S_z$ might contain $99\%$ of the probability mass but exclude the specific high-density peak, failing the containment condition $\operatorname{supp}_\epsilon \subseteq S_z$.
However, for risk management, the forward implication proved above is the critical one, as it ensures that geometric checks are a conservative proxy for probabilistic checks.
\end{remark}

\begin{theorem}[Early-Audit Safety (Modified Beck--Chevalley)]
\label{thm:early-audit}
Consider a commuting square of tight Feller kernels $P, Q$ and deterministic maps $g, h$. Assume the square commutes in distribution ($h_\# P(x) = Q(g(x))$). If we verify alignment upstream using a stricter risk budget $\delta \le \epsilon$, and if the solver satisfies the \textbf{pointwise stability condition}:
\[
\forall x \in X, \quad h(\operatorname{supp}_\delta(P(x))) \subseteq \operatorname{supp}_\epsilon(Q(g(x))),
\]
then:
\[
P_!^\delta (g^* R) \;\subseteq\; h^* (Q_!^\epsilon R).
\]
\end{theorem}

\begin{proof}
We proceed by element-wise verification. Let $(z, k)$ be an arbitrary element of the upstream pushforward $P_!^\delta (g^* R)$. By the definition of pushforward, there exists a hub portfolio $x \in X$ (or a sequence converging to such, in which case we apply the argument to the sequence) such that:
\begin{enumerate}
    \item The hub satisfies the pullback constraint: $(x, k) \in g^* R$.
    \item The spoke is in the safety region: $z \in \operatorname{supp}_\delta(P(x))$.
\end{enumerate}
From (1), we have $(g(x), k) \in R$. Let $y = g(x)$. Thus $(y, k) \in R$.

From (2), we apply the map $h$ to the element $z$. Since $z \in \operatorname{supp}_\delta(P(x))$, it follows that $h(z) \in h(\operatorname{supp}_\delta(P(x)))$. Recalling the fact that $g$ is deterministic, we now invoke the pointwise stability hypothesis:
\[
h(\operatorname{supp}_\delta(P(x))) \subseteq \operatorname{supp}_\epsilon(Q(g(x))) = \operatorname{supp}_\epsilon(Q(y)).
\]
Thus, $h(z) \in \operatorname{supp}_\epsilon(Q(y))$. Since $(y, k) \in R$, the pair $(h(z), k)$ is generated by a valid input $y$ and lies within the resulting $\epsilon$-support. Therefore $(h(z), k) \in Q_!^\epsilon R$. By definition of pullback, this implies $(z, k) \in h^* (Q_!^\epsilon R)$.
\end{proof}

\begin{theorem}[Conservative Filtering (Modified Frobenius)]
\label{thm:conservative-filtering}
Let $P: X \rightsquigarrow Y$ be a tight Feller kernel whose transition measures $P(x, \cdot)$ admit continuous densities.
For closed relations $R, S$, if we filter the hub to ensure certain compliance ($\delta=0$), the resulting safety region is contained in the intersection:
\[
P_!^\epsilon \left( R \cap P^{*,0} S \right) \;\subseteq\; P_!^\epsilon R \cap S
\]
\end{theorem}

\begin{proof}
Let $(y, z)$ be an element of the LHS set $P_!^\epsilon \left( R \cap P^{*,0} S \right)$.
This implies the existence of a hub $x$ such that:
\begin{enumerate}
    \item $(x, z) \in R$
    \item $(x, z) \in P^{*,0} S$
    \item $y \in \operatorname{supp}_\epsilon(P(x))$
\end{enumerate}
From (1) and (3), it is immediate that $(y, z) \in P_!^\epsilon R$ (by definition of pushforward).
It remains to show that $(y, z) \in S$.

From (2), the condition $(x, z) \in P^{*,0} S$ implies $P(x, S_z) \ge 1 - 0 = 1$.
Thus, the measure $P(x, \cdot)$ is concentrated entirely on the closed set $S_z$ (almost everywhere).
Let $\rho_x$ be the continuous density of $P(x, \cdot)$.
Since $\int_{S_z^c} \rho_x(u) du = 0$ and $S_z^c$ is open, it follows that $\rho_x(u) = 0$ for all $u \in S_z^c$.

Now consider the element $y$ from (3).
By definition of the $\epsilon$-HDR, $y \in \{ u \mid \rho_x(u) \ge \lambda_\epsilon \}$.
For any $\epsilon < 1$, the threshold $\lambda_\epsilon$ must be strictly positive. To see this, note that $\rho_x$ is continuous and $P(x, \cdot)$ is tight, so the superlevel sets $\{u \mid \rho_x(u) \ge \lambda\}$ are closed and, by tightness, contained in a compact set for $\lambda > 0$. If $\lambda_\epsilon = 0$, then $\operatorname{supp}_\epsilon(P(x)) = \operatorname{supp}(P(x))$, which has total mass $1 > 1 - \epsilon$; but the $\epsilon$-HDR is the \emph{smallest} superlevel set of mass $\ge 1 - \epsilon$, so $\lambda_\epsilon > 0$ whenever $\epsilon > 0$.
Since $\rho_x(y) \ge \lambda_\epsilon > 0$, $y$ cannot belong to the region where the density is zero.
Therefore, $y \notin S_z^c$, which implies $y \in S_z$.
Finally, $y \in S_z \iff (y, z) \in S$.
We have shown $(y, z) \in P_!^\epsilon R$ and $(y, z) \in S$, so the intersection holds.
\end{proof}

\begin{theorem}[Geometric Frobenius Reciprocity / Stress-Test Safety]
\label{thm:geometric-frobenius}
Let $P: X \rightsquigarrow Y$ be a tight Feller kernel. Define the \textbf{Geometric Pullback} (or Stress-Test Pullback) $P^{\dagger, \delta} S$ as the set of hubs where the entire $\delta$-safety region satisfies the constraint:
\[
P^{\dagger, \delta} S := \{ (x, z) \in X \times Z \mid \operatorname{supp}_\delta(P(x)) \subseteq S_z \}
\]
If the upstream verification risk budget $\delta$ is stricter than or equal to the downstream construction risk budget $\epsilon$ (i.e., $0 < \delta \le \epsilon$), then:
\[
P_!^\epsilon \left( R \cap P^{\dagger, \delta} S \right) \;\subseteq\; P_!^\epsilon R \cap S
\]
\end{theorem}

\begin{proof}
Let $(y, z)$ be an arbitrary element of the LHS. By definition of the pushforward, there exists a hub $x$ such that:
\begin{enumerate}
    \item $(x, z) \in R$
    \item $(x, z) \in P^{\dagger, \delta} S \implies \operatorname{supp}_\delta(P(x)) \subseteq S_z$
    \item $y \in \operatorname{supp}_\epsilon(P(x))$
\end{enumerate}

From (1) and (3), it follows immediately that $(y, z) \in P_!^\epsilon R$. We must now show that $(y, z) \in S$.

Recall that the $\epsilon$-HDR ($\operatorname{supp}_\epsilon$) is defined as the smallest closed set containing probability mass $1-\epsilon$. Since we assumed $\delta \le \epsilon$, we are requiring the upstream check to account for a larger proportion of the probability mass (or an equal proportion) than the downstream construction. By the nesting property of Highest Density Regions:
\[
\delta \le \epsilon \implies \operatorname{supp}_\epsilon(P(x)) \subseteq \operatorname{supp}_\delta(P(x))
\]
Combining this with condition (2):
\[
y \in \operatorname{supp}_\epsilon(P(x)) \subseteq \operatorname{supp}_\delta(P(x)) \subseteq S_z
\]
Therefore $y \in S_z$, which implies $(y, z) \in S$.
\end{proof}

\begin{remark}[Operational Interpretation]
This formulation validates the use of ``Stress Testing'' or ``Robust Optimization'' as a proxy for probabilistic safety. It replaces the integral check (``Is the probability of failure < 5\%?'') with a geometric check (``Do all outcomes within the 95\% confidence ellipsoid satisfy the rule?''). 

In practice, this justifies the use of \textbf{interval arithmetic} or \textbf{robust optimization} upstream. If we verify that the ``worst-case likely scenario'' (the boundary of the $\delta$-HDR) satisfies the constraint, we guarantee that the final portfolio is compliant, provided our verification budget $\delta$ is at least as conservative as our construction budget $\epsilon$.
\end{remark}

\begin{remark}[Handling atoms and discontinuous densities]
The proof of Theorem~\ref{thm:conservative-filtering} relies on the existence of a continuous density $\rho_x$ to argue via contradiction. In practical portfolio construction, however, probability distributions often contain \textbf{atoms} (Dirac deltas), particularly at the boundary values (e.g., a probability mass concentrated exactly at $0\%$ weight for a specific asset).

In such cases, the definition of the $\epsilon$-HDR (Definition~\ref{def:epsilon-supp}) should be generalized to the smallest closed set containing mass at least $1-\epsilon$. Under this generalization, the proof of Theorem~\ref{thm:conservative-filtering} adapts as follows: if $(x, z) \in P^{*,0}S$, then $P(x, S_z) = 1$, so the measure is entirely supported within the closed set $S_z$. Consequently, any closed set of mass $\ge 1-\epsilon$ (where $\epsilon < 1$) must be contained within $S_z$, since removing any point of $S_z^c$ costs zero mass. The safety guarantee therefore remains robust even for distributions with atoms. Theorem~\ref{thm:geometric-frobenius} (Geometric Frobenius) likewise generalizes, since its proof depends only on the nesting $\operatorname{supp}_\epsilon \subseteq \operatorname{supp}_\delta$, which holds for the generalized definition. The composition theorems (Theorems \ref{thm:pullback-comp} and \ref{thm:pushforward-comp}), however, rely on the density-based structure of the HDR---specifically, the Convolution Stability Condition---and do not automatically extend to distributions with atoms.
\end{remark}

\subsection{Compositional Safety of the Safety Radius and HDR Approaches}
\label{subsec:compositional-safety}

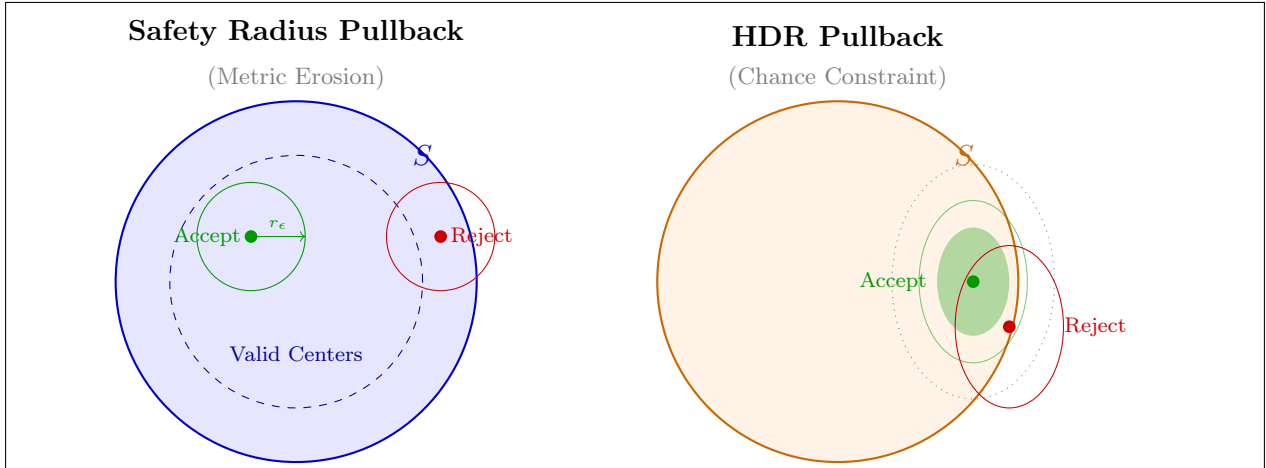
\begin{figure}
    \centering
\begin{tikzpicture}[scale=1.2]

    \def\constraintpath{(0,0) circle (2.0cm)}
    
    \begin{scope}[local bounding box=SafetyPanel]
        \node[above, font=\bfseries] at (0, 2.5) {Safety Radius Pullback};
        \node[below, font=\footnotesize, color=gray] at (0, 2.5) {(Metric Erosion)};

        \fill[blue!10] \constraintpath;
        \draw[blue!80!black, thick] \constraintpath;
        \node[blue!80!black] at (1.4, 1.4) {$S$};

        \draw[dashed, blue!60!black] (0,0) circle (1.4cm);
        \node[blue!60!black, font=\scriptsize] at (0, -0.8) {Valid Centers};

        \coordinate (Valid) at (-0.5, 0.5);
        \fill[green!60!black] (Valid) circle (2pt);
        \draw[green!60!black, thin] (Valid) circle (0.6cm); 
        \draw[->, green!60!black, ultra thin] (Valid) -- ++(0.6,0) node[midway, above=-1pt, font=\tiny] {$r_\epsilon$};
        \node[green!60!black, left, font=\scriptsize] at (Valid) {Accept};

        \coordinate (Invalid) at (1.6, 0.5);
        \fill[red!80!black] (Invalid) circle (2pt);
        \draw[red!80!black, thin] (Invalid) circle (0.6cm); 
        \node[red!80!black, right, font=\scriptsize] at (Invalid) {Reject};
    \end{scope}

    \begin{scope}[xshift=6cm, local bounding box=HDRPanel]
        \node[above, font=\bfseries] at (0, 2.5) {HDR Pullback};
        \node[below, font=\footnotesize, color=gray] at (0, 2.5) {(Chance Constraint)};

        \fill[orange!10] \constraintpath;
        \draw[orange!80!black, thick] \constraintpath;
        \node[orange!80!black] at (1.4, 1.4) {$S$};

        \coordinate (ProbValid) at (1.5, 0);
        \fill[green!60!black] (ProbValid) circle (2pt);
        
        \begin{scope}[shift={(ProbValid)}]
            \fill[green!60!black, opacity=0.3] (0,0) ellipse (0.4cm and 0.6cm);
            \draw[green!60!black, thin, opacity=0.5] (0,0) ellipse (0.6cm and 0.9cm);
            \draw[gray, dotted] (0,0) ellipse (0.9cm and 1.3cm);
        \end{scope}
        \node[green!60!black, left, font=\scriptsize] at (1.1, 0) {Accept};

        \coordinate (ProbInvalid) at (1.9, -0.5);
        \fill[red!80!black] (ProbInvalid) circle (2pt);
        \begin{scope}[shift={(ProbInvalid)}]
             \draw[red!80!black, thin] (0,0) ellipse (0.6cm and 0.9cm);
        \end{scope}
        \node[red!80!black, right, font=\scriptsize] at (2.4, -0.5) {Reject};
    \end{scope}

\end{tikzpicture}
    \caption{Comparison of probabilistic pullback in Safety Radius and HDR approaches. Left: shrinks the safety region by a uniform width. Right: a point is valid if the majority of its probability mass falls within $S$.}
    \label{fig:pullback-comparison}
\end{figure}

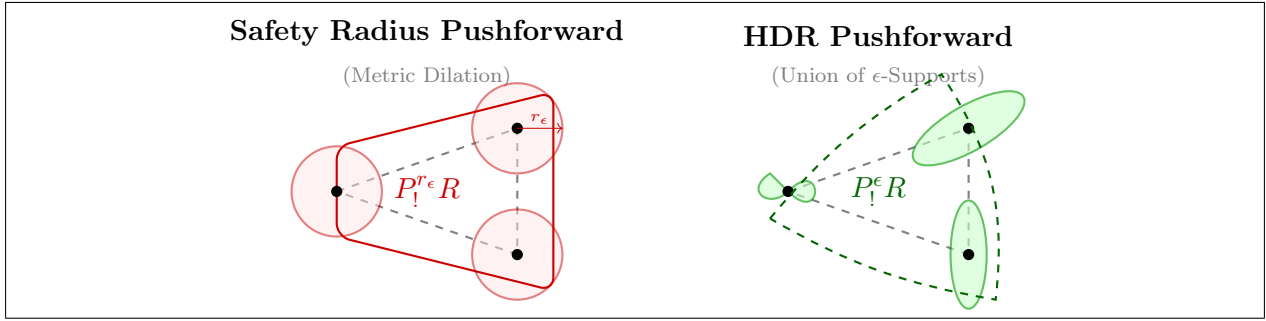
\begin{figure}
    \centering
\begin{tikzpicture}[scale=1.2]

    \tikzset{
        probcontour/.style={draw=green!60!black, fill=green!20, opacity=0.6, thick},
        safetyball/.style={draw=red!80!black, fill=red!10, opacity=0.5, thick},
        centerpt/.style={circle, fill=black, inner sep=1.5pt},
        maparrow/.style={->, gray, thick, shorten >= 2pt, shorten <= 2pt, -{Latex[length=2mm]}}
    }

    \def\setupcoords{
        \coordinate (y1) at (2.0, 2.2); 
        \coordinate (y2) at (2.0, 0.8); 
        \coordinate (y3) at (0.0, 1.5); 
    }

    \begin{scope}[local bounding box=SafetyPanel]
        \setupcoords
        
        \node[above, font=\bfseries] at (1.0, 3.0) {Safety Radius Pushforward};
        \node[below, font=\scriptsize, color=gray] at (1.0, 3.0) {(Metric Dilation)};

        \draw[thick, gray, dashed] (y3) -- (y1) -- (y2) -- cycle;
        
        \draw[safetyball] (y1) circle (0.5cm); 
        \draw[safetyball] (y2) circle (0.5cm); 
        \draw[safetyball] (y3) circle (0.5cm); 
        
        \draw[red!80!black, thick, rounded corners=5pt] 
            (0.0, 2.0) -- (2.4, 2.6) -- (2.4, 0.4) -- (0.0, 1.0) -- cycle;
            
        \node[centerpt] at (y1) {};
        \node[centerpt] at (y2) {};
        \node[centerpt] at (y3) {};
        
        \draw[->, red!80!black, ultra thin] (y1) -- ++(0.5, 0) node[midway, above=-2pt, font=\tiny] {$r_\epsilon$};
        
        \node[red!80!black, font=\bfseries] at (1.0, 1.5) {$P_!^{r_\epsilon} R$};
    \end{scope}

    \begin{scope}[xshift=5cm, local bounding box=HDRPanel]
        \setupcoords
        
        \node[above, font=\bfseries] at (1.0, 3.0) {HDR Pushforward};
        \node[below, font=\scriptsize, color=gray] at (1.0, 3.0) {(Union of $\epsilon$-Supports)};

        \draw[thick, gray, dashed] (y3) -- (y1) -- (y2) -- cycle;

        
        \begin{scope}[shift={(y1)}, rotate=30]
            \draw[probcontour] (0,0) ellipse (0.7cm and 0.25cm);
        \end{scope}
        \node[centerpt] at (y1) {};

        \begin{scope}[shift={(y2)}, rotate=90]
            \draw[probcontour] (0,0) ellipse (0.6cm and 0.2cm);
        \end{scope}
        \node[centerpt] at (y2) {};

        \begin{scope}[shift={(y3)}]
            \draw[probcontour] (0,0) .. controls (0.4,0.4) and (0.4,-0.4) .. (0,0) 
                .. controls (-0.2,-0.2) and (-0.5,0) .. (-0.2,0.2) -- cycle;
        \end{scope}
        \node[centerpt] at (y3) {};

        \draw[green!40!black, thick, dashed] 
            (-0.2, 1.2) to[bend left=10] (1.7, 2.8) 
            to[bend left=20] (2.3, 0.3) 
            to[bend left=10] (-0.2, 1.2);
            
        \node[green!40!black, font=\bfseries] at (1.0, 1.5) {$P_!^\epsilon R$};
    \end{scope}

\end{tikzpicture}
    \caption{Comparison of probabilistic pushforward in Safety Radius and HDR approaches. Left: expands the safety region by a uniform width. Right: respects the shape of the kernel's distribution at every point.}
    \label{fig:pushforward-comparison}
\end{figure}

\subsubsection*{The Asymmetry of Risk: Metric vs. Measure}

There is a structural asymmetry between the two approaches regarding compositional safety, specifically the Frobenius law.

\begin{itemize}
    \item \textbf{Safety Radius (Geometry):} The logic is symmetric.
    Verifying alignment upstream guarantees that the \textit{entire} error ball fits within the constraints.
    Consequently, the order of operations does not matter (up to the linearity of the radius), and risk budgets can be handled uniformly.

    \item \textbf{Highest Density Region (Probability):} The logic breaks down due to ``tail risk leakage.''
    In the deterministic theory (Theorem~\ref{thm:frobenius}), Frobenius Reciprocity is an \textbf{equality}: $f_!(R \cap f^*S) = f_!R \cap S$.
    However, in the probabilistic setting (Theorem~\ref{thm:conservative-filtering}), it is only an \textbf{inclusion}:
    \[
    P_!^\epsilon (R \cap P^{*,0} S) \;\subseteq\; P_!^\epsilon R \cap S
    \]
    The fact that this is a strict subset implies that the two workflows are not equivalent.
    The set on the right (filtering after construction) is larger because it accepts any portfolio where the probability of failure is merely small ($< \epsilon$).
    The set on the left (filtering before construction) is smaller because it conservatively rejects hubs that have \textit{any} non-zero failure probability.
    The difference between these two sets represents \textbf{leakage}: portfolios that pass the chance constraint despite having a 5\% tail risk.
    When composed with subsequent steps (especially those involving leverage), this ignored tail risk does not disappear; it propagates and can expand, invalidating the audit trail.
\end{itemize}

To close this gap and restore the tight audit trail required for the ``Action Menu'' $K \odot R$ (see Section~\ref{sec:DOTS}), we move from \textit{Chance Constraints} to \textit{Robust Constraints}.
By using the Geometric Frobenius law (Theorem~\ref{thm:geometric-frobenius}), we replace the probabilistic check (``Is failure unlikely?'') with a geometric one (``Is the entire $\epsilon$-HDR safe?'').
This treats the uncertainty set as a solid object, preventing the ``leakage'' of tail events into downstream optimization steps. 
For the intermediate stages of a multi-stage pipeline, \emph{Robust Optimization} is mathematically preferred over Chance Constraints.
\begin{itemize}
    \item A \textit{Chance Constraint} asks: ``Is the probability of failure less than 5\%?'' This allows probability to leak, breaking the audit trail.
    \item A \textit{Robust Constraint} asks: ``Do \emph{all} outcomes within the `95\% likely' region satisfy the rule?''
\end{itemize}
By using the latter (Theorem~\ref{thm:geometric-frobenius}), we convert the probability problem back into a geometry problem.
We define a ``shape'' (the $\epsilon$-HDR) that contains the bulk of the likely outcomes, and we treat that shape as a solid object that must not touch the constraint boundaries.
This ensures that the ``Action Menu'' $K \odot R$ preserves its integrity under composition.

\clearpage  \section{Transport-Based Safety}
\label{sec:wasserstein}

The Safety Radius and HDR approaches each simplify one side of the problem---the former assumes spherical error geometry, the latter reduces success to a mass threshold. Neither speaks the language that matters most to a portfolio manager: \emph{cost}. Practical portfolio transitions involve a \textbf{cost of compliance}.
\begin{itemize}
    \item If a solver outputs a distribution $P(x)$ that is technically non-compliant (e.g., $99\%$ mass inside, $1\%$ tail outside), the HDR approach might reject it entirely if the risk budget is $\epsilon < 0.01$.
    \item However, if that $1\%$ tail is only $1$ basis point away from the boundary, the economic cost of correcting it is negligible.
\end{itemize}

To capture this, we can define probabilistic operations using the \textbf{Wasserstein metric} (or Earth Mover's Distance) from Optimal Transport Theory (see \cite{Villani,Blanchet}). \textbf{Notational warning:} in the Wasserstein framework, the risk budget $\epsilon$ denotes \emph{transport cost} (a distance in the Wasserstein metric), not a \emph{probability threshold} as in the Safety Radius and HDR sections above. The budget now represents the maximum expected turnover required to force the portfolio into compliance.

The informal picture: imagine the solver has produced a distribution of outcomes, and a small fraction of them fall outside the constraint set. The HDR approach asks ``how much probability mass is outside?'' The Wasserstein approach asks ``how far outside is it, and what would it cost to move it back in?'' A distribution whose tail is one basis point beyond the boundary is, from a transport perspective, nearly compliant---even if the HDR approach, with a tight risk budget, would reject it. This reframing---from probability of failure to cost of correction---is what makes the Wasserstein approach natural for practitioners accustomed to thinking in terms of turnover and transaction costs.

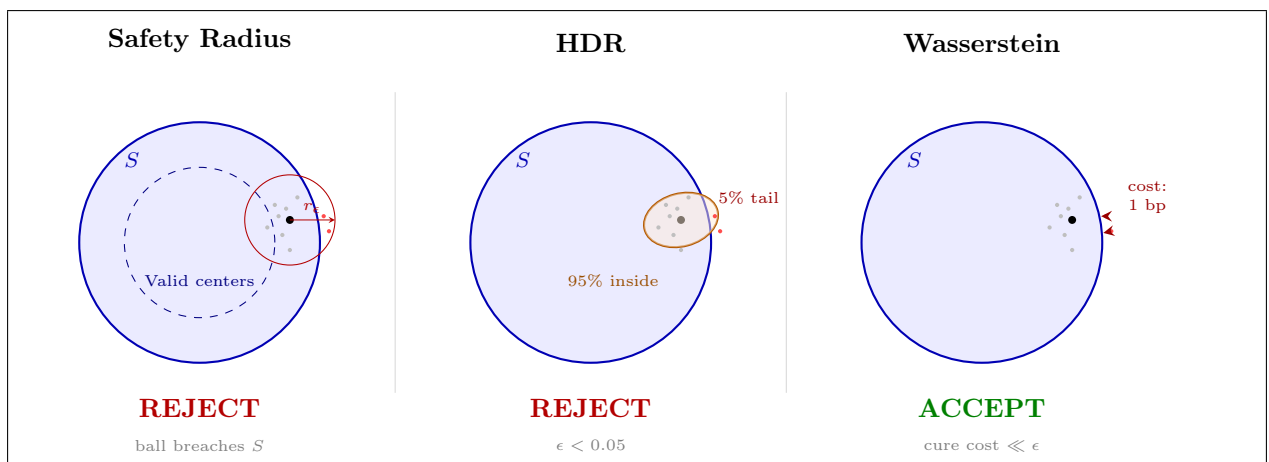
\begin{figure}[htbp]
\centering
\begin{tikzpicture}[scale=1.0, >=stealth]


    \def\constraintradius{1.6}
    \def\panelspacing{5.2}

    \begin{scope}[local bounding box=P1]
        \node[above, font=\bfseries\small] at (0, 2.4) {Safety Radius};

        \fill[blue!8] (0,0) circle (\constraintradius);
        \draw[blue!70!black, thick] (0,0) circle (\constraintradius);
        \node[blue!70!black, font=\scriptsize] at (-0.9, 1.1) {$S$};

        \draw[blue!50!black, dashed] (0,0) circle (1.0);
        \node[blue!50!black, font=\tiny] at (0, -0.5) {Valid centers};

        \coordinate (pt1) at (1.2, 0.3);
        \fill[black] (pt1) circle (1.5pt);

        \draw[red!70!black, thin] (pt1) circle (0.6);
        \draw[->, red!70!black, ultra thin] (pt1) -- ++(0.6, 0) node[midway, above=-1pt, font=\tiny] {$r_\epsilon$};

        \foreach \x/\y in {1.0/0.5, 1.1/0.1, 1.3/0.6, 0.9/0.2, 1.2/-0.1, 1.15/0.45, 1.05/0.35} {
            \fill[gray!50] (\x, \y) circle (0.8pt);
        }
        \fill[red!70] (1.65, 0.35) circle (0.8pt);
        \fill[red!70] (1.72, 0.15) circle (0.8pt);

        \node[red!70!black, font=\bfseries\small] at (0, -2.2) {REJECT};
        \node[font=\tiny, gray, align=center] at (0, -2.7) {ball breaches $S$};
    \end{scope}

    \begin{scope}[xshift=\panelspacing cm, local bounding box=P2]
        \node[above, font=\bfseries\small] at (0, 2.4) {HDR};

        \fill[blue!8] (0,0) circle (\constraintradius);
        \draw[blue!70!black, thick] (0,0) circle (\constraintradius);
        \node[blue!70!black, font=\scriptsize] at (-0.9, 1.1) {$S$};

        \coordinate (pt2) at (1.2, 0.3);
        \fill[black] (pt2) circle (1.5pt);

        \begin{scope}[shift={(pt2)}, rotate=15]
            \draw[orange!70!black, thick] (0,0) ellipse (0.5 and 0.35);
            \fill[orange!15, opacity=0.5] (0,0) ellipse (0.5 and 0.35);
        \end{scope}

        \foreach \x/\y in {1.0/0.5, 1.1/0.1, 1.3/0.6, 0.9/0.2, 1.2/-0.1, 1.15/0.45, 1.05/0.35} {
            \fill[gray!50] (\x, \y) circle (0.8pt);
        }
        \fill[red!70] (1.65, 0.35) circle (0.8pt);
        \fill[red!70] (1.72, 0.15) circle (0.8pt);

        \node[orange!60!black, font=\tiny, align=center] at (0.3, -0.5) {$95\%$ inside};

        \node[font=\tiny, red!60!black] at (2.1, 0.6) {$5\%$ tail};

        \node[red!70!black, font=\bfseries\small] at (0, -2.2) {REJECT};
        \node[font=\tiny, gray, align=center] at (0, -2.7) {$\epsilon < 0.05$};
    \end{scope}

    \begin{scope}[xshift=2*\panelspacing cm, local bounding box=P3]
        \node[above, font=\bfseries\small] at (0, 2.4) {Wasserstein};

        \fill[blue!8] (0,0) circle (\constraintradius);
        \draw[blue!70!black, thick] (0,0) circle (\constraintradius);
        \node[blue!70!black, font=\scriptsize] at (-0.9, 1.1) {$S$};

        \coordinate (pt3) at (1.2, 0.3);
        \fill[black] (pt3) circle (1.5pt);

        \foreach \x/\y in {1.0/0.5, 1.1/0.1, 1.3/0.6, 0.9/0.2, 1.2/-0.1, 1.15/0.45, 1.05/0.35} {
            \fill[gray!50] (\x, \y) circle (0.8pt);
        }

        \fill[red!70] (1.65, 0.35) circle (0.8pt);
        \fill[red!70] (1.72, 0.15) circle (0.8pt);

        \draw[->, red!60!black, thick, shorten >=1pt] (1.65, 0.35) -- (1.55, 0.35);
        \draw[->, red!60!black, thick, shorten >=1pt] (1.72, 0.15) -- (1.58, 0.12);

        \node[red!60!black, font=\tiny, align=center] at (2.2, 0.6) {cost:\\1 bp};

        \node[green!50!black, font=\bfseries\small] at (0, -2.2) {ACCEPT};
        \node[font=\tiny, gray, align=center] at (0, -2.7) {cure cost $\ll \epsilon$};
    \end{scope}

    \draw[gray!30] (0.5*\panelspacing, -2.0) -- (0.5*\panelspacing, 2.0);
    \draw[gray!30] (1.5*\panelspacing, -2.0) -- (1.5*\panelspacing, 2.0);

\end{tikzpicture}
\caption{The same non-compliant distribution evaluated by three approaches. The solver output has $95\%$ of its mass inside the constraint set $S$ and a $5\%$ tail barely outside the boundary. The Safety Radius rejects it (the metric ball breaches $S$). The HDR rejects it (the tail exceeds the risk budget $\epsilon < 0.05$). The Wasserstein approach accepts it: the cost of ``curing'' the two outliers---moving them back inside---is negligible.}
\label{fig:three-way-comparison}
\end{figure}

  \begin{definition}[Wasserstein Safety]
  \label{def:wasserstein-safety}
  Let $(Y, d)$ be a Polish space. Let
  $\mathcal{P}_1(Y) = \{\mu \in \mathcal{P}(Y) :
  \int_Y d(y, y_0)\, d\mu(y) < \infty \text{ for some (hence all) }
  y_0 \in Y\}$
  denote the space of probability measures with finite first moment.
  Let $\mathcal{P}_S = \{\nu \in \mathcal{P}_1(Y) :
  \operatorname{supp}(\nu) \subseteq S\}$ for a closed constraint set
  $S \subseteq Y$.
  For a kernel $P(x, \cdot) \in \mathcal{P}_1(Y)$ and a budget
  $\epsilon \ge 0$, the hub $x$ is \textbf{Wasserstein-compliant} if:
  \[
  W_1(P(x, \cdot),\, \mathcal{P}_S) :=
  \inf_{\nu \in \mathcal{P}_S} W_1(P(x, \cdot), \nu) \le \epsilon
  \]
  where $W_1$ is the $1$-Wasserstein distance.
  \end{definition}

Note that for tight Feller kernels on compact spaces (our setting), the finite-moment condition is automatic.

\noindent \textbf{Intuition:} This measures the ``distance to the nearest valid distribution.'' It answers: \textit{``How much trading (turnover) would be required, on average, to move the noisy solver output $P(x)$ onto the valid manifold $S$?''}

\subsection{Probabilistic Operations}

Unlike the previous approaches, the Wasserstein framework lifts the vertical category from relations on spaces ($R \subseteq X \times Y$) to relations on \emph{probability measures}. However, we can project this back to a ``Risk Profile'' on the base spaces:

\begin{definition}[Wasserstein Pullback]
\label{def:wasserstein-pullback}
Given a kernel $P: X \rightsquigarrow Y$ and a constraint $S \subseteq Y \times Z$, the \textbf{$\epsilon$-Wasserstein Pullback} is:
\[
P^{*, \epsilon}_W S = \{ (x, z) \in X \times Z \mid W_1(P(x), \mathcal{P}_{S_z}) \le \epsilon \}
\]
\end{definition}

\begin{definition}[Wasserstein Pushforward]
\label{def:wasserstein-pushforward}
Let $P: X \rightsquigarrow Y$ be a tight Feller kernel and $R \subseteq X \times Z$ be a closed relation. 
The \textbf{$\epsilon$-Wasserstein Pushforward} is defined as a subset of the product of the probability space $\mathcal{P}(Y)$ and $Z$:
\[
P_{!W}^{\epsilon} R = \{ (\nu, z) \in \mathcal{P}(Y) \times Z \mid \exists x \in X \text{ s.t. } (x, z) \in R \text{ and } W_1(P(x), \nu) \le \epsilon \}
\]
\end{definition}

Unlike the geometric approaches, the Wasserstein pushforward returns a set of \emph{distributions}. An element $(\nu, z) \in P_{!W}^{\epsilon} R$ represents a realized portfolio distribution $\nu$ that is a valid implementation of the hub requirement $z$, within a transport cost budget of $\epsilon$.

\begin{remark}[Serializing Pushforward]
    In a database implementation (see Appendix~\ref{app:architecture}), it is not feasible to store an infinite set of probability distributions. Instead, one would store a representative sample of distributions, each described by a finite set of parameters (e.g. mean vector and covariance matrix).
\end{remark}

\paragraph*{Why the Wasserstein Pushforward Must Return Distributions}

The structural divergence between the geometric approaches (Safety Radius, HDR) and the transport approach (Wasserstein) reflects the difference between \emph{Robust Safety} and \emph{Economic Safety}.

\begin{itemize}
    \item \textbf{Geometric Approaches (Robust Safety):} 
    The Safety Radius and HDR definitions (Section~\ref{sub:stoch_ops}) answer the question: \textit{``Is it possible for the solver to output a violation?''} 
    This is a binary worst-case check. Consequently, the pushforward projects the stochastic kernel $P(x)$ down to its \textbf{support} (a subset of $Y$):
    \[
    P_!^\epsilon R \subseteq Y \times Z
    \]
    We discard density information because a single point failure invalidates the set.

    \item \textbf{Wasserstein Approach (Economic Safety):}
    The Wasserstein definition (this section) answers the question: \textit{``What is the expected cost to fix the violation?''}
    This relies on the Kantorovich--Rubinstein integral:
    \[
    W_1(\mu, \nu) = \inf_{\pi} \int_{Y \times Y} d(y, y') \, d\pi(y, y')
    \]
    This calculation requires the full probability measure $\mu$. Two kernels with identical geometric supports can have vastly different transport costs. 
    
    Therefore, to preserve the data required for calculating ``Cure Costs'' across compositions, the pushforward must map into the space of measures $\mathcal{P}(Y)$:
    \[
    P_{!W}^\epsilon R \subseteq \mathcal{P}(Y) \times Z
    \]
    Collapsing this to a set of points would destroy the mass information, rendering the triangle inequality for transport costs (used in Theorem~\ref{thm:wasserstein-pushforward-comp}) impossible to apply.
\end{itemize}

\subsection{Coherence and Composition}

The Wasserstein approach has clean compositional properties, all following from the triangle inequality:

\begin{theorem}[Wasserstein Pullback Composition]
\label{thm:wasserstein-pullback-comp}
Let $P: X \rightsquigarrow Y$ and $Q: Y \rightsquigarrow Z$ be tight Feller kernels. Assume the mapping $y \mapsto Q(y)$ is $L$-Lipschitz with respect to the $W_1$ metric (i.e., $W_1(Q(y), Q(y')) \le L d_Y(y, y')$).
For any risk budgets $\delta, \epsilon \ge 0$ and any closed relation $S \subseteq Z \times M$ (where $M$ is a compact Hausdorff space; we write $M$ to avoid collision with the Wasserstein subscript $W$), we have \emph{lax functoriality} in the sense that the probabilistic pullbacks satisfy the inclusion:
\[
P^{*, \delta}_W \left( Q^{*, \epsilon}_W S \right) \;\subseteq\; (Q \circ P)^{*, L\delta + \epsilon}_W S
\]
\end{theorem}

\begin{proof}
Let $(x, m) \in P^{*, \delta}_W ( Q^{*, \epsilon}_W S )$.
By definition, $W_1(P(x), \mathcal{P}_{A}) \le \delta$, where $A = \{ y \in Y \mid (y, m) \in Q^{*, \epsilon}_W S \}$.
The set $A$ is the set of intermediate portfolios $y$ such that $W_1(Q(y), \mathcal{P}_{S_m}) \le \epsilon$.
Since $W_1(P(x), \mathcal{P}_{A}) \le \delta$, there exists a probability measure $\mu$ supported on the closure of $A$ such that $W_1(P(x), \mu) \le \delta$.

We construct a compliant distribution for the composite kernel.
For every $y \in \supp(\mu)$, since $y \in \bar{A}$, we have $W_1(Q(y), \mathcal{P}_{S_m}) \le \epsilon$.  (The set $A$ is closed: the map $y \mapsto W_1(Q(y), \mathcal{P}_{S_m})$ is continuous, since $Q$ is Feller and $\nu \mapsto W_1(\nu, \mathcal{P}_{S_m})$ is $1$-Lipschitz on $\mathcal{P}(Y)$, so $A$ is the preimage of $[0, \epsilon]$ under a continuous function.  Hence $\bar{A} = A$.)  This implies the existence of a measure $\nu_y$ supported on $S_m$ such that $W_1(Q(y), \nu_y) \le \epsilon$.
Define the composite target measure $\Gamma = \int_Y \nu_y \, \mu(dy)$, where the measurable selection $y \mapsto \nu_y$ exists by the Kuratowski--Ryll-Nardzewski theorem applied to the closed-valued correspondence $\Phi: y \mapsto \{\nu \in \mathcal{P}_{S_m} : W_1(Q(y), \nu) \le \epsilon\}$.  (This correspondence is closed-valued in the weak topology on $\mathcal{P}(Y)$ because $W_1$ is lower semi-continuous, and measurable because $Q$ is Feller and $(y, \nu) \mapsto W_1(Q(y), \nu)$ is Borel.)  Since each $\nu_y$ is supported on $S_m$, $\Gamma$ is supported on $S_m$.

We now bound the distance from the actual composite output $(Q \circ P)(x)$ to $\Gamma$:
\begin{align*}
W_1((Q \circ P)(x), \Gamma) &= W_1\left( \int Q(y) P(x, dy), \int \nu_y \mu(dy) \right) \\
&\le W_1\left( \int Q(y) P(x, dy), \int Q(y) \mu(dy) \right) + W_1\left( \int Q(y) \mu(dy), \int \nu_y \mu(dy) \right)
\end{align*}
Using the $L$-Lipschitz property of $Q$ and the joint convexity of $W_1$ (i.e., $W_1(\int \mu_y\, d\alpha(y),\, \int \nu_y\, d\alpha(y)) \le \int W_1(\mu_y, \nu_y)\, d\alpha(y)$):
1. The first term is bounded by $L \cdot W_1(P(x), \mu) \le L \delta$.
2. The second term is bounded by $\int W_1(Q(y), \nu_y) \mu(dy) \le \int \epsilon \, \mu(dy) = \epsilon$.

Thus, $W_1((Q \circ P)(x), \Gamma) \le L\delta + \epsilon$. Since $\Gamma \in \mathcal{P}_{S_m}$, this implies $(x, m) \in (Q \circ P)^{*, L\delta + \epsilon}_W S$.
\end{proof}

\begin{theorem}[Wasserstein Pushforward Composition]
\label{thm:wasserstein-pushforward-comp}
Let $P: X \rightsquigarrow Y$ and $Q: Y \rightsquigarrow Z$ be tight Feller kernels.
Assume the kernel $Q$ induces an $L$-Lipschitz mapping on the space of probability measures $\mathcal{P}(Y)$ with respect to the Wasserstein metric. Denoting by $\mu Q := \int Q(y, \cdot)\, d\mu(y)$ the pushforward of a measure $\mu$ through the kernel $Q$ (the Chapman--Kolmogorov integral), the Lipschitz condition reads $W_1(\mu Q, \nu Q) \le L W_1(\mu, \nu)$.
For any risk budgets $\delta, \epsilon \ge 0$ and any relation $R \subseteq X \times M$ (closedness of $R$ is not used in this proof; it is retained as a hypothesis for consistency with the framework), the sequential application of risk budgets satisfies the following inclusion within the linear composite budget:
\[
Q_{!W}^\epsilon \left( P_{!W}^\delta R \right) \;\subseteq\;
(Q \circ P)_{!W}^{L\delta + \epsilon} R
\]
where $Q_{!W}^\epsilon$ acts on the set of measures $\mathcal{M} \subseteq \mathcal{P}(Y)$ via the mapping $\mu \mapsto \mu Q$.
\end{theorem}

\begin{proof}
Let $(\gamma, m)$ be an element of the sequential pushforward $Q_{!W}^\epsilon ( P_{!W}^\delta R )$.
By the definition of the pushforward on sets of measures, there exists an intermediate distribution $\nu \in \mathcal{P}(Y)$ such that:
\begin{enumerate}
    \item $(\nu, m) \in P_{!W}^\delta R$.
    \item $W_1(\nu Q, \gamma) \le \epsilon$.
\end{enumerate}
From (1), by the definition of $P_{!W}^\delta R$, there exists a hub $x \in X$ such that:
\begin{enumerate}
    \item[(a)] $(x, m) \in R$.
    \item[(b)] $W_1(P(x), \nu) \le \delta$.
\end{enumerate}
We must show that $(\gamma, m)$ belongs to the composite pushforward $(Q \circ P)_{!W}^{L\delta + \epsilon} R$.
This requires showing that the distance between the composite kernel output $(Q \circ P)(x) = P(x)Q$ and the distribution $\gamma$ is bounded by $L\delta + \epsilon$.

We apply the Triangle Inequality to the Wasserstein metric:
\[
W_1(P(x)Q, \gamma) \le W_1(P(x)Q, \nu Q) + W_1(\nu Q, \gamma)
\]
Using the assumption that $Q$ is $L$-Lipschitz on measures (i.e., $W_1(\mu Q, \nu Q) \le L W_1(\mu, \nu)$), we bound the first term using (b):
\[
W_1(P(x)Q, \nu Q) \le L \cdot W_1(P(x), \nu) \le L\delta
\]
Using (2), we bound the second term:
\[
W_1(\nu Q, \gamma) \le \epsilon
\]
Combining these results:
\[
W_1((Q \circ P)(x), \gamma) \le L\delta + \epsilon
\]
Since $(x, m) \in R$, this confirms that $(\gamma, m) \in (Q \circ P)_{!W}^{L\delta + \epsilon} R$.
\end{proof}

\subsection{Adjunction, Lax Beck--Chevalley and Frobenius}
\label{subsec:wasserstein-properties}

\begin{theorem}[Wasserstein Safety Adjunction]
\label{thm:wasserstein-adjunction}
Let $P: X \rightsquigarrow Y$ be a tight Feller kernel representing a stochastic re-implementation.
For a closed set $S \subseteq Y$, let $\mathcal{Cure}_\delta(S) \subseteq \mathcal{P}(Y)$ be the set of distributions that can be transported to $S$ with expected turnover cost at most $\delta$:
\[
\mathcal{Cure}_\delta(S) := \{ \nu \in \mathcal{P}(Y) \mid W_1(\nu, \mathcal{P}_{S}) \le \delta \}
\]
When $S$ is a relation $S \subseteq Y \times Z$, we write $\mathcal{Cure}_\delta(S_z)$ for the cure set applied to the slice $S_z = \{y : (y,z) \in S\}$.
For any risk budgets $\epsilon, \delta \ge 0$, and any closed relations $R$ and $S$, if the $\epsilon$-pushforward of $R$ is contained within the slice-wise $\delta$-cure set of $S$, then $R$ satisfies the $\delta$-pullback of $S$:
\[
P_{!W}^{\epsilon} R \;\subseteq\; \{(\nu, z) \in \mathcal{P}(Y) \times Z : \nu \in \mathcal{Cure}_\delta(S_z)\} \implies R \;\subseteq\; P^{*, \delta}_W S
\]
\end{theorem}

\begin{proof}
Assume the slice-wise inclusion: for every $(\nu, z) \in P_{!W}^{\epsilon} R$, we have $\nu \in \mathcal{Cure}_\delta(S_z)$.
Let $(x, z)$ be an arbitrary element of $R$. We must show that $(x, z) \in P^{*, \delta}_W S$.

Consider the specific distribution output by the solver at $x$, denoted $\mu = P(x)$.
The Wasserstein distance is reflexive, so $W_1(P(x), \mu) = 0$. Since $0 \le \epsilon$, the pair $(\mu, z)$ satisfies the condition for membership in the Wasserstein pushforward $P_{!W}^{\epsilon} R$ (Definition~\ref{def:wasserstein-pushforward}).

By our inclusion hypothesis, $\mu \in \mathcal{Cure}_\delta(S_z)$, i.e.,
\[
W_1(\mu, \mathcal{P}_{S_z}) \le \delta
\]
Substituting $\mu = P(x)$, we obtain:
\[
W_1(P(x), \mathcal{P}_{S_z}) \le \delta
\]
This is precisely the definition of the Wasserstein pullback $P^{*, \delta}_W S$ (Definition~\ref{def:wasserstein-safety}). Thus, $(x, z) \in P^{*, \delta}_W S$.
\end{proof}

\begin{remark}[Interpretation: Economic Justification for Pre-Trade Compliance]
This theorem states that if we restrict our hub portfolios to those where the ``raw'' solver output is within distance $\delta$ of the rules ($P^{*, \delta}_W S$), we guarantee that the final realized portfolio can be fixed (``cured'') with a trading cost of at most $\delta$. This converts an abstract probability bound into a concrete budget for the trading desk.

We call this a ``safety adjunction'' by analogy with the deterministic case (Theorem~\ref{thm:adjunction}), where the pushforward-pullback relationship is a genuine Galois connection. Here the implication is one-directional: the converse fails because the hypothesis involves a \emph{uniform} bound over all distributions $\nu$ in the $\epsilon$-pushforward, not just the kernel output $P(x)$ itself. The forward direction suffices for pre-trade compliance: verifying the sufficient condition guarantees downstream curability.

Note that the parameter $\epsilon$ in the pushforward is vacuous in this implication (the proof uses only $\epsilon \ge 0$, which is always satisfied). The substantive content is the $\delta$-cure inclusion: it is the cure budget $\delta$ that controls the strength of the guarantee.
\end{remark}

\begin{theorem}[Wasserstein Lax Beck--Chevalley]
\label{thm:wasserstein-bc}
Consider a commuting square of kernels and maps where the aggregation map $h: Y \to W$ is $L$-Lipschitz.
Assume the square commutes up to a tolerance $\gamma$ in the Wasserstein metric:
\[
W_1(Q(g(x)), h_\# P(x)) \le \gamma \quad \text{for all } x.
\]
Then for any upstream risk budget $\delta \ge 0$, the pushforward of the upstream verification is contained in the downstream construction:
\[
h_{!W} \left( P_{!W}^{\delta} (g^* R) \right) \;\subseteq\; Q_{!W}^{L\delta + \gamma} R
\]
where $h_{!W}$ denotes the set-theoretic image under the measure pushforward $h_\#$: $h_{!W}(\mathcal{M}) = \{h_\# \mu : \mu \in \mathcal{M}\}$.
\end{theorem}

\begin{proof}
Let $(\nu, w)$ be an element of the set $P_{!W}^{\delta} (g^* R)$.
By definition of the Wasserstein pushforward $P_{!W}^\delta$ and pullback $g^*$, there exists a hub $x \in X$ such that:
\begin{enumerate}
    \item $(x, w) \in g^* R \implies (g(x), w) \in R$.
    \item $W_1(P(x), \nu) \le \delta$.
\end{enumerate}
We must show that the pushed-forward measure $\rho = h_\# \nu$ lies in $Q_{!W}^{L\delta + \gamma} R$.
Let $z = g(x)$. From (1), we have $(z, w) \in R$.
We now bound the Wasserstein distance from the target solver output $Q(z)$ to $\rho$:
\begin{align*}
W_1(Q(z), \rho) &= W_1(Q(g(x)), h_\# \nu) \\
&\le W_1(Q(g(x)), h_\# P(x)) + W_1(h_\# P(x), h_\# \nu) \quad \text{(Triangle Inequality)} \\
&\le \gamma + L \cdot W_1(P(x), \nu) \quad \text{(Commutativity + Lipschitz Property of } h) \\
&\le \gamma + L\delta
\end{align*}
Thus, $(\rho, w)$ is a valid implementation of $z$ within the cure budget $L\delta + \gamma$.
Consequently, $(\rho, w) \in Q_{!W}^{L\delta + \gamma} R$.
\end{proof}

\begin{remark}[Interpretation: Early Audit Safety]
This theorem governs the validity of auditing a strategy upstream. It tells us that verifying a strategy on the Hub is valid, but we must adjust our tolerance downstream. If the aggregation logic $h$ has leverage $L > 1$ (e.g., a 130/30 strategy), any ``slop'' or error $\delta$ allowed at the hub level is magnified. The term $\gamma$ represents the \textbf{irreducible noise} of the downstream solver. To ensure the downstream check doesn't falsely flag valid portfolios, its tolerance must be widened to $L\delta + \gamma$.
\end{remark}

\begin{theorem}[Wasserstein Frobenius Reciprocity]
\label{thm:wasserstein-frobenius}
Let $P: X \rightsquigarrow Y$ be a tight Feller kernel. For any closed relations $R, S$ and risk budget $\epsilon \ge 0$:
\[
P_{!W}^{\epsilon} (R \cap P^{*,0}_W S) \;\subseteq\; P_{!W}^{\epsilon} R \;\cap\; \{(\nu, z) \in \mathcal{P}(Y) \times Z : \nu \in \mathcal{Cure}_\epsilon(S_z)\}
\]
\end{theorem}

\begin{proof}
Let $(\nu, z)$ be an element of the LHS. There exists a hub $x$ satisfying:
\begin{enumerate}
    \item $(x, z) \in R$
    \item $(x, z) \in P^{*,0}_W S \implies P(x) \text{ is supported on } S_z$.
    \item $W_1(P(x), \nu) \le \epsilon$.
\end{enumerate}
From (1) and (3), $(\nu, z) \in P_{!W}^{\epsilon} R$.
From (2), since $P(x)$ is entirely supported on $S_z$, it belongs to the set of valid distributions $\mathcal{P}_{S_z}$.
The distance from $\nu$ to the valid set is bounded by the distance to $P(x)$:
\[
W_1(\nu, \mathcal{P}_{S_z}) \le W_1(\nu, P(x)) \le \epsilon
\]
Thus $\nu \in \mathcal{Cure}_\epsilon(S_z)$, proving the inclusion.
\end{proof}

The asymmetry---$\epsilon = 0$ on the pullback side, general $\epsilon$ on the pushforward---reflects the operational sequence: the pullback condition $P^{*,0}_W S$ requires exact compliance (the input portfolio must already satisfy the constraint with certainty), while the pushforward budget $\epsilon$ tolerates controlled deviation in the output.

\begin{remark}
The Wasserstein framework's main theorems rest on two assumptions that deserve scrutiny: the Lipschitz condition on kernels and the computability of the $W_1$ distance itself. The Lipschitz condition ($W_1(Q(y), Q(y')) \le L \cdot d(y, y')$) is what makes risk budgets compose linearly rather than combinatorially---without it, the $L\delta + \epsilon$ bound degrades. But verifying this condition for a given solver requires understanding the solver's sensitivity to perturbations, which is precisely the kind of information that commercial optimizers do not typically expose. The $W_1$ computation, meanwhile, is a linear program whose dimension scales with the support of the distributions involved; for continuous distributions on high-dimensional simplices, this is approximated rather than solved. The framework is correct; the question is whether its inputs can be obtained in practice. We suspect they often can be \emph{bounded} (upper bounds on $L$ via perturbation analysis, approximations of $W_1$ via sampling), and that these bounds suffice for the compositional guarantees to remain useful, if conservative.
\end{remark}

\begin{remark}[Interpretation: Robust Filtering]
This theorem validates the use of ``Compliance with Cure'' workflows: strictly filtering inputs allows us to accept outputs that are merely ``close enough.'' The set on the left represents ``Filtering before optimization'' (selecting hubs $x$ that are \emph{perfectly} safe). The set on the right represents the outcome: a portfolio distribution $\nu$ that might not be perfect (due to solver noise $\epsilon$), but is guaranteed to be \textbf{within distance $\epsilon$} of compliance. This allows the system to accept slight deviations in the solver output, knowing they are economically negligible and easily corrected.
\end{remark}

\subsection{Computing the Cure in Practice}
\label{sub:wasserstein-compute}

The Wasserstein formulation has a computational advantage over the density-based HDR approach. While computing the volume of a high-dimensional density slice is difficult ($O(e^d)$), computing the transport cost to a convex set reduces to standard optimization.

We exploit the following simplification (see \cite{Villani}, Chapter~6). If the transport cost is defined by a metric $d$ (e.g., $L_1$ turnover) and $S$ is a \emph{closed} set, the 1-Wasserstein distance from a distribution $\mu$ to the set of distributions supported on $S$ simplifies to the \textbf{expected distance to the set}:
\[
W_1(\mu, \mathcal{P}_S) = \mathbb{E}_{y \sim \mu} [ \inf_{z \in S} d(y, z) ]
\]
(The optimal coupling projects each point to its nearest point in $S$; intuitively, mass already inside $S$ costs nothing to move.  Existence of a nearest point is guaranteed by closedness, and a measurable selection exists by standard results.  When $S$ is additionally \emph{convex}, the projection is unique and non-expansive, which simplifies computation.)
This allows us to estimate the ``Cure'' (the expected turnover required to fix the portfolio) using a \textbf{Monte Carlo Projection} workflow rather than solving a full Optimal Transport problem.

\paragraph{Algorithm: The Cure Budget Estimator}
To verify if a hub portfolio $x$ satisfies the Wasserstein pullback $P^{*, \epsilon}_W S$ (i.e., ``Is the expected cure cost $\le \epsilon$?''):

\begin{enumerate}
    \item \textbf{Sampling:} Draw $N$ samples $\{y_1, \dots, y_N\}$ from the solver kernel $P(x)$. These represent potential raw spoke portfolios (before compliance filtering).
    \item \textbf{Projection (The Cure):} For each sample $y_i$, compute the minimum turnover required to reach the constraint set $S$. This is the projection problem:
    \[
    c_i = \min_{z \in \mathbb{R}^m} \| z - y_i \|_1 \quad \text{subject to } z \in S
    \]
    If $S$ is polyhedral (defined by linear constraints $Az \le b$), this is a standard \textbf{Linear Program}.  When $S$ is the intersection of multiple convex constraint sets, Dykstra's alternating projection algorithm \cite{BoyleDykstra} provides an efficient iterative approach.
    \item \textbf{Aggregation:} Compute the Monte Carlo estimate of the risk:
    \[
    \hat{W}_1 \approx \frac{1}{N} \sum_{i=1}^N c_i
    \]
    \item \textbf{Check:} If $\hat{W}_1 \le \epsilon$, the hub $x$ is approved.
\end{enumerate}

\begin{remark}[Operationalizing the Cure]
If the check fails ($\hat{W}_1 > \epsilon$), the projection vectors $z_i - y_i$ provide the \textbf{marginal sensitivities}. The gradient of the cure cost with respect to the hub positions $\nabla_x \hat{W}_1$ can be fed back into the upstream construction process (via the Chain Rule and the envelope theorem), allowing the hub manager to ``steer'' away from brittle constraints that induce high downstream turnover.
\end{remark}

\subsubsection*{The Choice of \texorpdfstring{$p=2$}{p=2} for Execution (Market Impact)}

Setting $p=1$ in the Wasserstein metric focuses on the portfolio turnover required to ``cure'' violations. For portfolios of modest size, such as usually arise in customization settings, it is reasonable to assume that trading costs scale linearly with turnover. This is no longer true for large portfolios, since large trades have disproportionate market impact. This can become relevant when a replicated ``spoke'' becomes very large relative to its ``hub''.

While any $p > 1$ provides the necessary convexity to penalize large trades, setting $p=2$ specifically is preferred for pre-trade execution analysis. This choice aligns the mathematical properties of the transport cost with standard financial theory and computational infrastructure in three distinct ways:

\paragraph{1. The Economic Model: Linear Price Impact}
In institutional finance, standard execution cost models (e.g. Almgren--Chriss \cite{AlmgrenChriss}, motivated by the linear price impact of Kyle \cite{Kyle}) assume that the price moves linearly with the size of the trade. If $\Delta P \approx \lambda \cdot \text{Size}$, then the total cost of execution scales quadratically:
\[
\text{Total Cost} \approx \text{Size} \times \Delta P \approx \lambda \cdot (\text{Size})^2
\]
Setting $p=2$ in the Wasserstein metric ($W_2$) implies a cost function proportional to $d(x,y)^2$, which matches this quadratic cost structure (cf. \cite{AlmgrenChriss}). Other values of $p$ would imply non-standard impact models (e.g., $p=1.5$ implies a square-root impact law) that are less commonly used in general asset classes.

\paragraph{2. Computational Architecture: Quadratic Programming (QP)}
Modern portfolio construction is built almost entirely on Mean-Variance Optimization, which is mathematically a Quadratic Program (minimizing $x^T \Sigma x$); cf. \cite{GrinoldKahn}.
\begin{itemize}
    \item \textbf{If $p=2$:} The transport cost is a quadratic term ($d^2$). Incorporating this into a standard risk model preserves the QP structure of the optimization problem. This allows the use of fast, stable, and standard industrial solvers (e.g., Gurobi, Mosek).
    \item \textbf{If $p \neq 2$:} The problem leaves the QP class. For $p = 1$, it reduces to a linear program; for other values of $p$, it becomes a general convex optimization problem requiring specialized conic solvers. In either case, the cost term no longer matches the quadratic structure of the risk model, breaking compatibility with standard mean-variance optimizer inputs. See \cite{BoydVandenberghe}.
\end{itemize}

\paragraph{3. Risk Compatibility (Variance)}
Financial risk is typically measured as Variance (the square of the standard deviation), which is an $L_2$ geometry. By setting $p=2$, we measure ``Execution Risk'' (Liquidity) in the same geometric units as ``Market Risk'' (Volatility). This allows for a unified utility function where risks can be summed directly:
\[
\text{Total Penalty} = \text{Risk Penalty} (\sigma^2) + \text{Trading Penalty} (\text{Impact}^2)
\]
Using $p \neq 2$ would result in mismatched units of sensitivity, complicating the trade-off between tracking error and execution cost.

\clearpage

\section{Estimation and Storage of Stochastic Kernels}

In industrial practice, the re-implementation map $f$ is rarely a closed-form function; it is typically defined by a numerical solver (e.g., Gurobi, Axioma, Mosek) acting on a complex constraint set. To operationalize the probabilistic framework (HSP-r), we need concrete procedures for estimating, storing, and validating the kernels $P(x)$.

This section details the methodologies for estimating the kernel $P(x)$ using standard industrial optimization tools, efficient schemas for serializing these kernels in the Evidence Ledger, and critical failure modes where Gaussian approximations are insufficient.

\subsection{Estimating the Kernel \texorpdfstring{$P(x)$}{P(x)}}

The primary computational challenge is to approximate the probability measure $\mu = P(x)$ that represents the solver's output distribution under noise. We distinguish three approaches based on the available computational budget and the nature of the solver access.

\subsubsection{Randomized Smoothing (The Monte Carlo Approach)}
\label{sub:randomized-smoothing}
This method is the direct numerical implementation of the Feller kernel definition. It treats the deterministic solver as a black box and induces a distribution via input perturbation.
\begin{itemize}
    \item \textbf{Technique:} Let $\mathcal{S}(x)$ denote the deterministic solver output for input $x$. We define the kernel by perturbing the input (e.g., the covariance matrix $\Sigma$ or the alpha vector $\alpha$) with noise $\epsilon \sim \mathcal{N}(0, \Sigma_{\text{noise}})$ and executing the solver $N$ times. The empirical measure is given by:
    \begin{equation}
        P(x) \approx \frac{1}{N} \sum_{i=1}^{N} \delta_{\mathcal{S}(x + \epsilon_i)}
    \end{equation}
    \item \textbf{Implementation in Mixed-Integer Programming (MIP):} For discrete solvers such as Gurobi, simple input perturbation may be insufficient to explore the solution space due to the discrete nature of the lattice. To properly induce the kernel, two specific strategies are employed:
    \begin{enumerate}
        \item \textbf{Objective Perturbation:} We add infinitesimal random noise to the linear objective vector $c$. This effectively breaks symmetries in the optimization landscape, forcing the solver to explore flat regions of the objective function that would otherwise be pruned.
        \item \textbf{Randomized Tie-Breaking:} By varying the internal random seed parameter (e.g., Gurobi \texttt{Seed}), we alter the tie-breaking rules within the branch-and-cut tree. This generates distinct valid leaves (optimal or near-optimal integer solutions) without modifying the problem definition.
    \end{enumerate}
    \item \textbf{Trade-off:} While mathematically rigorous and capable of capturing true non-convexities (such as tax cliffs), this method incurs a linear computational cost of $O(N)$ relative to the baseline optimization.
\end{itemize}

\subsubsection{The Solution Pool Method (Single-Pass Approximation)}
For Mixed-Integer Programming solvers, it is possible to approximate the kernel without repeated execution by exploiting the branch-and-bound tree structure.
\begin{itemize}
    \item \textbf{Configuration:} Modern solvers allow the retention of sub-optimal solutions found during the search. For example, configuring Gurobi with \texttt{PoolSolutions=N} and \texttt{PoolSearchMode=2} forces the solver to effectively exhaust the search tree to find the $N$ best distinct solutions.
    \item \textbf{Interpretation:} The set of solutions $\{y_1, \dots, y_N\}$ retrieved from the pool constitutes an empirical support for $P(x)$. This effectively delineates the ``Action Menu'' $K \odot R$ described in the DOTS framework (Section~\ref{sec:DOTS}), revealing the local geometry of the permissible region without multiple independent runs.
\end{itemize}

\subsubsection{Sensitivity Analysis (The Local Gaussian Approximation)}
For purely convex optimization problems (e.g., standard Mean-Variance Optimization without cardinality constraints), analytical derivatives provide a computationally inexpensive approximation.
\begin{itemize}
    \item \textbf{Technique:} By querying the solver's sensitivity report (specifically the Lagrange multipliers or dual values), we can construct a local quadratic approximation of the objective function.
    \item \textbf{Approximation:} The kernel $P(x)$ is approximated as a Gaussian distribution $\mathcal{N}(\mu, \Sigma_{H})$, where the covariance $\Sigma_{H}$ is derived from the inverse Hessian of the objective function at the optimum.
    \item \textbf{Trade-off:} This method offers $O(1)$ overhead but is strictly invalid for non-convex problems (MIPs), where it fails to capture the disconnected nature of the solution space.
\end{itemize}

\subsection{Efficient Representation and Storage}
\label{sub:storage}

Once the distribution $P(x)$ is estimated, it must be serialized into the Evidence Ledger (see Appendix~\ref{app:architecture}). Storing $N$ raw high-dimensional portfolio vectors is inefficient and opaque. We propose four storage schemas, selecting the appropriate format based on the topological complexity of the problem.

\begin{table}[h]
\centering
\begin{tabular}{|l|p{6cm}|l|}
\hline
\textbf{Storage Format} & \textbf{Description} & \textbf{Use Case} \\ \hline
\textbf{Gaussian} $(\mu, \Sigma)$ & Stores only the mean vector and covariance matrix. & Uni-modal Convex problems. \\ \hline
\textbf{Gaussian Mixture (GMM)} & Stores $k$ means and covariances: $\sum w_i \mathcal{N}(\mu_i, \Sigma_i)$. & Mixed-Integer Problems (MIPs). \\ \hline
\textbf{Quantile Vectors} & Stores the 5th, 50th, and 95th percentile weights per asset. & Reporting and Visualization. \\ \hline
\textbf{Convex Hull (Vertices)} & Stores the vertex portfolios defining the solution polytope. & Robust Safety (Safety Radius). \\ \hline
\end{tabular}
\caption{Efficient serialization formats for stochastic kernels.}
\label{tab:storage_formats}
\end{table}

\textbf{Recommendation:} For the general implementation of the HSP-r framework, the \textbf{Gaussian Mixture Model (GMM)} provides the optimal balance between compression and fidelity, as it is capable of capturing the multi-modal distributions characteristic of discrete portfolio constraints.

\subsection{Structural Failures of the Gaussian Assumption}

While it is tempting to simplify $P(x)$ to a unimodal Gaussian distribution for computational convenience, this assumption is structurally dangerous in financial engineering. Real-world constraints frequently ``clip,'' ``hollow,'' or ``warp'' the probability mass in ways that a single Gaussian cannot represent. The following sections detail three canonical failure modes---constraint hugging (Section~\ref{sub:hollow_shell}), the split peak from minimum position sizes (Section~\ref{sub:split_peak}), and the banana-shaped frontier from tax-efficient transitioning (Section~\ref{sub:banana})---each of which defeats the Gaussian assumption in a structurally distinct way.

\clearpage

\section{Handling Discrete Constraints}
\label{sec:discrete}

While the core theory of $\HSSimp$ relies on convexity and continuity to guarantee properness and stability, real-world portfolio construction is dominated by discrete constraints. These constraints—ranging from regulatory lot sizes to tax thresholds—introduce non-convexity and discontinuity into the re-implementation map $f$. This section details the specific failure modes these constraints create and how the probabilistic HSP-r framework adapts to handle them.

\subsection{The ``Split Peak'': Minimum Position Sizes}
\label{sub:split_peak}

A ubiquitous constraint in institutional asset management is the \textbf{Minimum Position Size} rule. This rule forbids holding ``dust'' positions that are operationally costly to trade or service.

\subsubsection*{Theoretical Formulation}
Let $w_i$ be the weight of asset $i$. The constraint is formally defined as a disjoint union of intervals:
\[
w_i \in \{0\} \cup [L, U] \quad \text{where } 0 < L \le U \le 1 \text{ (e.g., } L = 2\%, \; U = 10\%)
\]
This destroys the convexity of the permissible space $K$. If $x = 0$ (valid) and $y = 2\%$ (valid), the convex combination $0.5x + 0.5y = 1\%$ is strictly invalid.

\subsubsection*{The Solver Behavior (Bifurcation)}
Consider a scenario where the ``ideal'' (unconstrained) optimal weight for an asset is $w^* = 1.0\%$.
\begin{itemize}
    \item A deterministic solver must make a hard choice: either round down to $0\%$ or force up to $2\%$. This decision is often unstable, flipping based on microscopic changes in the input data.
    \item A stochastic solver (using randomized smoothing) will produce a \textbf{bimodal distribution}. With roughly 50\% probability, it outputs $0\%$; with 50\% probability, it outputs $\approx 2\%$.
\end{itemize}

\subsubsection*{Gaussian Failure and the Phantom Mean}
Approximating this bimodal distribution with a single Gaussian---a common shortcut in risk models---produces a qualitatively wrong answer.
\begin{itemize}
    \item \textbf{The Failure:} The mean of the distribution is $\approx 1\%$. A Gaussian model centers the probability mass at $1\%$.
    \item \textbf{The Phantom Portfolio:} This creates a ``Phantom Portfolio'' in the Evidence Ledger. The risk system records that the ``likely'' outcome is a 1
    \item \textbf{The Fix:} The system must store this kernel using a \textbf{Gaussian Mixture Model (GMM)} or a particle representation. This correctly identifies two distinct feasible regions (at 0 and 2) separated by the forbidden interval $(0, 2)$.
\end{itemize}

\subsection{The ``Hollow Shell'': Gross Exposure Limits}
\label{sub:hollow_shell}

Strategies that optimize for factors (e.g., Value or Momentum) often naturally seek maximum leverage. To control risk, a \textbf{Gross Exposure Limit} is applied.

\subsubsection*{Theoretical Formulation}
Let permissible leverage be bounded by $G$. The constraint is the $L_1$ ball:
\[
\sum_{i} |w_i| \le G \quad \text{(e.g., } G = 200\%)
\]
While this set is convex, the interaction with the objective function creates a topological problem known as ``Constraint Hugging.''

\subsubsection*{The Solver Behavior (Manifold Concentration)}
If the alpha signal is strong, the optimizer will push the portfolio to the absolute limit of the constraint to maximize exposure.
\begin{itemize}
    \item The solution $y^*$ does not lie in the interior of the feasible set; it lies on the boundary $\partial K$.
    \item Under perturbation (input noise), the solution does not move toward the center (de-leveraging); it slides along the ``skin'' of the constraint surface.
    \item \textbf{Dimensional Collapse:} Mathematically, the probability mass concentrates on a lower-dimensional manifold (a sphere of radius $G$).
\end{itemize}

\subsubsection*{Gaussian Failure and the Empty Interior}
\begin{itemize}
    \item \textbf{The Failure:} If we fit a high-dimensional Gaussian to points distributed on the surface of a sphere, the mean vector $\mu$ will lie at the center of the sphere.
    \item \textbf{The Implication:} The center represents a portfolio with \textit{lower} gross exposure than any actual realization. The Gaussian model implies that the solver occasionally de-leverages (moves internally).
    \item \textbf{The Risk:} This underestimates the ``crowding'' risk. In reality, the portfolio is always maximally levered. A Gaussian risk model would underestimate the probability of hitting a margin call because it assumes mass exists in the safer, lower-leverage interior.
    \item \textbf{The Fix:} The kernel must be represented using \textbf{manifold learning} techniques or simply raw particle storage to capture the fact that the interior is ``hollow.''
\end{itemize}

\subsection{The ``Banana'': Tax-Efficient Transitioning}
\label{sub:banana}

Transition management involves moving a portfolio from a legacy state $x$ to a target state $z$ while minimizing tracking error and keeping realized taxes below a budget $T$.

\subsubsection*{Theoretical Formulation}
The constraints combine an ellipsoid (tracking error) with a piecewise-linear function (taxes):
\[
(y - z)^T \Sigma (y - z) \le \epsilon^2 \quad \text{and} \quad \sum \text{RealizedGain}(y_i) \le T
\]
Crucially, the tax function is zero for losses (unless harvesting is allowed) and linear for gains, introducing a ``kink'' in the feasible space.

\subsubsection*{The Solver Behavior (Curvature)}
The solver attempts to sell legacy assets to buy the target.
\begin{itemize}
    \item \textbf{Phase 1:} It sells assets with losses or low gains.
    \item \textbf{Phase 2 (The Wall):} It hits the tax cap $T$. It cannot sell the high-gain legacy assets.
    \item \textbf{Phase 3 (Substitution):} To reduce tracking error further, it begins buying ``proxy'' assets—stocks highly correlated with the target but which are not currently held (thus incurring no tax).
    \item \textbf{Result:} As correlation structures shift, the ``efficient frontier'' of solutions curves around the expensive tax barrier, forming a shape often described in Monte Carlo literature as ``the banana.''
\end{itemize}

\subsubsection*{Gaussian Failure and False Compliance}
\begin{itemize}
    \item \textbf{The Failure:} A Gaussian distribution is defined by ellipsoid contours, which are inherently convex.
    \item \textbf{The Violation:} Fitting an ellipsoid to a banana shape forces the ellipsoid to cover the ``hole'' inside the curve.
    \item \textbf{The Reality:} The space inside the curve corresponds to portfolios that \textit{violate} the tax budget (they represent the direct trade, not the proxy trade).
    \item \textbf{The Fix:} The Gaussian model validates high-tax portfolios as ``likely outcomes.'' The system must use the \textbf{Highest Density Region (HDR)} approach (Section~\ref{sub:HDR}), which naturally wraps the safety region around the non-convex obstacle, excluding the tax-violation zone.
\end{itemize}

\subsection{DOTS: The Menu-Driven Alternative}
\label{subsec:dots_framework}

While the probabilistic HSP-r framework (Section~\ref{sec:Polish}) manages discrete constraints by quantifying the \textit{risk} of solver instability, the \textbf{Double Operadic Theory of Systems (DOTS)} \cite{LibkindMyers} manages them by quantifying the \textit{choice}. In the presence of non-convexities like minimum position sizes, a single "optimal" portfolio often does not exist or is unstable. DOTS resolves this by treating the re-implementation not as a function $f: K \to K'$, but as an \textbf{Action} that generates a \textbf{Menu} of valid options.

\subsubsection*{The Action Operator}

Instead of forcing a solver to pick a single point $y$ (which might be an invalid average), we define the re-implementation as the action of the alignment relation $R$ on the hub state $K$.

\begin{definition}[The Action Menu]
\label{def:action-menu}
Let $K \subseteq \Delta^n$ be the hub space and $R \subseteq \Delta^n \times \Delta^m$ be an alignment relation (e.g., a constraint). The \textbf{Action Menu}, denoted $K \odot R$, is the set of all permissible spoke portfolios:
\[
K \odot R := \{ y \in \Delta^m \mid \exists x \in K \text{ such that } (x, y) \in R \}
\]
\end{definition}

\subsubsection*{Application to Discrete Constraints}
Consider the \textbf{Minimum Position Size} constraint (Section~\ref{sub:split_peak}) where a holding must be either $0\%$ or $\ge 2\%$.
\begin{itemize}
    \item \textbf{Standard Optimizer:} Returns $1\%$ (invalid) or flips randomly between $0\%$ and $2\%$.
    \item \textbf{DOTS Action:} The set $K \odot R$ is explicitly \textbf{disconnected}. It contains two valid components:
    \[
    K \odot R = \{ \text{Portfolios with } w_i = 0 \} \cup \{ \text{Portfolios with } w_i \ge 2\% \}
    \]
    The "Action" does not make the decision; it \textit{presents} the valid decision space to the Portfolio Manager or the next stage of the algorithm.
\end{itemize}

\subsubsection*{Operadic Composition of Constraints}
\label{sub:dots_operad}

The power of DOTS lies in its ability to compose these menus safely. Complex constraints are built by ``wiring together'' simpler constraints using the \textbf{Operad of Wiring Diagrams}.

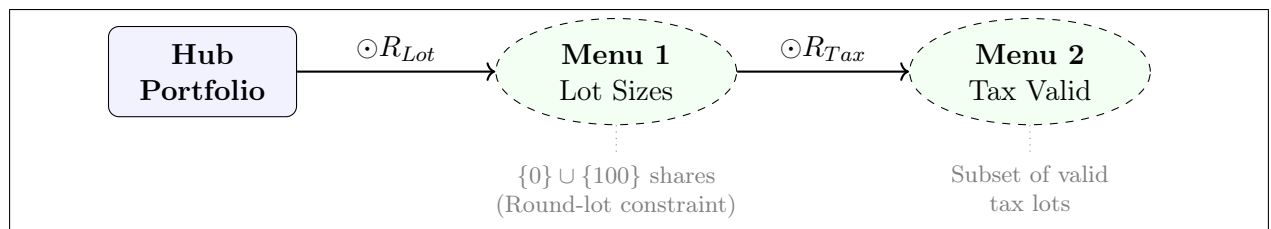
\begin{figure}[h]
    \centering
    \begin{tikzpicture}[scale=1.0, transform shape]
        \tikzstyle{block} = [draw, rectangle, minimum height=1.2cm, minimum width=2.5cm, rounded corners, fill=blue!5, align=center];
        \tikzstyle{menu} = [draw, ellipse, minimum height=1.2cm, minimum width=3.2cm, fill=green!5, dashed, align=center];
        
        \node[block] (Hub) at (0,0) {\textbf{Hub}\\\textbf{Portfolio}};
        \node[menu] (Menu1) at (5.5,0) {\textbf{Menu 1}\\Lot Sizes};
        \node[menu] (Menu2) at (11.0,0) {\textbf{Menu 2}\\Tax Valid};
        
        \draw[->, thick] (Hub) -- node[above] {$\odot R_{Lot}$} (Menu1);
        \draw[->, thick] (Menu1) -- node[above] {$\odot R_{Tax}$} (Menu2);
        
        \node[below=0.4cm of Menu1, font=\footnotesize, align=center, text=gray] (Viz1) {$\{0\} \cup \{100\}$ shares\\(Round-lot constraint)};
        \node[below=0.4cm of Menu2, font=\footnotesize, align=center, text=gray] (Viz2) {Subset of valid\\tax lots};
        
        \draw[dotted, gray] (Menu1) -- (Viz1);
        \draw[dotted, gray] (Menu2) -- (Viz2);
    \end{tikzpicture}
    \caption{The DOTS pipeline does not propagate a single portfolio; it propagates the \textit{set of valid possibilities}, narrowing the menu at each step.}
    \label{fig:dots_pipeline}
\end{figure}

\subsubsection*{The ``Design vs. Execution'' Split}
This framework formalizes the industrial separation of duties:
\begin{enumerate}
    \item \textbf{Design Phase (DOTS):} The ``Semantic Core'' calculates the Action Menu $K \odot R$. This is a set-valued operation that identifies \textit{all} theoretically valid portfolios (e.g., all valid integer combinations of round lots).
    \item \textbf{Selection Phase (Human/Solver):} A selection function $s: K \odot R \to \text{Spoke}$ picks a specific portfolio from the menu.
    \item \textbf{Execution Phase (HSP-r):} The selected portfolio is traded, with execution noise modeled by the probabilistic kernel.
\end{enumerate}

By separating the \textbf{Menu Generation} (which handles non-convexity via sets) from the \textbf{Selection} (which handles preference via optimization), we avoid the traps where a solver tries to average across two distinct valid options.

\clearpage  \section{Liquidity-Aware Portfolio Construction}
\label{sec:liquidity}

This section applies all three layers of the framework---the basic double category $\HSSimp$, the DOTS menu calculus, and the Wasserstein transport model---to the problem of adding illiquid or semi-liquid allocations to an existing portfolio.  The motivating setting is an investor who holds a liquid public-market portfolio and wishes to incorporate private assets (private equity, private credit, real assets, infrastructure) while respecting aggregate liquidity requirements, capacity constraints, and transaction cost budgets.  The foundational theory of the illiquidity premium \cite{Amihud,AmihudMendelson,PastorStambaugh,AcharyaPedersen}, optimal portfolio rebalancing with proportional transaction costs \cite{Constantinides,DavisNorman}, and portfolio choice under illiquidity constraints \cite{AngPapanikolaouWesterfield,Longstaff} all inform the constructions below.

The economic context draws on the analysis in \cite{Phoa}, which distinguishes liquidity as a \emph{property}---how costly it is to transact in an asset---from liquidity as a \emph{substance}: the portion of wealth that can be freely reallocated.  That analysis identifies eight components of the private asset premium (pure illiquidity premium, risk premium, guarantee value, information costs, transaction complexity, rent of balance sheet, control alpha, and offsets such as lower apparent volatility \cite{GetmanskyLoMakarov}).  The optimization aspects draw on the mean--variance framework in \cite{PhoaMVO}.

\subsection{Liquidity Constraints in the Basic Framework}
\label{sub:liq-basic}

Suppose the hub space is a public-market portfolio $K_{\mathrm{pub}}\subseteq\Delta^n$, and the spoke universe $\Delta^m$ includes both the original liquid assets and a set of illiquid or semi-liquid vehicles (private equity funds, private credit vehicles, interval funds, etc.).  The spoke permissible space $K_{\mathrm{blend}}\subseteq\Delta^m$ incorporates baseline feasibility constraints on the blended portfolio.

Liquidity requirements impose further structure.  We model three types of constraint, each defining a closed subset of $\Delta^m$:

\begin{enumerate}[label=(\roman*),leftmargin=*]
\item \textbf{Aggregate illiquidity cap.} Let $\mathcal{I}\subseteq\{0,\ldots,m\}$ index the non-daily-liquid assets.  The constraint
\[
L_{\alpha}:=\Big\{y\in\Delta^m:\sum_{i\in\mathcal{I}} y_i \le \alpha\Big\}
\]
caps the total weight in illiquid assets at $\alpha\in(0,1)$.  This is a closed half-space intersected with $\Delta^m$, hence a closed convex polytope.  Regulatory examples include the SEC's Rule~22e-4 classification of portfolio liquidity into buckets.

\item \textbf{Position-size and capacity constraints.}  For each illiquid asset $i\in\mathcal{I}$, let $c_i>0$ be a capacity limit reflecting the asset's available deal size or the investor's maximum commitment relative to fund size.  The constraint
\[
C:=\{y\in\Delta^m: y_i\le c_i\ \forall\, i\in\mathcal{I}\}
\]
is again a closed convex polytope.  In private markets, these bounds arise from fund subscription caps, co-investment limits, or internal concentration policies.

\item \textbf{Maintenance transaction cost cap.}  Let $\tau_i\ge 0$ be the expected per-unit annual rebalancing cost for asset~$i$ (incorporating bid--ask spread, market impact, and any fund-level redemption fees).  The constraint
\[
M_{\kappa}:=\Big\{y\in\Delta^m:\sum_{i=0}^m \tau_i\, y_i\le \kappa\Big\}
\]
bounds the expected annual maintenance cost at $\kappa$ basis points.  This is a weighted-fee projector of the same type as $R_{\mathrm{fee}\le\tau}$ in Section~\ref{sec:DOTS}.
\end{enumerate}

The permissible spoke space for the liquidity-constrained re-implementation is
\[
K_{\mathrm{liq}}:=K_{\mathrm{blend}}\cap L_{\alpha}\cap C\cap M_{\kappa},
\]
which is a closed (hence compact) convex subset of $\Delta^m$ and therefore a valid object in $\HSSimp$.

The liquidity-aware re-implementation is then constructed via the optimization framework of Section~\ref{sec:optimal-reimpl}.  Following that section's notation, let $(V,\|\cdot\|_V)$ be a finite-dimensional space of quantitative attributes (e.g., asset class exposures), and let $g_{\mathrm{pub}}\colon K_{\mathrm{pub}}\to V$ and $g_{\mathrm{liq}}\colon K_{\mathrm{liq}}\to V$ be continuous attribute maps that express each portfolio's characteristics in~$V$.  Let $u\colon V\to\R$ be a continuous objective function (e.g., expected return net of fees).  For a hub portfolio $x\in K_{\mathrm{pub}}$ and a candidate spoke portfolio $y\in K_{\mathrm{liq}}$, define the combined objective
\[
F_{\mathrm{liq}}(x,y):=\|g_{\mathrm{pub}}(x)-g_{\mathrm{liq}}(y)\|_V^p - \lambda\, u(g_{\mathrm{liq}}(y)),
\]
where $p>1$ controls the penalty for attribute mismatch and $\lambda\ge 0$ weights the objective against tracking.  If $\|\cdot\|_V$ is strictly convex (e.g., an inner-product norm), $g_{\mathrm{liq}}$ is affine and injective, and $u$ is strictly concave, then $F_{\mathrm{liq}}(x,\cdot)$ is strictly convex in $y$ and Theorem~\ref{thm:berge-optimal} guarantees a unique optimizer $y^*(x)$ for each $x\in K_{\mathrm{pub}}$, and the resulting map $f_{\mathrm{liq}}(x):=y^*(x)$ is a continuous (hence proper) horizontal morphism $f_{\mathrm{liq}}\colon K_{\mathrm{pub}}\to K_{\mathrm{liq}}$ in $\HSSimp$.

\begin{remark}[Alignment with the public-market portion]
In addition to re-implementing the hub, we typically require the liquid portion of the spoke to \emph{track} the hub.  Let $\pi_{\mathrm{liq}}\colon\Delta^m\to\Delta^n$ project out the liquid weights (renormalized by $1-\sum_{i\in\mathcal{I}}y_i$).  The tracking relation
\[
R_{\mathrm{track},\epsilon}:=\{(x,y)\in K_{\mathrm{pub}}\times K_{\mathrm{liq}}:\|g_{\mathrm{pub}}(x)-g_{\mathrm{pub}}(\pi_{\mathrm{liq}}(y))\|_2\le\epsilon\}
\]
is closed and can be imposed as an alignment constraint in the usual way (Section~\ref{sec:constraints}).  Pushforward and pullback along the re-implementation $f_{\mathrm{liq}}$ then propagate this tracking requirement through the hub--spoke tree.
\end{remark}

\subsection{DOTS: Menus with Liquidity and Capacity Screens}
\label{sub:liq-dots}

The DOTS framework (Section~\ref{sec:DOTS}) is particularly natural for illiquid allocation, for two reasons:
\begin{enumerate}[label=(\alph*),leftmargin=*]
\item Private asset allocations typically involve discrete, judgment-laden decisions---which fund, which vintage, how much to commit---that produce a \emph{menu} of feasible spoke portfolios rather than a unique optimum.
\item Liquidity and capacity constraints act as idempotent screens that narrow the menu, exactly the role of projectors in the DOTS calculus.
\end{enumerate}

\subsubsection*{Liquidity and capacity as projectors}
As in Section~\ref{sec:DOTS}, the aggregate illiquidity cap and capacity constraints define diagonal (projector) relations:
\begin{align*}
R_{\mathrm{liq},\alpha}&:=\{(y,y): y\in L_\alpha\},\\
R_{\mathrm{cap}}&:=\{(y,y): y\in C\},\\
R_{\mathrm{maint},\kappa}&:=\{(y,y): y\in M_\kappa\}.
\end{align*}
Each is idempotent ($K\odot R\odot R = K\odot R$) and order-insensitive with respect to itself (Proposition~\ref{prop:action}(e)).

\subsubsection*{The core--satellite template with liquidity gate}

The core--satellite architecture formalized by \cite{AmencMalaiseMartellini} extends naturally within the operadic template of Section~\ref{sec:DOTS}.  Let
\begin{itemize}
\item $K_{\mathrm{core}}\subseteq\Delta^n$ be the liquid public-market hub,
\item $K_{\mathrm{sat}}\subseteq\Delta^p$ be the illiquid satellite universe (e.g., a set of private credit funds and PE vintage pools),
\item $w\in[0,1]$ be the target allocation to the liquid core,
\item $R_{\mathrm{mix},w}$ be the blending relation that creates the composite portfolio $y=w\cdot y_{\mathrm{core}}+(1-w)\cdot y_{\mathrm{sat}}$ from a core allocation $y_{\mathrm{core}}\in K_{\mathrm{core}}$ and a satellite allocation $y_{\mathrm{sat}}\in K_{\mathrm{sat}}$.
\end{itemize}
The liquidity-aware menu is then
\[
K^{\star}=(K_{\mathrm{core}}\otimes K_{\mathrm{sat}})\;\odot\; R_{\mathrm{mix},w}\;\odot\; R_{\mathrm{liq},\alpha}\;\odot\; R_{\mathrm{cap}}\;\odot\; R_{\mathrm{maint},\kappa}.
\]
The projectors narrow the menu from both ends: $R_{\mathrm{liq},\alpha}$ enforces the aggregate cap, $R_{\mathrm{cap}}$ enforces position-level limits reflecting fund capacity, and $R_{\mathrm{maint},\kappa}$ enforces the ongoing cost budget.  By Proposition~\ref{prop:action}(a), $K^{\star}$ is closed.

\subsubsection*{Market impact awareness}
For assets where the investor's position is large relative to available capacity, market impact during implementation and ongoing rebalancing is a first-order concern \cite{Kyle,AlmgrenChriss}.  Let $\mathrm{ADV}_i$ denote the average daily volume (or, for private vehicles, the periodic subscription/redemption capacity) of asset~$i$.  The market-impact-aware constraint
\[
R_{\mathrm{impact}}:=\{(y,y): y_i\cdot W_{\mathrm{total}}\le \beta_i\cdot\mathrm{ADV}_i\ \forall\, i\in\mathcal{I}\}
\]
where $W_{\mathrm{total}}$ is the total portfolio value and $\beta_i$ is a participation-rate limit, is again a diagonal projector.  It ensures that the spoke portfolio can be \emph{implemented} without moving the market, and can be applied at any point in the DOTS pipeline.

\begin{remark}[Determinization]
If the blending relation $R_{\mathrm{mix},w}$ is combined with a strictly convex objective (e.g., minimizing tracking error plus a regularizer), it has a unique solution for each input and therefore determinizes to a graph of a continuous map.  By Frobenius reciprocity (Theorem~\ref{thm:DOTS-Frob}), the liquidity projectors then commute with this map: screening before or after optimization yields the same result.  This is practically important because it allows pre-screening the illiquid universe for capacity before running the optimizer.
\end{remark}

\begin{remark}[The private asset premium in the DOTS framework]
The eight components of the private asset premium identified in \cite{Phoa}---pure illiquidity premium, risk premium, guarantee value, information costs, transaction complexity, rent of balance sheet, control alpha, and offsets (cf.\ the empirical decompositions in \cite{KaplanSchoar,FranzoniNowakPhalippou})---can be viewed as attributes of the satellite objects in $K_{\mathrm{sat}}$.  The DOTS menu $K^{\star}$ reflects which combinations of these components are accessible to the investor after all screens have been applied.  In particular, control alpha (component~7) is only available to investors whose capacity constraints and governance structure permit concentrated, illiquid positions---a condition encoded by the interplay of $R_{\mathrm{cap}}$ and $R_{\mathrm{impact}}$.
\end{remark}

\subsection{Wasserstein: Transaction Cost Budgets}
\label{sub:liq-wasserstein}

The Wasserstein framework (Section~\ref{sec:wasserstein}) acquires a direct financial interpretation in the liquidity setting.  Recall that in the probabilistic theory (Section~\ref{sec:Polish}), a re-implementation is modeled as a stochastic kernel $P\colon K_{\mathrm{pub}}\rightsquigarrow K_{\mathrm{liq}}$ that maps each hub portfolio $x$ to a probability distribution $P(x)$ over spoke portfolios, capturing solver noise and execution uncertainty.  The Wasserstein safety condition $W_1(P(x),\mathcal{P}_S)\le\epsilon$ then measures the \emph{expected turnover required to move a realized portfolio into compliance} with a constraint set $S$.  For illiquid allocations, this turnover is expensive, making the Wasserstein budget the binding constraint in practice.

\subsubsection*{Asset-weighted transport cost}

The standard Wasserstein distance \cite{VillaniOT,Galichon} uses the Euclidean (or $\ell_1$) metric on the simplex.  In a liquidity-aware setting, this should be replaced by a \emph{transaction-cost-weighted} metric that reflects the heterogeneous cost of rebalancing across assets:
\[
d_{\tau}(y,y'):=\sum_{i=0}^{m}\tau_i\,|y_i-y_i'|,
\]
where $\tau_i\ge 0$ is the per-unit transaction cost for asset~$i$ (incorporating half-spread, market impact, and any redemption or subscription fees).  For liquid assets, $\tau_i$ is small (a few basis points for large-cap equities); for illiquid assets, $\tau_i$ can be large (tens or hundreds of basis points for private credit, or effectively infinite for locked-up PE commitments).

The resulting weighted Wasserstein distance
\[
W_{1,\tau}(\mu,\nu):=\inf_{\pi\in\Pi(\mu,\nu)}\int d_{\tau}(y,y')\,d\pi(y,y')
\]
inherits all the structural properties used in Section~\ref{sec:wasserstein}: it is a metric on $\mathcal{P}_1(\Delta^m)$, satisfies the triangle inequality, and admits the Kantorovich--Rubinstein dual.  The composition, adjunction, and coherence theorems of Section~\ref{sec:wasserstein} (Theorems~\ref{thm:wasserstein-pushforward-comp}--\ref{thm:wasserstein-frobenius}) carry over with $d_\tau$ replacing the unweighted metric throughout.

\subsubsection*{Pre-trade compliance: the cost of exit}

Before committing capital to illiquid allocations, a prudent investor asks: \emph{if conditions change, what would it cost to liquidate back to a fully liquid portfolio?}  (The illiquidity discount analysis of \cite{Longstaff} poses the same question in a single-asset setting.)  This is precisely the Wasserstein safety check.  Let $S_{\mathrm{liq}}:=\{y\in\Delta^m: y_i=0\ \forall\,i\in\mathcal{I}\}$ be the ``liquid-only'' face of the simplex.  For a stochastic re-implementation kernel $P\colon K_{\mathrm{pub}}\rightsquigarrow K_{\mathrm{liq}}$, the cost of exit from hub portfolio~$x$ is
\[
\mathrm{ExitCost}(x):=W_{1,\tau}(P(x),\,\mathcal{P}_{S_{\mathrm{liq}}}),
\]
the minimum expected transaction cost (under the weighted metric $d_\tau$) required to transport the realized spoke distribution onto the liquid-only face.  A pre-trade compliance check requires $\mathrm{ExitCost}(x)\le\epsilon_{\mathrm{exit}}$ for all $x\in K_{\mathrm{pub}}$.  In the notation of Definition~\ref{def:wasserstein-safety}, this is the $\epsilon_{\mathrm{exit}}$-Wasserstein pullback condition $K_{\mathrm{pub}}\subseteq P^{*,\epsilon_{\mathrm{exit}}}_{W}\, S_{\mathrm{liq}}$ (with the metric $d_\tau$).  The Wasserstein adjunction (Theorem~\ref{thm:wasserstein-adjunction}) then guarantees that verifying this at the hub level is equivalent to verifying that the pushforward of any hub constraint lands within the $\epsilon_{\mathrm{exit}}$-cure set of $S_{\mathrm{liq}}$.

\subsubsection*{Maintenance budgets over multiple periods}

Portfolio maintenance---rebalancing back toward target weights as illiquid positions are drawn down or marked up---incurs ongoing transaction costs.  The dynamic trading framework of \cite{GarleanuPedersen} shows that optimal rebalancing ``aims in front of the target,'' anticipating future costs; the Wasserstein composition theorem below provides a complementary budget-based approach.  Model a $T$-period rebalancing process as a chain of stochastic kernels $P_1,\ldots,P_T$, where each $P_t$ represents the solver that rebalances the portfolio at the end of period~$t$.  Let $\epsilon_t$ be the Wasserstein budget allotted to the rebalancing at period~$t$, and let $L_t$ be the Lipschitz constant of $P_t$ with respect to $d_\tau$ (measuring how much one period's execution noise is amplified by the next).  The Wasserstein composition theorem (Theorem~\ref{thm:wasserstein-pushforward-comp}) gives a single-step lax inclusion with budget $L\delta + \epsilon$.  Iterating: the budget after two stages is $L_2(L_1\epsilon_0 + \epsilon_1) + \epsilon_2 = L_2 L_1 \epsilon_0 + L_2 \epsilon_1 + \epsilon_2$; by induction, the cumulative Wasserstein budget for the entire chain is bounded by $\sum_{t=1}^{T}\bigl(\prod_{s=t+1}^{T}L_s\bigr)\cdot\epsilon_t$.  This gives the investor a principled way to set the per-period rebalancing tolerance $\epsilon_t$ so that the cumulative expected maintenance cost stays within a multi-period budget.

For semi-liquid vehicles (e.g., interval funds with quarterly redemption windows, or private credit funds with periodic liquidity), the Lipschitz constants $L_t$ vary with the redemption calendar: $L_t$ is small in periods when redemptions are available (rebalancing is relatively cheap) and large in lock-up periods (rebalancing requires secondary-market transactions at a discount).

\begin{remark}[Connection to the illiquidity premium]
The weighted transport cost $d_\tau$ endogenizes the illiquidity premium into the compliance framework.  An asset with high $\tau_i$ contributes disproportionately to the exit cost and maintenance budget, which constrains the feasible allocation to that asset.  In equilibrium, the investor requires a higher expected return to justify consuming a larger share of the Wasserstein budget---precisely the pure illiquidity premium of \cite{Phoa}; see \cite{AmihudMendelson} and \cite{Constantinides} for the foundational equilibrium analyses.  The mean--variance framework for analyzing this trade-off is developed in \cite{PhoaMVO}, where the illiquidity premium is derived as the marginal cost of relaxing the transaction cost constraint.
\end{remark}

\begin{remark}[Choosing the right framework layer]
The three layers address distinct aspects of liquidity-aware construction:
\begin{itemize}
\item \textbf{$\HSSimp$ (basic):} Static structure.  Defines the feasible space, ensures properness and coherent propagation of constraints.  Use for initial portfolio design and regulatory compliance verification.
\item \textbf{DOTS (menus):} Design-phase discretion.  Models the menu of valid illiquid allocations before a commitment is made.  Use when private asset selection involves judgment, vintage diversification, or manager selection.
\item \textbf{Wasserstein (transport):} Execution-phase cost.  Quantifies the ongoing cost of maintaining compliance and the cost of exiting positions.  Use for pre-trade analysis, maintenance budgeting, and stress testing.
\end{itemize}
\end{remark}

\clearpage

\part*{Conclusion}
\addcontentsline{toc}{part}{Conclusion}
\label{sec:conclusion}

We have shown that the compositional structure of portfolio re-implementation---the chains of transformations, alignment checks, and constraint propagations that constitute modern portfolio management---admits a rigorous formalization using elementary category theory and point-set topology. The key insight is architectural: the right objects are compact spaces (not open or arbitrary subsets of simplices), and once this is established, the required properness of morphisms follows automatically. Everything else---adjunction, coherence, path independence---is a consequence.

\begin{enumerate}
\item \textbf{Objects as Permissible Spaces}: Objects are permissible portfolio
spaces $K$, defined as \textbf{closed subsets} of ambient simplices
$\Delta^n$. This directly models baseline constraints.

\item \textbf{Automatic Properness}: Because $K$ is a closed subset of a
compact space $\Delta^n$, $K$ is \textbf{compact}. This ensures that all horizontal
morphisms $f: K_1 \to K_2$ (continuous maps) are \textbf{automatically
proper}.

\item \textbf{Coherent Framework $\HSSimp$}: This automatic properness is the
key: the pushforward operation $f_!$ preserves closedness
(Theorem~\ref{thm:proper-closed}). This allows for a well-defined double
category $\HSSimp$ where all compositional laws (Adjunction, Beck--Chevalley,
Frobenius) hold.

\item \textbf{Construction from Optimization}: We showed (Theorem
\ref{thm:berge-optimal}) that practical, optimization-based methods for
portfolio re-implementation (e.g., fee minimization) naturally construct a
continuous map $f$ on a \textbf{closed domain $K$}---precisely the
objects and morphisms of $\HSSimp$.

\item \textbf{Failure Analysis}: We proved that relaxing the closed-object
requirement (e.g., using open interiors) breaks properness and
catastrophically breaks all coherence properties (Theorem
\ref{thm:closure-fails}), justifying our definition.
\end{enumerate}

The technical insight is that by defining \textbf{objects as compact spaces}
(closed subsets of $\Delta^n$), the required \textbf{properness of morphisms}
becomes an automatic consequence, not an extra condition. This simplifies the
framework while strengthening its conclusions: it provides \textbf{compositional coherence} and \textbf{path independence};
\textbf{logical duality} (adjunction); and
\textbf{boundary stability} (properness).

The $\HSSimp$ framework enables:
\par
\begin{itemize}
\item \textbf{Multi-stage portfolio construction} with provable consistency.
\item \textbf{Modeling nonlinear customization}: Natively supports modeling
complex, optimization-based tasks like ``replace this active portfolio with a
low-fee/tax-efficient passive basket.''
\item \textbf{Flexible verification}: Check alignment at \textit{any} stage.
\end{itemize}

No new mathematics was required. The novelty is entirely in the model---in the claim that these particular objects and morphisms are the right ones for the domain, and in the observation that standard results from topology and category theory, once properly instantiated, yield the coherence properties that industrial portfolio management needs but has never had.

\begin{remark}
A reader who works through the proofs in Part~III will find that even the most technically demanding arguments---the measurable selection in the Wasserstein pullback composition, the integral splitting in the convolution stability proposition, the $\epsilon/2$ tightness argument for composed Feller kernels---are, once the definitions are in place, applications of standard machinery: the triangle inequality, Lipschitz contraction, Kuratowski--Ryll-Nardzewski. The one genuinely non-obvious modelling decision is the co-Lipschitz hypothesis (metric openness) in the safety radius composition theorems, where the insight is that Lipschitz continuity alone governs error \emph{amplification} but cannot guarantee \emph{coverage}. This is as it should be. If the theorems were not more or less inevitable consequences of the definitions, the definitions would be wrong. The difficulty in this paper---such as it is---lives in the architecture: in the choice of objects, morphisms, and 2-cells that make the standard results say what the application needs them to say.
\end{remark}

The extensions developed in Parts~II and~III address the gaps between the basic theory and industrial practice: a \textbf{menu-driven theory} (DOTS, Section~\ref{sec:DOTS}) for portfolio customization with discretionary overlays; a \textbf{probabilistic theory} (HSP-r, Section~\ref{sec:Polish}) for real-time replication at scale; and a \textbf{liquidity-aware framework} (Section~\ref{sec:liquidity}) integrating all three layers.

\subsection*{Critical Perspectives and Limitations}
\addcontentsline{toc}{subsection}{Critical Perspectives and Limitations}

The gap between a mathematical framework and an industrial implementation is real. It would be dishonest to present the theorems above without being explicit about where the difficulties lie.

\subsubsection*{Objections from Applied Category Theory}
\begin{enumerate}
    \item \textbf{Thinness Limits Expressivity:} The restriction to a \textbf{thin} double category (where 2-cells are unique inclusions) discards quantitative metadata. Treating alignment as a binary state cannot natively model \emph{degrees} of compliance (e.g., tracking error magnitude). Capturing this requires the richer \textbf{Span} model (Section~\ref{sec:2-cells}), which introduces the complexity of bicategorical coherence isomorphisms.
    \item \textbf{Loss of Companions in Probabilistic Settings:} In the stochastic extension (HSP-r), the construction of companions fails because the mapping from horizontal morphisms to vertical relations is not a strict functor (due to lax composition). This prevents treating functions as ``perfect'' relations and forces the system to treat $f_!$ and $f^*$ as primitive predicate transformers.
    \item \textbf{Restrictiveness of Strict Beck--Chevalley:} The strict path independence condition requires the commuting square to be \emph{pointwise cartesian} (Definition~\ref{def:pointwise-cartesian}). This implies surjectivity on points (no missing pre-images), a condition often violated in lossy aggregation pipelines. Consequently, the ``Lax'' condition (Audit Safety) is the only guarantee available in general settings.
\end{enumerate}

\subsubsection*{Objections from Quantitative Finance}
\begin{enumerate}
    \item \textbf{Computational Intractability of Exact Pushforwards:} While $f_!R$ is theoretically closed, calculating its explicit boundary requires quantifier elimination, which is intractable in high dimensions ($N \approx 500$). Practical implementations must rely on the Monte Carlo approximations or convex feasibility checks described in Section~\ref{sec:computation}.
    \item \textbf{Organizational Barriers to Commutativity:} The mathematical requirement for commuting squares (Bellman consistency) implies a unified utility function across the firm. In reality, different desks---Tax versus Model, say---operate with conflicting mandates, and the order of operations is often an irreducible business decision rather than a mathematical triviality. The framework can diagnose the non-commutativity; it cannot resolve the organizational politics that produce it.
    \item \textbf{Gaussian Failures on Structural Constraints:} As detailed in Section~\ref{sec:Polish}, approximating solver outputs with Gaussian distributions for storage efficiency fails in structurally important scenarios like ``constraint hugging'' (hollow shells) or ``split peaks'' (minimum position sizes). The framework requires particle representations or Gaussian Mixture Models to remain valid.
    \item \textbf{Discrete Constraints Break Continuity:} The topological foundation relies on the continuity of $f$. Real-world constraints involving discrete logic (e.g., lot sizes or fixed trade costs) challenge this assumption. While smoothed estimators (Example~\ref{ex:gurobi-feller}) restore continuity in the limit, the underlying discrete cliffs remain a source of residual operational risk.
\end{enumerate}

\subsection*{Future Directions: The Dual Calculus of Alignment}
\addcontentsline{toc}{subsection}{Future Directions: The Dual Calculus of Alignment}

The adjunction $f_! \dashv f^*$ established in Theorem~\ref{thm:adjunction} points toward a deeper structural insight: the logical duality between \emph{achievability} and \emph{compliance}. This relationship is formally captured by the theory of \textbf{Chu spaces} \cite{Pratt}.

In the $\HSSimp$ framework, we have prioritized the \emph{topological} perspective. We viewed the permissible space $K$ as the primary object and treated re-implementation $f: K_1 \to K_2$ as a forward transformation between asset universes.

A Chu-theoretic perspective invites the \emph{dual} view: treating the system of ``alignment relations'' as the primary object. Instead of asking how a portfolio moves through the pipeline, we ask how the \emph{conditions} of alignment transform under the operators $f_!$ and $f^*$.

\begin{itemize}
    \item \textbf{The Primal View (Forward Propagation):} We push an alignment relation forward via $f_!$ to determine the set of \emph{permissible} downstream relations (i.e., characterizing exactly which spoke portfolios are realizable from aligned hubs).
    \item \textbf{The Dual View (Backward Verification):} We pull an alignment relation backward via $f^*$ to determine the set of \emph{required} upstream relations (i.e., characterizing exactly which hub portfolios will satisfy the downstream constraints).
\end{itemize}

This mathematical duality mirrors---and, we suspect, partly explains---the practical \textbf{separation of duties} in institutional finance. Portfolio construction teams operate in the primal domain, designing the maps $f$ to transform assets; compliance and risk teams operate in the dual domain, defining the alignment relations $S$ that portfolios must respect. The two groups reason differently not because of organizational accident but because they inhabit genuinely dual mathematical structures. The tension between ``what can be built'' and ``what must be satisfied'' is not a defect of institutional design; it is the structural content of fiduciary obligation.

In this dual calculus, the re-implementation process functions as a \emph{predicate transformer}: it allows downstream rules $S$ applicable to spokes (defined by compliance) to be rigorously translated into upstream constraints $f^*S$ on the hub (executable by portfolio managers). This ensures that the logical requirements of the destination can be automatically and correctly imposed on the hub.

While the adjunction $f_! \dashv f^*$ in the basic theory already describes this situation adequately, it is possible that recasting the theory in the language of Chu spaces may enable a smoother formalization of the interaction between portfolio construction and compliance verification. The rigorous development of this idea must be left to future work.

\subsection*{Future Directions: Synthetic Probability and Support Functors}
\addcontentsline{toc}{subsection}{Future Directions: Synthetic Probability and Support Functors}

While this work grounds the probabilistic Hub-and-Spoke framework (HSP-r) in the explicit analysis of Feller kernels on Polish spaces (Section~\ref{sec:Polish}), recent developments in categorical probability suggest a natural generalization. The theory of \textit{Markov categories}, as developed by Fritz and others \cite{Fritz}, treats probability distributions and kernels purely as morphisms in a symmetric monoidal category, abstracting away the underlying measure-theoretic details.

Future work could explore replacing the analytic definitions of Section~\ref{sec:Polish} with this synthetic approach. Specifically, the operation of deriving a ``Safety Region'' or ``Action Menu'' from a stochastic solver can be formalized as a \textbf{Support Functor}
\[
\mathrm{Supp}_\epsilon: \mathbf{Stoch} \to \mathbf{Rel}
\]
that maps the category of stochastic kernels to the category of relations. This formalism would allow the coherence results for risk budgets (Theorems \ref{thm:pullback-comp} and \ref{thm:pushforward-comp}) to be derived as natural transformations, decoupling the audit logic from the specific topology of the asset space.

Note that this perspective offers a structural definition of ``determinism'' essential for audit. Rather than checking numerical variance, we may define a re-implementation morphism $f$ as deterministic if it commutes with the copy map (comultiplication) $\Delta$, i.e. $\Delta \circ f = (f \otimes f) \circ \Delta$.

Even remaining within the realm of Feller kernels on Polish spaces, a synthetic framework may shed conceptual light on an awkward aspect of Section~\ref{sec:Polish}: there are multiple definitions of ``pullback'' and ``pushforward'', each of which only satisfies functoriality and adjointness in a \emph{lax} sense (see Section~\ref{subsec:lax-axioms} and Appendix~\ref{app:lax_functors}). It may also enable us to interpret the ``Wasserstein pushforward'' (in Section~\ref{sec:wasserstein}) in terms of the Giry monad (see \cite{Giry}).

It is of course not possible to eliminate the engineering details entirely---the topology of the asset space, the regularity of the kernel, the choice of distance metric all matter. But the synthetic approach provides a way to encapsulate them within the definition of the support functor, enabling modular design of investment processes: the solver can be swapped without rewriting the audit logic.

\clearpage

\vspace*{\fill}

\begin{center}
\begin{tikzpicture}[
    scale=0.8,
    transform shape,
    node distance=1.6cm and 0.8cm,
    bus_input/.style={
        rectangle,
        rounded corners=8pt,
        draw=blue!50!black,
        fill=blue!5,
        thick,
        font=\sffamily,
        align=center,
        minimum width=3.2cm,
        minimum height=1.2cm
    },
    sketch/.style={
        decorate,
        decoration={random steps, segment length=3pt, amplitude=0.4pt},
        thick,
        draw=black,
        fill=white,
        line cap=round,
        font=\sffamily\bfseries,
        align=center,
        inner sep=6pt,
        minimum width=3.5cm,
        drop shadow
    },
    logic_step/.style={
        font=\sffamily\scriptsize,
        align=center,
        fill=white,
        inner sep=2pt,
        text=gray!40!black
    },
    engine_box/.style={
        draw=gray!40,
        dashed,
        thick,
        rounded corners=10pt,
        fill=gray!2
    },
    outcome/.style={
        rectangle,
        thick,
        draw=black,
        fill=white,
        font=\sffamily\bfseries,
        align=center,
        minimum width=3.5cm,
        minimum height=1.6cm,
        drop shadow
    },
    arr_main/.style={->, >=stealth, very thick, black},
    arr_logic/.style={->, >=stealth, thick, gray},
    arr_derive/.style={->, >=stealth, double, thick, gray!80}
]

    \node[bus_input] (Specs) {
        \textbf{Required Alignments}\\
        \scriptsize{Constraints, Guidance}
    };
    
    \node[bus_input, right=3.5cm of Specs] (Logic) {
        \textbf{Transition Logic}\\
        \scriptsize{Optimizers, Aggregators}
    };

    \node[sketch, below=1.0cm of Specs] (Verticals) {
        \textbf{Vertical Spans}\\
        \footnotesize{$A \leftarrow E \rightarrow B$}
    };
    
    \node[sketch, below=1.0cm of Logic] (Horizontals) {
        \textbf{Horizontal Maps}\\
        \footnotesize{Functions $f: A \to B$}
    };

    \draw[arr_main] (Specs) -- node[left, logic_step] {Defines} (Verticals);
    \draw[arr_main] (Logic) -- node[right, logic_step] {Defines} (Horizontals);

    
    \node[logic_step, below=1.2cm of Verticals] (Evidence) {Evidence Space\\(Quantitative Data)};
    \node[logic_step, below=1.2cm of Horizontals] (Companions) {Companions $\langle f \rangle$\\ (Graph Embedding)};
    
    \node[sketch, rectangle, minimum width=5cm, below=2.5cm of $(Verticals)!0.5!(Horizontals)$] (Inference) {
        \textbf{2-Cell Inference ($\Downarrow \alpha$)}\\
        \footnotesize{Evidence Transformation}\\
        \footnotesize{$\alpha: \text{Hub Cert} \Rightarrow \text{Spoke Cert}$}
    };
    
    \node[sketch, rectangle, minimum width=6cm, below=1.0cm of Inference] (Calculus) {
        \textbf{Operational Calculus}\\
        \footnotesize{Syntactic manipulation of}\\
        \footnotesize{Audit Trails}
    };
    
    \begin{scope}[on background layer]
        \node[engine_box, fit=(Evidence) (Companions) (Inference) (Calculus), label={[gray, font=\sffamily\bfseries]above:Abstract Equipment Syntax}] (Engine) {};
    \end{scope}

    \draw[arr_logic] (Verticals) -- (Evidence);
    \draw[arr_logic] (Horizontals) -- (Companions);
    \draw[arr_derive] (Evidence) -- (Inference);
    \draw[arr_derive] (Companions) -- (Inference);
    \draw[arr_main] (Inference) -- node[right, logic_step] {Generates} (Calculus);

    
    \node[outcome, below left=1.8cm and -1.0cm of Calculus] (Audit) {
        Quantitative\\Audit Trails\\
        \scriptsize{Full Provenance}\\
        \scriptsize{(Composite Spans)}
    };
    
    \node[outcome, below=1.8cm of Calculus] (Dual) {
        Primal / Dual\\
        Operations\\
        \scriptsize{Construction vs. Check}\\
        \scriptsize{($f_!$ vs $f^*$)}
    };
    
    \node[outcome, below right=1.8cm and -1.0cm of Calculus] (Safe) {
        Guaranteed\\
        Safety\\
        \scriptsize{Conservative Verification}\\
        \scriptsize{(Lax Beck--Chevalley)}
    };

    \draw[arr_main] (Calculus.south) -- (Audit.north);
    \draw[arr_main] (Calculus.south) -- (Dual.north);
    \draw[arr_main] (Calculus.south) -- (Safe.north);

    \node[logic_step, fill=white, inner sep=2pt, rotate=25] at ($(Calculus.south)!0.5!(Safe.north)$) {Syntactic\\Guarantee};

\end{tikzpicture}
\end{center}

\vspace*{\fill}

\clearpage

\enlargethispage{2\baselineskip}
\bibliographystyle{plain}
\bibliography{hub-spoke}
\addcontentsline{toc}{section}{References}

\clearpage

\appendix

\part*{Prerequisites and System Architecture}

\addcontentsline{toc}{part}{Prerequisites and System Architecture}

\clearpage

\section*{Mathematical Dependencies}
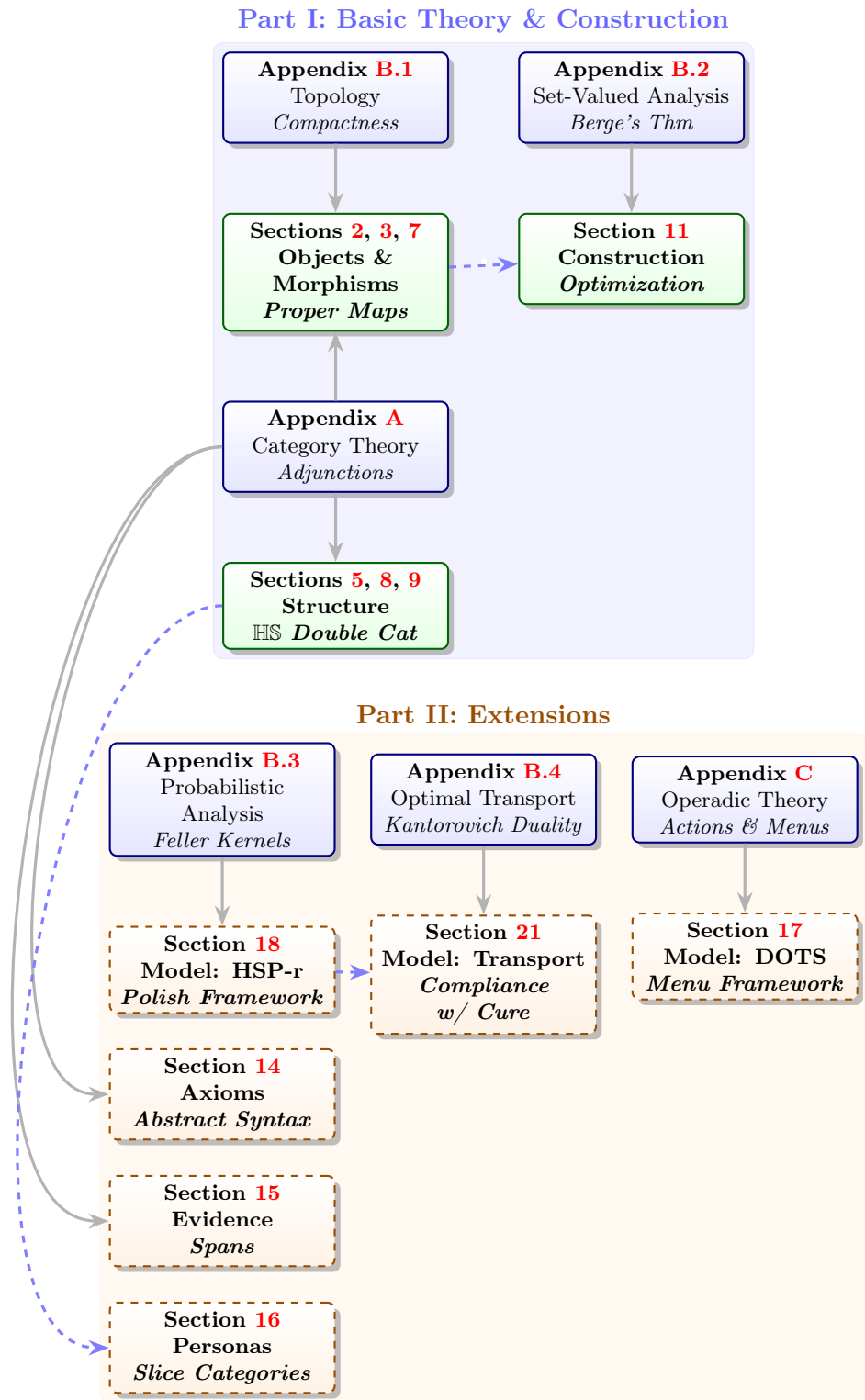
\addcontentsline{toc}{section}{Mathematical Dependencies}

\vspace{1cm}

\begin{center}
\begin{tikzpicture}[
    node distance=1.0cm and 0.5cm, 
    >={Stealth}, 
    math_tool/.style={
        rectangle, 
        draw=blue!50!black, 
        thick,
        top color=white, 
        bottom color=blue!10,
        rounded corners=3pt, 
        minimum width=3.2cm, 
        text width=3.0cm, 
        minimum height=1.0cm,
        align=center,
        drop shadow,
        font=\footnotesize
    },
    paper_sec/.style={
        rectangle, 
        draw=green!40!black, 
        thick,
        top color=white, 
        bottom color=green!10,
        rounded corners=3pt, 
        minimum width=3.2cm, 
        text width=3.0cm, 
        minimum height=1.0cm,
        align=center,
        drop shadow,
        font=\footnotesize\bfseries
    },
    extension_sec/.style={
        rectangle, 
        draw=orange!60!black, 
        thick,
        top color=white, 
        bottom color=orange!10,
        dashed,
        rounded corners=3pt, 
        minimum width=3.2cm, 
        text width=3.0cm, 
        minimum height=1.0cm,
        align=center,
        drop shadow,
        font=\footnotesize\bfseries
    },
    arrow/.style={->, very thick, draw=gray!60}, 
    gray_arrow/.style={->, very thick, draw=gray!60},
    label_text/.style={font=\tiny\itshape, text=blue!60!black, fill=white, inner sep=1pt}
]

    
    \node[math_tool] (Topo) {
        \textbf{Appendix~\ref{sub:topology}}\\
        Topology\\
        \textit{Compactness}
    };

    \node[paper_sec, below=1.0cm of Topo] (Sec_Basic) {
        \textbf{Sections \ref{sec:spaces}, \ref{sec:reimplement}, \ref{sec:proper}}\\
        Objects \& Morphisms\\
        \textit{Proper Maps}
    };

    \node[math_tool, below=1.0cm of Sec_Basic] (Cats) {
        \textbf{Appendix~\ref{app:categorical_prereqs}}\\
        Category Theory\\
        \textit{Adjunctions}
    };

    \node[paper_sec, below=1.0cm of Cats] (Sec_Struct) {
        \textbf{Sections \ref{sec:double}, \ref{sec:Beck--Chevalley}, \ref{sec:frobenius}}\\
        Structure\\
        \textit{$\mathbb{HS}$ Double Cat}
    };

    \node[math_tool, right=1.0cm of Topo] (SetValued) {
        \textbf{Appendix~\ref{sub:set-valued}}\\
        Set-Valued Analysis\\
        \textit{Berge's Thm}
    };

    \node[paper_sec, below=1.0cm of SetValued] (Sec_Opt) {
        \textbf{Section~\ref{sec:optimal-reimpl}}\\
        Construction\\
        \textit{Optimization}
    };


    \node[math_tool, below=4.0cm of $(Sec_Struct.south)!0.5!(Sec_Opt.south)$] (OptTrans) {
        \textbf{Appendix~\ref{app:optimal_transport}}\\
        Optimal Transport\\
        \textit{Kantorovich Duality}
    };

    \node[math_tool, left=0.5cm of OptTrans] (Stoch) {
        \textbf{Appendix~\ref{subsec:stochastic}}\\
        Probabilistic Analysis\\
        \textit{Feller Kernels}
    };

    \node[math_tool, right=0.5cm of OptTrans] (DOTS_Math) {
        \textbf{Appendix~\ref{app:dots_framework}}\\
        Operadic Theory\\
        \textit{Actions \& Menus}
    };

    \node[extension_sec, below=1.0cm of Stoch] (Ext_Polish) {
        \textbf{Section~\ref{sec:Polish}}\\
        Model: HSP-r\\
        \textit{Polish Framework}
    };

    \node[extension_sec, below=1.0cm of OptTrans] (Ext_Wass) {
        \textbf{Section~\ref{sec:wasserstein}}\\
        Model: Transport\\
        \textit{Compliance w/ Cure}
    };

    \node[extension_sec, below=1.0cm of DOTS_Math] (Ext_DOTS) {
        \textbf{Section~\ref{sec:DOTS}}\\
        Model: DOTS\\
        \textit{Menu Framework}
    };

    \node[extension_sec, below=0.5cm of Ext_Polish] (Sec_Axioms) {
        \textbf{Section~\ref{sec:abstract}}\\
        Axioms\\
        \textit{Abstract Syntax}
    };

    \node[extension_sec, below=0.5cm of Sec_Axioms] (Sec_Evidence) {
        \textbf{Section~\ref{sec:2-cells}}\\
        Evidence\\
        \textit{Spans}
    };

    \node[extension_sec, below=0.5cm of Sec_Evidence] (Sec_Personas) {
        \textbf{Section~\ref{sec:personas}}\\
        Personas\\
        \textit{Slice Categories}
    };


    \draw[arrow] (Topo) -- (Sec_Basic);
    \draw[arrow] (Cats) -- (Sec_Basic);
    \draw[arrow] (Cats) -- (Sec_Struct);
    \draw[arrow] (SetValued) -- (Sec_Opt);
    \draw[arrow, dashed, blue!50] (Sec_Basic) -- node[label_text, above] {} (Sec_Opt);

    \draw[arrow] (Stoch) -- (Ext_Polish);       
    \draw[gray_arrow] (OptTrans) -- (Ext_Wass); 
    \draw[arrow, dashed, blue!50] (Ext_Polish) -- (Ext_Wass); 
    \draw[arrow] (DOTS_Math) -- (Ext_DOTS);

    \draw[arrow] (Cats.west) to[out=180, in=180, looseness=0.6] (Sec_Axioms.west);
    \draw[arrow] (Cats.west) to[out=180, in=180, looseness=0.6] (Sec_Evidence.west);
    
    \draw[arrow, dashed, blue!50] (Sec_Struct.west) to[out=180, in=180, looseness=0.6] (Sec_Personas.west);


    \begin{scope}[on background layer]
        \node[fit=(Topo)(Cats)(Sec_Struct)(SetValued)(Sec_Opt), draw=blue!10, fill=blue!5, rounded corners, label={[blue!60, font=\bfseries]above:Part I: Basic Theory \& Construction}] {};

        \node[fit=(Stoch)(OptTrans)(DOTS_Math)(Ext_Wass)(Sec_Personas)(Ext_DOTS), draw=orange!10, fill=orange!5, rounded corners, label={[orange!60!black, font=\bfseries]above:Part II: Extensions}] {};
    \end{scope}

\end{tikzpicture}
\end{center}

\clearpage

\section{Categorical Prerequisites}
\label{app:categorical_prereqs}

This appendix provides the minimal category-theoretic background required to interpret the Hub-and-Spoke ($\HSSimp$) framework presented in Part I. We restrict our attention to the specific structures used in the main text: categories, functors, adjunctions and thin double categories (also known as framed bicategories or equipments).

\subsection{Categories and Functors}

\begin{definition}[Category]
A \textbf{category} $\mathcal{C}$ consists of:
\begin{enumerate}
    \item A collection of \textbf{objects}, denoted $A, B, C, \dots$
    \item For every pair of objects $A, B$, a set of \textbf{morphisms} (or arrows) $\text{Hom}_{\mathcal{C}}(A,B)$, denoted $f: A \to B$.
    \item A composition operation $\circ$: for $f: A \to B$ and $g: B \to C$, there exists a composite $g \circ f: A \to C$.
\end{enumerate}
Composition must be associative, and every object $A$ must possess an identity morphism $id_A$ such that $f \circ id_A = f$ and $id_B \circ f = f$.
\end{definition}

\begin{definition}[Functor]
A \textbf{functor} $F: \mathcal{C} \to \mathcal{D}$ is a structure-preserving map between categories. It assigns to every object $A \in \mathcal{C}$ an object $F(A) \in \mathcal{D}$, and to every morphism $f: A \to B$ a morphism $F(f): F(A) \to F(B)$, such that:
\begin{itemize}
    \item $F(id_A) = id_{F(A)}$
    \item $F(g \circ f) = F(g) \circ F(f)$
\end{itemize}
\end{definition}

An \textbf{adjunction} is widely considered the most important concept in category theory. It describes a specific relationship between two functors that act as ``best approximations'' or optimal solutions to each other's problems across different mathematical contexts.

\begin{definition}[Adjunction]
Let $\mathcal{C}$ and $\mathcal{D}$ be categories. Two functors $L: \mathcal{C} \to \mathcal{D}$ (the \textbf{Left Adjoint}) and $R: \mathcal{D} \to \mathcal{C}$ (the \textbf{Right Adjoint}) form an adjunction, denoted $L \dashv R$, if there is a natural isomorphism between their sets of morphisms (Hom-sets):
\[
\text{Hom}_{\mathcal{D}}(L(X), Y) \cong \text{Hom}_{\mathcal{C}}(X, R(Y))
\]
for all objects $X$ in $\mathcal{C}$ and $Y$ in $\mathcal{D}$.
\end{definition}

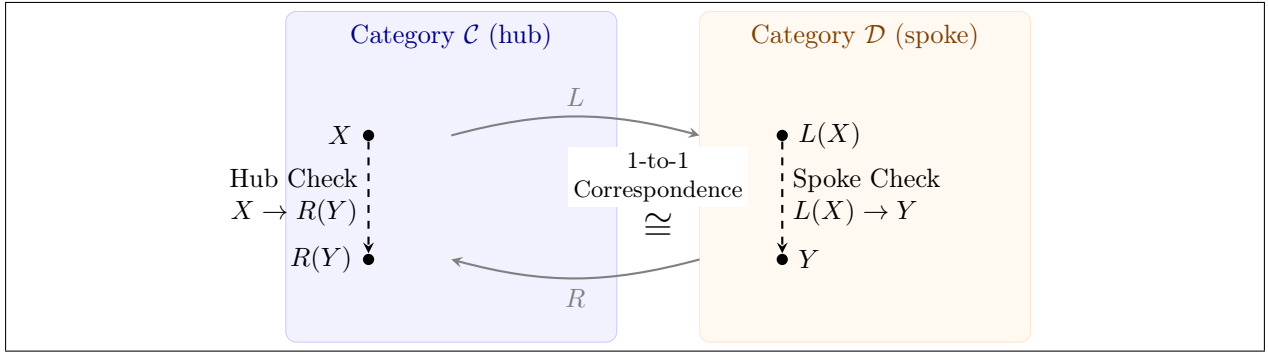
\begin{figure}[htbp]
\centering
\begin{tikzpicture}[>=stealth, scale=1.1, font=\small]

    \draw[fill=blue!5, draw=blue!20, rounded corners] (-4.5,-2) rectangle (-0.5,2);
    \node[blue!50!black] at (-2.5, 1.7) {Category $\mathcal{C}$ (hub)};

    \node[circle, fill=black, inner sep=1.5pt, label=left:$X$] (X) at (-3.5, 0.5) {};
    \node[circle, fill=black, inner sep=1.5pt, label=left:$R(Y)$] (RY) at (-3.5, -1.0) {};

    \draw[->, dashed, thick] (X) -- node[left, align=right] {Hub Check\\$X \to R(Y)$} (RY);

    \draw[fill=orange!5, draw=orange!20, rounded corners] (0.5,-2) rectangle (4.5,2);
    \node[orange!50!black] at (2.5, 1.7) {Category $\mathcal{D}$ (spoke)};

    \node[circle, fill=black, inner sep=1.5pt, label=right:$L(X)$] (LX) at (1.5, 0.5) {};
    \node[circle, fill=black, inner sep=1.5pt, label=right:$Y$] (Y) at (1.5, -1.0) {};

    \draw[->, dashed, thick] (LX) -- node[right, align=left] {Spoke Check\\$L(X) \to Y$} (Y);

    \draw[->, gray, thick] (-2.5, 0.5) to[bend left=15] node[midway, above] {$L$} (0.5, 0.5);

    \draw[->, gray, thick] (0.5, -1.0) to[bend left=15] node[midway, below] {$R$} (-2.5, -1.0);

    \node[font=\Large] at (0, -0.6) {$\cong$};
    \node[font=\footnotesize, align=center, fill=white, inner sep=2pt] at (0, 0) {1-to-1\\Correspondence};

\end{tikzpicture}
\caption{Visualizing an adjunction $L \dashv R$.}
\label{fig:adjunction_iso}
\end{figure}

\noindent \textbf{Interpretation:}
This isomorphism asserts a one-to-one correspondence between:
\begin{itemize}
    \item Morphisms from the approximation $L(X)$ to the spoke $Y$ in category $\mathcal{D}$.
    \item Morphisms from the hub $X$ to the approximation $R(Y)$ in category $\mathcal{C}$.
\end{itemize}

\noindent \textbf{Special Case: Galois Connections}
In the specific context of a \textbf{thin} double category (where categories are posets and morphisms are inequalities $\le$), an adjunction simplifies to a \textbf{Galois Connection}. The isomorphism condition becomes a pair of logical implications:
\[
L(x) \le y \iff x \le R(y)
\]
This is the form used throughout the $\HSSimp$ framework to relate alignment checks in different portfolio spaces.

\subsection{Double Categories}

A double category generalizes a category by allowing two distinct types of morphisms \cite{Benabou, GrandisPare}. This is essential for distinguishing between \textit{functional} transformations (re-implementations) and \textit{relational} constraints (alignments).

\begin{definition}[Double Category]
A \textbf{double category} $\mathbb{D}$ consists of:
\begin{enumerate}
    \item \textbf{Objects:} $A, B, C, \dots$
    \item \textbf{Horizontal Morphisms:} $f: A \to B$ (forming a category $\mathbb{D}_h$).
    \item \textbf{Vertical Morphisms:} $R: A \nrightarrow B$ (forming a category $\mathbb{D}_v$, often with weaker associativity).
    \item \textbf{2-Cells:} Squares of the form:
    \[
    \begin{tikzcd}
        A \arrow[r, "f"] \arrow[d, "R"'] & B \arrow[d, "S"] \\
        C \arrow[r, "g"'] & D
    \end{tikzcd}
    \]
    occupied by a 2-morphism $\alpha$.
\end{enumerate}
\end{definition}

The complexity of double categories is drastically reduced in our context because relations are either true or false.

\begin{definition}[Thin Double Category]
A double category is \textbf{thin} (or locally posetal) if, for any boundary of horizontal and vertical morphisms, there is at most one 2-cell filling the square. This structure is closely related to the concept of a \emph{framed bicategory} or \emph{equipment} \cite{Shulman}.
\end{definition}

\subsection{Lax Functors and Ordered Structures}
\label{app:lax_functors}

The probabilistic extensions in Part III (specifically Section~\ref{subsec:stochastic}) require a relaxation of standard category theory. When ``equality'' is replaced by ``approximation'' or ``inclusion,'' functors become \textbf{lax}.

\begin{definition}[Ordered Category]
An \textbf{ordered category} (or 2-poset) is a category $\mathcal{C}$ where each hom-set $\text{Hom}(A,B)$ is a partially ordered set $(P, \le)$, and composition is order-preserving. In our context, the vertical morphisms form ordered categories where the order represents inclusion of relations or risk profiles.
\end{definition}

\begin{definition}[Lax Functor]
Let $\mathcal{C}$ be a standard category and $\mathcal{D}$ be an ordered category. A \textbf{Lax Functor}\footnote{Our convention follows \cite{Shulman}: a lax functor satisfies $F(g \circ f) \le F(g) \circ F(f)$. Some authors call this \emph{oplax}; the choice depends on the direction of the ordering. In our setting, where $\le$ denotes set-theoretic inclusion of relations (or containment of risk envelopes), the ``lax'' label matches the interpretation that the composite pushforward is contained in the sequential one.} $F: \mathcal{C} \to \mathcal{D}$ assigns objects and morphisms such that:
\begin{enumerate}
    \item $F(id_A) \le id_{F(A)}$ (Sub-unitality)
    \item $F(g \circ f) \le F(g) \circ F(f)$ (Sub-compositionality)
\end{enumerate}
\end{definition}

\noindent \textbf{Relevance to Portfolio Construction:}
In the HSP-r model, the mapping that takes a re-implementation map $f$ to its pushforward operator $f_!$ is a Lax Functor.
\[
(Q \circ P)_!^\epsilon \le Q_!^\epsilon \circ P_!^\epsilon
\]
This inequality captures the \textbf{conservatism of sequential auditing}. The risk envelope of the end-to-end composite process $(Q \circ P)$ is \emph{contained in} the risk envelope obtained by composing the per-stage envelopes (the RHS). Equivalently, the sequential audit (checking each stage separately) is \emph{more conservative} than the composite audit: passing the sequential check implies passing the composite check, but not conversely. This is precisely the safety guarantee we need---the sequential audit may reject valid portfolios (type~I error) but never approves invalid ones (type~II error).

\clearpage  \section{Topology, Set-valued Analysis and Probability}
\label{app:analytical_foundations}

This appendix summarizes the background in general topology, set-valued analysis, and measure theory required to establish the well-definedness of the $\HSSimp$ objects and the convergence properties of the HSP-r (stochastic) extension.

\subsection{Topology and Proper Maps}
\label{sub:topology}

\begin{definition}[Compactness]
A topological space $X$ is \textbf{compact} if every open cover of $X$ admits a finite subcover. By the Heine--Borel theorem, a subset of Euclidean space $\mathbb{R}^n$ is compact if and only if it is closed and bounded.
\end{definition}

\begin{definition}[Proper Map]
A continuous map $f: X \to Y$ is \textbf{proper} if the preimage $f^{-1}(K)$ of every compact subset $K \subseteq Y$ is compact in $X$.
\end{definition}

\begin{theorem}[Closed Map Theorem for Proper Maps]
If $f: X \to Y$ is a proper map between locally compact Hausdorff spaces, then $f$ is a \textbf{closed map}: the image of any closed set in $X$ is closed in $Y$.
\end{theorem}

\subsection{Set-Valued Analysis and Optimization}
\label{sub:set-valued}

The study of correspondences and their continuity properties is central to optimization theory \cite{AubinFrankowska, RockafellarWets}.

\begin{definition}[Correspondence]
A \textbf{correspondence} (or multifunction) $\Phi: X \rightrightarrows Y$ assigns to every $x \in X$ a subset $\Phi(x) \subseteq Y$.
\end{definition}

\begin{definition}[Upper Hemicontinuity]
A correspondence $\Phi$ is \textbf{upper hemicontinuous} (UHC) at $x_0$ if, for every open set $V$ containing $\Phi(x_0)$, there exists a neighborhood $U$ of $x_0$ such that $\Phi(x) \subseteq V$ for all $x \in U$.
\end{definition}

\begin{theorem}[Berge's Maximum Theorem]
Let $X$ and $Y$ be topological spaces, $f: X \times Y \to \mathbb{R}$ a continuous objective function, and $C: X \rightrightarrows Y$ a correspondence of feasible sets that is continuous (both upper and lower hemicontinuous) with compact values. Define the value function $v(x) = \min_{y \in C(x)} f(x,y)$ and the solution correspondence $M(x) = \{y \in C(x) \mid f(x,y) = v(x)\}$. Then:
\begin{enumerate}
    \item The value function $v$ is continuous.
    \item The solution correspondence $M$ is upper hemicontinuous and compact-valued.
\end{enumerate}
See \cite{Berge} or \cite{AliprantisBorder} for detailed proofs.  For variable constraint sets, the analogous result requires the correspondence $x \mapsto C(x)$ to be continuous (upper and lower hemicontinuous); see Michael's Selection Theorem (Theorem~\ref{thm:michael-selection}) for the case where only a continuous selection, rather than continuity of the full solution set, is needed.
\end{theorem}

\begin{figure}[htbp]
\centering
\begin{tikzpicture}[
    >=Stealth,
    thick,
    font=\small, 
    space_style/.style={draw=gray!60, fill=gray!5, rounded corners=15pt, minimum height=5cm, minimum width=4cm},
    setU_style/.style={ellipse, draw=blue!80, dashed, fill=blue!10, minimum width=2.5cm, minimum height=1.5cm},
    setV_style/.style={ellipse, draw=orange!80, dashed, fill=orange!10, minimum width=3.2cm, minimum height=2.8cm},
    setCU_style/.style={ellipse, draw=red!70!blue!50, fill=red!30!blue!10, opacity=0.8, minimum width=2.2cm, minimum height=1.6cm},
    setCx_style/.style={rounded corners=4pt, draw=red!90, fill=red!50, minimum width=0.8cm, minimum height=0.5cm},
    point_style/.style={circle, fill=black, inner sep=1.5pt}
]

    
    \node[space_style] (X_space) at (0,0) {};
    \node[anchor=north] at (X_space.north) {\textbf{Hub Space} $X$};

    \node[space_style, right=2cm of X_space] (Y_space) {};
    \node[anchor=north] at (Y_space.north) {\textbf{Spoke Space} $Y$};

    \coordinate (x_center) at (X_space.center);
    
    \node[setU_style] (U_set) at (x_center) {};
    \node[blue!80, anchor=south east] at (U_set.north west) {$U$};
    
    \node[point_style, label={left:$x$}] (x_point) at (x_center) {};

    \coordinate (y_center) at (Y_space.center);

    \node[setV_style, rotate=-10] (V_set) at (y_center) {};
    \node[orange!80, anchor=south west] at (V_set.north east) {$V$};

    \node[setCU_style, rotate=-10] (CU_set) at (y_center) {};
    \node[red!60!blue!60, font=\footnotesize, anchor=south] at (CU_set.north) {$C(U)$};

    \node[setCx_style, rotate=-10] (Cx_set) at (y_center) {};
    \node[red!90, font=\scriptsize] at (Cx_set.center) {$C(x)$};

    \draw[->, ultra thick, gray!40, shorten >= 0.5cm, shorten <= 0.5cm] 
        (X_space.east) -- (Y_space.west) 
        node[midway, above, black] {$C$};

    \draw[->, red!80, dashed, shorten >= 2pt, shorten <= 2pt] 
        (x_point) to[out=20, in=160] (Cx_set);

    \draw[->, blue!50, dotted, shorten >= 2pt, shorten <= 2pt] 
        (U_set.north) to[out=30, in=150] (CU_set.north);

    \coordinate (middle_point) at ($(X_space.south)!0.5!(Y_space.south)$);

    \node[
        below=0.5cm of middle_point,
        draw=black!10, 
        fill=yellow!5, 
        rounded corners, 
        text width=9cm, 
        align=center,
        inner sep=10pt
    ] (def_box) {
        \textbf{Definition: Upper Hemicontinuity}\\
        A correspondence $C: X \rightrightarrows Y$ is upper hemicontinuous at $x$ if for any open neighborhood $V$ of the image set $C(x)$ (i.e., $C(x) \subset V$), there exists an open neighborhood $U$ of $x$ such that the image of $U$ is contained in $V$.
    };

\end{tikzpicture}
\caption{Visualizing Upper Hemicontinuity.}
\label{fig:uhc_stability}
\end{figure}
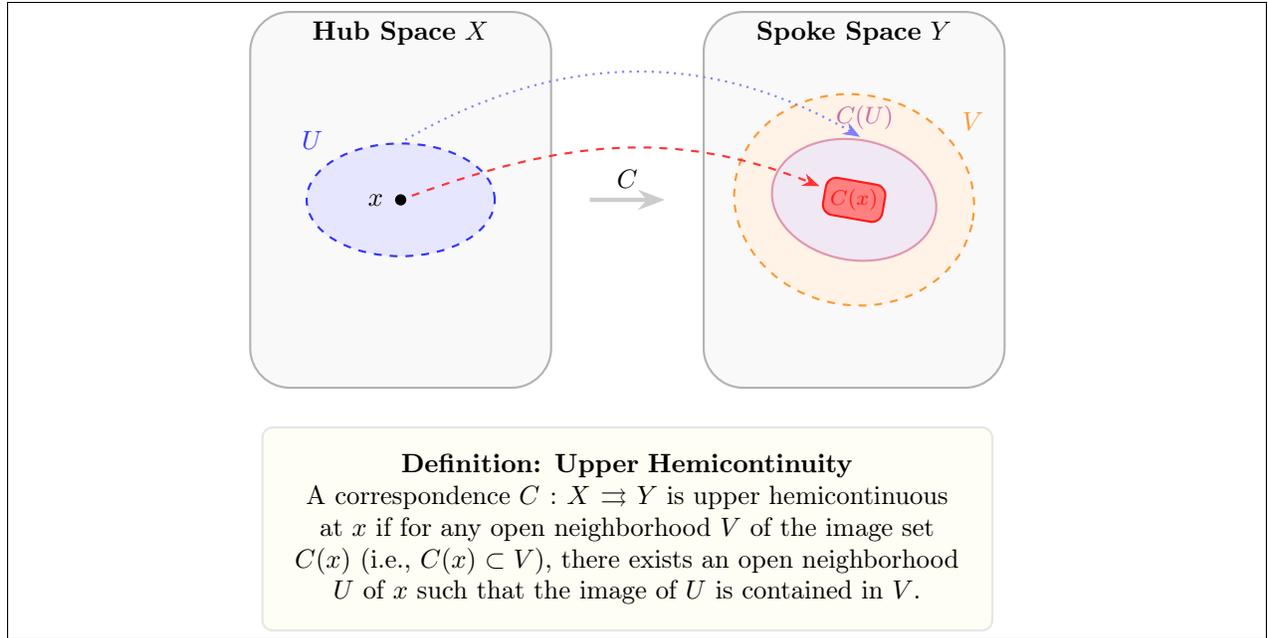

\subsection{Polish Spaces and Stochastic Kernels}
\label{subsec:stochastic}

For the probabilistic extensions of the framework, we rely on the properties of Polish spaces \cite{Billingsley, Bogachev, Kallenberg}.

\begin{definition}[Polish Space]
A \textbf{Polish space} is a separable, completely metrizable topological space. Examples include $\mathbb{R}^n$, separable Banach spaces, and the space of probability measures on a Polish space equipped with the weak topology \cite{AliprantisBorder}.
\end{definition}

\begin{definition}[Markov Kernel]
A \textbf{Markov kernel} (or probability kernel) from $X$ to $Y$ is a function $P: X \times \mathcal{B}(Y) \to [0,1]$ such that:
\begin{enumerate}
    \item For every fixed $x \in X$, the map $A \mapsto P(x, A)$ is a probability measure on $Y$.
    \item For every fixed measurable set $A$, the function $x \mapsto P(x, A)$ is measurable.
\end{enumerate}
\end{definition}

\begin{definition}[Feller Property]
A Markov kernel $P$ is \textbf{Feller} (or continuous) if, for every bounded continuous function $g: Y \to \mathbb{R}$, the function $x \mapsto \int_Y g(y) P{(x, dy)}$ is continuous on $X$.
\end{definition}

\begin{theorem}[Portmanteau Theorem]
\label{thm:portmanteau}
Let $Y$ be a Polish space and let $\mu_n, \mu$ be Borel probability measures on $Y$. The following conditions are equivalent to the weak convergence $\mu_n \Rightarrow \mu$:
\begin{enumerate}
    \item $\lim_{n \to \infty} \int_Y f \, d\mu_n = \int_Y f \, d\mu$ for all bounded continuous functions $f: Y \to \mathbb{R}$.
    \item $\limsup_{n \to \infty} \mu_n(F) \leq \mu(F)$ for all closed sets $F \subseteq Y$.
    \item $\liminf_{n \to \infty} \mu_n(G) \geq \mu(G)$ for all open sets $G \subseteq Y$.
\end{enumerate}
\end{theorem}

\begin{corollary}[Upper Semi-Continuity of Mass on Closed Sets]
\label{cor:USC}
Let $P: X \rightsquigarrow Y$ be a Feller kernel. For any closed set $F \subseteq Y$, the function $\phi(x) = P(x, F)$ is upper semi-continuous on $X$.

\begin{proof}
Since $P$ is Feller, if $x_n \to x$ in $X$, then $P(x_n, \cdot) \Rightarrow P(x, \cdot)$ weakly. By condition (2) of the Portmanteau Theorem, $\limsup_{n \to \infty} P(x_n, F) \leq P(x, F)$. This is the definition of upper semi-continuity. This property guarantees that the stochastic pullback of a closed relation is closed.
\end{proof}
\end{corollary}

\begin{definition}[Tight Feller Kernel]
\label{def:tight}
A Feller kernel $P: X \times \mathcal{B}(Y) \to [0,1]$ is said to be \textbf{tight} if for every compact subset $C \subseteq X$ and every $\epsilon > 0$, there exists a compact subset $P_{\epsilon, C} \subseteq Y$ such that
\[
P(x, P_{\epsilon, C}) \ge 1 - \epsilon \quad \text{for all } x \in C.
\]
\end{definition}

This condition ensures that the kernel does not ``leak mass to infinity'' when the input is constrained to a compact region, serving as the stochastic analogue to the properness of a continuous map. For the existence of measurable selectors in this context, see \cite{CastaingValadier}.

\begin{theorem}[Prokhorov's Theorem]
Let $Y$ be a Polish space and let $\mathcal{P}(Y)$ denote the space of Borel probability measures on $Y$ equipped with the topology of weak convergence. A subset $\Gamma \subseteq \mathcal{P}(Y)$ is relatively compact in the weak topology if and only if it is tight: see \cite{Billingsley}.
In the context of Feller kernels, this implies that if $P$ is a tight Feller kernel, then for any compact set $C \subseteq X$, the set of measures $\{P(x, \cdot) \mid x \in C\}$ is relatively compact in $\mathcal{P}(Y)$.
\end{theorem}

\begin{figure}[htbp]
    \centering
    \begin{tikzpicture}[>=latex, scale=0.85]

        \begin{scope}[local bounding box=LeftPanel]
            \node[font=\bfseries\large] at (3.5, 5.6) {Tight Kernel};
            \node[font=\small, text=gray] at (3.5, 4.6) {Input $x \in C \Rightarrow$ Mass bounded in $P_{\epsilon, C}$};

            \draw[->, thick] (-0.5,0) -- (7,0) node[right] {$Y$}; 
            \draw[->, thick] (0,-0.5) -- (0,3.5) node[above] {Density};

            \draw[very thick, orange] (2,-0.8) -- (4,-0.8);
            \draw[orange] (2, -0.9) -- (2, -0.7);
            \draw[orange] (4, -0.9) -- (4, -0.7);
            \node[orange, font=\bfseries] at (3, -1.2) {Input $C$};

            \draw[<->, blue!80!black] (1, 2.8) -- (6, 2.8) node[midway, fill=white] {Set $P_{\epsilon, C}$};
            \draw[dashed, blue!80!black] (1,0) -- (1, 2.8);
            \draw[dashed, blue!80!black] (6,0) -- (6, 2.8);

            \draw[thick, orange!80!black, smooth] plot[domain=-0.2:6.8, samples=40] (\x, {1.5*exp(-(\x-2.5)^2/0.4)});
            \draw[thick, orange!80!black, smooth] plot[domain=-0.2:6.8, samples=40] (\x, {1.5*exp(-(\x-4.5)^2/0.4)});
            \draw[thick, orange!80!black, smooth] plot[domain=-0.2:6.8, samples=40] (\x, {1.5*exp(-(\x-3.5)^2/0.4)});

            \draw[->, orange, dotted, thick] (2.2,-0.7) -- node[midway, sloped, above, font=\scriptsize] {Action of $P$} (2.5, 0.8);
            \draw[->, orange, dotted, thick] (3.8,-0.7) -- (4.2, 0.8);

            \node[blue!80!black, font=\footnotesize] at (6.5, 0.5) {Mass $< \epsilon$};
            \node[blue!80!black, font=\footnotesize] at (0.5, 0.5) {Mass $< \epsilon$};
        \end{scope}

        \begin{scope}[xshift=9cm, local bounding box=RightPanel]
            \node[font=\bfseries\large] at (3.5, 5.6) {Not a Tight Kernel};
            \node[font=\small, text=red!80!black] at (3.5, 4.6) {For any fixed $K$, mass escapes};

            \draw[->, thick] (-0.5,0) -- (7,0) node[right] {$Y$};
            \draw[->, thick] (0,-0.5) -- (0,3.5);

            \draw[very thick, orange] (2,-0.8) -- (4,-0.8);
            \draw[orange] (2, -0.9) -- (2, -0.7);
            \draw[orange] (4, -0.9) -- (4, -0.7);
            \node[orange, font=\bfseries] at (3, -1.2) {Input $C$};

            \draw[<->, blue!80!black] (1, 2.8) -- (6, 2.8) node[midway, fill=white] {Any fixed compact $K$};
            \draw[dashed, blue!80!black] (1,0) -- (1, 2.8);
            \draw[dashed, blue!80!black] (6,0) -- (6, 2.8);

            \draw[thick, gray, dashed, smooth] plot[domain=-0.2:6.8, samples=40] (\x, {1.5*exp(-(\x-2)^2/0.4)});
            \draw[thick, gray, dashed, smooth] plot[domain=-0.2:6.8, samples=40] (\x, {1.5*exp(-(\x-4)^2/0.4)});
            \draw[thick, red, smooth] plot[domain=-0.2:6.8, samples=40] (\x, {1.5*exp(-(\x-6.5)^2/0.4)});
            \node[red, font=\footnotesize] at (7, 1.8) {Mass escapes};

            \draw[->, orange, dotted, thick] (2.2,-0.7) -- (2.2, 0.8);
            \draw[->, orange, dotted, thick] (3.8,-0.7) -- node[midway, sloped, above, font=\scriptsize] {Action of $P$} (6.2, 0.8);
        \end{scope}

    \end{tikzpicture}
    \caption{Visualizing the definition of tight Feller kernel. Tightness means that for any compact input $C$, there exists a compact subset $P_{\epsilon, C}\subseteq Y$ containing most of the probability mass. The right panel shows failure: inputs stay in $C$, but mass drifts to infinity.}
    \label{fig:tightness_viz}
\end{figure}
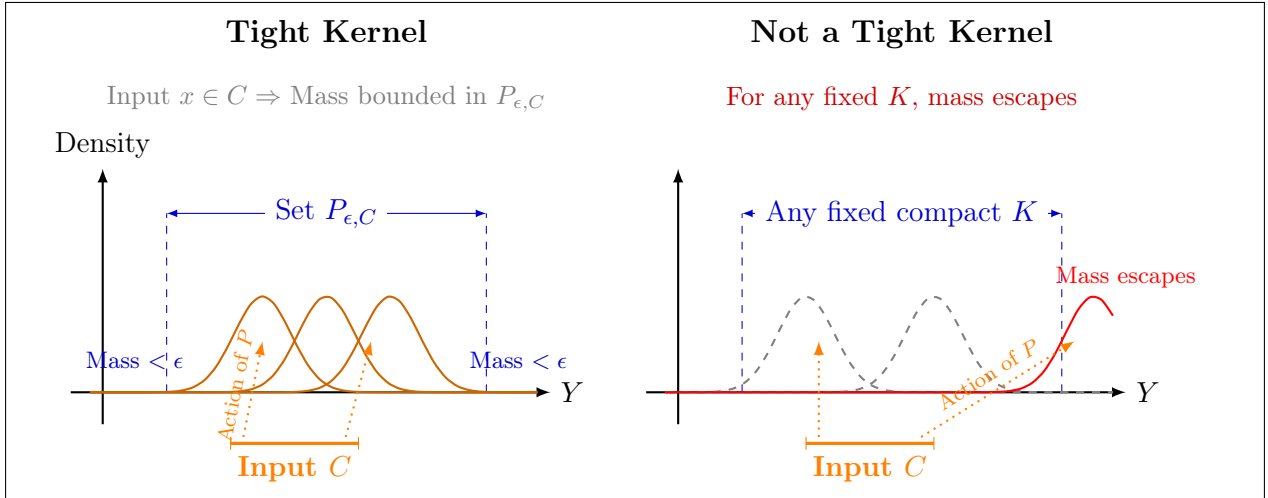

The composition of tight Feller kernels is defined by the standard composition of probability transition functions, known in stochastic analysis as the \emph{Chapman--Kolmogorov equation} \cite{Gardiner}. Given a kernel $P: X \rightsquigarrow Y$ and a kernel $Q: Y \rightsquigarrow Z$, their composite $Q \circ P: X \rightsquigarrow Z$ is defined by integrating the second kernel against the measure defined by the first:
\[
(Q \circ P)(x, C) = \int_Y Q(y, C) P(x, dy)
\]
for any measurable set $C \subseteq Z$. Intuitively, this represents a two-step random process: first, the system transitions from state $x$ to an intermediate state $y$ according to the distribution $P(x, \cdot)$; subsequently, it transitions from $y$ to the final set $C$ according to $Q(y, \cdot)$. We now establish that this operation preserves the structural properties required by our framework.

\begin{theorem}[Composition of Tight Feller Kernels]
Let $X, Y, Z$ be Polish spaces. Let $P: X \rightsquigarrow Y$ and $Q: Y \rightsquigarrow Z$ be tight Feller kernels. Then the composite kernel $Q \circ P$ is a tight Feller kernel.
\end{theorem}

\begin{proof}
\textbf{Feller Property}: We must show that for any $h \in C_b(Z)$, the function $x \mapsto \int_Z h(z) (Q \circ P)(x, dz)$ is continuous. By the definition of the composite kernel, this integral can be rewritten as iterated integrals:
\[
\int_Z h(z) (Q \circ P)(x, dz) = \int_Y \left( \int_Z h(z) Q(y, dz) \right) P(x, dy)
\]
Define $g(y) = \int_Z h(z) Q(y, dz)$. Since $Q$ is Feller and $h$ is bounded continuous, $g$ is a bounded continuous function on $Y$ (i.e., $g \in C_b(Y)$). The expression then simplifies to $\int_Y g(y) P(x, dy)$. Since $P$ is Feller and $g \in C_b(Y)$, this integral varies continuously with $x$. Thus $Q \circ P$ is Feller.

\textbf{Tightness} (Prokhorov-style argument): Let $K \subseteq X$ be compact and $\epsilon > 0$.
Since $P$ is tight, there exists a compact set $M \subseteq Y$ such that:
\[
\sup_{x \in K} P(x, M^c) < \frac{\epsilon}{2}
\]
Since $Q$ is tight and $M$ is compact, there exists a compact set $L \subseteq Z$ such that:
\[
\sup_{y \in M} Q(y, L^c) < \frac{\epsilon}{2}
\]
Now consider the mass assigned to $L^c$ by the composite kernel for any $x \in K$:
\begin{align*}
(Q \circ P)(x, L^c) &= \int_Y Q(y, L^c) P(x, dy) \\
&= \int_M Q(y, L^c) P(x, dy) + \int_{M^c} Q(y, L^c) P(x, dy)
\end{align*}
For the first term, $Q(y, L^c) < \epsilon/2$ for all $y \in M$, so the integral is bounded by $\epsilon/2 \cdot P(x, M) \le \epsilon/2$.
For the second term, the integral is bounded by $1 \cdot P(x, M^c) < \epsilon/2$.
Thus, $(Q \circ P)(x, L^c) < \epsilon$ for all $x \in K$, proving tightness.
\end{proof}

\subsection{Optimal Transport and Wasserstein Distances}
\label{app:optimal_transport}

The probabilistic framework in Section~\ref{sec:wasserstein} relies on the theory of Optimal Transport (see \cite{Villani}) to quantify the cost of moving a probability distribution into a compliant state.

\begin{definition}[Wasserstein Distance]
Let $(X, d)$ be a Polish space. For $p \in [1, \infty)$, the $p$-Wasserstein distance between two probability measures $\mu, \nu \in \mathcal{P}(X)$ is defined as:
\[
W_p(\mu, \nu) := \left( \inf_{\pi \in \Pi(\mu, \nu)} \int_{X \times X} d(x, y)^p \, d\pi(x, y) \right)^{1/p}
\]
where $\Pi(\mu, \nu)$ is the set of couplings (joint distributions) on $X \times X$ with marginals $\mu$ and $\nu$.
\end{definition}

The Wasserstein distance $W_1$ (where $p=1$) is often referred to as the \textbf{Earth Mover's Distance}.
\begin{itemize}
    \item \textbf{Physical Intuition:} Imagine the probability distribution $\mu$ as a pile of earth of unit mass, and $\nu$ as a target hole of unit capacity. $W_1(\mu, \nu)$ represents the \emph{minimum work} required to move the earth from $\mu$ to fill $\nu$, where Work = Mass $\times$ Distance.
    \item \textbf{Financial Intuition (Turnover Cost):} In the context of portfolio construction, if $\mu$ represents the noisy output of a stochastic solver, and $\nu$ represents a compliant portfolio distribution, $W_1(\mu, \nu)$ represents the \emph{minimum turnover volume} (in dollars or percent) required to rebalance the portfolio into compliance. Unlike the Euclidean distance (which sums squared errors), $W_1$ scales linearly with the size of the trades, making it the natural metric for transaction costs.
\end{itemize}

\begin{theorem}[Kantorovich--Rubinstein Duality]
For $p=1$, the Wasserstein distance admits a dual representation:
\[
W_1(\mu, \nu) = \sup \left\{ \int_X f \, d\mu - \int_X f \, d\nu : f \text{ is 1-Lipschitz} \right\}
\]
\end{theorem}

This duality is critical for computation. It implies that verifying Wasserstein safety is equivalent to checking the expectation of all Lipschitz ``test functions'' (risk factors). If every risk factor with unit sensitivity has an expected value close to the target, the portfolio is safe.

\begin{figure}[htbp]
\centering
\begin{tikzpicture}[scale=1.0, >=stealth]
    \draw[->] (-1,0) -- (6,0) node[right] {State Space $X$};
    \draw[->] (0,-1) -- (0,4) node[above] {Probability Density};

    \draw[thick, red!80!black, fill=red!10, fill opacity=0.5]
        plot [smooth, domain=0.5:2.5] (\x, {2.5*exp(-(\x-1.5)^2/0.2)}) -- cycle;
    \node[red!80!black] at (1.5, 3.0) {$\mu$ (Solver)};

    \draw[thick, blue!80!black, fill=blue!10, fill opacity=0.5]
        plot [smooth, domain=3.5:5.5] (\x, {2.5*exp(-(\x-4.5)^2/0.2)}) -- cycle;
    \node[blue!80!black] at (4.5, 3.0) {$\nu$ (Target)};

    \foreach \x in {1.2, 1.4, 1.5, 1.6, 1.8} {
        \draw[->, dashed, thick, gray] (\x, 0.5) to[bend left=20] (\x+3.0, 0.5);
    }
    \node[draw, fill=white, inner sep=3pt, rounded corners] at (3, 1.5) {Transport Cost $\int d(x,y) d\pi$};

    \draw[<->, thick] (1.5, -0.5) -- (4.5, -0.5) node[midway, below] {$W_1(\mu, \nu)$};
    \draw[dotted] (1.5, 0) -- (1.5, -0.5);
    \draw[dotted] (4.5, 0) -- (4.5, -0.5);

\end{tikzpicture}
\caption{Visualizing the Wasserstein Distance ($W_1$). It represents the minimum ``work'' required to move the probability mass from the noisy solver distribution $\mu$ to the compliant distribution $\nu$.}
\label{fig:wasserstein_transport}
\end{figure}
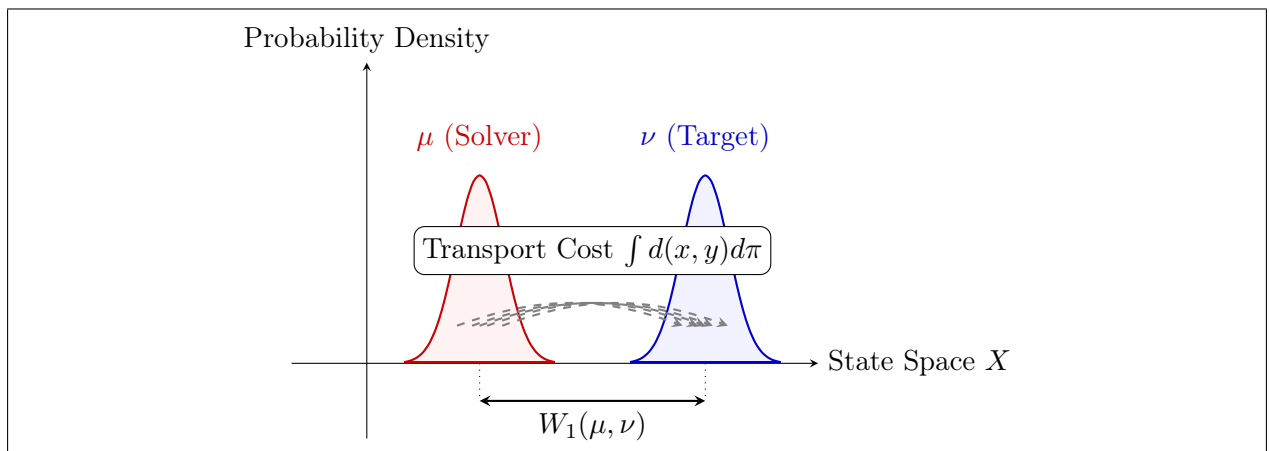

\begin{remark}[Computational Note: Solver Complexity]
The calculation of the Wasserstein safety check in Section~\ref{sub:wasserstein-compute} relies on projecting Monte Carlo samples onto the valid constraint set $S$. The choice of the order $p$ in the Wasserstein metric dictates the mathematical class of this projection problem and consequently the computational throughput.

\begin{itemize}
    \item \textbf{Case $p=1$ (Linear Transaction Costs):} 
    When modeling spreads or fixed commissions, the projection minimizes the $L_1$ norm $\|z - y\|_1$. If the constraint set $S$ is polyhedral (defined by linear inequalities $Ax \le b$), this constitutes a \emph{Linear Program (LP)}. Modern LP solvers (e.g., simplex or interior point methods) offer high throughput, capable of solving tens of thousands of projections per second for typical portfolio dimensions ($N \approx 500$).

    \item \textbf{Case $p=2$ (Execution Risk / Market Impact):} 
    When modeling quadratic market impact (following Kyle's $\lambda$, see \cite{Kyle}), the projection minimizes the Euclidean distance squared $\|z - y\|_2^2$. This constitutes a \emph{Quadratic Program (QP)}. This is the preferred choice for industrial implementation because it shares the exact mathematical structure of standard Mean-Variance Optimization. Consequently, it can use the same high-performance solvers (e.g., OSQP, Gurobi) already present in the portfolio construction stack, ensuring high throughput and numerical stability.

    \item \textbf{General Case $p \neq 1, 2$:} 
    For non-standard impact models (e.g., $p=1.5$ for square-root impact), the projection requires minimizing a general $L_p$ norm. This typically requires \emph{Second-Order Cone Programming (SOCP)} or general convex solvers. These methods incur significant computational overhead compared to LP or QP solvers and may create bottlenecks if used within the inner loops of a Monte Carlo simulation.
\end{itemize}
\end{remark}

\subsection{Proof of Convolution Stability (Proposition~\ref{prop:convolution-stability})}
\label{app:convolution-stability-proof}

\begin{proof}
For (1): we require that the kernel $Q(y, \cdot)$ acts as \emph{additive noise}---i.e., the output distribution has the form $y + \xi$ where $\xi$ is independent of $y$ with log-concave density $\rho_\xi$. (This is the natural model for solver noise: the optimizer targets $y$ and produces $y$ plus a random perturbation.) Under this assumption, writing $\rho_Y$ for the density of $P(x, \cdot)$, the composite density is a genuine convolution: $\rho_Z(z) = \int \rho_\xi(z - y) \rho_Y(y) \, dy$. The Pr\'ekopa--Leindler inequality \cite{Prekopa, Leindler} (see also \cite{Gardner} for a survey) implies this convolution is log-concave when both $\rho_Y$ and $\rho_\xi$ are log-concave. For a log-concave density, the $\alpha$-superlevel set $\{z : \rho(z) \ge \alpha\}$ is convex.

It remains to show that $\operatorname{supp}_{\delta+\epsilon}(\rho_{Q \circ P}) \subseteq A + B$, where $A = \operatorname{supp}_\delta(P(x))$ and $B = \operatorname{supp}_\epsilon(\rho_\xi)$, both convex and compact.

\emph{Step 1: Density bound outside $A+B$.}
Let $\lambda_\delta$ and $\lambda_\epsilon$ denote the HDR density thresholds for $\rho_Y$ and $\rho_\xi$ respectively, so that $A = \{\rho_Y \ge \lambda_\delta\}$ and $B = \{\rho_\xi \ge \lambda_\epsilon\}$.  For any $z \notin A + B$ and any $y$, we have $y \notin A$ or $z - y \notin B$.  Splitting the convolution integral:
\begin{align*}
\rho_Z(z) &= \int_{y \in A} \rho_Y(y)\, \rho_\xi(z - y)\, dy \;+\; \int_{y \notin A} \rho_Y(y)\, \rho_\xi(z - y)\, dy \\
&< \lambda_\epsilon \int_{y \in A} \rho_Y(y)\, dy \;+\; \lambda_\delta \int_{y \notin A} \rho_\xi(z - y)\, dy
\;\le\; \lambda_\delta + \lambda_\epsilon,
\end{align*}
where the first term uses $z - y \notin B$ (hence $\rho_\xi(z-y) < \lambda_\epsilon$) for $y \in A$, and the second uses $\rho_Y(y) < \lambda_\delta$ for $y \notin A$.

\emph{Step 2: Mass bound.}
Define the ``escape event'' $E = \{z \notin A + B\}$. For any $z \in E$ and any decomposition $z = y + (z-y)$, we must have $y \notin A$ or $(z-y) \notin B$. By the union bound:
\[
\Pr(Z \in E) \le \Pr(Y \notin A) + \Pr(\xi \notin B) \le \delta + \epsilon,
\]
where $Y \sim P(x, \cdot)$ and $\xi$ is the independent noise. Therefore $A + B$ has probability mass at least $1 - (\delta + \epsilon)$.

\emph{Step 3: Convolution stability.}
From Step~1, the superlevel set $C := \{z : \rho_Z(z) \ge \lambda_\delta + \lambda_\epsilon\} \subseteq A + B$.  From Step~2, $\Pr(Z \in A+B) \ge 1 - (\delta + \epsilon)$.  In particular, $(A+B)^c$ has mass at most $\delta + \epsilon$ and, by Step~1, the density is below $\lambda_\delta + \lambda_\epsilon$ throughout $(A+B)^c$.

The HDR containment $\operatorname{supp}_{\delta+\epsilon}(\rho_Z) \subseteq A+B$ follows if $\Pr(C) \ge 1 - (\delta+\epsilon)$ (for then the HDR threshold satisfies $\lambda_{\delta+\epsilon} \ge \lambda_\delta + \lambda_\epsilon$, giving $\operatorname{supp}_{\delta+\epsilon}(\rho_Z) \subseteq C \subseteq A+B$).  This mass condition holds in the principal applications:
\begin{itemize}
\item \emph{Gaussian noise:} when $\rho_Y$ and $\rho_\xi$ are Gaussian, $\rho_Z$ is Gaussian and $\lambda_\delta + \lambda_\epsilon$ is small relative to the peak density, so $C$ captures essentially all of the mass.
\item \emph{Concentrated noise:} when the noise variance is small relative to the HDR scale, $A + B$ is only slightly larger than $A$, and the density is near-constant inside $A + B$, so $\Pr(C) \approx \Pr(A+B) \ge 1 - (\delta+\epsilon)$.
\end{itemize}
In general, the Convolution Stability Condition should be verified for the specific noise model in use.  The mass bound $\Pr(Z \in A+B) \ge 1 - (\delta+\epsilon)$ (Step~2) always holds and provides the weaker guarantee that $A + B$ is a valid $(\delta+\epsilon)$-confidence region; the stronger HDR containment requires the additional regularity above.

For general (non-additive) kernels, the Pr\'ekopa--Leindler argument does not directly apply; one would need the joint density $(z, y) \mapsto \rho_Q(z \mid y) \rho_Y(y)$ to be log-concave, which is a stronger condition.

For (2): the inclusion is immediate from the support containment.
\end{proof}

\clearpage  \section{The General DOTS Framework}
\label{app:dots_framework}

This appendix outlines the general mathematical syntax of the \textit{Double Operadic Theory of Systems} (DOTS) \cite{LibkindMyers, MyersDOTS}. While the main text instantiates this theory using portfolio spaces and alignment relations, the underlying structure is a general formalism for modeling systems, states, and constraints using the language of framed bicategories (also known as proarrow equipments) \cite{Shulman}.

\subsection{The Ambient Structure: Framed Bicategories}

The DOTS framework operates within a specific type of double category that allows horizontal morphisms (functions) to be treated as special cases of vertical morphisms (relations).

\begin{definition}[Framed Bicategory / Equipment]
A \textbf{framed bicategory} (or equipment) is a double category satisfying the following structural axioms:
\begin{enumerate}
    \item \textbf{Vertical Structure:} The vertical morphisms form a bicategory (or a 2-category if the composition is strictly associative) \cite{Benabou}.
    \item \textbf{Horizontal Structure:} The horizontal morphisms form a category.
    \item \textbf{Companions and Conjoints:} For every horizontal morphism $f: A \to B$, there exist vertical morphisms $f_*: A \nrightarrow B$ (the companion) and $f^*: B \nrightarrow A$ (the conjoint) that represent $f$ in the vertical domain \cite{Shulman}. These essentially allow one to ``turn'' a function into a relation (its graph) and its transpose.
\end{enumerate}
\end{definition}

In a \textbf{thin} equipment (like the one used in this paper), the vertical structure is locally posetal, meaning parallel 2-cells are unique. This simplifies the bicategorical coherence laws to simple inclusions.

\subsection{States as Generalized Elements}

In general systems theory, a ``state'' is not just a property of an object, but a selection of valid configurations. We formalize this using generalized elements.

\begin{definition}[State]
Let $I$ be the monoidal unit object of the horizontal category (representing the trivial system). A \textbf{state} on an object $A$ is a horizontal morphism:
\[
\sigma: I \to A
\]
\end{definition}

In the context of set-based theories, $\sigma$ identifies a subset (or subspace) of $A$. In the context of logic, it represents a predicate on $A$.

\subsection{The Action Operator}

The core logic of DOTS is the interaction between states (horizontal data) and constraints (vertical data). This is formalized as an \textbf{action}.

\begin{definition}[The Action $\odot$]
Let $\sigma: I \to A$ be a state on $A$, and let $R: A \nrightarrow B$ be a vertical morphism (a system transformation or constraint). The \textbf{action} of $R$ on $\sigma$, denoted $\sigma \odot R$, is defined as the composite in the vertical category:
\[
\sigma \odot R := R \circ \sigma_*
\]
where $\sigma_*: I \nrightarrow A$ is the companion of the state $\sigma$.
\end{definition}

The action $\sigma \odot R$ represents the \textit{image} of the state $\sigma$ under the relation $R$. It is the set of all possible output configurations in $B$ that are reachable from valid input configurations in $A$ via the system $R$. This is analogous to functorial data migration concepts found in database theory \cite{Spivak}.

\subsection{Monoidal Structure and Parallelism}

To model systems composed of parallel subsystems, the equipment is equipped with a symmetric monoidal structure.

\begin{definition}[Symmetric Monoidal Double Category]
A double category is \textbf{symmetric monoidal} if it is equipped with a tensor product functor $\otimes: \mathbb{D} \times \mathbb{D} \to \mathbb{D}$ and a unit object $I$, satisfying the standard coherence axioms for both horizontal and vertical compositions.
\end{definition}

This structure implies two key properties for system modeling:
\begin{itemize}
    \item \textbf{Parallel States:} If $\sigma_A$ is a state on $A$ and $\sigma_B$ is a state on $B$, then $\sigma_A \otimes \sigma_B$ is a state on $A \otimes B$.
    \item \textbf{Parallel Actions:} The action operator distributes over the tensor product:
    \[
    (\sigma_A \otimes \sigma_B) \odot (R_A \otimes R_B) \cong (\sigma_A \odot R_A) \otimes (\sigma_B \odot R_B)
    \]
\end{itemize}

\subsection{Operadic Composition}

Finally, complex systems are built by wiring simpler systems together. This is captured formally by the \textbf{operad of wiring diagrams} \cite{Leinster, LibkindMyers}.

\begin{definition}[Wiring Diagram]
A wiring diagram is a morphism in the operad underlying the symmetric monoidal category. It specifies a method for constructing a global relation $R_{global}$ from a tensor product of local relations $R_1 \otimes \dots \otimes R_k$ via composition and contraction of indices.
\end{definition}

In the DOTS framework, these diagrams allow for the template-based construction of complex system constraints from reusable components, ensuring that local actions (on subsystems) compose coherently to form global actions.

\clearpage

\section{Conceptual System Architecture and Workflows}
\label{app:architecture}

This appendix translates the mathematical theory of the $\HSSimp$ double category and its extensions into a high-level software architecture. We describe a modular system designed to support multi-stage portfolio construction, continuous alignment verification, and the coherent propagation of changes.

\subsection{Architectural Principles}

The system is designed around the separation of \emph{topology} (the permissible spaces and maps) from \emph{computation} (the numerical solvers). This separation allows the system to support multiple semantic models (as described in Part II) within a single infrastructure.

\begin{enumerate}
    \item \textbf{Categorical Schema:} The database schema directly mirrors the double category structure. It stores \textit{Objects} (space definitions), \textit{Horizontal Morphisms} (re-implementation configs), and \textit{Vertical Morphisms} (alignment logic) as first-class entities.
    \item \textbf{Solver Agnosticism:} The core logic engine operates on abstract morphisms. The actual execution (computing $f(x)$ or checking $(x,z) \in R$) is delegated to a pluggable solver layer (supporting convex optimization, Monte Carlo, or discrete search).
    \item \textbf{Evidence Persistence:} Following the Span model (Section~\ref{sec:2-cells}), the system does not just store boolean pass/fail results. It stores \textit{Evidence Objects} (audit trails) that link specific hub portfolios to spoke portfolios.
\end{enumerate}

\subsection{Component Diagram}

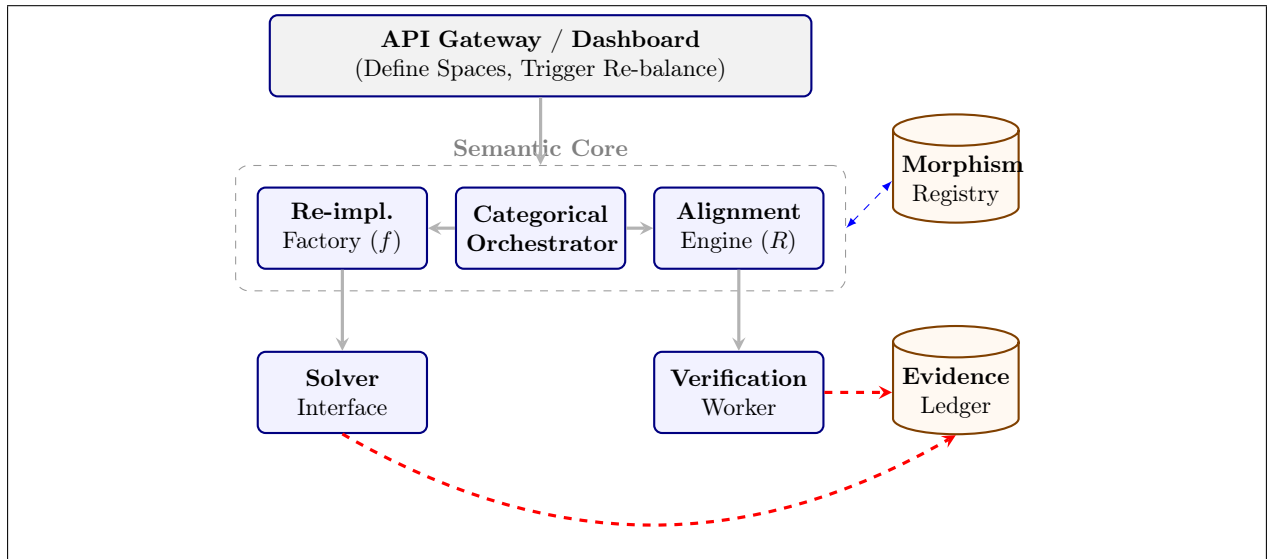
\begin{figure}[htbp]
\centering
\begin{tikzpicture}[
    scale=0.9, transform shape,
    >=stealth,
    font=\small, 
    component/.style={
        draw=blue!50!black, 
        fill=blue!5, 
        rounded corners=3pt, 
        minimum width=2.5cm, 
        minimum height=1.2cm,
        text width=2.2cm,
        align=center,
        thick
    },
    store/.style={
        draw=orange!50!black, 
        fill=orange!5, 
        cylinder, 
        shape border rotate=90, 
        aspect=0.25, 
        minimum width=1.8cm, 
        minimum height=1.6cm,
        text width=1.6cm,
        align=center,
        thick
    },
    engine/.style={
        draw=gray!80,
        fill=white,
        dashed,
        rounded corners=5pt,
        inner sep=8pt
    },
    arrow/.style={->, very thick, gray!60}
]

    \node[component] (Orch) at (0, 0) {\textbf{Categorical}\\ \textbf{Orchestrator}};
    
    \node[component, left=0.4cm of Orch] (Reimpl) {\textbf{Re-impl.}\\Factory ($f$)};
    \node[component, right=0.4cm of Orch] (Align) {\textbf{Alignment}\\Engine ($R$)};

    \begin{scope}[on background layer]
        \node[engine, fit=(Orch) (Reimpl) (Align), label={[gray]above:\textbf{Semantic Core}}] (CoreGroup) {};
    \end{scope}

    \node[component, fill=gray!10, minimum width=8cm, text width=7cm, above=1.0cm of CoreGroup] (API) {\textbf{API Gateway / Dashboard}\\(Define Spaces, Trigger Re-balance)};

    \node[component, below=1.2cm of Reimpl] (Solver) {\textbf{Solver}\\Interface};
    \node[component, below=1.2cm of Align] (Verifier) {\textbf{Verification}\\Worker};

    \node[store, right=1.0cm of Verifier] (Ledger) {\textbf{Evidence}\\Ledger};
    \node[store, above=1.5cm of Ledger] (DB) {\textbf{Morphism}\\Registry};

    \draw[arrow] (API) -- (CoreGroup);
    
    \draw[arrow] (Orch) -- (Reimpl);
    \draw[arrow] (Orch) -- (Align);
    
    \draw[arrow] (Reimpl) -- (Solver);
    \draw[arrow] (Align) -- (Verifier);

    \draw[<->, >=latex, dashed, blue] (CoreGroup.east) -- (DB.west);
    
    \draw[arrow, dashed, red] (Verifier) -- (Ledger);
    \draw[arrow, dashed, red] (Solver.south) to[out=-30, in=210] (Ledger.south);

\end{tikzpicture}
\caption{System Architecture for the $\HSSimp$ Framework. The Semantic Core abstracts the mathematical logic, using specialized solvers for execution and an Evidence Ledger for persistent audit trails. The Categorical Orchestrator coordinates the Re-implementation Factory and the Alignment Engine; it operates on abstract morphisms, delegating numerical execution to the layer below.}
\label{fig:system_arch}
\end{figure}

Figure~\ref{fig:system_arch} illustrates the modular architecture. The layout groups the semantic logic centrally, with the execution layer below and persistence to the right, ensuring a clear flow of data and evidence.

\subsection{Execution Flows}

The following diagrams illustrate how specific portfolio events traverse these architectural components for a single stage of the portfolio construction (re-implementation) process. Compositionality guarantees the validity of a multi-stage process, when each stage is executed in turn.

These workflows are simple, almost trivial, to specify. The complications come when executing them at an industrial scale---building hundreds or thousands of portfolios, often in multi-stage processes, that are required to satisfy equally numerous alignment relations. The operational complexity of \emph{ad hoc} processes can grow geometrically. This is where the category-theoretic framework, and associated system architecture, can play a valuable organizing role.

\subsubsection*{Workflow A: Portfolio Change Propagation}
When a Hub portfolio changes ($x \leadsto x'$), the system executes parallel branches through the factory and alignment engines to ensure safety before committing.

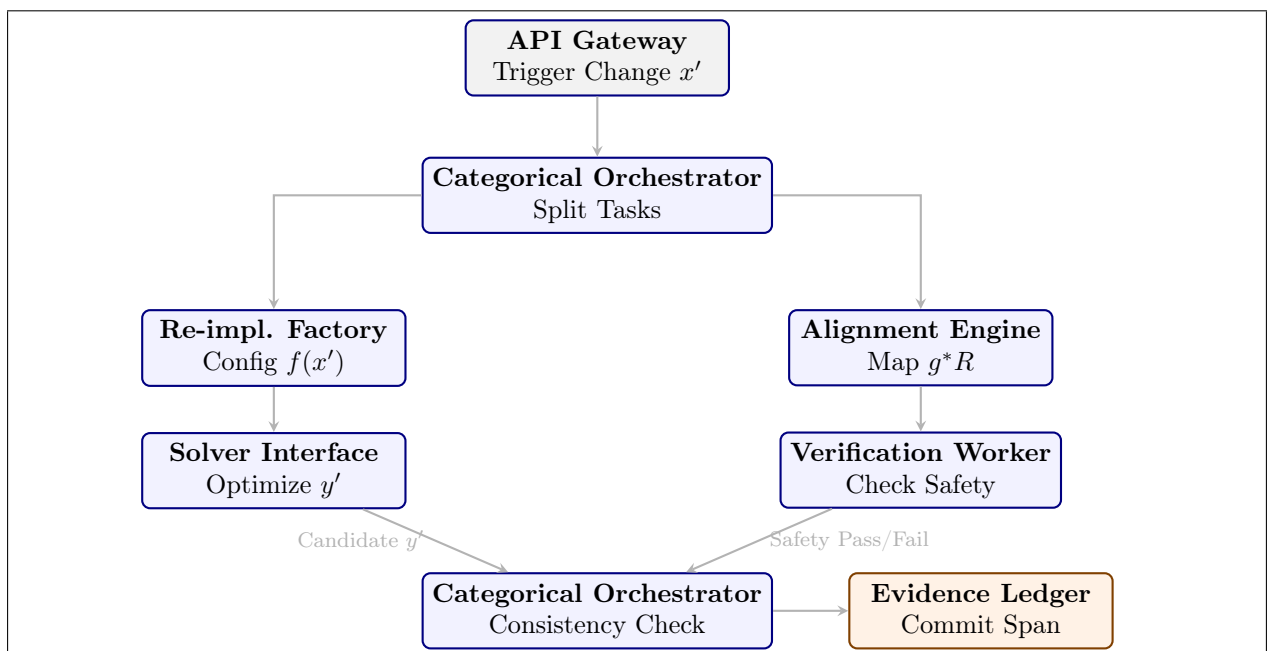
\begin{figure}[htbp]
\centering
\begin{tikzpicture}[
    node distance=1.0cm and 0.5cm,
    >=stealth,
    font=\small,
    sys_node/.style={
        rectangle, 
        draw=blue!50!black, 
        fill=blue!5, 
        rounded corners=3pt, 
        minimum width=3.5cm, 
        minimum height=1.0cm, 
        align=center,
        thick
    },
    data_node/.style={
        rectangle, 
        draw=gray!60, 
        fill=white, 
        dashed,
        minimum width=2.5cm, 
        minimum height=0.8cm, 
        align=center,
        font=\footnotesize
    },
    arrow/.style={->, thick, gray!60}
]

    \node[sys_node, fill=gray!10] (API) {\textbf{API Gateway}\\Trigger Change $x'$};

    \node[sys_node, below=0.8cm of API] (Orch) {\textbf{Categorical Orchestrator}\\Split Tasks};

    \node[sys_node, below left=1.0cm and 0.2cm of Orch] (Factory) {\textbf{Re-impl. Factory}\\Config $f(x')$};
    \node[sys_node, below=0.6cm of Factory] (Solver) {\textbf{Solver Interface}\\Optimize $y'$};

    \node[sys_node, below right=1.0cm and 0.2cm of Orch] (Align) {\textbf{Alignment Engine}\\Map $g^*R$};
    \node[sys_node, below=0.6cm of Align] (Verifier) {\textbf{Verification Worker}\\Check Safety};

    \node[sys_node, below=4.5cm of Orch] (Decide) {\textbf{Categorical Orchestrator}\\Consistency Check};

    \node[sys_node, fill=orange!10, draw=orange!50!black, right=1.0cm of Decide] (Ledger) {\textbf{Evidence Ledger}\\Commit Span};

    \draw[arrow] (API) -- (Orch);
    
    \draw[arrow] (Orch) -| (Factory);
    \draw[arrow] (Factory) -- (Solver);
    \draw[arrow] (Solver) -- node[left, font=\scriptsize] {Candidate $y'$} (Decide);

    \draw[arrow] (Orch) -| (Align);
    \draw[arrow] (Align) -- (Verifier);
    \draw[arrow] (Verifier) -- node[right, font=\scriptsize] {Safety Pass/Fail} (Decide);

    \draw[arrow] (Decide) -- (Ledger);

\end{tikzpicture}
\caption{Workflow A mapped to Architecture. The Orchestrator delegates construction to the Solver and safety checks to the Verifier. The Spoke is only committed if the Verifier confirms the Hub change is safe via the Lax Beck--Chevalley condition.}
\label{fig:workflow_arch_a}
\end{figure}

\subsubsection*{Workflow B: Alignment Change Propagation}
When a new constraint is defined ($R \leadsto R'$), the system uses the Alignment Engine to compute both the primal implications (new menus) and dual implications (compliance checks).

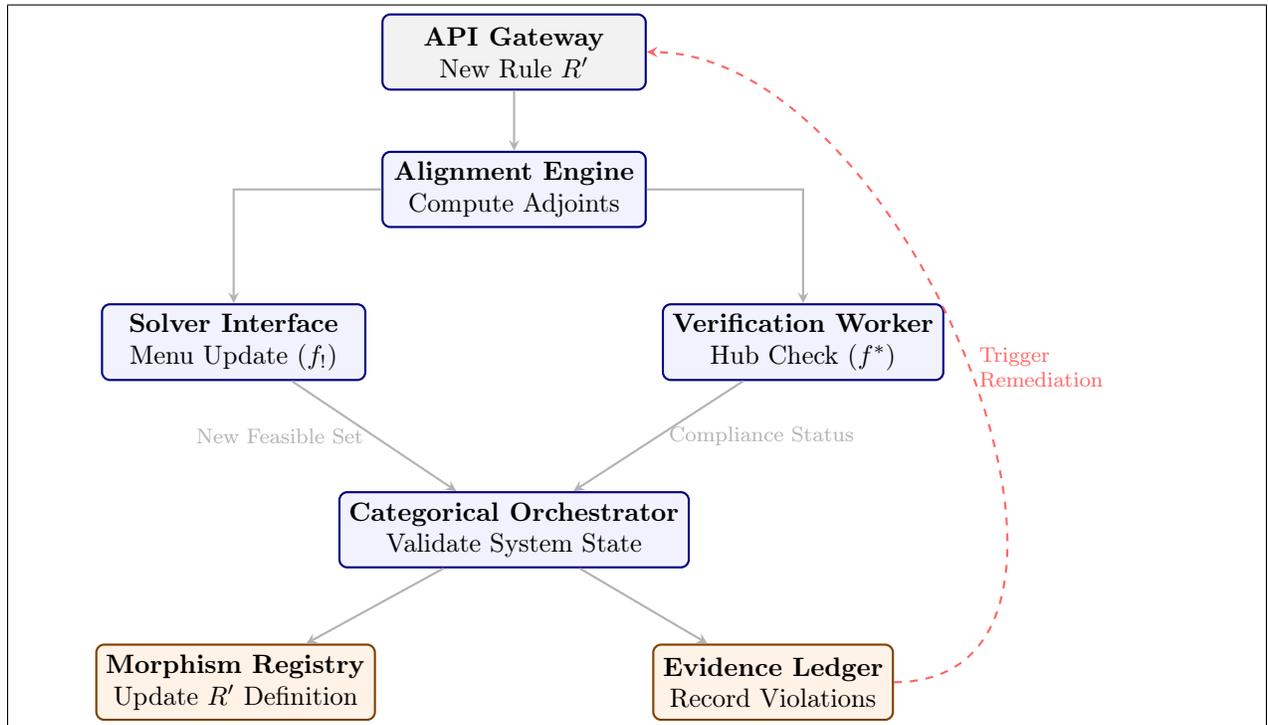
\begin{figure}[htbp]
\centering
\begin{tikzpicture}[
    node distance=1.0cm and 0.5cm,
    >=stealth,
    font=\small,
    sys_node/.style={
        rectangle, 
        draw=blue!50!black, 
        fill=blue!5, 
        rounded corners=3pt, 
        minimum width=3.5cm, 
        minimum height=1.0cm, 
        align=center,
        thick
    },
    persist_node/.style={
        rectangle, 
        draw=orange!50!black, 
        fill=orange!10, 
        rounded corners=3pt, 
        minimum width=3.0cm, 
        minimum height=1.0cm, 
        align=center,
        thick
    },
    arrow/.style={->, thick, gray!60}
]

    \node[sys_node, fill=gray!10] (API) {\textbf{API Gateway}\\New Rule $R'$};

    \node[sys_node, below=0.8cm of API] (Align) {\textbf{Alignment Engine}\\Compute Adjoints};

    \node[sys_node, below left=1.0cm and 0.2cm of Align] (Solver) {\textbf{Solver Interface}\\Menu Update ($f_!$)};
    \node[sys_node, below right=1.0cm and 0.2cm of Align] (Verifier) {\textbf{Verification Worker}\\Hub Check ($f^*$)};

    \node[sys_node, below=3.5cm of Align] (Orch) {\textbf{Categorical Orchestrator}\\Validate System State};

    \node[persist_node, below left=1.0cm and -0.5cm of Orch] (Registry) {\textbf{Morphism Registry}\\Update $R'$ Definition};
    \node[persist_node, below right=1.0cm and -0.5cm of Orch] (Ledger) {\textbf{Evidence Ledger}\\Record Violations};

    \draw[->, dashed, thick, red!60] (Ledger.east) to[out=0, in=0] node[right, font=\scriptsize, align=left] {Trigger\\Remediation} (API.east);

    \draw[arrow] (API) -- (Align);
    
    \draw[arrow] (Align) -| (Solver);
    \draw[arrow] (Align) -| (Verifier);

    \draw[arrow] (Solver) -- node[left, font=\scriptsize] {New Feasible Set} (Orch);
    \draw[arrow] (Verifier) -- node[right, font=\scriptsize] {Compliance Status} (Orch);

    \draw[arrow] (Orch) -- (Registry);
    \draw[arrow] (Orch) -- (Ledger);

\end{tikzpicture}
\caption{Workflow B mapped to Architecture. The Alignment Engine drives the dual process: The Solver defines the new feasible region for spokes, while the Verifier identifies broken hubs.}
\label{fig:workflow_arch_b}
\end{figure}

\subsubsection*{Workflow C: Build New Spokes}
When a series of new investment vehicles is required ($K_{new}$), the system constructs the re-implementation map $f$ by optimizing an objective $u$ subject to a given structural alignment $R$.

\begin{figure}[htbp]
\centering
\begin{tikzpicture}[
    node distance=1.0cm and 0.5cm,
    >=stealth,
    font=\small,
    sys_node/.style={
        rectangle, 
        draw=blue!50!black, 
        fill=blue!5, 
        rounded corners=3pt, 
        minimum width=3.5cm, 
        minimum height=1.0cm, 
        align=center,
        thick
    },
    persist_node/.style={
        rectangle, 
        draw=orange!50!black, 
        fill=orange!10, 
        rounded corners=3pt, 
        minimum width=3.0cm, 
        minimum height=1.0cm, 
        align=center,
        thick
    },
    arrow/.style={->, thick, gray!60}
]

    \node[sys_node, fill=gray!10] 
    (API) {\textbf{API Gateway}\\Request Spokes ($R, u$)};

    \node[sys_node, below=0.8cm of API] (Align) {\textbf{Alignment Engine}\\Define Constraints $F_R(x)$};

    \node[sys_node, below left=1.0cm and 0.2cm of Align] (Solver) {\textbf{Solver Interface}\\Construct $f(x)$\\(Maximize $u$)};
    \node[sys_node, below right=1.0cm and 0.2cm of Align] (Verifier) {\textbf{Verification Worker}\\Validate Object $K_{new}$\\ (Closedness)};

    \node[sys_node, below=3.5cm of Align] (Orch) {\textbf{Categorical Orchestrator}\\Finalize Morphism $f$};

    \node[persist_node, below=1.0cm of Orch] (Registry) {\textbf{Morphism Registry}\\Register $f$ and $K_{new}$};

    \draw[arrow] (API) -- (Align);
    
    \draw[arrow] (Align) -| (Solver);
    \draw[arrow] (Align) -| (Verifier);

    \draw[arrow] (Solver) -- node[left, font=\scriptsize] {Optimal Map} (Orch);
    \draw[arrow] (Verifier) -- node[right, font=\scriptsize] {Object Valid} (Orch);

    \draw[arrow] (Orch) -- (Registry);

\end{tikzpicture}
\caption{Workflow C mapped to Architecture. The Alignment Engine defines the search space based on $R$, allowing the Solver to construct a new horizontal morphism $f$ that maximizes the objective $u$ (Section~\ref{subsec:opt-align}).}
\label{fig:workflow_arch_c}
\end{figure}
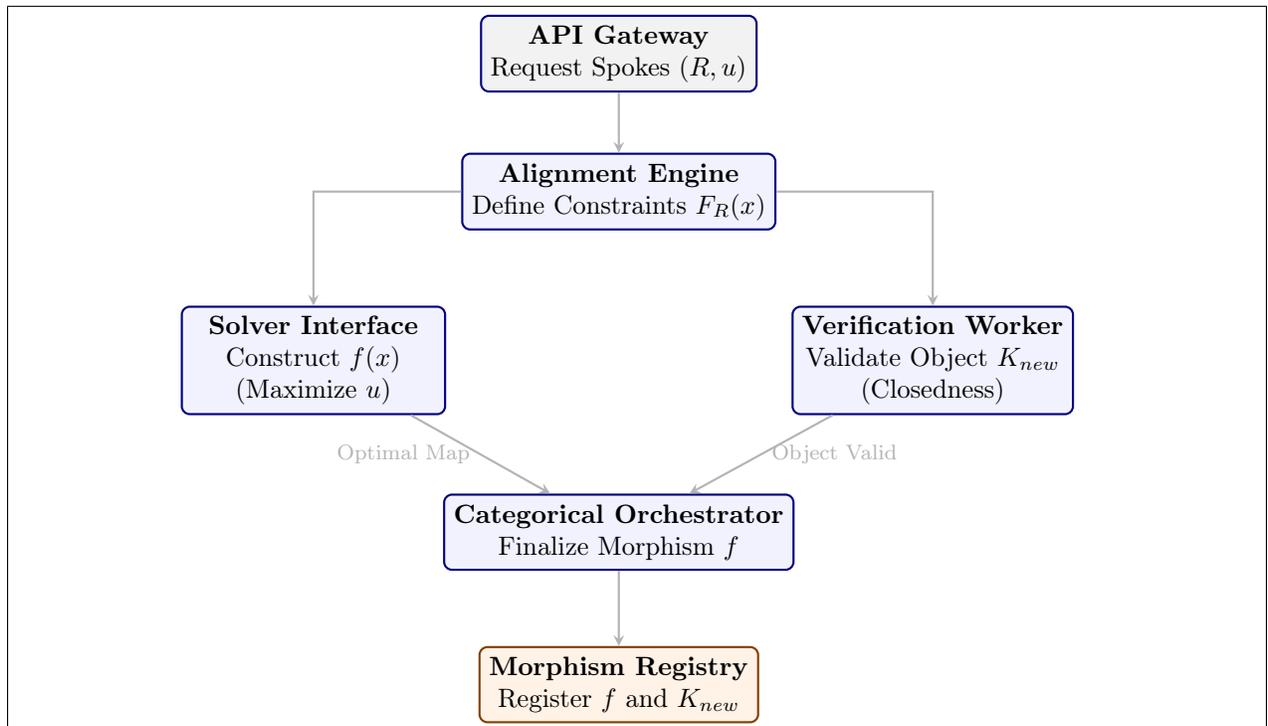

\section{Proof-of-Concept Implementations}
\label{app:poc}

This appendix describes two lightweight software implementations that instantiate, respectively, the deterministic DOTS framework of Section~\ref{sec:DOTS} with the system architecture of Appendix~\ref{app:architecture}, and the probabilistic extensions of Part~III.  Both are Python applications with interactive browser-based interfaces, accompanied by automated test suites.

The implementations serve three purposes: to verify numerically that the coherence laws of the framework hold in both the deterministic and probabilistic settings; to demonstrate that the architectures described in the paper can be instantiated with modest effort; and to provide interactive environments in which the reader can build geometric intuition for the operations involved.

\subsection{Application~I: The DOTS Framework}
\label{subsec:poc-dots}

The first application consists of a mathematical engine (\texttt{dots\_core.py}, approximately 500 lines) and an interactive interface (\texttt{dots\_app.py}, approximately 400 lines), with 28 automated tests.  It operates on the $2$-simplex $\Delta^2$ (three assets), using lattice enumeration at step~$0.01$ (5,151 points) for all set-valued computations.

\medskip
\noindent\textbf{Included.}\quad The following components are implemented:
\begin{itemize}
    \item \textbf{Simplex and permissible spaces} (Definition~\ref{def:ambient-portfolio-space}, Definition~\ref{def:permissible-portfolio-space}): lattice enumeration, constraint predicates, barycentric-to-Cartesian visualization.
    \item \textbf{Three concrete relations} (Section~\ref{subsec:frequent-relations}): the tracking relation $R_{\mathrm{track},\epsilon}$, the fee-cap projector $R_{\mathrm{fee}\le\tau}$, and the turnover relation $R_{\mathrm{turn}\le\kappa}$.
    \item \textbf{The DOTS action} $K \odot R$ (Proposition~\ref{prop:action}): lattice-based computation of the menu, including single-point fibers~$F_R(x)$.
    \item \textbf{Pullback} $\pullback{f}S$ (Definition~\ref{def:pullback}) and \textbf{pushforward} $f_! R$ (Definition~\ref{def:pushforward}).
    \item \textbf{Vertical composition} of relations.
    \item \textbf{Determinization} via $\ell^2$-regularization (Section~\ref{subsec:worked-example}): given a regularization parameter $\alpha > 0$, selects the minimum-norm element of each fiber.
    \item \textbf{Coherence verification}: automated checks of all five action laws (Proposition~\ref{prop:action}), Frobenius reciprocity (Theorem~\ref{thm:DOTS-Frob}), and Beck--Chevalley (Theorem~\ref{thm:DOTS-BC}).
    \item \textbf{Appendix~\ref{app:architecture} components}: a Morphism Registry, an Evidence Ledger, and a Categorical Orchestrator implementing Workflow~A.
\end{itemize}

\subsubsection*{The worked example}

The central demonstration implements the worked example of Section~\ref{subsec:worked-example} with the following parameters:
\begin{itemize}
    \item \textbf{Hub:} $\Delta^2$ with constraint $x_1 \le 0.6$.
    \item \textbf{Spoke:} $\Delta^2$ (unconstrained).
    \item \textbf{Fee map:} $\Fee(y) = 10\,y_1 + 5\,y_2 + 0\,y_3$ (basis points).
    \item \textbf{Tracking:} $R_{\mathrm{track},\epsilon}$ with identity exposure maps ($g_A = g_B = \mathrm{id}$).
    \item \textbf{Fee cap:} $R_{\mathrm{fee}\le\tau}$ with adjustable cap $\tau$.
\end{itemize}

With the default parameters $\epsilon = 0.05$ and $\tau = 6$\,bps, the hub permissible space contains 4,331 lattice points (at step~$0.01$).  The tracking menu $K \odot R_{\mathrm{track},0.05}$ comprises 4,485 spoke-space points (slightly larger than the hub, since spokes outside the hub constraint boundary can still track hub points near it); intersecting with the fee-cap projector reduces this to 3,511 points.  The menu visibly shrinks toward the Asset~3 vertex, confirming the projector law (Proposition~\ref{prop:action}(e)).

\subsubsection*{Coherence verification}

The test suite verifies the five action laws and Frobenius reciprocity numerically at step~$0.05$:
\begin{enumerate}[label=(\alph*),leftmargin=*]
    \item \emph{Closedness}: the menu $K \odot R$ is non-empty and well-defined.
    \item \emph{Unitality}: $K \odot \Delta_{\Delta^2} = K$ (195 lattice points on both sides).
    \item \emph{Associativity}: $(K \odot R) \odot S = K \odot (S \vcomp R)$, verified with $R = R_{\mathrm{track},0.10}$ and $S = R_{\mathrm{turn}\le 0.3}$.
    \item \emph{Isotonicity}: $K \subseteq K'$ implies $K \odot R \subseteq K' \odot R$, verified with $K' = \Delta^2$.
    \item \emph{Projectors}: $K \odot R_{\mathrm{fee}\le 6} = K \cap \{y : \Fee(y) \le 6\}$, with matching point counts.
\end{enumerate}

\noindent Frobenius reciprocity (Theorem~\ref{thm:DOTS-Frob}) is verified for $f(x) = 0.8\,x + 0.2\,(\tfrac{1}{3}, \tfrac{1}{3}, \tfrac{1}{3})$, $R = R_{\mathrm{track},0.10}$, and $S = R_{\mathrm{turn}\le 0.3}$ at step~$0.10$; both sides of $f_!(R \cap \pullback{f}S) = f_! R \cap S$ produce the same set.

\subsubsection*{Workflow~A demonstration}

The seven components of Appendix~\ref{app:architecture} are instantiated as Python classes:

\begin{center}
\begin{tabular}{ll}
\toprule
\textbf{Paper Component} & \textbf{Implementation} \\
\midrule
Morphism Registry & \texttt{MorphismRegistry}: objects, h-morphisms, v-morphisms \\
Evidence Ledger & \texttt{EvidenceLedger}: timestamped audit entries \\
Categorical Orchestrator & \texttt{CategoricalOrchestrator.workflow\_a()} \\
Re-implementation Factory & Configurable callable $f: \Delta^n \to \Delta^m$ \\
Alignment Engine & Relation membership test via \texttt{R.contains\_pairs()} \\
Solver Interface & Direct evaluation of $f(x')$ \\
Verification Worker & Check $(x', y') \in R$ \\
\bottomrule
\end{tabular}
\end{center}

Workflow~A proceeds as described in Appendix~\ref{app:architecture}: on receiving a hub change $x \mapsto x'$, the orchestrator computes a candidate spoke, checks safety, commits or rejects, and records an evidence span.  The interactive interface lets the user step through this workflow and inspect the audit trail.

\subsection{Application~II: The Probabilistic Framework}
\label{subsec:poc-prob}

The second application implements the three probabilistic approaches of Part~III: the safety radius of Section~\ref{sec:safety-radius}, the $\epsilon$-HDR of Definition~\ref{def:epsilon-supp}, and the Wasserstein transport-based approach of Section~\ref{sec:wasserstein}.  It consists of a mathematical engine (\texttt{prob\_core.py}, approximately 350 lines) and an interactive interface (\texttt{prob\_app.py}, approximately 400 lines), with 28 automated tests.  It uses the same $2$-simplex $\Delta^2$ and lattice-based approach as Application~I.

\subsubsection*{Scope}

The implementation instantiates:
\begin{itemize}
    \item \textbf{Stochastic kernels} on $\Delta^2$ (Definition~\ref{def:HSP-r}): three kernel types---Gaussian, bimodal (``split peak''), and curved (``banana'')---each modeled as a sample generator that projects to the simplex.
    \item \textbf{Safety radius} (Section~\ref{sec:safety-radius}): empirical computation of $r_\epsilon$ as the $(1-\epsilon)$-quantile of distances, metric erosion, Minkowski dilation, and composition of radii under both the linear bound $L r_P + r_Q$ and the quadratic bound $\sqrt{(Lr_P)^2 + r_Q^2}$.
    \item \textbf{HDR} (Definition~\ref{def:epsilon-supp}): kernel density estimation on the simplex lattice, density threshold computation, and $\epsilon$-HDR set extraction.
    \item \textbf{Wasserstein transport} (Definition~\ref{def:wasserstein-safety}): per-sample $L^1$ cure costs, mean $W_1$ distance to the constraint set, and budget-based compliance checks.
    \item \textbf{Failure mode scenarios} (Section~\ref{sec:discrete}): the split peak (Section~\ref{sub:split_peak}) and the banana (Section~\ref{sub:banana}).
    \item \textbf{Three-way comparison}: side-by-side evaluation of all three approaches on the same scenario.
    \item \textbf{Composition}: sequential kernels with linear and quadratic radius accumulation.
\end{itemize}

\noindent\textbf{Not included.}\quad The full categorical structure of $\textbf{HSP-r}$ (2-cells, lax functors), persona-indexed categories, convolution stability proofs, the hollow shell scenario, and production optimal-transport solvers.

\subsubsection*{Stochastic kernels}

Each kernel $P(x, \cdot)$ is represented as a Monte Carlo sampler: given a hub point $x \in \Delta^2$, the kernel returns $N$ samples on $\Delta^2$ by adding noise to $x$ and projecting back to the simplex.  Three shapes are implemented:
\begin{enumerate}[label=(\alph*),leftmargin=*]
    \item \emph{Gaussian}: isotropic noise $x + \sigma\,\xi$, $\xi \sim \mathcal{N}(0, I)$, projected to $\Delta^2$.  This produces a unimodal, approximately symmetric distribution.
    \item \emph{Bimodal}: a 50/50 mixture of two Gaussians offset along the first coordinate, modeling the split peak of Section~\ref{sub:split_peak}.  The center~$f(x)$ lies in the gap between the two modes.
    \item \emph{Banana}: a curved distribution generated by $x + (t,\, \kappa t^2 + \eta,\, -(t + \kappa t^2 + \eta))$ for $t \sim \mathcal{N}(0, 3\sigma)$, modeling the tax-aware transitioning of Section~\ref{sub:banana}.
\end{enumerate}

\noindent The three shapes are designed to illustrate, concretely, why the choice among the three compliance approaches matters: for (a) they all agree, while for (b) and (c) the safety radius approach becomes unnecessarily conservative.

\subsubsection*{Safety radius}

The safety radius $r_\epsilon$ is computed empirically as the $(1-\epsilon)$-quantile of $\{\|y_i - f(x)\|_2 : y_i \sim P(x)\}$.  With Gaussian noise at $\sigma = 0.03$ and $\epsilon = 0.05$, the implementation yields $r_{0.05} \approx 0.074$.  For the bimodal kernel the radius inflates to $r_{0.05} \approx 0.098$, and for the banana kernel to $r_{0.05} \approx 0.124$---the ball must grow to encompass the extended tails even though the high-density regions are compact.

Metric erosion implements the inner parallel set: a constraint point~$z$ survives if every lattice point within distance~$r$ of $z$ is also in the constraint set.  With a weight constraint $x_1 \le 0.4$ (861 lattice points at step~$0.02$), the eroded set shrinks to 798 points for the Gaussian kernel and to 700 for the banana kernel.  In the banana case, the center~$f(x)$ falls outside the eroded set: the safety radius approach rejects the hub even though the high-density region of the distribution lies entirely within the constraints.

Composition is verified by constructing sequential kernels $P$ (Gaussian, $\sigma = 0.025$) and $Q$ (Gaussian, $\sigma = 0.020$) with $L = 1$.  The individual radii are $r_P \approx 0.060$ and $r_Q \approx 0.049$.  The linear bound gives $r_{\mathrm{comp}} = Lr_P + r_Q \approx 0.109$; the quadratic bound gives $\sqrt{(Lr_P)^2 + r_Q^2} \approx 0.078$.  The actual composed radius (measured from the composed samples) is approximately $0.080$, confirming that the quadratic bound is tighter and that both bounds are conservative.

\subsubsection*{$\epsilon$-HDR}

Density estimation uses a Gaussian KDE on the simplex lattice (step~$0.02$, 1,326 points) with bandwidth matching the kernel parameter.  The density threshold~$\lambda_\epsilon$ is the $\epsilon$-quantile of density values at the sample points; the $\epsilon$-HDR is $\{y : \hat\rho(y) \ge \lambda_\epsilon\}$.

For the bimodal kernel, the HDR correctly identifies two disjoint clusters, one around each mode.  The ball, by contrast, covers the gap between the modes---a region of near-zero density.  Shrinking the risk budget from $\epsilon = 0.20$ (the 20\%-HDR) to $\epsilon = 0.05$ (the 5\%-HDR) reduces the HDR from 40 lattice points to 10, capturing only the peak cores.  This contraction is the geometric content of the monotonicity property of the $\epsilon$-support.

\subsubsection*{Wasserstein transport}

The Wasserstein-1 cure cost is computed as $W_1 = \frac{1}{N}\sum_i \min_{z \in S} \|y_i - z\|_1$: for each sample, the nearest point in the constraint set~$S$ is found by brute-force search over the lattice.  With a Gaussian kernel at $\sigma = 0.03$ and a weight constraint $x_1 \le 0.5$, the mean cure cost is $W_1 \approx 0.015$ with a violation rate of approximately 3\%.

The per-sample cost map reveals where violations concentrate: samples overshooting the weight cap require the largest adjustments, while samples landing inside~$S$ have zero cure cost.  The transport arrows---visualized in the interface as directed segments from each violating sample to its nearest compliant point---make the ``economic cost of non-compliance'' interpretation of Section~\ref{sec:wasserstein} visually tangible.

\subsubsection*{Three-way comparison}

The comparison tab runs all three approaches on the same scenario with the same $\epsilon$.  The results illuminate the structural differences:

\begin{center}
\begin{tabular}{lccc}
\toprule
\textbf{Scenario} & \textbf{Safety Radius} & \textbf{HDR} & \textbf{Wasserstein} \\
\midrule
Gaussian & Safe ($r = 0.074$) & Safe (21 pts) & Approved ($W_1 = 0.015$) \\
Split Peak & Safe ($r = 0.100$) & Safe (10 pts) & Approved ($W_1 = 0.016$) \\
Banana ($x_1 \le 0.4$) & \textbf{Rejected} ($r = 0.130$) & Safe (27 pts) & Approved ($W_1 = 0.020$) \\
\bottomrule
\end{tabular}
\end{center}

\noindent In the banana scenario with a tight weight constraint ($x_1 \le 0.4$), the safety radius ball extends beyond the eroded constraint set and the hub is rejected.  The HDR, which traces the actual density contour of the curved distribution rather than imposing a symmetric ball, correctly identifies the high-probability region as compliant.  The Wasserstein approach, measuring the average cost of curing the few boundary violations, finds the cost well within budget.  This is the failure mode described in Section~\ref{sub:banana}: the safety radius approach is sound but overly conservative for non-spherical distributions.

\subsection{Observations}
\label{subsec:poc-observations}

Four observations emerge from the two implementations.

First, the lattice-enumeration approach, while na\"ive, is sufficient to verify every algebraic identity in the deterministic framework and to demonstrate every probabilistic concept.  The discretization error is bounded by the lattice step and does not accumulate across composed operations---a consequence of the closedness guarantees that the theory provides.  For production use, the lattice would be replaced by convex optimization and optimal transport solvers, but the algebraic structure that governs how solvers compose remains the same.

Second, the separation of concerns in Appendix~\ref{app:architecture} translates directly into code modularity.  In Application~I, the \texttt{Relation} base class captures the abstract vertical morphism interface; concrete relations are subclasses.  In Application~II, the \texttt{StochasticKernel} dataclass captures the abstract kernel interface (Definition~\ref{def:HSP-r}); concrete kernels differ only in their sampling functions.  In both cases, the compliance checks operate on the abstract interface without knowing which specific type is in use.

Third, the probabilistic implementation makes concrete the paper's claim that the three approaches occupy distinct points in a trade-off space.  The safety radius is fastest to compute (a single quantile) and easiest to compose (radii add), but fails on non-spherical distributions.  The HDR is shape-adaptive but requires density estimation and a bandwidth parameter.  The Wasserstein approach provides a financially interpretable cost metric but requires nearest-neighbor computation against the constraint set.  Running all three on the same scenario, as the comparison tab does, allows the user to calibrate when the simpler approach suffices and when it does not.

Fourth, the visual interfaces make the geometric content of the theory tangible.  Watching the menu shrink as the fee cap tightens (Application~I), or seeing the HDR trace a curved contour while the safety ball extends into forbidden territory (Application~II), provides immediate intuition that complements the formal proofs.  This could be particularly useful for communicating the framework to portfolio managers who are not conversant with category theory or measure theory.

\clearpage

\section*{Glossary}
\addcontentsline{toc}{section}{Glossary and Key References}

\begin{description}

\item[Action ($\odot$)]
The fundamental operator in the DOTS framework that applies a relational constraint to a state. It generates the set of all possible outputs reachable from valid inputs via the relation.
(see Proposition~\ref{prop:action}, Section~\ref{sec:DOTS}, Appendix~\ref{app:dots_framework})

\item[Adjunction ($f_! \dashv f^*$)]
A relationship between two functors (or operators) where one is the optimal approximation of the other. In this framework, it formally links the pushforward ($f_!$) and pullback ($f^*$) operations, asserting that verifying a constraint $S$ in the spoke space is logically equivalent to verifying the pulled-back constraint $f^*S$ in the hub space.
(see Theorem~\ref{thm:adjunction}, Theorem~\ref{thm:safety-adjunction}, Theorem~\ref{thm:wasserstein-adjunction})

\item[Alignment Relation]
A vertical morphism in the double category $\mathbb{HS}$: a closed subset $R \subseteq K_1 \times K_2$ of the product of two permissible spaces, expressing a constraint or compatibility condition between portfolios. Examples include tracking error bounds, asset class exposure limits, factor exposure tolerances, and strategy identity requirements. Unlike a function, a single hub portfolio may be aligned with many spoke portfolios, and vice versa.
(see Section~\ref{sec:constraints}, Section~\ref{sec:double})

\item[Beck--Chevalley Condition]
A commutation law relating the composition of pushforward and pullback operations across a commuting square. The strict condition implies complete path independence (order of operations does not matter), while the lax version serves as a ``Safety Guarantee,'' ensuring upstream verification is conservative.
(see Proposition~\ref{prop:lax-BC}, Theorem~\ref{thm:Beck--Chevalley}, Theorem~\ref{thm:DOTS-BC}, Theorem~\ref{thm:early-audit})

\item[Berge's Maximum Theorem]
A foundational result in optimization theory providing conditions under which the set of optimal solutions varies continuously with parameters. It is used to prove that optimizing a re-implementation problem constructs a well-defined, continuous horizontal morphism.
(see Theorem~\ref{thm:berge-optimal}, Appendix~\ref{app:analytical_foundations})

\item[Convolution Stability]
A condition required for the coherent composition of probabilistic pushforwards. It ensures that the ``likely'' outcomes of a composite process are contained within the set of outcomes generated by chaining the likely intermediate steps, preventing mass from leaking into the tails.
(see Theorem~\ref{thm:pushforward-comp})

\item[Cost of Exit]
The Wasserstein distance from a realized portfolio distribution to the nearest fully liquid portfolio. It measures the expected turnover required to liquidate all illiquid positions.
(see Section~\ref{sub:liq-wasserstein})

\item[Correspondence (Multifunction)]
A mapping $\Phi: X \rightrightarrows Y$ that assigns to every input $x$ a subset of $Y$. In this framework, correspondences model feasible sets of portfolios or ``menus'' of options before a specific selection is made.
(see Appendix~\ref{app:analytical_foundations})

\item[DOTS (Double Operadic Theory of Systems)]
An extension of the framework that models systems where a single hub portfolio maps to a ``menu'' of permissible spoke portfolios rather than a single deterministic result.
(see Section~\ref{sec:DOTS}, Appendix~\ref{app:dots_framework})

\item[Double Category ($\mathbb{HS}$)]
A categorical structure consisting of objects (permissible spaces), horizontal morphisms (re-implementations), and vertical morphisms (alignment relations). It allows simultaneous modeling of functional transformations and relational constraints.
(see Section~\ref{sec:double}, Appendix~\ref{app:categorical_prereqs})

\item[Feller Kernel]
A stochastic kernel that maps continuous bounded functions to continuous bounded functions. In the probabilistic extension (HSP-r), horizontal morphisms are modeled as tight Feller kernels to ensure stability under perturbation.
(see Section~\ref{sub:HSP-r}, Appendix~\ref{subsec:stochastic})

\item[Frobenius Reciprocity]
A logical law stating that filtering a universe of assets before optimization yields the same valid set as optimizing first and filtering the results afterwards ($f_!(R \cap f^*S) = f_!R \cap S$). In probabilistic settings, this becomes an inclusion that guarantees conservative filtering.
(see Theorem~\ref{thm:frobenius}, Theorem~\ref{thm:DOTS-Frob}, Theorem~\ref{thm:conservative-filtering})

\item[Highest Density Region (HDR)]
A method for defining the ``effective support'' of a probability distribution by selecting the smallest closed set containing a specified mass (e.g., $1-\epsilon$) based on probability density. It defines the ``safety region'' of a stochastic solver.
(see Definition~\ref{def:epsilon-supp}, Section~\ref{sub:HDR})

\item[Hub / Spoke]
The organizing metaphor of the framework. A \emph{hub} portfolio is the source of investment intent (e.g., a model portfolio defined by a CIO); a \emph{spoke} portfolio is a concrete implementation of that intent in a different vehicle, universe, or constraint environment (e.g., an ETF, a tax-managed SMA, a target-date fund). The mathematical structure formalizes the translation between them and guarantees that the translation is coherent.
(see Section~\ref{sec:double}, Section~\ref{sec:propagation})

\item[Kantorovich--Rubinstein Duality]
A duality theorem expressing the Wasserstein distance as the supremum of expectations over 1-Lipschitz functions. It allows the ``cost of cure'' to be estimated efficiently using risk factors.
(see Appendix~\ref{app:optimal_transport})

\item[Markov Kernel]
A function $P: X \times \mathcal{B}(Y) \to [0,1]$ that acts as a stochastic mapping, assigning a probability measure on $Y$ to every $x \in X$.
(see Appendix~\ref{subsec:stochastic})

\item[Michael's Selection Theorem]
A result stating that a lower hemicontinuous correspondence with closed convex values admits a continuous selection. It ensures a continuous re-implementation map can be chosen even from a ``menu'' of options.
(see Theorem~\ref{thm:michael-selection})

\item[Path Independence]
The operational consequence of the Beck--Chevalley condition: if compliance is verified at an upstream stage, it need not be re-verified downstream. Equivalently, the order in which re-implementations and compliance checks are performed does not affect the final set of valid portfolios. The lax version (Audit Safety) guarantees that upstream verification is conservative; the strict version requires a pointwise cartesian condition.
(see Proposition~\ref{prop:lax-BC}, Theorem~\ref{thm:Beck--Chevalley}, Section~\ref{sec:propagation})

\item[Permissible Space]
An object in the double category $\mathbb{HS}$: a closed subset $K \subseteq \Delta^n$ of a simplex, representing the universe of portfolios that satisfy a given set of constraints. Closedness is not a technicality but the condition that prevents phantom portfolios; compactness (automatic, since $\Delta^n$ is compact) is what makes properness free.
(see Section~\ref{sec:spaces})

\item[Persona]
A parameter indexing a family of permissible spaces---e.g., the age cohort in a target-date glide path, or a client segment in a demographic conditioning scheme. The persona framework ($\mathbb{HS}/B$) fibers the entire theory over a base space $B$ of personas, so that all results (adjunction, Beck--Chevalley, Frobenius) hold persona by persona.
(see Section~\ref{sec:personas})

\item[Phantom Portfolio]
A spoke portfolio that satisfies every local compliance check but has no valid parent in the hub space. Phantoms arise when permissible spaces are not closed: a sequence of valid portfolios converges to a limit outside the permissible set, and the optimizer, with impeccable numerical precision, produces it. The closedness requirement on all objects in $\mathbb{HS}$ exists to prevent this pathology.
(see Section~\ref{sec:spaces}, Section~\ref{sec:proper})

\item[Pointwise Cartesian Square]
A condition on a commuting square of morphisms ensuring that for every matching pair of outputs, there exists a valid common pre-image in the domain. It is necessary for the strict Beck--Chevalley law.
(see Definition~\ref{def:pointwise-cartesian}, Section~\ref{sec:personas})

\item[Polish Space]
A separable, completely metrizable topological space. These spaces provide the necessary topological foundation for the probabilistic extension (HSP-r) to ensure measures are well-behaved.
(see Section~\ref{sec:Polish}, Appendix~\ref{subsec:stochastic})

\item[Private Asset Premium]
The total excess return demanded by investors for holding private (illiquid) assets, decomposed into eight components: pure illiquidity premium, risk premium, guarantee value, information costs, transaction complexity, rent of balance sheet, control alpha, and offsets.
(see Section~\ref{sec:liquidity}, \cite{Phoa})

\item[Portmanteau Theorem]
A fundamental theorem establishing equivalent conditions for the weak convergence of probability measures. It relates solver stability to the upper semi-continuity of safety regions.
(see Theorem~\ref{thm:portmanteau})

\item[Prokhorov's Theorem]
A theorem characterizing the relative compactness of probability measures in terms of tightness. It underpins the definition of tight Feller kernels, ensuring mass does not escape to infinity.
(see Appendix~\ref{subsec:stochastic})

\item[Proper Map]
A continuous map where the preimage of every compact set is compact. In $\mathbb{HS}$, this guarantees that the pushforward of a closed alignment relation remains closed, preventing ``phantom portfolios.''
(see Proposition~\ref{prop:reimpl-are-proper}, Appendix~\ref{app:analytical_foundations})

\item[Pullback ($f^*$)]
The operation that translates a spoke-side constraint into a hub-side constraint. Given a re-implementation $f: K_1 \to K_2$ and an alignment relation $S \subseteq K_2 \times K_3$ on the spoke, the pullback $f^*S = \{(x, z) : (f(x), z) \in S\}$ is the set of hub portfolios whose images under $f$ satisfy $S$. Informally: ``which hub portfolios will automatically produce compliant spokes?''
(see Section~\ref{sec:proper}, Theorem~\ref{thm:adjunction})

\item[Pushforward ($f_!$)]
The operation that propagates a hub-side constraint to the spoke side. Given a re-implementation $f: K_1 \to K_2$ and an alignment relation $R \subseteq K_1 \times K_3$ on the hub, the pushforward $f_!R = \{(f(x), z) : (x, z) \in R\}$ is the set of spoke portfolios that arise from aligned hub portfolios. Properness guarantees that $f_!R$ is closed whenever $R$ is.
(see Section~\ref{sec:proper}, Theorem~\ref{thm:proper-closed}, Theorem~\ref{thm:adjunction})

\item[Re-implementation]
A horizontal morphism in the double category $\mathbb{HS}$: a continuous map $f: K_1 \to K_2$ between permissible spaces, representing a deterministic process that converts a hub portfolio into a spoke portfolio. Typical re-implementations include multi-objective optimization, vehicle replacement (minimum tracking error or maximum return), and mechanical substitution (e.g., swapping foreign securities for ADR equivalents). In the DOTS extension, re-implementations are generalized to set-valued ``menus''; in HSP-r, to stochastic kernels.
(see Section~\ref{sec:double}, Section~\ref{sec:DOTS}, Section~\ref{sec:Polish})

\item[Safety Radius ($r_\epsilon$)]
A scalar quantifying the maximum deviation of a stochastic solver from its deterministic center, such that the solver output falls within a metric ball of radius $r_\epsilon$ with probability at least $1 - \epsilon$. The Safety Radius approach models solver error as symmetric and uniform, enabling fast geometric verification at the cost of conservatism. It is contrasted with the HDR approach, which tracks exact density contours.
(see Section~\ref{sub:stoch_ops}, Section~\ref{sub:concentration})

\item[Set-Valued Analysis]
The mathematical study of correspondences and their analytic properties, providing the rigorous basis for handling ``menus'' of portfolios and stability of constraints.
(see Appendix~\ref{app:analytical_foundations})

\item[Tightness]
A property of a stochastic kernel indicating that probability mass effectively remains within compact regions. It is the stochastic analogue to the properness of deterministic maps.
(see Definition~\ref{def:tight}, Section~\ref{sub:tightness})

\item[Upper Hemicontinuity (UHC)]
A continuity property for correspondences. It ensures that the set of valid portfolios does not ``explode'' or drastically change shape with small changes in the input.
(see Appendix~\ref{app:analytical_foundations}, Remark~\ref{rem:UHC})

\item[Wasserstein Distance ($W_1$)]
A metric defined as the minimum cost to transport mass from one distribution to another. In the Transport-Based Safety model, it quantifies the expected turnover required to move a non-compliant portfolio into compliance.
(see Section~\ref{sec:wasserstein}, Appendix~\ref{app:optimal_transport})

\end{description}

\clearpage

\begin{theindex}
\addcontentsline{toc}{section}{Index}

  \item Action ($\odot$), Section~\ref{sec:DOTS}, Prop.~\ref{prop:action}
  \item Action Menu, Def.~\ref{def:action-menu}
  \item Adjunction ($f_! \dashv f^*$)
    \subitem in $\mathbb{HS}$, Theorem~\ref{thm:adjunction}
    \subitem in $\mathbb{HS}/B$ (Personas), Section~\ref{sec:personas}
    \subitem in $\mathbb{S}\mathbf{pan}$, Section~\ref{sec:2-cells}
    \subitem Probabilistic (Safety Guarantee), Theorem~\ref{thm:safety-adjunction}
    \subitem Wasserstein, Theorem~\ref{thm:wasserstein-adjunction}
  \item Alignment Relation
    \subitem Definition, Section~\ref{sec:constraints}
    \subitem Vertical Morphism, Section~\ref{sec:double}
  \item Ambient Portfolio Space ($\Delta^n$), Def. in Section~\ref{sec:spaces}
  \item Architecture (System), Appendix~\ref{app:architecture}
  \item Asset Set, Def. in Section~\ref{sec:spaces}
  \item Attribute Map, Section~\ref{sec:optimal-reimpl}
  \item Audit Safety, Prop.~\ref{prop:lax-BC}, Section~\ref{sec:propagation}
  \item Audit Trail, Section~\ref{sec:operations}, Section~\ref{sec:2-cells}
  \item Axioma, Ex.~\ref{ex:gurobi-feller}

  \indexspace

  \item Banana (Distribution Shape), Section~\ref{sub:banana}, Section~\ref{sec:Polish}
  \item Beck--Chevalley Condition
    \subitem DOTS, Theorem~\ref{thm:DOTS-BC}
    \subitem Lax (Audit Safety), Prop.~\ref{prop:lax-BC}
    \subitem Probabilistic (Early-Audit), Theorem~\ref{thm:early-audit}
    \subitem Strict, Theorem~\ref{thm:Beck--Chevalley}, Theorem~\ref{thm:BC}
    \subitem Wasserstein, Theorem~\ref{thm:wasserstein-bc}
  \item Bellman Consistency, Theorem~\ref{thm:bellman_commutativity}
  \item Berge's Maximum Theorem, Theorem~\ref{thm:berge-optimal}, Appendix~\ref{app:analytical_foundations}
  \item Bicategory (Framed), Section~\ref{sec:abstract}, Appendix~\ref{app:dots_framework}

  \indexspace

  \item Cartesian Square
    \subitem Pointwise (in $\mathbb{HS}$), Def.~\ref{def:pointwise-cartesian}
    \subitem over Personas, Def. in Section~\ref{sec:personas}
  \item Category Theory
    \subitem Double Category, Section~\ref{sec:double}, Appendix~\ref{app:categorical_prereqs}
    \subitem Thin Category, Section~\ref{sec:double}
  \item Closedness
    \subitem of Objects, Section~\ref{sec:spaces}
    \subitem of Pushforward, Theorem~\ref{thm:proper-closed}
    \subitem of Vertical Composition, Lemma~\ref{lem:vertical-closed}
  \item Compactness
    \subitem of Permissible Spaces, Section~\ref{sec:spaces}
    \subitem Role in Properness, Section~\ref{sec:proper}, Appendix~\ref{app:analytical_foundations}
  \item Companion and Conjoint, Def.~\ref{def:equipment}
  \item Compliance (as Alignment), Remark~\ref{rem:compliance}
  \item Compliance with Cure, Section~\ref{sec:wasserstein}
  \item Composition
    \subitem Horizontal, Section~\ref{sec:double}
    \subitem Operadic, Section~\ref{sec:DOTS}
    \subitem Probabilistic (Accumulation of Risk), Theorem~\ref{thm:pullback-comp}
    \subitem Vertical (Relational), Def.~\ref{def:vertical-composition}
    \subitem Wasserstein (Lax), Theorem~\ref{thm:wasserstein-pushforward-comp}
  \item Concentration of Measure, Section~\ref{sub:concentration}
  \item Constraint
    \subitem Cross-vintage, Section~\ref{sec:personas}
    \subitem Global (Operadic), Section~\ref{sec:DOTS}
  \item Constraint Hugging, Section~\ref{sub:hollow_shell}, Section~\ref{sec:Polish}
  \item Convolution Stability, Theorem~\ref{thm:pushforward-comp}
  \item Core--Satellite Strategy, Section~\ref{sec:DOTS}, Section~\ref{sub:liq-dots}
  \item Correspondence (Multifunction), Appendix~\ref{app:analytical_foundations}
  \item Cure Budget, Section~\ref{sub:wasserstein-compute}

  \indexspace

  \item DOTS (Double Operadic Theory of Systems), Section~\ref{sec:DOTS}
  \item Dynamic Trading, Section~\ref{sub:liq-wasserstein}
  \item Double Category
    \subitem $\mathbb{HS}$ (Hub-and-Spoke), Def. in Section~\ref{sec:double}
    \subitem $\mathbb{HS}/B$ (Persona-indexed), Def. in Section~\ref{sec:personas}
    \subitem $\mathbf{HSP-r}$ (Risk), Def. in Section~\ref{sub:HSP-r}
    \subitem $\mathbb{S}\mathbf{pan}$ (Evidence), Section~\ref{sec:2-cells}

  \indexspace

  \item Equipment (Framed Bicategory), Def.~\ref{def:equipment}, Appendix~\ref{app:dots_framework}
  \item $\epsilon$-Support (HDR), Def.~\ref{def:epsilon-supp}
  \item Evidence Ledger, Appendix~\ref{app:architecture}, Section~\ref{sec:Polish}
  \item Evidence Space, Section~\ref{sec:2-cells}

  \indexspace

  \item Factor Decomposition, Example~\ref{ex:factor}
  \item Feller Kernel
    \subitem Definition, Def. in Section~\ref{sub:HSP-r}, Def. in Section~\ref{subsec:stochastic}
    \subitem Tightness, Def.~\ref{def:tight}, Def. in Section~\ref{sub:tightness}
  \item Frobenius Reciprocity
    \subitem DOTS, Theorem~\ref{thm:DOTS-Frob}
    \subitem Geometric (Stress-Test), Theorem~\ref{thm:geometric-frobenius}
    \subitem in $\mathbb{HS}$, Theorem~\ref{thm:frobenius}
    \subitem Probabilistic (Conservative Filtering), Theorem~\ref{thm:conservative-filtering}
    \subitem Wasserstein, Theorem~\ref{thm:wasserstein-frobenius}

  \indexspace

  \item Galois Connection, Appendix~\ref{app:categorical_prereqs}
  \item Gaussian Mixture Model (GMM), Section~\ref{sub:storage}
  \item Glide Path, Section~\ref{sec:personas}
  \item Graph (of a map), Section~\ref{sec:double}
  \item Gross Exposure Limit, Section~\ref{sub:hollow_shell}
  \item Gurobi, Example~\ref{ex:gurobi-feller}

  \indexspace

  \item Hausdorff Space, Prop.~\ref{prop:simplex-properties}
  \item Highest Density Region (HDR), Def.~\ref{def:epsilon-supp}, Section~\ref{sub:HDR}
  \item Hollow Shell (Distribution Shape), Section~\ref{sub:hollow_shell}, Section~\ref{sec:Polish}
  \item Horizontal Morphism
    \subitem Definition, Def. in Section~\ref{sec:reimplement}
    \subitem over Personas, Def. in Section~\ref{sec:personas}
  \item Hub-and-Spoke ($\mathbb{HS}$), \textit{See Double Category}

  \indexspace

  \item Illiquidity Discount, Section~\ref{sub:liq-wasserstein}
  \item Illiquidity Premium, Section~\ref{sec:liquidity}, Section~\ref{sub:liq-wasserstein}
  \item Interchange Law, Section~\ref{sec:abstract}

  \indexspace

  \item Kantorovich--Rubinstein Duality, Theorem in Appendix~\ref{app:optimal_transport}

  \indexspace

  \item Lax Beck--Chevalley, Prop.~\ref{prop:lax-BC}
  \item Lax Functor, Appendix~\ref{app:lax_functors}
  \item Leverage (Lipschitz Constant), Remark~\ref{rem:lipschitz-leverage}
  \item Liquidity Constraints, Section~\ref{sec:DOTS}, Section~\ref{sec:liquidity}
  \item Liquidity-Aware Portfolio Construction, Section~\ref{sec:liquidity}
    \subitem Aggregate Illiquidity Cap, Section~\ref{sub:liq-basic}
    \subitem Core--Satellite with Liquidity Gate, Section~\ref{sub:liq-dots}
    \subitem Cost of Exit, Section~\ref{sub:liq-wasserstein}
    \subitem Dynamic Rebalancing, Section~\ref{sub:liq-wasserstein}
    \subitem Maintenance Budget, Section~\ref{sub:liq-wasserstein}
    \subitem Market Impact Constraint, Section~\ref{sub:liq-dots}
    \subitem Weighted Transport Cost, Section~\ref{sub:liq-wasserstein}

  \indexspace

  \item Markov Kernel, Def. in Section~\ref{subsec:stochastic}
  \item Menu (DOTS), Section~\ref{sec:DOTS}
  \item Michael's Selection Theorem, Theorem~\ref{thm:michael-selection}
  \item Minimum Position Size, Section~\ref{sub:split_peak}
  \item Minkowski Sum, Section~\ref{sub:concentration}
  \item Multi-hub Alignment, Section~\ref{sec:DOTS}
  \item Monoidal Structure, Appendix~\ref{app:dots_framework}
  \item Monte Carlo Projection, Section~\ref{sub:wasserstein-compute}

  \indexspace

  \item Object (Permissible Portfolio Space)
    \subitem Definition, Def. in Section~\ref{sec:spaces}
    \subitem Invalid examples, Section~\ref{sec:non-proper-examples}
    \subitem Persona-indexed, Def. in Section~\ref{sec:personas}
  \item Objective Function, Section~\ref{sec:optimal-reimpl}
  \item Operad, Section~\ref{sec:DOTS}
  \item Optimal Execution, Section~\ref{sub:liq-dots}
  \item Optimal Transport, Appendix~\ref{app:optimal_transport}, Section~\ref{sec:wasserstein}
  \item Optimization
    \subitem Convex, Section~\ref{sec:computation}
    \subitem for Re-implementation, Section~\ref{sec:optimal-reimpl}
    \subitem Stochastic, Section~\ref{sec:Polish}

  \indexspace

  \item Path Independence, Section~\ref{sec:Beck--Chevalley}, Section~\ref{sec:propagation}
  \item Permissible Portfolio Space ($K$), Def. in Section~\ref{sec:spaces}
  \item Persona (Base Space $B$), Section~\ref{sec:personas}
  \item Phantom Portfolio, Example~\ref{ex:pushforward-not-closed}, Section~\ref{sec:failure}, Section~\ref{sub:split_peak}
  \item Pointwise Cartesian, Def.~\ref{def:pointwise-cartesian}
  \item Polish Space, Section~\ref{sec:Polish}, Def. in Section~\ref{subsec:stochastic}
  \item Portmanteau Theorem, Theorem~\ref{thm:portmanteau}
  \item Portfolio Space of Quantitative Attributes, Def. in Section~\ref{sec:optimal-reimpl}
  \item Pre-Trade Compliance, Section~\ref{subsec:wasserstein-properties}, Section~\ref{sub:liq-wasserstein}
  \item Private Asset Premium, Section~\ref{sec:liquidity}
  \item Prokhorov's Theorem, Theorem in Section~\ref{subsec:stochastic}
  \item Proper Map
    \subitem Definition, Def. in Section~\ref{sec:proper}, Def. in Appendix~\ref{app:analytical_foundations}
    \subitem Horizontal Morphisms as, Prop.~\ref{prop:reimpl-are-proper}
  \item Pullback ($f^*$)
    \subitem Definition, Def.~\ref{def:pullback}
    \subitem Diagnostic vs.\ Prescriptive Use, Section~\ref{sec:propagation}
    \subitem HDR ($\epsilon$-Pullback), Section~\ref{sub:HDR}
    \subitem Safety Radius ($\epsilon$-Pullback), Def. in Section~\ref{sub:concentration}
    \subitem Wasserstein, Def. in Section~\ref{sec:wasserstein}
  \item Pushforward ($f_!$)
    \subitem Definition, Def.~\ref{def:pushforward}
    \subitem HDR ($\epsilon$-Pushforward), Section~\ref{sub:HDR}
    \subitem Safety Radius ($\epsilon$-Pushforward), Def. in Section~\ref{sub:concentration}
    \subitem Wasserstein, Def.~\ref{def:wasserstein-pushforward}

  \indexspace

  \item Quantitative Attributes ($V$), Section~\ref{sec:optimal-reimpl}

  \indexspace

  \item Randomized Smoothing, Example~\ref{ex:gurobi-feller}, Section~\ref{sub:randomized-smoothing}
  \item Ratio Constraint (Reformulation), Section~\ref{sec:non-proper-examples}
  \item Re-implementation
    \subitem Definition, Def. in Section~\ref{sec:reimplement}
    \subitem Optimal Construction, Theorem~\ref{thm:berge-optimal}
  \item Relation, \textit{See Alignment Relation}
  \item Risk Budget, Section~\ref{sec:Polish}

  \indexspace

  \item Safety Radius, Def. in Section~\ref{sub:concentration}
  \item Simplex ($\Delta^n$), Def. in Section~\ref{sec:spaces}
  \item Slice Category, Section~\ref{sec:personas}
  \item Solver (Interface), Appendix~\ref{app:architecture}
  \item Span (Evidence), Section~\ref{sec:2-cells}
  \item Split Peak (Distribution Shape), Section~\ref{sub:split_peak}, Section~\ref{sec:Polish}
  \item Spoke, Remark in Section~\ref{sec:reimplement}
  \item State (DOTS), Section~\ref{sec:DOTS}, Appendix~\ref{app:dots_framework}
  \item Stochastic Kernel, Section~\ref{sec:Polish}, Appendix~\ref{subsec:stochastic}
  \item Support Functor, Section~\ref{sec:conclusion}

  \indexspace

  \item Target Date Fund, Section~\ref{sec:personas}
  \item Thinness (of Double Category), Section~\ref{sec:double}
  \item Tightness (of Kernels), Def.~\ref{def:tight}, Def. in Section~\ref{sub:tightness}
  \item Tracking Error, Ex. in Section~\ref{sec:2-cells}, Section~\ref{sec:DOTS}
  \item Transport Cost, Section~\ref{sec:wasserstein}
  \item Two-Cell (2-cell)
    \subitem Evidence Transformation, Section~\ref{sec:2-cells}
    \subitem Inclusion In $\mathbb{HS}$, Section~\ref{sec:double}

  \indexspace

  \item Value-at-Risk (VaR), Section~\ref{sec:Polish}
  \item Vertical Morphism
    \subitem Definition, Def. in Section~\ref{sec:constraints}
    \subitem over Personas, Def. in Section~\ref{sec:personas}

  \indexspace

  \item Virtual Hub (Interpolation), Section~\ref{sec:personas}

  \indexspace

  \item Wasserstein Distance ($W_1$), Section~\ref{sec:wasserstein}, Def. in Appendix~\ref{app:optimal_transport}
  \item Wiring Diagram, Def. in Section~\ref{app:dots_framework}, Fig.~\ref{fig:operad-wiring}

\end{theindex}

\end{document}